\documentclass[11pt,a4paper]{amsart}
\usepackage[margin=2.5cm]{geometry}
\usepackage{contactgmt-macros}
\usepackage{xcolor,needspace}
\usepackage[colorlinks=true,linkcolor=blue!45!black,citecolor=green!35!black,urlcolor=blue!45!black]{hyperref}
\usepackage[expansion=false]{microtype}
\theoremstyle{plain}
\newtheorem*{claim}{Claim}
\theoremstyle{definition}
\newtheorem{fact}[theorem]{Fact}

\providecommand{\cN}{\mathcal N}

\providecommand{\cV}{\mathcal V}
\providecommand{\cB}{\mathcal B}
\providecommand{\cW}{\mathcal W}

\providecommand{\cG}{\mathcal G}

\newcommand{\bB}{\mathbf B}
\newcommand{\bb}{\mathsf b}
\newcommand{\Kreg}{\mathcal K}

\hypersetup{pdftitle={Local regularity and the first singular time of Legendrian mean curvature flow},pdfauthor={Shu-Cheng Chang, Hongbing Qiu, Liuyang Zhang}}
\begin{document}
\title[Local regularity for Legendrian mean curvature flow]{Local regularity and the first singular time of Legendrian mean curvature flow}
\author[Chang]{Shu-Cheng Chang}
\author[Qiu]{Hongbing Qiu}
\author[Zhang]{Liuyang Zhang}
\date{23 November 2026}
\begin{abstract}
We prove a local regularity theorem for Legendrian mean curvature flow with bounded
angle in a compact Sasakian manifold. We assume an evolution equation for the angle;
it holds in the $\eta$-Einstein case with exact mean curvature form.
If the height-localised Gaussian density is close to one at all flow points and at all scales
below the truncated boundary distance, the flow is locally a single Legendrian graph with
$C^{2,\alpha}$ estimates at a comparable scale. No global projection along the Reeb field
is needed. For a closed flow, let $\eps_*$ be the density excess in the local theorem
at $\alpha=1/2$. At a point of the first singular time where the limiting density is
less than $1+\eps_*$, there is either an empty neighbourhood or one embedded sheet
smooth up to that time.
At a Type~I curvature singularity in the $\eta$-Einstein case the density is an integer
at least two.
\end{abstract}
\maketitle
\enlargethispage{2pt}
\tableofcontents
\section{Introduction}\label{sec:intro}

White's local regularity theorem gives a way to recognise a smooth sheet of a mean curvature
flow from its Gaussian density ratios. If the ratios are sufficiently close to one
at all flow points and all scales below the distance to the boundary, the flow has
$C^{2,\alpha}$ estimates at a comparable scale \cite[Theorem~3.1]{wh}.
We prove a local regularity theorem of this kind for the Legendrian mean curvature flow.
The choice of density matters here. Two Legendrian sheets may project to the same
transverse point and lie at different heights along the Reeb field.

The first and third authors, together with Wu, used a transverse Gaussian
that counts sheets approaching the same Reeb orbit at any height
\cite[Remark~4.13]{ChangWuZhang2026ContactBlowup}. In the regular case a small terminal density gives
a curvature bound above a neighbourhood of the projected point \cite[Lemma~4.12]{ChangWuZhang2026ContactBlowup}.
To study the flow in the contact manifold one needs to distinguish the heights.
In \cite{PaperI} we used the angle-corrected Reeb height to localise the transverse Gaussian.
Here we use that density to prove single-sheet regularity. Our estimates are local
in the contact manifold and require no global projection along the Reeb field.

The static contact problem has other regularity theories. For Hamiltonian stationary
Legendrian surfaces there are monotonicity formulae sensitive to the Reeb direction
\cite{riv24,pr24}; for Legendrian area-minimising currents a partial regularity
theory is available without such a formula \cite{or26}.

\subsection{Local regularity}
Let $M^{2n+1}$ be a compact Sasakian manifold. We consider a smooth local Legendrian
mean curvature flow with bounded real angle $\theta$, in the sense of
Definition~\ref{def:localflow}. We assume the flow and angle equations \eqref{eq:F};
these hold in the $\eta$-Einstein case with exact mean curvature form
\cite[\S2.2]{PaperI}. The flow moves in its mean curvature direction and
along the Reeb field with velocity $-2\theta\xi$ \cite{s1,ChangWuZhang2026ContactBlowup}.

For a flow point $X$ in an open space-time region $\cU$, let $d_1(X,\cU)$ be the
distance to the complement truncated at the geometric radius $r_1$.
The density $\Theta_X^\chi(\sigma)$ is measured at backward time $\sigma$.
In the Heisenberg model the weight selects heights within order $\sigma$ of
$4\sigma\theta$; the kernel in contact normal coordinates is recalled in
Section~\ref{subsec:inputs}. Theorem~\ref{thm:A} shows that if
\[
  \Theta_X^\chi(\sigma)\le 1+\eps
       \quad\text{for every flow point }X
       \text{ and every }0<\sigma<d_1(X,\cU)^2 ,
\]
then near each such point the flow is a single Legendrian graph with
$C^{2,\alpha}$ control at a scale comparable to $d_1(X,\cU)$.
More precisely, $d_1(X,\cU)\Kreg(\cM,X)\le C$. The constants $\eps$ and $C$
depend only on the dimension, the weight parameter, the exponent and the angle
bound. The ambient geometry enters through $r_1$. To estimate curvature at time $t$, one can truncate the flow at $t$.
Corollary~\ref{cor:past} then gives the curvature bound under the density
hypothesis up to the truncated distance to the past boundary.

\subsection{The first singular time}
Let the flow be closed and satisfy the hypotheses (S) and (F) of
\cite[\S2.2]{PaperI}, and let $T$ be its first singular time. The limiting density
$\Theta^\chi(p,T)$ exists by \cite[Theorem~A]{PaperI}. Set
$\eps_*=\eps(n,a,\tfrac12,C_\theta)$, the threshold in Theorem~\ref{thm:A}
at $\alpha=1/2$. If $\Theta^\chi(p,T)<1+\eps_*$, Theorem~\ref{thm:B}
gives two alternatives. Either the flow
is absent from a neighbourhood of $(p,T)$ near time $T$, or it consists there of
one embedded sheet, covered once, which extends smoothly up to $T$. In particular
$\Theta^\chi(p,T)\ge 1+\eps_*$ at every curvature singularity. Under the Type~I assumption in
an $\eta$-Einstein manifold, Corollary~\ref{cor:Bprime} identifies the density at a
curvature singularity with the total multiplicity of the lifted tangent planes;
it is an integer at least two. The transverse-density result for regular ambient
manifolds is \cite[Corollary~1.3]{ChangWuZhang2026ContactBlowup}.

\subsection{Outline of the proof}
We argue by contradiction, following \cite[Theorem~3.1]{wh}.
After point selection (Lemma~\ref{lem:pick}), the rescaled flows have regularity
scale one at the origin. Proposition~\ref{prop:compact} gives an ancient
Heisenberg limit with bounded regularity scale. Its density is at most one
at every point and scale, so
Proposition~\ref{prop:rigidity} identifies it as a single plane translating
vertically. For sufficiently large indices, Lemma~\ref{lem:onesheet} gives a
single graph near the origin on the rescaled flow. Lemma~\ref{lem:schauder}
improves the graph estimates, giving an admissible radius greater than one. This
contradicts the normalisation (Section~\ref{subsec:proofA}). The model theorem yields a gap for ancient flows
(Theorem~\ref{thm:gap}). Theorem~\ref{thm:alpha} gives a density threshold independent
of the H\"older exponent in both the compact ambient and Heisenberg settings.

At time $T$ we follow White's truncation argument \cite[Remark~3.3 and Theorem~3.5]{wh}.
A density bound at $(p,T)$ must first be transferred to nearby centres. The charts need
not vary continuously with their centres; Lemma~\ref{lem:compare-centres} compares the
densities with a continuous auxiliary density in one fixed chart. If the flow accumulates at $(p,T)$, uniform regularity on the truncations
gives one embedded sheet smooth up to $T$ (Theorem~\ref{thm:B}).
Appendix~\ref{app:schauder} supplies the parabolic estimates up to a closed upper end.

\section{Preliminaries}\label{sec:setting}

Throughout, $(M^{2n+1},\eta,\xi,\Phi,g)$ is a compact Sasakian manifold, and we use the
conventions of \cite[Section~2]{PaperI}, which are those of
\cite[\S2]{ChangWuZhang2026ContactBlowup}. In a contact normal chart we measure two quantities: the
height-localised Gaussian density and the radius on which the flow is a single
Legendrian graph. Theorem~\ref{thm:A} gives a lower bound for this radius from an
upper bound for the density. We first recall the charts and kernel of \cite{PaperI}
(Section~\ref{subsec:inputs}). The local flows considered here may include their upper
time slice (Section~\ref{subsec:localflows}). Sections~\ref{subsec:rM}
and~\ref{subsec:scale} define the geometric scale and the graph scale, respectively.
The estimates specific to closed flows enter in Section~\ref{sec:terminal}.

\subsection{Contact normal charts and the height-localised kernel}\label{subsec:inputs}

\emph{The model.} The Heisenberg group is $\Hn=\C^n\times\R$ with coordinates
$(w,z)$, $w=(x,y)$, and the structure
\[
\overline\eta=\tfrac12dz+\overline\beta,\quad
\overline\beta=-\tfrac14\sum_j(y_j\,dx_j-x_j\,dy_j),\quad
\overline\xi=2\partial_z,\quad
\overline g=g_0+\overline\eta\otimes\overline\eta,\quad
g_0=\tfrac14|dw|^2 .
\]
The left translations $\ell_{(a,b,c)}(x,y,z)=(x+a,y+b,z+c+\frac12(b\cdot x-a\cdot y))$
preserve the structure, and $\omega_0:=d\overline\eta=d\overline\beta=\frac12\sum_jdx_j\wedge dy_j$
is the corresponding symplectic form of $\C^n$. We put $\tau_b=\ell_{(0,0,b)}$. The contact dilations are
$\delta_\lambda(w,z)=(\lambda w,\lambda^2z)$, with generator
$X=\sum_j(x_j\partial_{x_j}+y_j\partial_{y_j})+2z\partial_z=X^{\cH}+z\overline\xi$. The
horizontal potential is $f=\frac14|X^{\cH}|^2=|w|^2/16$, with $\overline\nabla
f=\frac12X^{\cH}$. For $\hat y\in\Hn$ we write $f_{\hat y}=f\circ\ell_{\hat y}^{-1}$ and
$z_{\hat y}=z\circ\ell_{\hat y}^{-1}$. The complex structure $J_0$ on $\C^n$ is induced by
$\overline\Phi$, the functions $\zeta_j=\frac12(x_j-iy_j)$ are holomorphic, and in the
coordinates $\zeta$, the metric $g_0$ is the Euclidean metric. If $F=(\varphi,z)$ is a
Legendrian immersion into $\Hn$, then $dz=-2\varphi^*\overline\beta$ along $F$, and $\varphi$
is a Lagrangian immersion into $(\C^n,g_0,J_0)$ (not necessarily an embedding, even when
$F$ is embedded). Conversely, a Lagrangian immersion $\varphi$ of a connected manifold has
a horizontal lift $(\varphi,z)$ into $\Hn$ if and only if the closed form
$\varphi^*\overline\beta$ is exact. Exactness always holds locally and on simply connected
domains (in particular for Lagrangian planes). When a lift exists, it is unique up to a
vertical translation $\tau_b$.

\emph{The flow.} A Legendrian mean curvature flow is understood in the sense of hypothesis
(F) of \cite[Section~2.2]{PaperI}: a smooth normally parametrised family of Legendrian
immersions $F_t$ together with a smooth function $\theta$, the Legendrian angle, with
\begin{equation}\label{eq:F}
H=\Phi\nabla\theta,\qquad \partial_tF=H-2\theta\xi,\qquad
(\partial_t-\Delta)\theta=(K+2)\theta
\end{equation}
for a constant $K$. In the model, $K=-2$. We always assume a bound $|\theta|\le C_\theta$,
which for closed flows is the maximum principle bound of \cite[Section~2.2]{PaperI}.
Along \eqref{eq:F}, the velocity $\partial_tF=H-2\theta\xi$ is normal, and the mean
curvature vector $H=\Phi\nabla\theta$ is horizontal, so $\langle H,\xi\rangle=0$ and
$\partial_td\mu=-\langle H,H-2\theta\xi\rangle d\mu=-|H|^2d\mu$.

\Needspace{6\baselineskip}
\emph{Charts.} The first group of inputs concerns the contact normal coordinates.

\begin{fact}[{\cite[Lemma~2.3, Lemma~2.5, Lemma~2.7]{PaperI}}]\label{fact:charts}
There are $r_0,z_0>0$ and constants $C_k$, depending only on $(M,\eta,g)$, such that every
$p\in M$ has a chart $\Upsilon_p\colon U_p\to\{|w|<r_0,\ |z|<z_0\}$, $\Upsilon_p(p)=0$,
with the following properties.
\begin{enumerate}[(i)]
\item $(\Upsilon_p)_*\xi=2\partial_z$ and $\eta=\frac12dz+\beta_p$ with $\beta_p$ a
$1$-form in $w$ alone in radial gauge; $w$ are holomorphic coordinates of a local leaf
space in which the transverse metric $h_p$ satisfies $h_p=g_0+O(|w|^2)$, $\partial
h_p=O(|w|)$, $h_p\ge\frac12g_0$; $|\beta_p-\overline\beta|\le C|w|^3$,
$|\partial(\beta_p-\overline\beta)|\le C|w|^2$; all derivatives of $h_p,\beta_p$ are
bounded by $C_k$ uniformly in $p$.
\item (Model identities with errors.) With $f_p=|w|^2/16$ and $z_p=z$ in $U_p$: $\xi
f_p=0$, $\xi z_p=2$; $|\overline\nabla^2f_p(U,W)-\frac12g(U,W)|\le C|w|^2|U||W|$ and
$|\overline\nabla^2z_p(U,W)|\le C|w|^2|U||W|$ for horizontal $U,W$; the horizontal part of
$\overline\nabla z_p$ is $4\Phi\overline\nabla f_p+e_H$ with $|e_H|\le C|w|^3$; and along a
flow \eqref{eq:F}, $(\partial_t-\Delta)z_p=-4\theta+e_z$ with $|e_z|\le C|w|^2$, and
$\partial_tz_p=dz_p(H)-4\theta$ with $dz_p(V)=-2\beta_p(V)$ for horizontal $V$.
\item (Rescaling.) For $\lambda\ge1$ the map $\delta_\lambda\circ\Upsilon_p$ is a chart
of the same kind for the $\cD$-homothetic structure with $\kappa=\lambda^2$, on
$\{|w|<\lambda r_0,\ |z|<\lambda^2z_0\}$. For $R\ge1$ and $\lambda\ge R\max\{r_0^{-1},z_0^{-1/2}\}$
the rescaled structure
$\sigma_\lambda=(\delta_\lambda)_*(\eta_\kappa,\xi_\kappa,\Phi,g_\kappa)$ differs from the
model structure by at most $C_k(R)\lambda^{-2}$ in $C^k(\{|w|<R,|z|<R^2\})$. For $R\ge1$
there is $\lambda_*(R)$ such that, for $q\in U_p$ with $\hat y=(\hat w,\hat z):=\delta_\lambda\Upsilon_p(q)$,
$|\hat w|+|\hat z|^{1/2}\le R$, and $\lambda\ge\lambda_*(R)$, the transition
$\delta_\lambda\Upsilon_q\Upsilon_p^{-1}\delta_\lambda^{-1}$ is defined on $\{|w|<R,|z|<R^2\}$
and $C_k(R)\lambda^{-2}$-close there in $C^k$ to $\tilde U\circ\ell_{\hat y}^{-1}$ for some $U\in
U(n)$, $\tilde U(w,z)=(Uw,z)$; consequently $\lambda^2f_q$, $\lambda^2z_q$, written in
the rescaled $p$-chart, converge in $C^k_{\rm loc}$ to $f_{\hat y}$, $z_{\hat y}$,
uniformly for $\hat y$ in compact sets.
\end{enumerate}
\end{fact}

Fix once and for all a cut-off $\zeta_*\in C^\infty_c(\{|w|<r_0,|z|<z_0\})$ with
$0\le\zeta_*\le1$ and $\zeta_*=1$ on $\{|w|<\frac{r_0}2,|z|<\frac{z_0}2\}$, and, for
$q\in M$, put $\zeta_q:=\zeta_*\circ\Upsilon_q$ on $U_q$, extended by $0$, and
$\psi_q:=\zeta_q^2$, as in \cite[Section~2.3]{PaperI}. In the chart $\Upsilon_q$, the
metric $g$ is given by $h_q$ on horizontal vectors and by $\eta=\frac12dz+\beta_q$. Hence, by
Fact~\ref{fact:charts}(i) (the uniform bounds on $h_q$, $\beta_q$ and their derivatives,
and $h_q\ge\frac12g_0$), the coefficients of the metric, those of its inverse and all their derivatives
are bounded independently of $q$. Since $\zeta_*$ is fixed, $|\overline\nabla\zeta_q|$ and
$|\overline\nabla^2\zeta_q|$ are bounded independently of $q$.

\emph{The kernel.} The second group of inputs concerns the weighted kernel. Fix
$0<a<\frac1{32}$ and $\chi(v)=e^{-av^2}$, so that $\chi+16\chi''\ge c_\chi\chi$ with
$c_\chi=1-32a>0$. For a centre $(q,t_*)$ and $\tau=t_*-t>0$ put
\begin{gather*}
a_K(\tau)=\tfrac4{K+2}\big(e^{(K+2)\tau}-1\big)\ (=4\tau\text{ if }K=-2),\qquad
u=z_q-a_K(\tau)\theta,\qquad v=\frac u\tau,\\
\rho_{q,t_*}=(4\pi\tau)^{-n/2}e^{-f_q/\tau},\qquad \rho^\chi_{q,t_*}=\rho_{q,t_*}\chi(v),
\qquad \cS_q=H+\frac{(\overline\nabla f_q)^\perp}{\tau},\\
D^\chi=\frac{|v\chi'(v)|}\tau+(\chi+16\chi'')(v)|\cS_q|^2 .
\end{gather*}
In the model, $u=z-4\tau\theta$ is the \emph{caloric height}, which satisfies $(\partial_t-\Delta)u=0$
and $\nabla u=4\tau\overline\Phi\cS$ \cite[Lemma~3.1]{PaperI}.

\begin{fact}[{\cite[Proposition~3.2, Lemma~3.11]{PaperI}}]\label{fact:kernel}
\leavevmode
\begin{enumerate}[(i)]
\item (Exact monotonicity in the model under an area bound.) Let $\Sigma$ be an
$n$-manifold without boundary and $F\colon\Sigma\times[t_0,T)\to\Hn$ a smooth family of
Legendrian immersions satisfying \eqref{eq:F} with $K=-2$, with $(x,t)\mapsto(F(x,t),t)$
proper on $\Sigma\times[t_0,T']$ for every $T'<T$, $\sup|\theta|<\infty$ and
$\mu_t(F_t^{-1}\{|w|<R\})\le D_0R^n$ for all $R>0$. Then for $t_0\le t_1<t_2<T$ and the
centre $(0,T)$,
\[
\int\rho^\chi d\mu_{t_2}-\int\rho^\chi d\mu_{t_1}=-\int_{t_1}^{t_2}\!\!\int D^\chi\rho\,d\mu_t\,dt ,
\]
and $t\mapsto\int\rho^\chi d\mu_t$ is constant on $[t_1,t_2]$ if and only if $z=4\tau\theta$
and $H=-(\overline\nabla f)^\perp/\tau$ on $\Sigma_t$ for $t\in[t_1,t_2]$. The same holds
for every centre $(\hat y,T)$.
\item (Scaling.) Under the rescaling $(w',z',s)=(\lambda w,\lambda^2z_q,\lambda^2(t-t_*))$
with the $\cD$-homothetic structure, $\rho\,d\mu_t=\rho'd\mu'$, $v=v'$, $|H|^2dt=|H'|^2ds$,
$|\cS_q|^2dt=|\cS'|^2ds$ and $D^\chi dt=D^{\chi\prime}ds$, where the primed quantities are
computed in the rescaled structure with $\hat a(s)=a_K(\lambda^{-2}|s|)/(\lambda^{-2}|s|)\to4$.
\end{enumerate}
\end{fact}

\subsection{Local flows with a closed upper end}\label{subsec:localflows}

The local regularity theorem concerns flows defined in an open subset of space-time, as in
\cite[\S2.2--2.3]{wh}. We allow the flow to stop at a common closed upper end, which is
needed for the truncation argument of \cite[Remark~3.3]{wh}.

\begin{definition}\label{def:localflow}
Let $\cU\subset M\times\R$ be open and $T_*\in(-\infty,\infty]$. A \emph{smooth local
Legendrian flow in $\cU$ with upper end $T_*$} is a triple $(\Sigma,\cN,F)$ together with
a function $\theta$, where $\Sigma$ is an $n$-manifold without boundary,
$\cN\subset\Sigma\times(-\infty,T_*]$ is relatively open, and $F\colon\cN\to M$ is smooth
(with one-sided derivatives at $t=T_*$), such that
\begin{enumerate}[(a)]
\item $\hat F:=(F,t)\colon\cN\to\cU$ is a proper map;
\item for every $t$ the map $F_t=F(\cdot,t)$ on $\cN_t=\{x:(x,t)\in\cN\}$ is a normally
parametrised Legendrian immersion;
\item $\theta\colon\cN\to\R$ is smooth, $|\theta|\le C_\theta$, and \eqref{eq:F} holds on
$\cN$ for a constant $K$.
\end{enumerate}
The \emph{trace} is $\cM=\hat F(\cN)\subset\cU$; by (a) it is relatively closed in
$\cU$. We write $\cM(t)$ for the time-$t$ slice and $\mu_t$ for the area measure of
$F_t$ on $\cN_t$. When $T_*=\infty$ the flow is called \emph{fully smooth}.
\end{definition}

\begin{remark}[examples and truncation]\label{rem:localflow}\leavevmode
\begin{enumerate}[(1)]
\item A closed flow $F\colon\Sigma\times[0,T)\to M$ and an open set $\cU\subset
M\times(0,T)$ give the fully smooth local flow $(\Sigma,\hat F^{-1}(\cU),F|)$. The map
$\hat F$ is proper because $\Sigma$ is closed. The initial slice $t=0$ is excluded from
$\cU$.
\item (Truncation.) If $(\Sigma,\cN,F)$ is a smooth local flow in $\cU$ with upper end
$T_*$ and $t_*\le T_*$, then $F^{\le t_*}:=F|_{\cN\cap\{t\le t_*\}}$ is a smooth local flow
with upper end $t_*$ in
\[
\cU_{t_*}:=\cU\cup\{(q,t):t>t_*\}.
\]
Indeed, $\cN\cap\{t\le t_*\}$ is relatively open in $\Sigma\times(-\infty,t_*]$. For a
compact $\cK\subset\cU_{t_*}$, the set $\cK\cap\{t\le t_*\}$ is a compact subset of $\cU$,
whose preimage under $\hat F$ is compact. Hence $\hat F^{\le t_*}$ is proper. The trace of
the truncation is $\cM\cap\{t\le t_*\}$.
\item A flow with upper end $T_*$ restricted to an open set $\cU'\subset\cU$ is a flow
in $\cU'$ with the same upper end.
\end{enumerate}
\end{remark}

The limits produced in Section~\ref{sec:graphs} are given by their traces rather than by a
global parametrisation. We therefore also use the following trace form of the
definition, which is that of \cite[\S2.2]{wh}.

\begin{definition}[flows in trace form]\label{def:traceflow}
A \emph{smooth Legendrian flow in trace form} in $\cU$ with upper end $T_*$ is a
relatively closed subset $\cM\subset\cU\cap\{t\le T_*\}$, together with a smooth function
$\theta\colon\cM\to\R$, such that $\cM$ is an embedded $(n+1)$-dimensional submanifold
with boundary of $M\times\R$ (the boundary being $\cM\cap\{t=T_*\}$), the time function
has no critical points on $\cM$, and every point of $\cM$ has a neighbourhood $\cU'$ in
which $\cM\cap\cU'$ is the trace of a smooth local Legendrian flow in $\cU'$ in the sense
of Definition~\ref{def:localflow}, with the given $\theta$. Its slices $\cM(t)$ are
embedded Legendrian submanifolds, and $\mu_t$ denotes their area measure (multiplicity
one). A flow in the sense of Definition~\ref{def:localflow} whose map $\hat F$ is
injective is a flow in trace form; conversely, for a non-empty model flow in trace form
on $(t_0,T_*]$, that is, a flow in trace form in $\cU=\Hn\times(t_0,\infty)$ for the
model structure with $K=-2$ (Definition~\ref{def:modelflow}), whose regularity scale
$\Kreg$ (Definition~\ref{def:K}, for some exponent $\alpha\in(0,1)$) is bounded on
$\cM\cap(\Hn\times J)$ for every compact $J\subset(t_0,T_*]$, Lemma~\ref{lem:param} below
gives a global normal parametrisation which is proper, so that such a flow is a flow in the
sense of Definition~\ref{def:localflow}.
\end{definition}

\emph{Contact space-time quasi-distance.} Put
\begin{equation}\label{eq:bb}
\bb^2:=8C_\theta+4 .
\end{equation}
For $X=(q,t)$ and $Y=(q',t')$ with $q'\in U_q$ let
\[
\|Y-X\|_X:=\max\big\{|w_q(q')|,\ \bb^{-1}|z_q(q')|^{1/2},\ |t'-t|^{1/2}\big\},
\]
where $(w_q,z_q)=\Upsilon_q(q')$, and $\|Y-X\|_X:=+\infty$ if $q'\notin U_q$. In the
model, $\|Y-X\|=\max\{|w|,\bb^{-1}|z|^{1/2},|t'-t|^{1/2}\}$ with $(w,z)=\ell_q^{-1}(q')$.
This quasi-distance is left-invariant and homogeneous of degree one under $\delta_\lambda$ combined with
$t\mapsto\lambda^2t$. The \emph{cylinder} of radius $r$ about $X$ is
\[
\cC_r(X):=\Upsilon_q^{-1}\{|w|<r,\ |z|<\bb^2r^2\}\times(t-r^2,t+r^2)=\{Y:\|Y-X\|_X<r\}.
\]
For an open set $\cU\subset M\times\R$ the \emph{two-sided distance} and the \emph{past
distance} are
\[
d(X,\cU):=\inf\{\|Y-X\|_X:\ Y\notin\cU\},\qquad
d^-(X,\cU):=\inf\{\|Y-X\|_X:\ Y\notin\cU,\ t(Y)\le t(X)\},
\]
and we put $d_1:=\min\{d,r_1\}$ and $d_1^-:=\min\{d^-,r_1\}$, where $r_1$ is the constant of
Section~\ref{subsec:rM}. Note that $d\le d^-$. Thus $\cC_r(X)\subset\cU$ for $r<d(X,\cU)$,
and $\cC_r(X)\cap\{t\le t(X)\}\subset\cU$ for $r<d^-(X,\cU)$. For a truncation,
\begin{equation}\label{eq:dpast}
d\big(Y,\cU_{t_*}\big)\le d^-(Y,\cU)\qquad\text{whenever }t(Y)\le t_* ,
\end{equation}
because every $Z\notin\cU$ with $t(Z)\le t(Y)$ also satisfies $Z\notin\cU_{t_*}$.
Moreover, $d(X,\cU_{t(X)})=d^-(X,\cU)$.

\begin{lemma}[quasi-triangle inequality]\label{lem:quasi}
In the model, $\|Z-X\|\le2(\|Y-X\|+\|Z-Y\|)$ for all $X,Y,Z$. In $M$ there is
$r_\triangle>0$, depending only on $(M,\eta,g)$ and $C_\theta$, such that
$\|Z-X\|_X\le3(\|Y-X\|_X+\|Z-Y\|_Y)$ whenever the right-hand side is at most
$r_\triangle$.
\end{lemma}
\begin{proof}
In the model, write $Y=\ell_X(w,z,s)$ and $Z=\ell_Y(w',z',s')$. By the group law
$\ell_{(a,b,c)}(x',y',z')=(x'+a,y'+b,z'+c+\frac12(b\cdot x'-a\cdot y'))$, the model
coordinates of $Z$ relative to $X$ are $(w+w',\ z+z'+\frac12(y\cdot a'-x\cdot b'),\ s+s')$,
with $w=(x,y)$ and $w'=(a',b')$. Hence $|w+w'|\le\|Y-X\|+\|Z-Y\|$,
$|s+s'|^{1/2}\le\|Y-X\|+\|Z-Y\|$, and
\begin{multline*}
\bb^{-1}|z+z'+\tfrac12(y\cdot a'-x\cdot b')|^{1/2}\le\bb^{-1}\big(|z|^{1/2}+|z'|^{1/2}\big)+\bb^{-1}\big(\tfrac12|w||w'|\big)^{1/2}\\
\le(1+\bb^{-1})\big(\|Y-X\|+\|Z-Y\|\big),
\end{multline*}
using $\sqrt{ab}\le\frac12(a+b)$ and $\bb\ge1$.

In the curved case, let $X=(q,t)$, $Y=(q',t')$, $Z=(q'',t'')$ and
$r:=\|Y-X\|_X+\|Z-Y\|_Y\le r_\triangle:=r_M$. Then $q'\in U_q$ and $q''\in U_{q'}$. If $r=0$,
then $\Upsilon_q(q')=0$ and $t'=t$, so $Y=X$ since $\Upsilon_q$ is injective with
$\Upsilon_q(q)=0$. Likewise $Z=Y$, and there is nothing to prove. Now let $r>0$ and put
$a_1:=\|Y-X\|_X/r$, $a_2:=\|Z-Y\|_Y/r$, $\hat y:=\delta_{1/r}\Upsilon_q(q')$ and
$P:=\delta_{1/r}\Upsilon_{q'}(q'')$. Then $|\hat y_w|\le1$, $|\hat y_z|\le\bb^2$, $|P_w|\le1$,
$|P_z|\le\bb^2$, so $\hat y\in\overline{\cB_1}\subset\cB_8$. By Corollary~\ref{cor:rM}(b) with
$\lambda=1$ (whose proof uses only Lemma~\ref{lem:rM}(b) and Definition~\ref{def:rM}, not the
present lemma), the map $S:=\delta_{1/r}\Upsilon_{q'}\Upsilon_q^{-1}\delta_r$ is defined on $\cB_8$, and
$S=A+E$ with $A:=\tilde U\ell_{\hat y}^{-1}$ and both components of $E$ bounded by $\delta_*$ on
$\cB_8$. The map $A$ is affine: $\ell_{\hat y}^{-1}(w,z)=(w-\hat y_w,\ z-\hat y_z+\mathsf c\cdot w)$
with the vector $\mathsf c:=\frac12(-\hat y_y,\hat y_x)\in\R^{2n}$, where $\hat y_w=(\hat y_x,\hat y_y)$. Hence
$|\mathsf c|=\frac12|\hat y_w|\le\frac12$, and $L_A^{-1}(v_w,v_z)=(U^{-1}v_w,\ v_z-\mathsf c\cdot U^{-1}v_w)$
for the linear part $L_A$ of the affine map. The point $P_0:=A^{-1}(P)=\ell_{\hat y}\tilde U^{-1}P$ satisfies
$|P_{0,w}|\le2$ and $|P_{0,z}|\le2\bb^2+\frac12$. On the compact convex set
$\mathsf K:=\{P':\ |P'_w-P_{0,w}|\le1,\ |P'_z-P_{0,z}|\le2\}\subset\{|w|\le3,\ |z|\le2\bb^2+\frac52\}\subset\cB_8$,
the continuous map $\mathsf G(P'):=P_0-L_A^{-1}E(P')$ satisfies $|\mathsf G(P')_w-P_{0,w}|\le\delta_*$ and
$|\mathsf G(P')_z-P_{0,z}|\le\frac32\delta_*$, so $\mathsf G(\mathsf K)\subset\mathsf K$. By Brouwer's theorem
there is $P_Z\in\mathsf K$ with $\mathsf G(P_Z)=P_Z$, i.e.\ $S(P_Z)=P$. As $S$ is defined at $P_Z$,
the point $\Upsilon_q^{-1}\delta_rP_Z$ lies in $U_q\cap U_{q'}$, and $\Upsilon_{q'}$ maps it to
$\Upsilon_{q'}(q'')$. By injectivity of $\Upsilon_{q'}$, this point is $q''$. Hence $q''\in U_q$ and
$\delta_{1/r}\Upsilon_q(q'')=P_Z$. Let $N(w,z):=\max\{|w|,\bb^{-1}|z|^{1/2}\}$. By the model
case (the coordinates of a composition) and $N(\tilde U^{-1}P)=N(P)$, we have
$N(P_0)\le(1+\bb^{-1})(N(\hat y)+N(P))\le(1+\bb^{-1})(a_1+a_2)\le\frac32$. Then
$N(P_Z)\le N(P_0)+\delta_*+\bb^{-1}(\frac32\delta_*)^{1/2}<\frac74$, using $\delta_*\le\frac1{16}$ and
$\bb\ge2$. The time part satisfies $|t''-t|^{1/2}\le|t'-t|^{1/2}+|t''-t'|^{1/2}\le r$. Rescaling
back gives $\|Z-X\|_X\le\frac74r\le3r$.
\end{proof}

\subsection{The normalised bounded geometry scale}\label{subsec:rM}

Fix the integer $k_0:=8$ and the constant
\[
\delta_*:=\min\big\{\delta_0(n),\ \tfrac1{16}\big\},
\]
where $\delta_0(n)\in(0,1)$ is the dimensional constant of Section~\ref{sec:graphs} for which
Lemmas~\ref{lem:grapheq} and~\ref{lem:heightangle} hold. Thus $\delta_*$ depends only on
$n$, not on $M$, the flow, the centre or the scale.

For $r\in(0,1]$ and $p\in M$, let $\varsigma_{1/r}$ be the structure
$(\delta_{1/r})_*(\eta_\kappa,\xi_\kappa,\Phi,g_\kappa)$, $\kappa=r^{-2}$, written in the
chart $\delta_{1/r}\Upsilon_p$ (Fact~\ref{fact:charts}(iii)). This structure is defined on
$\{|w|<r_0/r,\ |z|<z_0/r^2\}$ and has the form $(\frac12dz+\beta^{(1/r)}_p,\ 2\partial_z,\
\Phi^{(1/r)}_p,\ \pi^*h^{(1/r)}_p+\eta\otimes\eta)$ with
\begin{equation}\label{eq:rescaled-coeff}
h^{(1/r)}_p(w)=h_p(rw),\qquad \beta^{(1/r)}_p(w)=r^{-1}\beta_p(rw),\qquad
\Phi^{(1/r)}_p\partial_{w_a}=J_0\partial_{w_a}-2\beta^{(1/r)}_p(J_0\partial_{w_a})\partial_z ,
\end{equation}
the last identity because $\Phi$ maps the horizontal lift $\partial_{w_a}-2\beta(\partial_{w_a})\partial_z$
to the horizontal lift of $J_0\partial_{w_a}$ and kills $\partial_z$.

For $R\ge1$ put
$\cB_R:=\{|w|<R,\ |z|<\bb^2R^2\}$. The $C^k$ norm on $\cB_R$ of a tensor field or of a map,
written in the coordinates $(w,z)=(x,y,z)$ of $\R^{2n+1}$, is
\[
\|\cdot\|_{C^k(\cB_R)}:=\sum_{j=0}^k\,\sup_{\cB_R}|\partial^j\,\cdot\,| ,
\]
where $|\partial^j\,\cdot\,|$ is the Euclidean (Hilbert--Schmidt) norm of the full array of
the components of the $j$-th coordinate derivatives: the square root of the sum of the
squares of all partial derivatives of order $j$, taken along all ordered $j$-tuples of
the $2n+1$ coordinate directions, of all components in these coordinates. The $C^k$ norm
on any other open subset of the coordinate space is defined in the same way. The norm $\|\varsigma-\varsigma_{\rm Heis}\|_{C^k(\cB_R)}$
is the sum of the $C^k(\cB_R)$ norms of $h-g_0$, $\beta-\overline\beta$ and
$\Phi-\overline\Phi$, written in the coordinates $(w,z)$. By
\eqref{eq:rescaled-coeff}, this norm is bounded by a constant times
$\|h-g_0\|_{C^k}+\|\beta-\overline\beta\|_{C^k}$, and the coefficients of $g$ are
polynomial in those of $h$ and $\beta$.

These norms are invariant under the rotations $\tilde U(w,z)=(Uw,z)$, $U\in U(n)$, which
are linear isometries of the Euclidean coordinate space $\R^{2n+1}$ and map every $\cB_R$
onto itself. Indeed, let $V,W\in U(n)$. Suppose that a tensor field on $\cB_R$ is written in
the coordinates rotated by $\tilde V$, that is, replaced by its push-forward under
$\tilde V$, or that a map $E$ on $\cB_R$ is replaced by $\tilde W\circ E\circ\tilde V^{-1}$.
Then the array of the $j$-th derivatives undergoes an orthogonal transformation in each of
its slots (those of the components and those of the derivatives) and is evaluated at the
rotated point. This changes neither its Euclidean norm nor the supremum over $\cB_R$. Moreover,
the push-forward $(\tilde V)_*\varsigma_{\rm Heis}$ of the model structure is
$\varsigma_{\rm Heis}$ itself, because $V$ preserves $g_0$,
$J_0$ and $\overline\beta$. Hence, for every structure $\varsigma$ of the above form and
every map $E$ on $\cB_R$,
\[
\|(\tilde V)_*\varsigma-\varsigma_{\rm Heis}\|_{C^k(\cB_R)}=\|\varsigma-\varsigma_{\rm Heis}\|_{C^k(\cB_R)},\qquad
\|\tilde W\circ E\circ\tilde V^{-1}\|_{C^k(\cB_R)}=\|E\|_{C^k(\cB_R)} .
\]
Consequently, the bounds of Lemma~\ref{lem:rM} and Corollary~\ref{cor:rM} below hold, with
the same constants, when the charts are composed with such rotations. In a chart rotated
by $\tilde V$, the structure is $(\tilde V)_*\varsigma_{1/r}$, for which we keep the
notation $\varsigma_{1/r}$. If the charts at $p$ and at $q$ are rotated by $\tilde V$
and $\tilde W$, then a transition of the form $\tilde U\ell_{\hat y}^{-1}+E$ as in
\eqref{eq:transition}, with $\hat y$ the image of $q$ in the chart at $p$, becomes
\[
\tilde W\big(\tilde U\ell_{\hat y}^{-1}+E\big)\tilde V^{-1}=\tilde W\tilde U\tilde V^{-1}\,\ell_{\tilde V\hat y}^{-1}+\tilde WE\tilde V^{-1},
\]
where $\tilde V\hat y$ is the image of $q$ in the rotated chart at $p$. Here we used
$\tilde V\ell_{\hat y}\tilde V^{-1}=\ell_{\tilde V\hat y}$, which holds because the
bilinear form $b\cdot x-a\cdot y$ in the definition of the left translations is invariant
under $U(n)$.

The $C^k$ norms of structures and of transitions between charts are always taken in the
above sense. The constants $C_k$ and $C_k(R)$ of Fact~\ref{fact:charts}, the constants
$C_\sharp$ and $C^{\rm tr}$ of Lemma~\ref{lem:rM} and the constants $\delta_0$ and
$\delta_*$ are understood with respect to these norms. For Fact~\ref{fact:charts}, this is
a normalisation only: measuring the arrays by the maximum or by the sum of the absolute
values of their entries, or taking the maximum over $j$ instead of the sum, changes $C_k$
and $C_k(R)$ at most by factors depending only on $n$ and $k$.

\begin{lemma}[finite-order closeness at small scales]\label{lem:rM}
\leavevmode
\begin{enumerate}[(a)]
\item There is a constant $C_\sharp=C_\sharp(M,k_0)$, depending only on the chart
constants of Fact~\ref{fact:charts}, such that for every $p\in M$, every $0<r\le1$ and
every $R\ge1$ with $rR\le r_0/2$ and $\bb^2r^2R^2\le z_0/2$,
\begin{equation}\label{eq:closeness}
\|\varsigma_{1/r}-\varsigma_{\rm Heis}\|_{C^{k_0}(\cB_R)}\le C_\sharp\,r^2(1+R)^3 .
\end{equation}
\item (Transitions.) There are $r_*=r_*(M,C_\theta)\in(0,\min\{1,r_0/64\}]$ and $C^{\rm tr}=C^{\rm tr}(M,k_0,C_\theta)$
such that for every $p\in M$, every $0<r\le r_*$ and every $q\in U_p$ with
$\hat y:=\delta_{1/r}\Upsilon_p(q)\in\cB_8$ there is $U\in U(n)$ with
\begin{equation}\label{eq:transition}
\big\|\delta_{1/r}\Upsilon_q\Upsilon_p^{-1}\delta_r-\tilde U\circ\ell_{\hat y}^{-1}\big\|_{C^{k_0}(\cB_8)}\le C^{\rm tr}\,r^2,
\qquad \tilde U(w,z)=(Uw,z).
\end{equation}
\end{enumerate}
\end{lemma}
\begin{proof}
(b) Put $R_\bb:=8(1+\bb)$. For $\hat y\in\cB_8$ we have $|\hat y_w|+|\hat y_z|^{1/2}<8+8\bb=R_\bb$,
and $\cB_8=\{|w|<8,\ |z|<64\bb^2\}\subset\{|w|<R_\bb,\ |z|<R_\bb^2\}$, since $64\bb^2<64(1+\bb)^2$.
So (b) is Fact~\ref{fact:charts}(iii) with $\lambda=1/r$, $R=R_\bb$ and $k=k_0$: it holds for
$\lambda\ge\max\{1,\lambda_*(R_\bb)\}$, i.e.\ for $r\le r_*:=\min\{1,1/\lambda_*(R_\bb),r_0/64\}$,
with the constant $C^{\rm tr}:=C_{k_0}(R_\bb)$ of that statement. In particular, the transition
is defined on $\cB_8$. Through $\bb^2=8C_\theta+4$ both constants depend on $C_\theta$ besides
$M$ and $k_0$; so does $r_M$, through $(z_0/(128\bb^2))^{1/2}$.

(a) By Fact~\ref{fact:charts}(i), $h_p-g_0$ vanishes to second order at $w=0$, with all
derivatives up to order $k_0$ bounded by $C_{k_0}$. Hence, by Taylor's theorem,
$|\partial^j(h_p-g_0)(w)|\le C|w|^{2-j}$ for $j\le2$ and $\le C$ for $2\le j\le k_0$.
Likewise, $|\partial^j(\beta_p-\overline\beta)(w)|\le C|w|^{3-j}$ for $j\le3$ and $\le C$ for
$3\le j\le k_0$. By \eqref{eq:rescaled-coeff}, $\partial^j_w[h^{(1/r)}_p-g_0](w)=r^j(\partial^jh_p)(rw)$
for $j\ge1$, and $\partial^j_w[\beta^{(1/r)}_p-\overline\beta](w)=r^{j-1}(\partial^j(\beta_p-\overline\beta))(rw)$,
where we used that $\overline\beta$ is homogeneous of degree one, $r^{-1}\overline\beta(rw)=\overline\beta(w)$.
Hence on $\cB_R$, for $r\le1$:
\[
|\partial^j_w(h^{(1/r)}_p-g_0)|\le Cr^2(1+R)^2,\qquad
|\partial^j_w(\beta^{(1/r)}_p-\overline\beta)|\le Cr^2(1+R)^3\qquad(0\le j\le k_0),
\]
using $r^j\le r^2$ for $j\ge2$ and $r^{j-1}\le r^2$ for $j\ge3$. The conditions on $rR$
and $\bb^2r^2R^2$ ensure that $\cB_R$ lies in the domain of $\varsigma_{1/r}$ with margin. The
bound for $\Phi^{(1/r)}_p-\overline\Phi$ follows from \eqref{eq:rescaled-coeff}, since
$(\Phi^{(1/r)}_p-\overline\Phi)\partial_{w_a}=-2(\beta^{(1/r)}_p-\overline\beta)(J_0\partial_{w_a})\partial_z$
and $(\Phi^{(1/r)}_p-\overline\Phi)\partial_z=0$. The metric $g$ is not part of the norm.
\end{proof}

\begin{definition}\label{def:rM}
The \emph{normalised bounded geometry scale} of $M$ is
\[
r_M:=\min\Big\{\frac{r_0}{64},\ r_*,\ \Big(\frac{z_0}{128\,\bb^2}\Big)^{1/2},\ 1,\
\Big(\frac{\delta_*}{729\,C_\sharp}\Big)^{1/2},\ \Big(\frac{\delta_*}{C^{\rm tr}}\Big)^{1/2}\Big\},
\]
and
\begin{equation}\label{eq:r1}
r_1:=\min\big\{r_M,\ |K+2|^{-1/2}\big\}\qquad(\text{the second term omitted if }K=-2).
\end{equation}
Here $\delta_*=\min\{\delta_0(n),\frac1{16}\}$ is the constant fixed at the beginning of
Section~\ref{subsec:rM}; it depends only on $n$.
\end{definition}

\begin{corollary}\label{cor:rM}
Let $r\le r_M$ and $\lambda\ge1$.
\begin{enumerate}[(a)]
\item For $1\le R\le8\lambda$,
$\|\varsigma_{\lambda/r}-\varsigma_{\rm Heis}\|_{C^{k_0}(\cB_R)}\le\delta_*\lambda^{-2}(1+R)^3/729$;
in particular $\le\delta_*\lambda^{-2}$ for $R\le8$.
\item For every $p$ and every $q$ with $\hat y=\delta_{\lambda/r}\Upsilon_p(q)\in\cB_8$,
the transition $\delta_{\lambda/r}\Upsilon_q\Upsilon_p^{-1}\delta_{r/\lambda}$ is
$\delta_*\lambda^{-2}$-close in $C^{k_0}(\cB_8)$ to $\tilde U\ell_{\hat y}^{-1}$ for some
$U\in U(n)$.
\item (Fixed larger domains.) Let $R\ge8$ and $A:=R/8$. For $\lambda\ge A$ and every
$q$ with $\hat y=\delta_{\lambda/r}\Upsilon_p(q)\in\cB_R$, the transition
$\delta_{\lambda/r}\Upsilon_q\Upsilon_p^{-1}\delta_{r/\lambda}$ is $\delta_*A^4\lambda^{-2}$-close
in $C^{k_0}(\cB_R)$ to $\tilde U\ell_{\hat y}^{-1}$ for some $U\in U(n)$; the constant
$A^4=(R/8)^4$ depends on $R$ but not on $M$.
\end{enumerate}
Thus every chart $\Lambda_{X,r}$ with $r\le r_1$ is $\delta_*$-close to the model on
$\cB_8$, its transitions to the charts of nearby centres are $\delta_*$-close to
rotations composed with left translations, and for fixed $R$ both errors are
$O(\lambda^{-2})$ with constants independent of $M$; these geometric and transition
bounds use only $r\le r_M$. Moreover, for $r\le r_1$ (by \eqref{eq:r1}, not only $r\le r_M$),
$|K+2|(r/\lambda)^2\le\lambda^{-2}$. Lemma~\ref{lem:quasi} holds with the constant $3$ at
all scales $\le r_M$; its proof uses only (b), which rests on Lemma~\ref{lem:rM}(b).
\end{corollary}
\begin{proof}
We have $\varsigma_{\lambda/r}=\varsigma_{1/(r/\lambda)}$ with $r/\lambda\le r_M$. For $R\le8\lambda$,
we have $(r/\lambda)R\le8r\le r_0/8$ and $\bb^2(r/\lambda)^2R^2\le64\bb^2r^2\le z_0/2$,
so Lemma~\ref{lem:rM}(a) applies with $r/\lambda$ and gives $C_\sharp(r/\lambda)^2(1+R)^3\le
\delta_*\lambda^{-2}(1+R)^3/729$. For (b), Lemma~\ref{lem:rM}(b) at the scale
$r/\lambda\le r_*$ gives the error $C^{\rm tr}(r/\lambda)^2\le\delta_*\lambda^{-2}$.

For (c), put $\lambda':=\lambda/A\ge1$ and write the transition at the scale $r/\lambda$
as a conjugate of the transition at the scale $r/\lambda'$:
\[
S_\lambda:=\delta_{\lambda/r}\Upsilon_q\Upsilon_p^{-1}\delta_{r/\lambda}=\delta_A\circ S_{\lambda'}\circ\delta_{1/A},\qquad
S_{\lambda'}=\delta_{\lambda'/r}\Upsilon_q\Upsilon_p^{-1}\delta_{r/\lambda'} .
\]
The centre $\hat y':=\delta_{1/A}\hat y$ lies in $\cB_{R/A}=\cB_8$, and $\delta_{1/A}(\cB_R)=\cB_8$,
so (b) at the scale $r/\lambda'$ gives $S_{\lambda'}=\tilde U\ell_{\hat y'}^{-1}+E'$ with
$\|E'\|_{C^{k_0}(\cB_8)}\le\delta_*\lambda'^{-2}=\delta_*A^2\lambda^{-2}$. Dilations
commute with $\tilde U$ and conjugate left translations: $\delta_A\ell_{\hat y'}^{-1}\delta_{1/A}=\ell_{\delta_A\hat y'}^{-1}=\ell_{\hat y}^{-1}$.
Hence $S_\lambda=\tilde U\ell_{\hat y}^{-1}+\delta_AE'\delta_{1/A}$ on $\cB_R$. The conjugated
error has components $(AE'_w,A^2E'_z)\circ\delta_{1/A}$. A partial derivative of order $j$ of
$E'_\bullet\circ\delta_{1/A}$, $\bullet\in\{w,z\}$, is the same partial derivative of $E'_\bullet$, taken at the
image point and multiplied by $A^{-m}$ with $m\ge j$ (each $w$-derivative
brings a factor $A^{-1}$ and each $z$-derivative a factor $A^{-2}$). As $A\ge1$, every entry
of the array of the $j$-th derivatives of $\delta_AE'\delta_{1/A}$ is thus, in absolute
value, at most $A^{2-j}\le A^2$ times the corresponding entry for $E'$. Hence
$|\partial^j(\delta_AE'\delta_{1/A})|\le A^2\,|\partial^jE'|\circ\delta_{1/A}$ on $\cB_R$ for $0\le j\le k_0$, and
$\|\delta_AE'\delta_{1/A}\|_{C^{k_0}(\cB_R)}\le A^2\|E'\|_{C^{k_0}(\cB_8)}\le\delta_*A^4\lambda^{-2}$.
The remaining assertions follow from $|K+2|r_1^2\le1$ and from the proof of Lemma~\ref{lem:quasi}.
\end{proof}

\subsection{The contact regularity scale}\label{subsec:scale}

\emph{Coordinates and norms.} Throughout, a Lagrangian subspace $P\subset\C^n$ is
brought to $\{y=0\}$ by an element of $U(n)$. Such an element preserves $g_0$, $J_0$, $\omega_0$ and
$\overline\beta$ and leaves the bounds of Lemma~\ref{lem:rM} and Corollary~\ref{cor:rM}
unchanged (Section~\ref{subsec:rM}). We always use the Euclidean coordinates $w=(x,y)$ of the chart, so
that $P=\R^n=\{y=0\}$, $J_0P=\{x=0\}$, and a function $u$ on a subset $\Omega$ of
$\R^n\times\R$ has the graph $\Gamma_u=\{(x,u(x,s),s)\}$. In these coordinates, the metric
$g_0$ is $\frac14|dw|^2$, so that the induced metric of a graph is
$G_{ij}=\frac14(\delta_{ij}+\partial_iu\cdot\partial_ju)$ in the model, and the area of the
graph of $u\equiv0$ over a set $A$ is $2^{-n}|A|$. For a function $u$ on a set
$\Omega\subset\R^n\times\R$ we use the parabolic H\"older norm of \cite[\S8]{wh}
\emph{without its zero-order term},
\[
\|u\|_{2,\alpha;\Omega}:=\|Du\|_{0;\Omega}+\|D^2u\|_{0;\Omega}+\|\partial_su\|_{0;\Omega}+[Du]_{\alpha;\Omega}+[D^2u]_{\alpha;\Omega}+[\partial_su]_{\alpha;\Omega},
\]
where $[\cdot]_\alpha$ is the parabolic H\"older seminorm with exponent $\alpha$ (so that
$[Du]_\alpha$ is the $\alpha$-H\"older seminorm of $Du$ in space-time, not the
$(1+\alpha)/2$-seminorm in time). This norm is bounded by the parabolic $C^{2,\alpha}$
norm of \cite{wh}, and it controls every quantity used below. On the base domains
$\Omega_{T'}$, the zero-order term is controlled by
$|u(x,s)|\le|s|\,\|\partial_su\|_0+|x|\,\|Du\|_0$ for $u(0,0)=0$, along the path from
$(0,0)$ to $(0,s)$ to $(x,s)$.

The reason for this choice is scaling: the rescaled function
$u_\mu(x,s):=\mu^{-1}u(\mu x,\mu^2s)$ satisfies $\|u_\mu\|_{2,\alpha;\Omega}\le\|u\|_{2,\alpha;\mu\Omega}$
for $0<\mu\le1$, where $\mu\Omega=\{(\mu x,\mu^2s)\}$, since every term is multiplied by
a power $\mu^j$ with $j\ge0$ ($j=0$ for $\|Du\|_0$, $j=\alpha$ for $[Du]_\alpha$, $j=1$ for
$\|D^2u\|_0$ and $\|\partial_su\|_0$, $j=1+\alpha$ for the top seminorms), whereas a
zero-order term would be multiplied by $\mu^{-1}$.

The \emph{unit base domain} with
upper end $T'\ge0$ is
\[
\Omega_{T'}:=\{(x,s)\in\R^n\times\R:\ |x|<1,\ |s|<1,\ s\le T'\},
\]
and the \emph{unit box} in the frame of $P$ is
\[
\bB:=\{(x,y,z)\in\Hn:\ |x|<1,\ |y|<2,\ |z|<\bb^2\}\times(-1,1),
\]
which satisfies $\{|w|<1,|z|<\bb^2\}\times(-1,1)\subset\bB\subset\{|w|<\sqrt5,|z|<\bb^2\}\times(-1,1)$.
In particular, for a centre $X$, $r>0$ and $U\in U(n)$, the closure $\overline\bB$ lies in
the $\tilde U$-invariant set $\{|w|\le\sqrt5,\ |z|\le\bb^2\}\times[-1,1]$, which lies in
$\{|w|<3,\ |z|<9\bb^2\}\times(-9,9)$. Hence $(\tilde U\Lambda_{X,r})^{-1}(\overline\bB)\subset\cC_{3r}(X)$.
If moreover $0<r\le r_1$ (any $r>0$ in the model, where the charts are global), the box
$\{|w|<3r,\ |z|<9\bb^2r^2\}$ lies in the range of the chart of the centre, since
$3r\le3r_0/64<r_0$ and $9\bb^2r^2\le9z_0/128<z_0$ by Definition~\ref{def:rM}. In that case
$\tilde U\Lambda_{X,r}(\cC_{3r}(X))=\{|w|<3,\ |z|<9\bb^2\}\times(-9,9)$. Thus every point of
$\overline\bB$ has a preimage in $\cC_{3r}(X)$, and $\overline\bB\subset\tilde U\Lambda_{X,r}(\cU)$
whenever $\cC_{3r}(X)\subset\cU$.

\begin{definition}[Legendrian graphs]\label{def:leggraph}
Let $\varsigma=(\frac12dz+\beta_\varsigma,2\partial_z,\Phi_\varsigma,g_\varsigma)$ be a
structure on $\cB_8$ as in Section~\ref{subsec:rM}, with transverse metric $h_\varsigma$,
and let $u\colon\Omega\to\R^n$ be $C^1$, $\Omega\subset\R^n\times\R$. The graph
$\Gamma_u$ is \emph{Lagrangian} if $\Gamma_{u(\cdot,s)}^*d\beta_\varsigma=0$ for every
$s$. Then on every simply connected $\Omega_s=\{x:(x,s)\in\Omega\}$ there is a height
function $\zeta(\cdot,s)$, unique up to an additive constant, with
\begin{equation}\label{eq:height}
d_x\zeta=-2\,\Gamma_{u(\cdot,s)}^*\beta_\varsigma ,
\end{equation}
and the \emph{Legendrian graph} $\hat\Gamma_{u,\zeta}=\{(x,u(x,s),\zeta(x,s),s)\}$ is
Legendrian for $\eta_\varsigma$ at every time. A Legendrian graph carries exactly one
height above each point of $\Gamma_u$; the additive constants $\zeta(0,s)$ are additional
data.
\end{definition}

\begin{definition}[admissible radii and the regularity scale]\label{def:K}
Let $\cM$ be the trace of a smooth local Legendrian flow in $\cU\subset M\times\R$ with
upper end $T_*$ (Definition~\ref{def:localflow}, or Definition~\ref{def:traceflow}), let
$X=(q,t)\in\cM$ and $0<r\le r_1$, and let $\Lambda_{X,r}\colon(q',t')\mapsto(\delta_{1/r}\Upsilon_q(q'),(t'-t)/r^2)$
be the rescaling centred at $X$. The radius $r$ is \emph{admissible} for $(\cM,X)$ if
there are $U\in U(n)$ and a function $u\colon\Omega_{T'}\to\R^n$, $T':=(T_*-t)/r^2\ge0$,
with
\[
\|u\|_{2,\alpha;\Omega_{T'}}\le1,\qquad u(0,0)=0,\qquad Du(0,0)=0,
\]
such that, in the coordinates $\tilde U\circ\Lambda_{X,r}$, $\tilde U(w,z)=(Uw,z)$,
the graph $\Gamma_u$ is Lagrangian for the structure $\varsigma_{1/r}$, and
\begin{enumerate}[(i)]
\item (containment) $\tilde U\Lambda_{X,r}(\cM)\cap\bB\subset\hat\Gamma_{u,\zeta}$, where
$\zeta\in C^1(\Omega_{T'})$ (one-sided at $s=T'$) is a height function in the sense of
\eqref{eq:height} with $\zeta(0,0)=0$;
\item (injectivity) the map $\hat F$ is injective on the preimage of
$(\tilde U\Lambda_{X,r})^{-1}(\bB)$; for a flow in trace form this condition is empty,
the trace having multiplicity one by definition.
\end{enumerate}
The regularity required of $\zeta$ concerns the data $\zeta(0,s)$, which are not
determined by \eqref{eq:height}; for the witnesses produced from a smooth flow
(Lemma~\ref{lem:finite}, Step~1; cf.\ Lemma~\ref{lem:heightangle}(b)) the height is smooth.
The \emph{regularity scale} is $\Kreg(\cM,X):=1/\sup\{r\le r_1:\ r\text{ admissible}\}\in[1/r_1,\infty]$,
with $\Kreg=\infty$ if no radius is admissible. For $\cU$ as above
\[
\Kreg_{\cU}(\cM):=\sup_{X\in\cM}d_1(X,\cU)\,\Kreg(\cM,X),\qquad
\Kreg^-(\cM,X):=\Kreg\big(\cM\cap\{t'\le t\},X\big),
\]
the latter being the scale of the truncation at $X$. In the model we use the same
definition with $r_1=\infty$, with the charts $\Upsilon_q:=\ell_q^{-1}$, defined on all of
$\Hn$ (in these charts $h_q=g_0$ and $\beta_q=\overline\beta$, so that the identities and
expansions of Fact~\ref{fact:charts}(i) hold with vanishing errors; the coefficients of
$\overline\beta$, of $\overline g$ and of its inverse are polynomials in $w$ which do not
depend on $q$, for instance $\overline g(\partial_{x_j},\partial_{x_j})=\frac14+\frac1{16}y_j^2$,
so they are bounded, together with their derivatives, on bounded subsets of the charts,
but not on all of $\Hn$; see Section~\ref{subsec:gap}), and with
$\Kreg(\cM,X)\in[0,\infty]$.
\end{definition}

\begin{lemma}[admissible radii]\label{lem:admissible}
\leavevmode
\begin{enumerate}[(a)]
\item If $r$ is admissible for $(\cM,X)$ and $0<r'<r$, then $r'$ is admissible, with the
witness $u'(x,s)=\mu^{-1}u(\mu x,\mu^2s)$, $\mu=r'/r$, and the same $U$. Consequently
every $r<1/\Kreg(\cM,X)$ is admissible.
\item If $\cM'\subset\cM$ is the trace of a restriction or of a truncation, then every
radius admissible for $(\cM,X)$ is admissible for $(\cM',X)$; hence
$\Kreg(\cM',X)\le\Kreg(\cM,X)$ and $\Kreg^-\le\Kreg$.
\item In the model, $\Kreg(D_\lambda\cM,D_\lambda X)=\lambda^{-1}\Kreg(\cM,X)$ for
$D_\lambda(q,t)=(\delta_\lambda q,\lambda^2t)$. In $M$ the same identity holds for the
$\cD$-homothetic structure with $\kappa=\lambda^2$, written in the rescaled chart, with
$r_1$ replaced by $\lambda r_1$.
\end{enumerate}
\end{lemma}
\begin{proof}
(a) The chart $\Lambda_{X,r'}=D_{1/\mu}\circ\Lambda_{X,r}$ with $D_{1/\mu}(w,z,s)=(w/\mu,z/\mu^2,s/\mu^2)$,
and $T'/\mu^2$ is the upper end belonging to the radius $r'$. The set
$\mu\Omega_{T'/\mu^2}=\{(\mu x,\mu^2s):(x,s)\in\Omega_{T'/\mu^2}\}=\{|x|<\mu,\ |s|<\mu^2,\ s\le T'\}$
is contained in $\Omega_{T'}$, as $\mu<1$. So $u$ is defined there, and $D_{1/\mu}$ maps
the Legendrian graph of $(u,\zeta)$ over $\mu\Omega_{T'/\mu^2}$
to the Legendrian graph of $(u',\zeta')$, $\zeta'(x,s)=\mu^{-2}\zeta(\mu x,\mu^2s)$, over
$\Omega_{T'/\mu^2}$ (a dilation preserves the Legendrian condition, and the structure
$\varsigma_{1/r'}$ is the dilation of $\varsigma_{1/r}$). The function $\zeta'$ is $C^1$ with
$\zeta'(0,0)=0$. The scaling property of the norm with $\Omega=\Omega_{T'/\mu^2}$ and its
monotonicity under restriction of the domain give
$\|u'\|_{2,\alpha;\Omega_{T'/\mu^2}}\le\|u\|_{2,\alpha;\mu\Omega_{T'/\mu^2}}\le\|u\|_{2,\alpha;\Omega_{T'}}\le1$,
and $u'(0,0)=0$, $Du'(0,0)=0$. The box $\bB$ in the chart $\Lambda_{X,r'}$ corresponds to the smaller set
$D_\mu\bB\subset\bB$ in the chart $\Lambda_{X,r}$, so (i) and (ii) are inherited.

(b) The witness of $(\cM,X)$ serves for
the subset. For a truncation at $t_*$ (with $t\le t_*\le T_*$), it is restricted to
$\Omega_{(t_*-t)/r^2}\subset\Omega_{T'}$.

(c) This follows from the identity $\Lambda_{D_\lambda X,\lambda r}\circ D_\lambda=\Lambda_{X,r}$.
\end{proof}

\begin{remark}\label{rem:K}\leavevmode
\begin{enumerate}[(1)]
\item The definition follows \cite[\S2.5]{wh}: the trace is required to be
\emph{contained} in a graph, not to fill it. Lemma~\ref{lem:filling} shows that the full
graph over $\Omega_{T'}$ is realised by the trace when $\cC_{3r}(X)\subset\cU$. The
function $u$ is only constrained on $\Omega_{T'}$, the part of the base domain before the
upper end, and no extension across $s=T'$ is required.
\item The height $\zeta$ and the Legendrian angle $\theta$ are not part of the norm. By
\eqref{eq:height}, $u$ and the structure determine $D_x\zeta$ and hence control
$\zeta(\cdot,s)-\zeta(0,s)$ (Lemma~\ref{lem:heightangle}(a)). The values $\zeta(0,s)$ are
additional data. Where the Legendrian graph is the trace of the flow (on all of
$\Omega_{T'}$ when $\cC_{3r}(X)\subset\cU$, by Lemma~\ref{lem:filling}),
Lemma~\ref{lem:heightangle}(b),(c) express $\theta$ through $\partial_s\zeta$, $\partial_su$
and the structure, and bound $D_x\theta$. Lemma~\ref{lem:heightangle}(d), using
$|\theta|\le C_\theta$ and $\zeta(0,0)=0$, gives $|\zeta(x,s)|<\bb^2/2$ whenever the graph is
the trace of the flow near every point of the time segment from $(0,0)$ to $(0,s)$.
Condition (i) excludes several heights above one projected point inside $\bB$, and (ii)
excludes double points in $\bB$.
\item The supremum in the definition of $\Kreg$ need not be attained. By
Lemma~\ref{lem:admissible}(a), we always work with a radius $r<1/\Kreg(\cM,X)$, typically
$r=1/(2\Kreg)$.
\end{enumerate}
\end{remark}

\subsection{The density hypothesis and the main theorem}\label{subsec:statements}

For a smooth local flow in $\cU$, a point $X=(q,t)\in\cM$ and $0<\sigma\le r_1^2$, put
\begin{equation}\label{eq:density-local}
\Theta^\chi_X(\sigma):=\int_{\cN_{t-\sigma}}\psi_q\,\rho^\chi_{q,t}\,d\mu_{t-\sigma},
\end{equation}
where the integral is taken over the time-$(t-\sigma)$ slice of the parameter domain. This
integral counts every branch of the flow in $\cU$ with multiplicity, and $\theta$ enters through
$v$ as a function on $\cN_{t-\sigma}$. For a truncation and $t(X)\le t_*$, the slices
before $t(X)$ are unchanged, so $\Theta^\chi_X(\sigma;F^{\le t_*})=\Theta^\chi_X(\sigma;F)$.
For a restriction to $\cU'\subset\cU$, the density can only decrease, since the integrand is
non-negative.

\begin{definition}\label{def:Heps}
A smooth local flow in $\cU$ satisfies $(\mathrm H_\eps)$ if
\[
\Theta^\chi_X(\sigma)\le1+\eps\qquad\text{for all }X\in\cM\text{ and all }0<\sigma<d_1(X,\cU)^2 ,
\]
and $(\mathrm H^-_\eps)$ if the same holds for all $0<\sigma<d_1^-(X,\cU)^2$. Since $d\le
d^-$, $(\mathrm H^-_\eps)$ implies $(\mathrm H_\eps)$.
\end{definition}

\begin{theorem}[local regularity; two-sided form]\label{thm:A}
Let $0<\alpha<1$. There are $\eps=\eps(n,a,\alpha,C_\theta)\in(0,\frac12]$ and
$C=C(n,a,\alpha,C_\theta)\in[3,\infty)$ with the following property (any smaller $\eps$ and
any larger $C$ then have it as well; the normalisations $\eps\le\frac12$ and $C\ge3$ are
used in Sections~\ref{sec:proof} and~\ref{sec:terminal}). Let $M$ be a compact
Sasakian manifold with the constant $r_1$ of \eqref{eq:r1}, let $\cU\subset M\times\R$ be
open and let $F$ be a smooth local Legendrian flow in $\cU$ (with any upper end) satisfying
$(\mathrm H_\eps)$. Then
\[
\Kreg_{\cU}(\cM)=\sup_{X\in\cM}d_1(X,\cU)\,\Kreg(\cM,X)\le C .
\]
\end{theorem}

Section~\ref{subsec:gap} fixes a common pair of constants $(\eps,C)$ for this
theorem and its Heisenberg model version. We use this choice throughout.

Theorem~\ref{thm:A} is proved in Section~\ref{sec:proof}, following the argument of
\cite[Theorem~3.1]{wh}; the proof rests on Sections~\ref{sec:graphs} and~\ref{sec:model} and
Appendix~\ref{app:schauder}. If the theorem fails, there are flows satisfying
$(\mathrm H_{\eps_i})$ with $\eps_i\to0$ and $\Kreg_{\cU_i}(\cM_i)\to\infty$.
After point selection and rescaling, the quantity $\Kreg$ equals one at the origin and
is uniformly bounded on expanding cylinders. Proposition~\ref{prop:compact} gives an
ancient model limit with bounded $\Kreg$.
By Lemma~\ref{lem:lsc} the density of this limit is at most one at every point and scale,
so the limit is a plane by Proposition~\ref{prop:rigidity}. The height enters at this
point: density one at every scale forces $z=4\tau\theta$ together with the self-shrinker
equation, and the bound on the regularity scale then forces the second fundamental form to
vanish. Finally, Lemmas~\ref{lem:onesheet} and~\ref{lem:schauder} show that the rescaled
flows are graphs with small norm over a box larger than the unit box, contradicting the
normalisation.

Theorem~\ref{thm:B} applies the local estimate to truncations of a closed flow.

\begin{corollary}[past form]\label{cor:past}
If a smooth local Legendrian flow in $\cU$ satisfies
$(\mathrm H^-_\eps)$, then for every $X\in\cM$
\[
d_1^-(X,\cU)\,\Kreg^-(\cM,X)\le C,\qquad\text{and}\qquad |B|(X)\le\frac{C'}{d_1^-(X,\cU)},
\]
with $C'=C'(n,a,\alpha,C_\theta)$.
\end{corollary}
\begin{proof}
Fix $X=(q,t)\in\cM$ and consider the truncation $F^{\le t}$, a smooth local flow in
$\cU_t=\cU\cup\{t'>t\}$ with trace $\cM^{\le t}=\cM\cap\{t'\le t\}$
(Remark~\ref{rem:localflow}(2)). For $Y\in\cM^{\le t}$ and $0<\sigma<d_1(Y,\cU_t)^2$ we
have $\sigma<d_1^-(Y,\cU)^2$ by \eqref{eq:dpast}, hence
$\Theta^\chi_Y(\sigma;F^{\le t})=\Theta^\chi_Y(\sigma;F)\le1+\eps$ by $(\mathrm H^-_\eps)$.
So $F^{\le t}$ satisfies $(\mathrm H_\eps)$ in $\cU_t$, and Theorem~\ref{thm:A} gives
$d_1(X,\cU_t)\Kreg(\cM^{\le t},X)\le C$. Since $d(X,\cU_t)=d^-(X,\cU)$ and
$\Kreg(\cM^{\le t},X)=\Kreg^-(\cM,X)$, this is the first inequality. The second inequality follows
from Lemma~\ref{lem:subord}(iii) below, applied at the admissible scale
$r=\min\{d_1^-(X,\cU)/3,\,1/(2\Kreg^-(\cM,X))\}$ to the truncated flow in $\cU_t$. The
trace of this flow coincides with that of $F$ at times $\le t$, and its second fundamental form at
$X$ is that of $F$.
\end{proof}

\begin{remark}
In White's theory, the two-sided distance controls the graph norm \cite[Theorem~3.1]{wh},
and the past distance enters through the truncation of \cite[Remark~3.3]{wh}, which yields
the curvature bound of \cite[Corollary~3.4]{wh}. Corollary~\ref{cor:past} follows the same
route: Theorem~\ref{thm:A} applied to the truncation bounds the truncated scale
$\Kreg^-(\cM,X)$ in terms of $d_1^-(X,\cU)$ (the first inequality), and the curvature bound
(the second inequality) follows from this bound by Lemma~\ref{lem:subord}(iii). The
hypothesis $(\mathrm H^-_\eps)$ only involves times before $t(X)$ and the scales
$\sigma<d_1^-(X,\cU)^2$, at which the time-$(t-\sigma)$ slice of the cylinder
$\cC_{d_1^-(X,\cU)}(X)$ lies in $\cU$. Indeed, a point $Y$ of that slice has $t(Y)=t-\sigma<t(X)$
and $\|Y-X\|_X<d_1^-(X,\cU)\le d^-(X,\cU)$, so $Y\in\cU$ by the definition of $d^-$ as an
infimum.
\end{remark}

\section{Legendrian graphs, realisation and compactness}\label{sec:graphs}

The regularity scale of Section~\ref{subsec:scale} is defined using Legendrian
graphs in contact normal charts. Within each chart the transverse projection
moves by mean curvature
\cite[Lemma~2.3]{ChangWuZhang2026ContactBlowup}, so its graph function satisfies a
quasilinear parabolic system (Lemma~\ref{lem:grapheq}). Lemma~\ref{lem:heightangle} relates
the spatial derivatives of the height to the graph function and the Legendrian angle to
the height velocity. With these identities we show that an
admissible radius yields a full graph realised by the flow, with interior estimates
(Lemmas~\ref{lem:filling} and~\ref{lem:subord}). We then prove the compactness of
sequences of local flows with locally bounded regularity scale
(Proposition~\ref{prop:compact}), and show that the graph bounds improve near a plane
(Lemmas~\ref{lem:onesheet} and~\ref{lem:schauder}). Throughout, the charts of nearby
centres differ, at the scales used, by left translations composed with unitary rotations,
up to small errors (Corollary~\ref{cor:rM}).

Throughout the section we work in a rescaled chart $\Lambda_{X,r}$, $r\le r_1$, with structure
$\varsigma=\varsigma_{1/r}$ or, after a further rescaling by $\lambda\ge1$,
$\varsigma=\varsigma_{\lambda/r}$, so that by Corollary~\ref{cor:rM}
\begin{equation}\label{eq:delta}
\|\varsigma-\varsigma_{\rm Heis}\|_{C^{k_0}(\cB_8)}\le\delta:=\delta_*\lambda^{-2}\le\delta_* ,
\end{equation}
and all tensors are written in the Euclidean coordinates $(w,z)=(x,y,z)$ of the chart,
with the frame of Section~\ref{subsec:scale}. Constants $C$ depend only on $n$, $\alpha$
and $C_\theta$ unless stated otherwise. The dimensional constant $\delta_0=\delta_0(n)\in(0,1)$ is
the one provided by Lemma~\ref{lem:grapheq}, decreased if necessary in
Lemma~\ref{lem:heightangle}(d). The constant $\delta_*=\min\{\delta_0,\frac1{16}\}$ was fixed in
Section~\ref{subsec:rM}, so that $\delta\le\delta_*\le\delta_0$. We write
$\kappa_\theta:=(K+2)r^2\lambda^{-2}$ for the zero-order constant in the rescaled angle
equation \eqref{eq:angleeq} below. It is the constant $K+2$ of \eqref{eq:F} transformed by the
$\cD$-homothety with parameter $(\lambda/r)^2$, and $|\kappa_\theta|\le\lambda^{-2}$ by \eqref{eq:r1}.
(The letter $\kappa$ without subscript denotes the $\cD$-homothety parameter itself, as in
Fact~\ref{fact:charts}(iii) and Lemma~\ref{lem:admissible}(c), and is not used for this
constant.)

\subsection{The graph equation and the angle equation}\label{subsec:grapheq}

Let $h=h_\varsigma$ be the transverse metric of the chart, with Christoffel symbols
$\Gamma^A_{BC}$, $A,B,C=1,\dots,2n$, and let $u\colon\Omega\to\R^n$ be $C^{2,\alpha}$ with
$|x|<2$, $|u|<2$, $|Du|\le1$ on $\Omega$, so that the graph lies in $\{|w|<4\}$. Put
\[
F(x,s)=(x,u(x,s))\in\R^{2n},\qquad Q^A_i=\partial_iF^A,\qquad
G_{ij}=h_{AB}Q^A_iQ^B_j ,
\]
so that $G$ is the induced metric of $\Gamma_{u(\cdot,s)}$ and $G_{ij}=\frac14(\delta_{ij}+\partial_iu\cdot\partial_ju)$
in the model.

\begin{lemma}[graph equation]\label{lem:grapheq}
There is $\delta_0=\delta_0(n)>0$ such that for $\delta\le\delta_0$ the following holds.
Suppose that the graphs $\Gamma_{u(\cdot,s)}$ are Lagrangian in $(\{|w|<8\},h)$. Then
the family of graphs moves by the mean curvature flow of $(\{|w|<8\},h)$, as a family of
sets, if and only if $u$ satisfies
\begin{equation}\label{eq:grapheq}
\partial_su^\gamma=a^{ij}(x,u,Du)\,D_{ij}u^\gamma+b^\gamma(x,u,Du),\qquad \gamma=1,\dots,n,
\end{equation}
with
\begin{equation}\label{eq:coeff-def}
a^{ij}=G^{ij},\qquad
b^\gamma=G^{ij}\big(\Gamma^{n+\gamma}_{AB}-u^\gamma_k\Gamma^k_{AB}\big)Q^A_iQ^B_j ,\qquad
e^{ij}:=G^{ij}-4\big((I+Du^TDu)^{-1}\big)_{ij},
\end{equation}
where $u^\gamma_k=\partial_ku^\gamma$ and $\Gamma$ is evaluated at $(x,u)$. The
coefficients are $C^{k_0-2}$ functions of $(x,u,p)=(x,u,Du)$ on
$\{|x|<2,\ |u|<2,\ |p|\le1\}$, and there
\begin{equation}\label{eq:coeff}
|e|+|\partial_{(x,u,p)}e|+|b|+|\partial_{(x,u,p)}b|\le C(n)\,\delta,\qquad
a^{ij}\xi_i\xi_j\ge(2-C(n)\delta)|\xi|^2\ge|\xi|^2 .
\end{equation}
\end{lemma}
\begin{proof}
\emph{The equation.} Let $\vec h:=G^{ij}\big(\partial_{ij}F+\Gamma(Q_i,Q_j)\big)$, whose
components are $\vec h^A=G^{ij}(\partial_{ij}F^A+\Gamma^A_{BC}Q^B_iQ^C_j)$. Its
$h$-orthogonal projection onto the normal space of the graph is the mean curvature vector
$H$, and $\vec h-H$ is tangent. Consider the linear map
$L\colon\R^{2n}\to\R^n$, $L(V)^\gamma:=V^{n+\gamma}-u^\gamma_kV^k$. It vanishes on the
tangent vectors $Q_k=(e_k,\partial_ku)$, $k=1,\dots,n$, which span the tangent space.
Since $L$ is onto, its kernel is exactly the tangent space. The graph point above $x$
moves with velocity $(0,\partial_su)$, and the family moves by mean curvature flow as sets
if and only if $(0,\partial_su)-H$ is tangent, i.e.\ if and only if $L((0,\partial_su))=L(H)=L(\vec h)$.
Now $L((0,\partial_su))=\partial_su$, and since $\partial_{ij}F^k=0$ for $k\le n$ and
$\partial_{ij}F^{n+\gamma}=\partial_{ij}u^\gamma$,
\[
L(\vec h)^\gamma=G^{ij}\big(\partial_{ij}u^\gamma+\Gamma^{n+\gamma}_{BC}Q^B_iQ^C_j\big)-u^\gamma_kG^{ij}\Gamma^k_{BC}Q^B_iQ^C_j ,
\]
which is the right-hand side of \eqref{eq:grapheq} with \eqref{eq:coeff-def}. In the
model ($\Gamma=0$, $h=g_0$) this is the equation of \cite[Theorem~8.4]{wh} with
$b=0$, the factor $4$ coming from $g_0=\frac14|dw|^2$, that is, $(g_0)_{AB}=\frac14\delta_{AB}$.

\emph{The bounds.} By \eqref{eq:delta}, $\|h-g_0\|_{C^2(\cB_8)}\le\delta$ and
hence $|\Gamma|+|\partial\Gamma|\le C\delta$ on $\{|w|<8\}$. As functions of $(x,u,p)$ on
$\{|x|<2,|u|<2,|p|\le1\}$, $G_{ij}=\frac14(\delta_{ij}+p_i\cdot p_j)+(h-g_0)_{AB}Q^A_iQ^B_j$,
whose second term and its first derivatives are bounded by $C\delta$. Since
$\frac14(I+p^Tp)\ge\frac14I$, the inverse $G^{ij}$ differs from $4((I+p^Tp)^{-1})_{ij}$ by
$C\delta$ in $C^1$, which is the bound on $e$, and $G^{ij}\ge\frac{4}{(1+C\delta)(1+|p|^2)}I\ge(2-C\delta)I$
for $|p|\le1$. The term $b$ is a polynomial in $G^{ij}$, $p$ and the values of $\Gamma$ at
$(x,u)$, hence bounded by $C\delta$ in $C^1$. The regularity $C^{k_0-2}$ follows from
$h\in C^{k_0}$.
\end{proof}

The vertical velocity $-2\theta\xi$ of \eqref{eq:F} does not appear in
\eqref{eq:grapheq}: by \cite[Lemma~2.3]{ChangWuZhang2026ContactBlowup}, whose proof is
pointwise and applies in the chart, the projection of a Legendrian flow \eqref{eq:F} moves
by the mean curvature flow of the transverse K\"ahler metric. The height and the angle
are controlled by the graph as follows.

\begin{lemma}[height and angle of a Legendrian graph]\label{lem:heightangle}
Let $\hat\Gamma_{u,\zeta}$ be a Legendrian graph over $\Omega_{T'}$ with
$\|u\|_{2,\alpha}\le1$, $u(0,0)=0$, $Du(0,0)=0$ and $\zeta(0,0)=0$; then
$|u(x,s)|\le\frac12|x|^2+|s|<\frac32$ and $|w|<2$ on $\Omega_{T'}$. In (a), (b) and (d) the
constants $C$ depend only on $n$.
\begin{enumerate}[(a)]
\item $|D_x\zeta|\le\frac54+C\delta$ and $|D^2_x\zeta|\le C(1+|D^2u|)$ on
$\Omega_{T'}$; hence $|\zeta(x,s)-\zeta(0,s)|\le(\frac54+C\delta)|x|$.
\item Suppose that the Legendrian graph is the trace of a flow \eqref{eq:F} for the
structure $\varsigma$ (satisfying \eqref{eq:F} for $\varsigma$ with $K+2$ replaced by the
constant $\kappa_\theta$ in the angle equation), with Legendrian angle $\theta$, regarded
as a function of $(x,s)$. Then
\begin{equation}\label{eq:theta-from-graph}
\partial_s\zeta=-2\beta_\varsigma\big((0,\partial_su)\big)-4\theta,\qquad\text{i.e.}\qquad
\theta=-\tfrac14\Big(\partial_s\zeta+2\beta_\varsigma(0,\partial_su)\Big),
\end{equation}
with $|2\beta_\varsigma((0,\partial_su))|\le\frac12|x||\partial_su|+C\delta$. In
particular $|\partial_s\zeta|\le\frac12+4C_\theta+C\delta$, and at $x=0$,
$|\partial_s\zeta(0,s)+4\theta(0,s)|\le C\delta$.
\item Under the hypothesis of (b), $|\nabla\theta|=|H|\le\sqrt n|B|$, and
$|D_x\theta|\le C(1+|D^2u|)$.
\item Let $(x,s)\in\Omega_{T'}$, and suppose that the hypothesis of (b) holds locally at
the points of the time segment from $(0,0)$ to $(0,s)$: every point $(0,\sigma)$ of this
segment has a neighbourhood in $\Omega_{T'}$ over which the Legendrian graph is the trace
of a flow \eqref{eq:F} for the structure $\varsigma$ (with $K+2$ replaced by
$\kappa_\theta$) whose Legendrian angle satisfies $|\theta|\le C_\theta$. Nothing is assumed
about the graph away from this segment. Then
\[
|\zeta(x,s)|\le(4C_\theta+C\delta)|s|+(\tfrac54+C\delta)|x|<4C_\theta+2=\tfrac{\bb^2}2
\]
for $\delta\le\delta_0$; here $\delta_0=\delta_0(n)$ of Lemma~\ref{lem:grapheq} is
decreased, if necessary, so that $2C\delta_0\le\frac34$ for the constant $C$ of (a)
and (b).
\end{enumerate}
\end{lemma}
\begin{proof}
The bound on $u$ follows from $|D^2u|\le1$, $|\partial_su|\le1$ and the normalisation,
and $|w|^2=|x|^2+|u|^2<1+\frac94$.

(a) By \eqref{eq:height}, $D_x\zeta\cdot e=-2\beta_\varsigma((e,Du\,e))$ for a unit
vector $e$. In the model, $\overline\beta((V_x,V_y))=-\frac14(y\cdot V_x-x\cdot V_y)$, so
$|\overline\beta((e,Du\,e))|\le\frac14(|u|+|x||Du|)\le\frac14(\frac32+1)=\frac58$;
and $|\beta_\varsigma-\overline\beta|\le C\delta$ on $\{|w|<4\}$ by \eqref{eq:delta}.
Differentiating \eqref{eq:height} once more gives the bound on $D^2_x\zeta$.

(b) Let $\hat F$ be the normal parametrisation of the flow near the graph point above
$(x,s)$, and let the graph points above the nearby $(x',s')$ be $\hat F(\gamma(x',s'),s')$.
The computation is local and uses the flow only near this point. Then
$\partial_s[\hat F(\gamma,s)]=\partial_s\hat F+d\hat F(\partial_s\gamma)=H-2\theta\xi+V$
with $V$ tangent to the Legendrian, hence horizontal. Apply $dz$: by Fact~\ref{fact:charts}(ii), $dz(H)=-2\beta_\varsigma(\pi_*H)$,
$dz(V)=-2\beta_\varsigma(\pi_*V)$ and $dz(\xi)=2$, so
$\partial_s\zeta=-2\beta_\varsigma(\pi_*(H+V))-4\theta$. Here $\pi_*(H+V)=(0,\partial_su)$
is the velocity of the projected graph point. In the model
$\overline\beta((0,\partial_su))=\frac14x\cdot\partial_su$, which gives the bound and
vanishes at $x=0$.

(c) The identity $|\nabla\theta|=|H|$ follows from \eqref{eq:F} and the fact that $\Phi$ is an
isometry on the contact distribution, and $|H|\le\sqrt n|B|$ by Cauchy--Schwarz. The second
fundamental form of the Legendrian lift equals that of the projected graph
\cite[Lemma~2.3]{ChangWuZhang2026ContactBlowup}, which is bounded by $C(1+|D^2u|)$, and
$|D_x\theta|\le|Q||\nabla\theta|\le C|H|$.

(d) The proof of (b) is local, so under the hypothesis of (d) the identity of (b) holds at
every point $(0,\sigma)$ of the time segment, where it gives
$|\partial_s\zeta(0,\sigma)|\le4C_\theta+C\delta$. Integrate this from $\zeta(0,0)=0$
along the time segment, and integrate (a), which holds on all of $\Omega_{T'}$ by
\eqref{eq:height} alone, along the space segment. Since $|s|<1$ and $|x|<1$, the sum is
less than $4C_\theta+\frac54+2C\delta\le4C_\theta+2$ for $\delta\le\delta_0$. The
constants in (a) and (b) bound $|\beta_\varsigma-\overline\beta|$ on $\{|w|<4\}$ by
\eqref{eq:delta} and depend only on $n$.
\end{proof}

\begin{lemma}[angle equation in graph coordinates]\label{lem:angleeq}
Under the hypothesis of Lemma~\ref{lem:heightangle}(b), the function $\theta(x,s)$
satisfies the linear equation
\begin{equation}\label{eq:angleeq}
\partial_s\theta=G^{ij}\big(\partial_{ij}\theta-{}^G\Gamma^k_{ij}\partial_k\theta\big)-H^k\partial_k\theta+\kappa_\theta\theta
\qquad\text{on }\Omega_{T'},
\end{equation}
where ${}^G\Gamma$ are the Christoffel symbols of the induced metric $G$ and $H^k$,
$k\le n$, are the $x$-components of the mean curvature vector of the projected graph. The
coefficients $G^{ij}$, $G^{ij}\,{}^G\Gamma^k_{ij}$ and $H^k$ are bounded in
$C^{0,\alpha}(\Omega_{T'})$ by a constant $C$, and $G^{ij}\ge I$.
\end{lemma}
\begin{proof}
In the normal parametrisation, $(\partial_t-\Delta)\theta=\kappa_\theta\theta$, and
$\Delta\theta=\Delta_\Gamma\theta$ is the Laplacian of the projected graph, since the
projection is an isometric immersion \cite[Lemma~2.3]{ChangWuZhang2026ContactBlowup}. In
graph coordinates, $\Delta_\Gamma\theta=G^{ij}(\partial_{ij}\theta-{}^G\Gamma^k_{ij}\partial_k\theta)$.
Relative to the normal parametrisation, the projected graph point has tangential velocity
$\pi_*V=c^k(e_k,\partial_ku)$ by the local calculation in the proof of
Lemma~\ref{lem:heightangle}(b), and $c^k=-H^k$ by
comparing $x$-components in $(0,\partial_su)=\pi_*H+\pi_*V$. Hence
$\partial_s\theta(x,s)=(\partial_t\theta)_{\rm normal}+c^k\partial_k\theta$. The
coefficients involve $h$, $\partial h$ at $(x,u)$, $Du$ and $D^2u$, all bounded in
$C^{0,\alpha}$ by $\|u\|_{2,\alpha}\le1$ and \eqref{eq:delta}. Finally, $G^{ij}\ge I$ by
\eqref{eq:coeff}.
\end{proof}

\subsection{Realisation of the graph and interior estimates}\label{subsec:subord}

\begin{lemma}[realisation of the full graph]\label{lem:filling}
Let $\cM$ be the trace of a smooth local Legendrian flow in $\cU$ (in the sense of
Definition~\ref{def:localflow}, or of Definition~\ref{def:traceflow}), $X\in\cM$, and let
$r$ be an admissible radius for $(\cM,X)$ with $\cC_{3r}(X)\subset\cU$, with witness
$(U,u,\zeta)$. Then, in the coordinates $\tilde U\Lambda_{X,r}$, for every
$(x,s)\in\Omega_{T'}$ the point $\hat\Gamma_{u,\zeta}(x,s)$ lies in $\bB$ and belongs to
the rescaled trace. Consequently
\[
\tilde U\Lambda_{X,r}(\cM)\cap\bB=\hat\Gamma_{u,\zeta}(\Omega_{T'}),\qquad
|\zeta|<\tfrac{\bb^2}2\ \text{on }\Omega_{T'},
\]
$u$ solves \eqref{eq:grapheq} on $\Omega_{T'}$ for the structure $\varsigma_{1/r}$, and
$\theta$ solves \eqref{eq:angleeq} there. Moreover $u$, $\zeta$ and $\theta$, as functions
of $(x,s)$, are smooth on $\Omega_{T'}$ up to the upper end $s=T'$; in particular their
parabolic H\"older norms of every order are finite on each
$\{|x|<\rho,\ |s|<\rho^2,\ s\le T'\}$ with $\rho<1$, whose closure is a compact subset of
$\Omega_{T'}$.
\end{lemma}
\begin{proof}
Write $\cM'$ for the rescaled trace in the coordinates $\tilde U\Lambda_{X,r}$. It is
relatively closed in $\tilde U\Lambda_{X,r}(\cU)$ (for a flow in the sense of
Definition~\ref{def:localflow} by the properness of $\hat F$; for a flow in trace form
because $\cM$ is relatively closed in $\cU\cap\{t\le T_*\}$, which is relatively closed in
$\cU$). The set $\tilde U\Lambda_{X,r}(\cU)$ contains $\overline\bB$ by Section~\ref{subsec:scale},
since the admissible radius $r$ is at most $r_1$ and $\cC_{3r}(X)\subset\cU$ (the closed box and its rotations lie
in $\{|w|\le\sqrt5,\ |z|\le\bb^2\}\times[-1,1]$, the time end points included). Throughout, the flow in trace form
is handled with the data of Definition~\ref{def:traceflow} only: no global map $\hat F$
or parameter domain $\cN$ is available for it, and Lemma~\ref{lem:param} is not used.

Fix $(x,s)\in\Omega_{T'}$ and let
$\gamma\colon[0,1]\to\Omega_{T'}$ be the path from $(0,0)$ to $(0,s)$ along the time axis
followed by the straight segment from $(0,s)$ to $(x,s)$. This path stays in $\Omega_{T'}$, which
is star-shaped with respect to this family of paths. Let $A\subset[0,1]$ be the set of
$\lambda$ such that $\hat\Gamma_{u,\zeta}(\gamma(\lambda'))\in\bB\cap\cM'$ for all
$\lambda'\le\lambda$. Then $0\in A$, since $\hat\Gamma(0,0)=0\in\cM'$.

\emph{The local structure of $\cM'$ in $\bB$.} We first prove the following claim: for every
$Z\in\bB\cap\cM'$ there are a neighbourhood $V$ of $Z$ in $\cM'$, contained in $\bB$, and
a neighbourhood $N$ of $\pi(Z)$ in $\Omega_{T'}$ (relative to $\Omega_{T'}$), where
$\pi\colon(x,y,z,s)\mapsto(x,s)$, such that $\pi|_V$ is a diffeomorphism of $V$ onto $N$
with inverse $\hat\Gamma_{u,\zeta}|_N$; in particular $\hat\Gamma_{u,\zeta}(N)=V\subset\cM'$.

Near $Z$ the trace $\cM'$ is a smooth embedded $(n+1)$-dimensional submanifold of
space-time with boundary $\cM'\cap\{s=T'\}$, on which the time function $s$ has no
critical points. For a flow in the sense of Definition~\ref{def:localflow} this holds
because $\hat F$ is proper (Definition~\ref{def:localflow}(a)), injective on the preimage
of $\bB$ by Definition~\ref{def:K}(ii), and an immersion of space-time since $\partial_t$
is transverse to the slices. Thus near $Z$ the trace is the image of an open subset of $\cN$
under the embedding $\hat F$, whose boundary points are those with $s=T'$, and $s$ is a
coordinate of $\cN\subset\Sigma\times(-\infty,T_*]$. For a flow in trace form it is part
of Definition~\ref{def:traceflow}: $\cM$ is an embedded $(n+1)$-dimensional submanifold
with boundary $\cM\cap\{t=T_*\}$ on which $t$ has no critical points. The
diffeomorphism $\tilde U\Lambda_{X,r}$ carries this submanifold to $\cM'$ with boundary
$\cM'\cap\{s=T'\}$, as $T'=(T_*-t(X))/r^2$. By Definition~\ref{def:K}(i), $\cM'\cap\bB$ is
contained in $\hat\Gamma_{u,\zeta}(\Omega_{T'})$. The map $\hat\Gamma_{u,\zeta}\colon(x,s)\mapsto(x,u(x,s),\zeta(x,s),s)$
is a $C^1$ embedding of $\Omega_{T'}$ (one-sided at $s=T'$), as $u\in C^{2,\alpha}$ and
$\zeta\in C^1$, with $\pi\circ\hat\Gamma_{u,\zeta}=\mathrm{id}$.

\emph{Tangent spaces.} Let $v=c'(0)$ for a $C^1$ curve $c$ in $\cM'\cap\bB$ with
$c(0)=Z$, allowing one-sided curves at a boundary point. Such velocities span $T_Z\cM'$.
Since
$c(\sigma)\in\cM'\cap\bB\subset\hat\Gamma_{u,\zeta}(\Omega_{T'})$ for all $\sigma$, we have
$c=\hat\Gamma_{u,\zeta}\circ(\pi\circ c)$ with $\pi\circ c$ a $C^1$ curve in $\Omega_{T'}$,
so $v=d\hat\Gamma_{u,\zeta}\big((\pi\circ c)'(0)\big)\in T_Z\hat\Gamma_{u,\zeta}$. Hence
$T_Z\cM'\subset T_Z\hat\Gamma_{u,\zeta}$, and the two spaces are equal, both having
dimension $n+1$. From $\pi\circ\hat\Gamma_{u,\zeta}=\mathrm{id}$ we get
$d\pi\circ d\hat\Gamma_{u,\zeta}=\mathrm{id}$, so $d\pi$ is injective on
$T_Z\hat\Gamma_{u,\zeta}=T_Z\cM'$, hence an isomorphism $T_Z\cM'\to\R^{n+1}$.

\emph{Interior points.} If $s(Z)<T'$, then $Z$ is an interior point of $\cM'$. By the
inverse function theorem, $\pi$ maps a neighbourhood $V$ of $Z$ in $\cM'$, which we take
inside $\bB$, diffeomorphically onto an open set $\widetilde N\ni\pi(Z)$ of $\R^{n+1}$,
which we take inside $\{|x|<1,\ |s|<1\}$. Then $N:=\widetilde N\cap\{s\le T'\}=\pi(V)$ is
a neighbourhood of $\pi(Z)$ in $\Omega_{T'}$ (as $V\subset\{s\le T'\}$ and
$\pi$ preserves $s$). For $p\in N$ we have $p=\pi(v)$ with $v\in V\subset\hat\Gamma_{u,\zeta}(\Omega_{T'})$,
so $v=\hat\Gamma_{u,\zeta}(\pi(v))=\hat\Gamma_{u,\zeta}(p)$. Hence the inverse of $\pi|_V$ is
$\hat\Gamma_{u,\zeta}|_N$.

\emph{Points at the upper end.} If $s(Z)=T'$ (which requires $T'<1$), then $Z$ is a
boundary point of $\cM'$. A smooth embedded submanifold with boundary extends, near a
boundary point, to a smooth embedded submanifold $\widetilde{\cM}$ without boundary of the
same dimension (extend a boundary chart smoothly across the boundary; on a small enough
ball the extension is still an embedding). Near $Z$ we have
$\cM'=\widetilde{\cM}\cap\{s\le T'\}$. Indeed, the function $s-T'$ on $\widetilde{\cM}$ has
non-zero differential at $Z$, vanishes on the boundary $\cM'\cap\{s=T'\}$ of $\cM'$, which is
a hypersurface of $\widetilde{\cM}$ through $Z$, and is negative on the interior of $\cM'$
near $Z$, so $\{s<T'\}$ is the side of this hypersurface containing the interior of $\cM'$.
Now $T_Z\widetilde{\cM}=T_Z\cM'$, so $d\pi$ is an isomorphism on $T_Z\widetilde{\cM}$.
The inverse function theorem gives a neighbourhood $\widetilde V$ of $Z$ in
$\widetilde{\cM}$, inside $\bB$ and inside the neighbourhood where
$\cM'=\widetilde{\cM}\cap\{s\le T'\}$, mapped by $\pi$ diffeomorphically onto an open set
$\widetilde N\ni\pi(Z)$, inside $\{|x|<1,|s|<1\}$. As $\pi$ preserves $s$,
$\pi(\widetilde V\cap\{s\le T'\})=\widetilde N\cap\{s\le T'\}$. Put
$V:=\widetilde V\cap\{s\le T'\}=\widetilde V\cap\cM'$, a neighbourhood of $Z$ in $\cM'$,
and $N:=\widetilde N\cap\{s\le T'\}$, a neighbourhood of $\pi(Z)$ in $\Omega_{T'}$. Then the
map $\pi|_V$ is a diffeomorphism onto $N$ (smooth up to the boundary), and as before
its inverse is $\hat\Gamma_{u,\zeta}|_N$. This proves the claim.

\emph{$A$ is open in $[0,1]$.} Let $\lambda\in A$ and $Z=\hat\Gamma(\gamma(\lambda))\in\bB\cap\cM'$.
By the claim, $\hat\Gamma_{u,\zeta}(N)=V\subset\cM'\cap\bB$ for a neighbourhood $N$ of
$\gamma(\lambda)=\pi(Z)$ in $\Omega_{T'}$, and $\gamma(\lambda')\in N$ for $\lambda'$ near
$\lambda$. Hence $\hat\Gamma(\gamma(\lambda'))\in\bB\cap\cM'$ for $\lambda'$ near $\lambda$.

\emph{$A$ is closed in $[0,1]$.} Let $\lambda_k\uparrow\lambda_*$ with $\lambda_k\in A$.
For every $\lambda<\lambda_*$ we have $\lambda\in A$, and by the claim the graph over a
neighbourhood $N$ of $\gamma(\lambda)$ in $\Omega_{T'}$ coincides with the trace near the
point $\hat\Gamma(\gamma(\lambda))$: $\hat\Gamma_{u,\zeta}(N)=V\subset\cM'$. For a flow in
the sense of Definition~\ref{def:localflow} this set is the trace of the (rescaled) flow
\eqref{eq:F}, for the structure $\varsigma_{1/r}$ with the constant $\kappa_\theta$ and
with $|\theta|\le C_\theta$. For a flow in trace form, Definition~\ref{def:traceflow}
provides a neighbourhood of the point in which $\cM'$ is the trace of a smooth local
Legendrian flow \eqref{eq:F} with the same structure, the same constant and the given
$\theta$, $|\theta|\le C_\theta$. So the local hypothesis of Lemma~\ref{lem:heightangle}(d)
is satisfied at every point $\gamma(\lambda)$, $\lambda<\lambda_*$. Indeed, if $\gamma(\lambda)$
lies on the time segment, the time segment from $(0,0)$ to $\gamma(\lambda)$ consists of
the points $\gamma(\lambda')$ with $\lambda'\le\lambda$. If $\gamma(\lambda)=(x',s)$ lies
on the space segment, the time segment from $(0,0)$ to $(0,s)$ is $\gamma([0,\lambda_1])$
for the $\lambda_1<\lambda$ with $\gamma(\lambda_1)=(0,s)$. In both cases every point of
the time segment is some $\gamma(\lambda')$ with $\lambda'<\lambda_*$. Hence
Lemma~\ref{lem:heightangle}(d) gives $|\zeta(\gamma(\lambda))|<\bb^2/2$ for
$\lambda<\lambda_*$ (\eqref{eq:height} holds for the whole witness graph, which is all
that part (a) of that lemma uses). By continuity,
$|\zeta(\gamma(\lambda_*))|\le\bb^2/2<\bb^2$. Together with $|x|<1$, $|u|<\frac32<2$,
$|s|<1$ on $\Omega_{T'}$, this shows $\hat\Gamma(\gamma(\lambda_*))\in\bB$. Being a limit of
points of $\cM'$ inside $\overline\bB\subset\tilde U\Lambda_{X,r}(\cU)$, this point belongs to
$\cM'$ by the relative closedness of $\cM'$. Hence $\lambda_*\in A$, and $A=[0,1]$.

Thus the whole graph over $\Omega_{T'}$ lies in $\bB\cap\cM'$, and the reverse inclusion
is Definition~\ref{def:K}(i). The equations hold because the graph is the trace of the
flow (Lemma~\ref{lem:grapheq}, Lemma~\ref{lem:angleeq}). Both lemmas are local
statements: Lemma~\ref{lem:grapheq} characterises the motion of the slices as sets near
each point, and the proof of Lemma~\ref{lem:angleeq} uses the normal parametrisation only
near the point considered. For a flow in trace form they are therefore applied with the local
flows of Definition~\ref{def:traceflow}, whose Legendrian angle is the given $\theta$ and
whose traces cover the graph. The height bound is Lemma~\ref{lem:heightangle}(d), applied
at every $(x,s)\in\Omega_{T'}$ as in the closedness step (now with $A=[0,1]$).

\emph{Regularity of the witness.} By the above, $\cM'\cap\bB=\hat\Gamma_{u,\zeta}(\Omega_{T'})$,
an open subset of $\cM'$, hence a smooth embedded $(n+1)$-dimensional submanifold with
boundary at the upper end $s=T'$ when $T'<1$ (one-sided derivatives there;
Definition~\ref{def:localflow}, respectively Definition~\ref{def:traceflow}). The map
$\hat\Gamma_{u,\zeta}\colon(x,s)\mapsto(x,u(x,s),\zeta(x,s),s)$ is a $C^1$ bijection of
$\Omega_{T'}$ onto this set, as $u\in C^{2,\alpha}$ and $\zeta\in C^1$, and
$\pi\circ\hat\Gamma_{u,\zeta}=\mathrm{id}$. Hence the restriction of $\pi$ to $\cM'\cap\bB$ is a
bijection onto $\Omega_{T'}$ with inverse $\hat\Gamma_{u,\zeta}$. By the claim, this restriction is a
local diffeomorphism at every point, at interior points and at points over $s=T'$ alike,
with local inverse $\hat\Gamma_{u,\zeta}$. A bijective local diffeomorphism is a
diffeomorphism; hence $\hat\Gamma_{u,\zeta}$ is smooth on $\Omega_{T'}$ up to $s=T'$, and
so are $u$ and $\zeta$. Consider next $\theta(x,s)$, the Legendrian angle at the graph point
$\hat\Gamma_{u,\zeta}(x,s)$. For a flow in the sense of Definition~\ref{def:localflow},
$\theta$ is smooth on $\cN$ and $\hat F^{-1}$ is smooth on $\cM'\cap\bB$, the image of an
open subset of $\cN$ under the embedding $\hat F$, so
$\theta(x,s)=\theta\circ\hat F^{-1}\circ\hat\Gamma_{u,\zeta}(x,s)$ is smooth up to $s=T'$.
For a flow in trace form, $\theta$ is by Definition~\ref{def:traceflow} a smooth function
on $\cM$, hence (in the rescaled coordinates) on $\cM'$, and
$\theta(x,s)=\theta\circ\hat\Gamma_{u,\zeta}(x,s)$ is smooth up to $s=T'$; no inverse of a
global parametrisation is used. All derivatives of these three functions are therefore continuous on
$\Omega_{T'}$, and their norms are finite on the closed subcylinders named in the
statement. This qualitative finiteness is used in
Lemma~\ref{lem:subord} and in Step~4 of the proof of Theorem~\ref{thm:B}; no bound on these
norms enters any estimate.
\end{proof}

\begin{lemma}[interior estimates of finite order]\label{lem:subord}
In the situation of Lemma~\ref{lem:filling} (for a flow in the sense of
Definition~\ref{def:localflow} or of Definition~\ref{def:traceflow}), put
$\Omega'_{T'}:=\{|x|<\frac12,\ -\frac14<s\le\min\{\frac14,T'\}\}$ and let
$\Omega^{(\rho)}_{T'}:=\{|x|<\rho,\ |s|<\rho^2,\ s\le T'\}$. Then:
\begin{enumerate}[(i)]
\item $\|u\|_{C^{m,\alpha}(\Omega^{(9/16)}_{T'})}\le C_m$ for $2\le m\le4$, where $C_m$ depends
only on $n,\alpha,m$ and $\delta_*$ (parabolic H\"older norms, one-sided at $s=T'$);
\item $|\zeta|<\bb^2/2$ on $\Omega_{T'}$, and $\|\zeta\|_{C^{1,\alpha}(\Omega'_{T'})}\le C$;
\item $|B|\le C_1$ and $|\nabla B|\le C'_1$ on the rescaled trace over $\Omega'_{T'}$;
in the original scale $|B|(X)\le C_1/r$, so that
\[
|B|(X)\le2C_1\,\Kreg(\cM,X)\qquad\text{whenever }\cC_{3/(2\Kreg(\cM,X))}(X)\subset\cU ;
\]
for $\Kreg(\cM,X)=0$, which occurs only in the model, the condition means
$\cC_r(X)\subset\cU$ for every $r>0$, and for $\Kreg(\cM,X)=\infty$ the bound is trivial;
\item $\|\theta\|_{C^{2,\alpha}(\Omega^{(9/16)}_{T'})}\le C$; in particular
$|D\theta|+|\partial_s\theta|\le C$ on $\Omega^{(9/16)}_{T'}$, hence
\[
|\theta(x,s)-\theta(0,0)|\le C(|x|+|s|)\le C(|x|+|s|^{1/2})\qquad\text{and}\qquad|\partial_s\theta|\le C
\qquad\text{on }\Omega'_{T'}.
\]
\end{enumerate}
\end{lemma}
\begin{proof}
All Schauder estimates below are the interior estimate of Theorem~\ref{thm:schauder-var}
(Appendix~\ref{app:schauder}) for linear parabolic systems with scalar principal part and
H\"older coefficients on one-sided cylinders, whose constants do not depend on the upper
end $T'$.

(i) By Lemma~\ref{lem:filling}, $u$ solves \eqref{eq:grapheq} on $\Omega_{T'}$ with
$\|u\|_{2,\alpha}\le1$. Put $\bar a(x,s):=a(x,u,Du)$ and $\bar b(x,s):=b(x,u,Du)$. Then
$u$ solves the linear system $\partial_su^\gamma=\bar a^{ij}D_{ij}u^\gamma+\bar b^\gamma$,
which is uniformly parabolic with $\bar a\ge I$ by \eqref{eq:coeff} and has the same principal
coefficients for every component. Since $u\in C^{2,\alpha}$ and $a,b$ are $C^{k_0-2}$ in
their arguments, $\bar a,\bar b\in C^{1,\alpha}(\Omega_{T'})$ with norms bounded by
$C(n,\delta_*)$. Differentiating the system in a direction $e_k$ gives, for
$w=D_ku$,
\[
\partial_sw^\gamma=\bar a^{ij}D_{ij}w^\gamma+(D_k\bar a^{ij})D_{ij}u^\gamma+D_k\bar b^\gamma ,
\]
a system of the same type with $C^{0,\alpha}$ coefficients and right-hand side bounded in
$C^{0,\alpha}$ by $C$ on $\Omega_{T'}$. By Lemma~\ref{lem:filling}, $u$ is smooth on
$\Omega_{T'}$ up to the upper end, so $w$ has finite parabolic $C^{2,\alpha}$ norm on
$\Omega^{(7/8)}_{T'}$, whose closure is a compact subset of $\Omega_{T'}$ (this finiteness is
the only use of the smoothness; the constants below do not involve it). Theorem~\ref{thm:schauder-var} with the outer domain
$\Omega^{(7/8)}_{T'}$ and the inner domain $\Omega^{(3/4)}_{T'}$ gives
$\|Du\|_{2,\alpha;\Omega^{(3/4)}}\le C$, i.e.\ $u\in C^{3,\alpha}(\Omega^{(3/4)})$ with a bound.
Then $\bar a,\bar b\in C^{2,\alpha}(\Omega^{(3/4)})$ with bounds, since $D_ku$ and $D^2_{kl}u$
are $C^{1,\alpha}$ there and $a,b$ are $C^{k_0-2}$ in $(x,u,p)$. Differentiating the system twice, we find that $w'=D_{kl}u$ solves a system of the same type on $\Omega^{(3/4)}_{T'}$, with finite norm on $\Omega^{(5/8)}_{T'}$ by Lemma~\ref{lem:filling} again. The theorem with the outer domain $\Omega^{(5/8)}_{T'}$
and the inner domain $\Omega^{(9/16)}_{T'}$ gives $\|D^2u\|_{2,\alpha;\Omega^{(9/16)}}\le C$,
i.e.\ $u\in C^{4,\alpha}$ with a bound on $\Omega^{(9/16)}_{T'}\supset\Omega'_{T'}$ (note that
$\Omega'_{T'}$ contains the time $s=\frac14$ when $T'\ge\frac14$, while $\Omega^{(1/2)}_{T'}$
does not; $\Omega^{(9/16)}_{T'}$, with $|s|<\frac{81}{256}$, does). The time derivative
$\partial_su$ is $C^{2,\alpha}$ by the equation. This uses $h\in C^{k_0}$ with $k_0=8\ge5$,
through the regularity of $a,b$ (Lemma~\ref{lem:grapheq}); higher orders are not claimed.

(ii) The first assertion is Lemma~\ref{lem:filling}. By \eqref{eq:height},
$D_x\zeta\cdot e=-2\beta_\varsigma((e,Du\,e))$ for unit vectors $e\in\R^n$, so that
$D_x\zeta$ is $C^{1,\alpha}$ in $(x,s)$. By
\eqref{eq:theta-from-graph}, $\partial_s\zeta=-2\beta_\varsigma((0,\partial_su))-4\theta$, which is
$C^{0,\alpha}$ by (i) and (iv). Together with $|\zeta|<\bb^2/2$ this gives the
$C^{1,\alpha}$ bound. No higher time derivatives of $\zeta$ are claimed.

(iii) The second fundamental form of the Legendrian lift equals that of the projected
graph \cite[Lemma~2.3]{ChangWuZhang2026ContactBlowup}, which is bounded by $C(1+|D^2u|)$
and whose first derivative is bounded by $C(1+|D^3u|)(1+|D^2u|)$ plus terms with
$\partial h$. Rescaling back multiplies $|B|$ by $1/r$. With $r=1/(2\Kreg(\cM,X))$, which
is admissible by Lemma~\ref{lem:admissible}(a), this gives the last bound when
$0<\Kreg(\cM,X)<\infty$. If $\Kreg(\cM,X)=0$, then every $r>0$ is admissible and
$\cC_{3r}(X)\subset\cU$, so that $|B|(X)\le C_1/r$ for every $r>0$, with $C_1$ independent of
$r$; hence $|B|(X)=0$.

(iv) By Lemma~\ref{lem:filling}, $\theta$ solves the linear uniformly parabolic equation
\eqref{eq:angleeq} on $\Omega_{T'}$, whose coefficients $G^{ij}$,
$-G^{ij}\,{}^G\Gamma^k_{ij}-H^k$ and $\kappa_\theta$ are bounded in $C^{0,\alpha}$ by $C$, with
zero right-hand side, and $|\theta|\le C_\theta$. By Lemma~\ref{lem:filling} again, $\theta(x,s)$
is smooth up to the upper end, so its parabolic $C^{2,\alpha}$ norm is finite on
$\Omega^{(3/4)}_{T'}$.
Theorem~\ref{thm:schauder-var}, which allows first- and zero-order terms with
$C^{0,\alpha}$ coefficients, with the outer domain $\Omega^{(3/4)}_{T'}$ and the inner
domain $\Omega^{(9/16)}_{T'}\supset\Omega'_{T'}$, gives
$\|\theta\|_{2,\alpha;\Omega^{(9/16)}_{T'}}\le C\|\theta\|_{0;\Omega^{(3/4)}_{T'}}\le CC_\theta$. No
H\"older bound on $\theta$ is assumed. This norm contains $\|D\theta\|_0$ and
$\|\partial_s\theta\|_0$ (Section~\ref{subsec:scale}, Appendix~\ref{app:schauder}), which
gives $|D\theta|+|\partial_s\theta|\le C$ on $\Omega^{(9/16)}_{T'}$.

For the modulus of
continuity, let $(x,s)\in\Omega'_{T'}$. Since $T'\ge0$ (Definition~\ref{def:K}), the point
$(0,0)$ lies in $\Omega'_{T'}$. The points $(0,\sigma s)$, $0\le\sigma\le1$, have times
between $0$ and $s$, hence in $(-\frac14,\min\{\frac14,T'\}]$, and the points $(\sigma x,s)$
have $|\sigma x|<\frac12$. So the path from $(0,0)$ to $(0,s)$ to $(x,s)$ lies in
$\Omega'_{T'}\subset\Omega^{(9/16)}_{T'}$, and integrating along it gives
$|\theta(x,s)-\theta(0,0)|\le|s|\,\|\partial_s\theta\|_0+|x|\,\|D\theta\|_0\le C(|x|+|s|)$.
Finally $|s|\le\frac14<1$ gives $|s|\le|s|^{1/2}$.
\end{proof}

\begin{remark}
The finite order in (i) is dictated by the finite-order closeness \eqref{eq:delta} of the
chart to the model and by the orders actually used below (the $C^{4,\alpha'}$ convergence in
Proposition~\ref{prop:compact}(ii), the start $m=0$ of the induction in Step~4 of the proof of
Theorem~\ref{thm:B}, $C^{2,\alpha}$ bounds for $\theta$). Lemma~\ref{lem:param} uses only the
bound on $|B|$. Limits of flows in the model are $C^\infty$ by a difference-quotient bootstrap
for the model equations, whose coefficients are smooth (Step~3 of the proof of
Proposition~\ref{prop:compact}). The a priori estimates of this lemma alone do not give this smoothness.
\end{remark}

\subsection{Compactness of local flows}\label{subsec:compactness}

We now consider a sequence of local flows in rescaled charts converging to the model. The
convergence of sets is understood \emph{relative to the limit region}, as in
\cite[\S2.6]{wh} with $U_i^c\to U^c$: for an open $\cU\subset\Hn\times\R$ and subsets
$S_i,S\subset\Hn\times\R$ we say that $S_i\to S$ \emph{locally in $\cU$} if
\begin{enumerate}[(a)]
\item every point of $S$ is the limit of a sequence of points $X_i\in S_i$, and
\item whenever $(i_k)$ is a subsequence and $X_{i_k}\in S_{i_k}$ converge to a point
$X\in\cU$, then $X\in S$.
\end{enumerate}
Condition (b) constrains every subsequence, not only the full sequence (for the full
sequence alone, $S_i=\{a\}$ for even $i$ and $S_i=\{b\}$ for odd $i$, $a\ne b$, would
converge to the empty set). Consequently the notion passes to subsequences: if
$S_i\to S$ locally in $\cU$ and $(i_k)$ is a subsequence, then $S_{i_k}\to S$ locally in
$\cU$. Indeed, a subsequence of $(i_k)$ is a subsequence of the original sequence, and the
sequences of (a) converge along $(i_k)$. (Points of $S_i$ accumulating on $\partial\cU$
are not constrained; this is the appropriate notion for proper flows in $\cU_i$ with
$\cU_i^c\to\cU^c$, since such a flow may approach $\partial\cU_i$.)

On $\Hn\times\R$ we
use the metric $\mathsf m$ induced by the Euclidean norm of $\R^{2n+2}$. By the formulas for the gauge and the group law in the proof of
Lemma~\ref{lem:quasi}, the quasi-distance $\|\cdot\|$ and $\mathsf m$ define the same
bounded sets and the same convergence.

\begin{proposition}[compactness]\label{prop:compact}
For $i=1,2,\dots$ let $F_i$ be smooth local Legendrian flows with upper ends $T_i$ in open
sets $\cU_i\subset\Hn\times\R$, for structures $\varsigma_i$ with
$\|\varsigma_i-\varsigma_{\rm Heis}\|_{C^{k_0}(\cB_R)}\le\delta_*\lambda_i^{-2}(1+R)^3/729$
for $1\le R\le8\lambda_i$ (Corollary~\ref{cor:rM}), $\lambda_i\to\infty$, the ambient
coordinates being the rescaled charts $\delta_{\lambda_i/r_i}\Upsilon_{p_i}$ of centres
$p_i\in M_i$ at scales $r_i/\lambda_i$ with $r_i\le r_{M_i}$ (so that
Corollary~\ref{cor:rM}, with $r=r_i$ and $\lambda=\lambda_i$, gives the stated bound and
its part (c) applies to the charts of all centres in compact sets; no lower bound on the
scales $r_i$ is assumed), with constants $\kappa_{\theta,i}\to0$ in the angle equation
\eqref{eq:angleeq} (the constants $\kappa_\theta$ of the flows $F_i$) and
$|\theta_i|\le C_\theta$. Assume that
$\cU_i^c\to\cU^c$ locally in $\Hn\times\R$, in the sense above with $S_i=\cU_i^c$,
$S=\cU^c$ and $\Hn\times\R$ in place of $\cU$ (every point of $\cU^c$ is a limit of points
of $\cU_i^c$, and every subsequential limit of points $X_{i_k}\in\cU_{i_k}^c$ lies in
$\cU^c$), for an open $\cU\subset\Hn\times\R$, that $T_i\to T\in(-\infty,\infty]$,
and that for every compact $\cK\subset\cU$ there are $L=L(\cK)$ and $i_0$ with
$\Kreg(\cM_i,X)\le L$ for all $X\in\cM_i\cap\cK$ and $i\ge i_0$. Then, after passing to
a subsequence:
\begin{enumerate}[(i)]
\item $\cM_i\to\cM$ locally in $\cU$, where $\cM\subset\cU$ is a smooth Legendrian flow
in trace form (Definition~\ref{def:traceflow}) for the model structure, with upper end
$T$, $K=-2$, and a Legendrian angle $\theta$ with $|\theta|\le C_\theta$; $\cM$ is
$C^\infty$ up to its upper end;
\item (local graphs at a uniform scale) for every $X\in\cM$, every sequence $X_i\in\cM_i$
with $X_i\to X$, and every $r>0$ such that $\cC_{4r}(X)\subset\cU$ and $r$ is admissible
for $(\cM_i,X_i)$ for all large $i$ (for instance $r\le1/(2L(\cK))$ for a compact
$\cK\subset\cU$ containing $\cC_{4r}(X)$, by Lemma~\ref{lem:admissible}(a), since $X_i\in\cK$
for $i$ large), there is a further subsequence along which
$\cM_i\cap\bB_{X_i,r}$ and $\cM\cap\bB_{X,r}$ are Legendrian graphs
$\hat\Gamma_{u_i,\zeta_i}$, $\hat\Gamma_{u,\zeta}$ over the base domains
$\Omega_{T'_i}$, $\Omega_{T'}$ in charts $\tilde U_i\Lambda_{X_i,r}$, $\tilde W\Lambda_{X,r}$,
where $W=\lim U_iV_i\in U(n)$ along the subsequence and $V_i\in U(n)$ are the rotations
of the transitions from the ambient chart to the charts of $X_i$
(Corollary~\ref{cor:rM}(c)), and $\bB_{Y,r}$ denotes the preimage of $\bB$ under the
respective chart; and $u_i\to u$ in $C^{2,\alpha'}$ on $\Omega_{T'}$ and in $C^{4,\alpha'}$ on
$\Omega'_{T'}$, $\zeta_i\to\zeta$ in $C^0(\Omega_{T'})$ and in $C^{1,\alpha'}(\Omega'_{T'})$,
$\theta_i\to\theta$ in $C^{2,\alpha'}(\Omega'_{T'})$, for every $\alpha'<\alpha$ (on
compact subsets of these domains, up to the upper end); in particular $r$ is admissible
for $(\cM,X)$;
\item[(ii$'$)] (lower semicontinuity with margin) for every $X\in\cM$ and every sequence
$X_i\in\cM_i$ with $X_i\to X$,
\[
\Kreg(\cM,X)\le\max\Big\{\liminf_i\Kreg(\cM_i,X_i),\ \frac4{d(X,\cU)}\Big\},
\]
and $\Kreg(\cM,X)\le\liminf_i\Kreg(\cM_i,X_i)$ when $\cU=\Hn\times\R$;
\item for every $s<T$ and every $s_i\to s$, the spatial area measures
$(F_i(\cdot,s_i))_\#\mu^i_{s_i}$ converge weakly to $\mu_s$ on $\cU(s)$, and the joint
measures $\fm^i_{s_i}:=(F_i,\theta_i)_\#\mu^i_{s_i}$ converge weakly to
$\fm_s:=(\mathrm{id},\theta)_\#\mu_s$ on $\cU(s)\times\R$. Here $\cU(s)$ is the
time-$s$ slice of $\cU$ and $\mu_s$ is the area measure of the embedded slice $\cM(s)$
(multiplicity one).
\end{enumerate}
\end{proposition}
\begin{proof}
\emph{Step 1: the limit set.} Let $d_i(Z):=\inf\{\mathsf m(Z,Y):Y\in\cM_i\}$ for
$Z\in\Hn\times\R$. We may assume $\cM_i\ne\emptyset$. Fix $Z_0\in\Hn\times\R$. If
$d_i(Z_0)\to\infty$ along a subsequence, then, since the $d_i$ are $1$-Lipschitz,
$d_i\to\infty$ uniformly on compact sets along it. The flows then leave every compact set
(no further subsequence of points $X_{i_k}\in\cM_{i_k}$ converges), and the limit
$\cM=\emptyset$ satisfies the proposition trivially, conditions (a) and (b) of local
convergence being empty. Otherwise, pass to a
subsequence with $d_i(Z_0)$ bounded. Then the $d_i$ are locally uniformly bounded and
$1$-Lipschitz, so by Arzel\`a--Ascoli a further subsequence converges locally uniformly
on $\Hn\times\R$ to a $1$-Lipschitz function $d_\infty\ge0$. We work along this
subsequence from now on. Put
$\widetilde\cM:=\{d_\infty=0\}$, a closed subset of $\Hn\times\R$, and
$\cM:=\widetilde\cM\cap\cU$, relatively closed in $\cU$. Condition (b): let $(i_k)$ be a
subsequence and $X_{i_k}\in\cM_{i_k}$ with $X_{i_k}\to X\in\cU$. Since $d_{i_k}\to d_\infty$
locally uniformly, $d_{i_k}$ is $1$-Lipschitz and $d_{i_k}(X_{i_k})=0$, we have
$d_\infty(X)=\lim_kd_{i_k}(X)\le\lim_k\big(d_{i_k}(X_{i_k})+\mathsf m(X,X_{i_k})\big)=0$,
so $X\in\cM$. Condition (a): if $X\in\cM$, then $d_i(X)\to0$, so there are $X_i\in\cM_i$
with $\mathsf m(X_i,X)\to0$. Hence $\cM_i\to\cM$ locally in $\cU$. (Points of
$\widetilde\cM\setminus\cU$ lie in $\cU^c$, possibly in its interior, and play no role:
only $\cM=\widetilde\cM\cap\cU$ is used below.)

\emph{Step 2: local graphs.} Let $X\in\cM$, $X_i\in\cM_i$ with $X_i\to X$, and let
$r>0$ satisfy $\cC_{4r}(X)\subset\cU$ and be admissible for $(\cM_i,X_i)$ for all $i\ge i_0$.
The hypothesis of the proposition provides such $r$. Indeed, fix a compact neighbourhood
$\cK_0\subset\cU$ of $X$; then $X_i\in\cK_0$ for $i$ large, and every $r>0$ with
$\cC_{4r}(X)\subset\cK_0$ and $r\le1/(2L(\cK_0))$ has both properties, by
Lemma~\ref{lem:admissible}(a). For the given $r$ no compact subset of $\cU$ containing
$\cC_{4r}(X)$ is used below.

First, $\cC_{3r}(X_i)\subset\cU_i$ for $i$ large. We begin by placing $\cC_{3r}(X_i)$ in a
bounded part of the ambient chart. Let $P_i$ be the point with ambient coordinates
$(0,t(X_i))$, that is, the centre $p_i$ at the time of $X_i$. Fix $\rho_*\ge1$ with
$\|X_i-P_i\|_{P_i}\le\rho_*$ for all $i$, computed in the ambient coordinates; such
$\rho_*$ exists since $X_i\to X$. In $M_i$, at the original scale $r_i/\lambda_i$,
this reads $\|X_i-P_i\|_{P_i}\le\rho_*r_i/\lambda_i$, and every $Z\in\cC_{3r}(X_i)$ satisfies
$\|Z-X_i\|_{X_i}<3r\,r_i/\lambda_i$. Let $i$ be so large that $\lambda_i\ge3(\rho_*+3r)$. Then the
right-hand side $3(\|X_i-P_i\|_{P_i}+\|Z-X_i\|_{X_i})$ of Lemma~\ref{lem:quasi} for the
points $P_i$, $X_i$, $Z$ is less than $3(\rho_*+3r)r_i/\lambda_i\le r_i\le r_{M_i}$.
Lemma~\ref{lem:quasi}, which holds with the constant $3$ at all scales $\le r_{M_i}$
(Corollary~\ref{cor:rM}), therefore gives $\|Z-P_i\|_{P_i}\le3(\rho_*+3r)r_i/\lambda_i<\infty$. Hence the
spatial point of $Z$ lies in the domain $U_{p_i}$ of the ambient chart, and in the
ambient coordinates $|w|\le3(\rho_*+3r)$, $|z|\le9\bb^2(\rho_*+3r)^2$ and
$|t(Z)-t(X_i)|<9r^2$, for every $Z\in\cC_{3r}(X_i)$. (In the model the charts are global,
and the model case of Lemma~\ref{lem:quasi} gives the same conclusion without a condition
on $\lambda_i$.) Fix $R\ge8$ with
$\{|w|\le3(\rho_*+3r),\ |z|\le9\bb^2(\rho_*+3r)^2\}\subset\cB_R$ (this box also contains the
spatial points of the $X_i$, as $\|X_i-P_i\|_{P_i}\le\rho_*$), and let $i$ be so large
that also $\lambda_i\ge R/8$.

Now suppose that there are $Z_i\in\cC_{3r}(X_i)\setminus\cU_i$
along a subsequence. By the bound just proved, they lie in $\cB_R\times\R$ with bounded
times, so $Z_i\to Z$ along a further subsequence. The limit satisfies
$\|Z-X\|_X\le3r$, because the quasi-distance $\|\cdot\|_{X_i}$, computed in the chart of
$X_i$, converges to $\|\cdot\|_X$ locally uniformly as $X_i\to X$. Indeed, by Corollary~\ref{cor:rM}(c) with this $R$,
in the ambient coordinates the chart of
$X_i$ is $\tilde V_i\ell_{X_i}^{-1}+E_i$ with $V_i\in U(n)$ and
$\|E_i\|_{C^{k_0}(\cB_R)}\le\delta_*(R/8)^4\lambda_i^{-2}\to0$, while in the model $\Upsilon_q=\ell_q^{-1}$. Moreover, the gauge of Lemma~\ref{lem:quasi} is
continuous and $U(n)$-invariant. Hence $Z\in\overline{\cC_{3r}(X)}\subset\cC_{4r}(X)\subset\cU$.
On the other hand, $Z_i\notin\cU_i$ along the subsequence and condition (b) of $\cU_i^c\to\cU^c$ (a
subsequential limit of points of $\cU_i^c$ lies in $\cU^c$) give $Z\in\cU^c$, a
contradiction.

Let
$(U_i,u_i,\zeta_i)$ be witnesses. In the chart $\tilde U_i\Lambda_{X_i,r}$ the structure is
$\hat\varsigma_i$, the structure of $M_i$ rescaled at the point of $X_i$ to the scale
$r\,r_i/\lambda_i$ and rotated by $\tilde U_i$. By Corollary~\ref{cor:rM}(a) at that point, with
$r_i$ in place of $r$, $\lambda_i/r\ge1$ in place of $\lambda$ and $R=8$, and since the model
structure is $U(n)$-invariant, we have
$\|\hat\varsigma_i-\varsigma_{\rm Heis}\|_{C^{k_0}(\cB_8)}\le\delta_*(r/\lambda_i)^2\to0$. The
witness graphs, and $\beta_{\hat\varsigma_i}$, $\Phi_{\hat\varsigma_i}$ below, refer to this
structure, whereas $\varsigma_i$ is the structure in the ambient coordinates. By
Lemma~\ref{lem:filling}, $\tilde U_i\Lambda_{X_i,r}(\cM_i)\cap\bB=\hat\Gamma_{u_i,\zeta_i}(\Omega_{T'_i})$
with $T'_i=(T_i-t(X_i))/r^2\to T':=(T-t(X))/r^2$. Lemmas~\ref{lem:subord}
and~\ref{lem:heightangle} give the bounds $\|u_i\|_{2,\alpha;\Omega_{T'_i}}\le1$,
$\|u_i\|_{C^{4,\alpha}(\Omega^{(9/16)}_{T'_i})}\le C_4$, $\|\theta_i\|_{C^{2,\alpha}(\Omega^{(9/16)}_{T'_i})}\le C$,
$\|\zeta_i\|_{C^{1,\alpha}(\Omega'_{T'_i})}\le C$, $|\zeta_i|<\bb^2/2$, $|D_x\zeta_i|\le2$,
$|\partial_s\zeta_i|\le1+4C_\theta$ on $\Omega_{T'_i}$.

Pass to a subsequence with
$W_i:=U_iV_i\to W\in U(n)$, $V_i$ as above. By Arzel\`a--Ascoli, along a further subsequence $u_i\to u$ in $C^{2,\alpha'}$ on compact
subsets of $\Omega_{T'}$ (including the upper end, the one-sided bounds holding up to
$s=T'_i$) and in $C^{4,\alpha'}$ on compact subsets of $\Omega^{(9/16)}_{T'}$, $\theta_i\to\theta$
in $C^{2,\alpha'}$ on compact subsets of $\Omega^{(9/16)}_{T'}$, $\zeta_i\to\zeta$ in $C^{1,\alpha'}$ on
compact subsets of $\Omega'_{T'}$, and $\zeta_i\to\zeta$ uniformly on compact subsets of
$\Omega_{T'}$. All these convergences hold along this one subsequence, with the same
one-sided convention at the upper end. (When $T'_i<T'$, a compact subset $\cK'$ of
$\Omega_{T'}$ need not lie in $\Omega_{T'_i}$, and convergence on $\cK'$ means convergence of
the stated norms of the differences on $\cK'\cap\Omega_{T'_i}$. The limits are first obtained
on the sets $\cK'\cap\{s\le T'-\epsilon\}$, which lie in $\Omega_{T'_i}$ for $i$ large, with a
diagonal sequence $\epsilon\to0$. The uniform one-sided H\"older bounds extend the limits, with
the same bounds, to $s=T'$, and through their time moduli they give the convergence on the
remaining slabs $\cK'\cap\{T'-\epsilon<s\le T'_i\}$.) The H\"older bounds pass to the limit:
$\|u\|_{C^{4,\alpha}(\Omega^{(9/16)}_{T'})}\le C_4$ and $\|\theta\|_{C^{2,\alpha}(\Omega^{(9/16)}_{T'})}\le C$. The
limit satisfies $\|u\|_{2,\alpha;\Omega_{T'}}\le1$ (the norm being lower semicontinuous
under $C^2$ convergence), $u(0,0)=0$ and $Du(0,0)=0$. Moreover, $\Gamma_u$ is Lagrangian for the
model, since $\Gamma_{u_i}^*d\beta_{\hat\varsigma_i}=0$ and $d\beta_{\hat\varsigma_i}\to\omega_0$
in $C^{k_0-1}$ on $\cB_8$. The height converges with its derivatives on $\Omega'_{T'}$:
$D_x\zeta_i\to D_x\zeta$ by \eqref{eq:height}, and $\partial_s\zeta_i=-2\beta_{\hat\varsigma_i}((0,\partial_su_i))-4\theta_i\to-2\overline\beta((0,\partial_su))-4\theta$
uniformly by \eqref{eq:theta-from-graph}. Hence $\zeta\in C^1(\Omega'_{T'})$ with
$\partial_s\zeta=-2\overline\beta((0,\partial_su))-4\theta$, and $\zeta$ is the height
function of $u$ with $\zeta(0,0)=0$.

In the ambient coordinates,
$\tilde U_i\Lambda_{X_i,r}=D_{1/r}\circ(\tilde W_i\ell_{X_i}^{-1}+\tilde U_iE_i)$
(composed with the time shift $t'\mapsto(t'-t(X_i))/r^2$), which converges in
$C^{k_0}_{\rm loc}$ to $\tilde W\Lambda_{X,r}=D_{1/r}\circ\tilde W\ell_X^{-1}$ as
$X_i\to X$, $W_i\to W$, $E_i\to0$. So the sets $\hat\Gamma_{u_i,\zeta_i}(\Omega_{T'_i})$,
mapped to $\Hn\times\R$ by the inverse charts, converge as sets to the image of
$\hat\Gamma_{u,\zeta}(\Omega_{T'})$ under $(\tilde W\Lambda_{X,r})^{-1}$. By Step 1 (local
convergence passing to the subsequence extracted here),
this image is $\cM\cap\bB_{X,r}$, i.e.\ $\hat\Gamma_{u,\zeta}(\Omega_{T'})=\tilde W\Lambda_{X,r}(\cM)\cap\bB$.
(The rotation $V_i$ of the chart transition and the rotation $U_i$ of the witness are
not separately controlled; only their product $W_i$ enters the limit.) In particular $r$ is admissible for
$(\cM,X)$ once we know that $\cM$ is a flow (Step 3). This is (ii). The subsequence
extracted here depends on $X$; the limit set $\cM$ of Step 1 does not, and all
statements about $\cM$ below are insensitive to further extractions.

For (ii$'$) let $\ell:=\liminf_i\Kreg(\cM_i,X_i)$ and let $r<\min\{d(X,\cU)/4,\ 1/\ell\}$.
Choose a subsequence along which $\Kreg(\cM_i,X_i)\to\ell$. Along it $r<1/\Kreg(\cM_i,X_i)$
for all large $i$, so $r$ is admissible for $(\cM_i,X_i)$ by Lemma~\ref{lem:admissible}(a).
The argument of Step 2 along this subsequence then shows that $r$ is admissible for
$(\cM,X)$. Letting $r\uparrow\min\{d(X,\cU)/4,1/\ell\}$ gives the bound. The margin
$\cC_{4r}(X)\subset\cU$ is needed because outside $\cU_i$ the witness graph is not known
to be the trace, and the estimates of Lemma~\ref{lem:subord} are not available there.
(We do not claim that $r$ is admissible for all large $i$ of the original sequence: if
$\Kreg(\cM_i,X_i)$ alternates between $1$ and $3$, then $\ell=1$ and $r=\frac12$ is
admissible only along the subsequence with $\Kreg=1$.)

\emph{Step 3: the limit is a flow in trace form.} We pass to the limit in
\eqref{eq:grapheq} (with $e_i,b_i\to0$ by \eqref{eq:coeff}), in \eqref{eq:angleeq}
(whose constant in the chart $\tilde U_i\Lambda_{X_i,r}$ is $r^2\kappa_{\theta,i}\to0$, with
coefficients converging in $C^{0,\alpha'}$) and in
\eqref{eq:theta-from-graph}. Then the limit graph on $\Omega^{(9/16)}_{T'}$ satisfies the model graph
equation, the model angle equation and the vertical velocity identity. These are the
three equations of \eqref{eq:F} with $K=-2$ for the Legendrian graph $\hat\Gamma_{u,\zeta}$
with angle $\theta$. The first holds because $H=\lim H_i=\lim\Phi_{\hat\varsigma_i}\nabla\theta_i=\overline\Phi\nabla\theta$
(uniform convergence of $D^2u_i$ and of $D\theta_i$).

The second holds by the following
computation, in which (in this paragraph only) $G(x,s):=(x,u(x,s),\zeta(x,s))$ denotes the
graph point, with velocity $\partial_sG=(0,\partial_su,\partial_s\zeta)$. By
Lemma~\ref{lem:grapheq}, applied to the limit graphs (which satisfy the model graph
equation), the horizontal velocity of the graph point is $H$ up to a tangential term: there
is a vector field $V$ tangent to the slices with $\pi_*(H+V)=(0,\partial_su)$, the velocity of
the projected graph point. Hence $\partial_sG-H-V$ has zero projection under $\pi$, so
$\partial_sG-H-V=c\,\overline\xi$ for a function $c$, where $\overline\xi=2\partial_z$. Both
$H=\overline\Phi\nabla\theta$ and $V$ (tangent to a Legendrian slice) are contact-horizontal,
that is, $\overline\eta(H)=\overline\eta(V)=0$ with $\overline\eta=\frac12dz+\overline\beta$, so
$dz(H)=-2\overline\beta(\pi_*H)$ and $dz(V)=-2\overline\beta(\pi_*V)$. Here $dz(H)$ need not
vanish. The limit of \eqref{eq:theta-from-graph} established in Step~2,
$\partial_s\zeta=-2\overline\beta((0,\partial_su))-4\theta$, therefore reads
\[
\begin{aligned}
\partial_s\zeta&=dz(H)-4\theta+dz(V),\qquad dz(H)+dz(V)=-2\overline\beta\big((0,\partial_su)\big),\\
2c&=dz(\partial_sG-H-V)=\partial_s\zeta-dz(H)-dz(V)=-4\theta .
\end{aligned}
\]
Thus $c=-2\theta$ and $\partial_sG=H-2\theta\overline\xi+V$: up to the tangential
reparametrisation term $V$, the graph point moves with the normal velocity
$H-2\theta\overline\xi$. The third equation is \eqref{eq:angleeq}.

The model graph equation is $\partial_su^\gamma=a^{ij}(Du)D_{ij}u^\gamma$ with $a$
smooth and $a\ge I$ for $|Du|\le1$ ($b=0$, the transverse metric being flat). We show
$u\in C^\infty$ on $\Omega'_{T'}$ by a difference-quotient bootstrap, which does not
presuppose the finiteness of unknown higher norms. Put $\rho_m:=\frac12+2^{-m-5}$, $m\ge0$.
Then $\rho_0=\frac{17}{32}<\frac9{16}$, $\rho_m\downarrow\frac12$, and $\Omega'_{T'}\subset\Omega^{(\rho_m)}_{T'}$
for all $m$, since $\rho_m>\frac12$ and $\rho_m^2>\frac14$. We claim that for $m\ge2$ and every
multi-index $|\beta|=m$, $v:=\partial^\beta_xu\in C^{2,\alpha}(\Omega^{(\rho_{m-2})}_{T'})$ with finite
norm. For $m=2$ this is the bound of Step~2. Given the claim for $m$, differentiation of the
equation gives $\partial_sv=\bar a^{ij}D_{ij}v+F_\beta$ with $\bar a:=a(Du)$ and
$F_\beta=\sum_{0<\beta'\le\beta}\binom{\beta}{\beta'}\partial^{\beta'}_x[a^{ij}(Du)]\,D_{ij}\partial^{\beta-\beta'}_xu$,
which involves derivatives of $u$ of order at most $m+1$, while $D_xF_\beta$ involves order at
most $m+2$. Hence $F_\beta$, $D_xF_\beta$ and $D_x\bar a$ are bounded in $C^{0,\alpha}(\Omega^{(\rho_{m-2})}_{T'})$.
For a spatial unit vector $e$ and $0<|h|<\eta:=\frac12(\rho_{m-2}-\rho_{m-1})$ the difference
quotient $\Delta_hv:=(v(\cdot+he,\cdot)-v)/h$ solves $\partial_s\Delta_hv=\bar a^{ij}D_{ij}\Delta_hv+f_h$
on $\Omega^{(\rho_{m-2}-\eta)}_{T'}$, where $f_h:=(\Delta_h\bar a^{ij})\,(D_{ij}v)(\cdot+he,\cdot)+\Delta_hF_\beta$. The inequality
$\|\Delta_hg\|_{0,\alpha}\le\|D_eg\|_{0,\alpha}$ bounds $\|f_h\|_{0,\alpha}$ independently of $h$. For
fixed $h$, $\Delta_hv$ has finite $C^{2,\alpha}$ norm, being a difference quotient of a
$C^{2,\alpha}$ function, which is all that Theorem~\ref{thm:schauder-var} requires. Spatial
difference quotients do not move the time, so the one-sided upper end is respected. The
theorem, with the outer domain $\Omega^{(\rho_{m-2}-\eta)}_{T'}$ and the inner domain
$\Omega^{(\rho_{m-1})}_{T'}$, gives
$|\Delta_hv|_{2,\alpha;\Omega^{(\rho_{m-1})}_{T'}}\le C(\|D_ev\|_0+\sup_h\|f_h\|_{0,\alpha})$, uniformly in
$h$. Letting $h\to0$ (Arzel\`a--Ascoli, the limit being $D_ev$) gives
$D_ev\in C^{2,\alpha}(\Omega^{(\rho_{m-1})}_{T'})$, i.e.\ the claim for $m+1$. Time derivatives
follow from the equation, so $u\in C^\infty$ on $\Omega'_{T'}$ (one-sided at the upper end).
The coefficients of the limit angle equation (with $\kappa_\theta=0$) are then smooth functions of
$Du,D^2u$, and the same bootstrap, started from $\theta\in C^{2,\alpha}(\Omega^{(9/16)}_{T'})$, gives
$\theta\in C^\infty$ on $\Omega'_{T'}$. Finally, $D_x\zeta\cdot e=-2\overline\beta((e,Du\,e))$ for unit vectors $e$ by
\eqref{eq:height}, and $\partial_s\zeta=-2\overline\beta((0,\partial_su))-4\theta$ by
\eqref{eq:theta-from-graph}. These derivatives are smooth, so $\zeta\in C^\infty$ on $\Omega'_{T'}$.

Every point of $\cM$ is the centre of such a chart
(with a margin $r$ depending on the point), and the inner domains $\Omega'$ of these
charts cover $\cM$. On the overlap of two charts both graphs equal $\cM$, so the local
Legendrian graphs patch to an embedded
$(n+1)$-dimensional $C^\infty$ submanifold with boundary $\cM\cap\{t=T\}$, on which
$\partial_t$ is transverse to the slices and $\theta$ is a well-defined smooth function
(it is determined by the trace through \eqref{eq:theta-from-graph}). Each slice
$\cM(s)$ is an embedded Legendrian submanifold: in a box it is one graph, with one point
above each base point. Near each point the normal parametrisation is obtained by
integrating the normal velocity $H-2\theta\overline\xi$, which is smooth; this gives the
local flows of Definition~\ref{def:traceflow}. Hence $\cM$ is a smooth Legendrian flow in
trace form in $\cU$ with upper end $T$.

Finally, in the chart of Step~2 the projection
$(x,u,z,s)\mapsto(x,s)$, restricted to $\cM\cap\bB=\hat\Gamma_{u,\zeta}(\Omega_{T'})$, is a bijection
onto $\Omega_{T'}$ whose differential is everywhere an isomorphism. Indeed, each slice is the graph of
the $C^2$ function $u(\cdot,s)$, and the time function has no critical points on $\cM$. By the
inverse function theorem (up to the boundary $s=T'$), this projection is a diffeomorphism, so $\zeta$ is
smooth on $\Omega_{T'}$; in particular $\zeta\in C^1(\Omega_{T'})$, and $r$ is admissible for
$(\cM,X)$. This proves (i) and (ii).

\emph{Step 4: measures, (iii).} Fix $s<T$, $s_i\to s$ and $\gamma\in C_c(\cU(s)\times\R)$,
the second variable being the angle. Let $K_\gamma:=\pi_{\Hn}(\spt\gamma)\subset\cU(s)$ be
the compact projection of the support of $\gamma$ to the space variable. Let $A\subset\cU$
be a compact neighbourhood of $K_\gamma\times\{s\}$. Cover $\cM\cap A$ by finitely many
inner boxes $\frac12\bB_{X^{(m)},r_m}$ of the charts of Step 2 (the preimages of
$\{|x|<\frac12,|y|<2,|z|<\bb^2\}\times(-\frac14,\frac14)$), $m=1,\dots,N$.

The inner box
of a centre $X$ at an admissible radius $r$ does not depend on the subsequence extracted
in Step~2. Indeed, let $W$ and $W'$ be limit rotations at $X$ obtained along two subsequences,
with limit witnesses $(u,\zeta)$ and $(u',\zeta')$. Near the origin of the chart
$\Lambda_{X,r}$ the slice $\Lambda_{X,r}(\cM)(0)$ is both
$\tilde W^{-1}\hat\Gamma_{u(\cdot,0),\zeta(\cdot,0)}$ and
$\tilde W'^{-1}\hat\Gamma_{u'(\cdot,0),\zeta'(\cdot,0)}$ (by
$\hat\Gamma_{u,\zeta}(\Omega_{T'})=\tilde W\Lambda_{X,r}(\cM)\cap\bB$ and its primed
version; $0\le T'$). The tangent plane at the origin of the horizontal projection of
this embedded slice is $W^{-1}P$ as well as $W'^{-1}P$, because $u(0,0)=0$, $Du(0,0)=0$,
and likewise for $u'$. Hence $A_0:=W'W^{-1}\in U(n)$ satisfies $A_0P=P$. As $A_0$
commutes with $J_0$ and preserves the real inner product, $O:=A_0|_P\in O(n)$ and
$A_0(x+J_0y)=Ox+J_0Oy$, i.e.\ $A_0(x,y)=(Ox,Oy)$, which preserves
$\{|x|<\frac12,|y|<2,|z|<\bb^2\}\times(-\frac14,\frac14)$ (with $\tilde A_0$ fixing $z$ and
$s$). Thus
$(\tilde W'\Lambda_{X,r})^{-1}(\tfrac12\bB)=\Lambda_{X,r}^{-1}\tilde W^{-1}\tilde A_0^{-1}(\tfrac12\bB)=(\tilde W\Lambda_{X,r})^{-1}(\tfrac12\bB)$.

We extract the subsequences of Step~2 at the finitely many centres $X^{(m)}$
successively, so that one common subsequence serves for all $m$. Along it the charts of
$X^{(m)}_i$ converge to a limit chart of $X^{(m)}$ whose inner box is the one just used
to form the cover. By local set
convergence in $\cU$, $\cM_i\cap A$ lies in the union of these boxes for $i$ large
(otherwise a subsequence of points of $\cM_i\cap A$ outside the union would accumulate at
a point of the compact set $A\subset\cU$; by condition (b) of
local convergence this point lies in $\cM\cap A$ and hence in the union, which is open: a contradiction). The time-$s$
slices $\cW_m\subset\Hn$ of the inner boxes are open and cover
$\cM(s)\cap K_\gamma$, and $\cU(s)\setminus\cM(s)$ is open in $\Hn$ ($\cM$ being relatively
closed in $\cU$). Hence there is a partition of unity on the compact set $K_\gamma$
subordinate to the open cover $\{\cW_1,\dots,\cW_N,\ \cU(s)\setminus\cM(s)\}$: functions
$\chi_0,\chi_1,\dots,\chi_N\in C_c(\Hn)$, $0\le\chi_m\le1$, with $\sum_{m=0}^N\chi_m=1$ on
$K_\gamma$, $\spt\chi_m\subset\cW_m$ for $m\ge1$ (so that $\spt\chi_m\times\{s\}$ is compactly
contained in the inner box of $X^{(m)}$), and $\spt\chi_0\subset\cU(s)\setminus\cM(s)$. Since
the charts of $X^{(m)}_i$ converge to the chart of $X^{(m)}$ and $s_i\to s$, for $i$ large
$\spt\chi_m\times\{s_i\}$ lies in the inner box of $X^{(m)}_i$ ($m\ge1$). Moreover, by local set
convergence in $\cU$, $\cM_i(s_i)\cap\spt\chi_0=\emptyset$ for $i$ large (otherwise a
subsequence of points of $\cM_i$ in the compact sets $\spt\chi_0\times\{s_i\}$ would
accumulate at a point of $\spt\chi_0\times\{s\}\subset\cU\setminus\cM$, contradicting
condition (b)), so the term with $\chi_0$ vanishes for $i$
large and in the limit. Hence
$\int\gamma\,d\fm^i_{s_i}=\sum_{m=1}^N\int(\chi_m\circ\pi_{\Hn})\,\gamma\,d\fm^i_{s_i}$ for
$i$ large.

Let $\Phi^{(m)}_i:=\tilde U^{(m)}_i\Lambda_{X^{(m)}_i,r_m}$ be the chart of $X^{(m)}_i$ of
Step~2, with witness $(u^{(m)}_i,\zeta^{(m)}_i)$ and angle $\theta^{(m)}_i$, and let
$\Phi^{(m)}:=\tilde W^{(m)}\Lambda_{X^{(m)},r_m}$ be the limit chart. In the box of
$X^{(m)}_i$ the slice of $\cM_i$ at the time $s_i$ is the Legendrian graph of
$u^{(m)}_i(\cdot,s'_i)$, where $s'_i:=(s_i-t(X^{(m)}_i))/r_m^2\to s':=(s-t(X^{(m)}))/r_m^2$ is the
rescaled time. By Definition~\ref{def:K}(ii), the map $\hat F_i$ is injective on the
preimage of the box, so the parameter measure $\mu^i_{s_i}$ restricted to the box is the
area measure, for the ambient structure $\varsigma_i$, of this graph slice (multiplicity
one). Let $G^{(m)}_i$ be the induced metric of the projected graph
$\Gamma_{u^{(m)}_i(\cdot,s'_i)}$ for the transverse metric of the chart structure of
$\Phi^{(m)}_i$, as in Section~\ref{subsec:grapheq}, so that the area of the graph slice
over a set $E\subset\{|x|<1\}$ in the chart structure is $\int_E\sqrt{\det G^{(m)}_i}\,dx$,
the area of a Legendrian graph being that of its projection. The chart structure is the
ambient structure written in the coordinates of $\Phi^{(m)}_i$: it is transported by the
transition to the chart of $X^{(m)}_i$ at the ambient scale (the same structure in other
coordinates) and then rescaled by $\delta_{1/r_m}$, the $\cD$-homothety with parameter
$r_m^{-2}$ of \eqref{eq:rescaled-coeff}. The transverse metric of the chart structure assigns to a horizontal
vector $r_m^{-1}$ times its $\varsigma_i$-length (in the chart coordinates $w'=w/r_m$ the
transverse metric is $h(r_mw')$, while the ambient one is $r_m^2h(r_mw')$). Hence
$n$-dimensional $\varsigma_i$-areas are $r_m^n$ times chart areas. Writing $q^{(m)}_i(x)$
for the spatial component of $(\Phi^{(m)}_i)^{-1}\big(\hat\Gamma^{(m)}_i(x,s'_i)\big)$, the
ambient point of the graph point over $x$, the $m$-th term equals
\[
r_m^n\int_{\{|x|<1/2\}}(\chi_m\gamma)\big(q^{(m)}_i(x),\theta^{(m)}_i(x,s'_i)\big)\sqrt{\det G^{(m)}_i(x,s'_i)}\,dx ,
\]
the integrand being pulled back by the inverse chart to the fixed inner base domain
$\{|x|<\frac12\}$ of the $m$-th chart. (For the plane $u\equiv0$ in a chart of radius $r$
the chart area over $\{|x|<\frac12\}$ is $2^{-n}|\{|x|<\frac12\}|$ by
Section~\ref{subsec:scale}, and the ambient area is $r^n$ times this, in accordance with
the factor.) Since $\spt\chi_m\times\{s_i\}$ is compactly contained in the inner box, the
pulled-back integrands are supported in a fixed compact subset $\{|x|\le\frac12-\epsilon_m\}$
of that domain, for some margin $\epsilon_m\in(0,\frac12)$, for all large $i$. There the
integrands converge uniformly ($u^{(m)}_i\to u^{(m)}$ in $C^1$,
$\zeta^{(m)}_i\to\zeta^{(m)}$ and $\theta^{(m)}_i\to\theta^{(m)}$ uniformly on the inner
domain; $G^{(m)}_i\to G^{(m)}$ uniformly, where
$G^{(m)}_{jk}=\frac14(\delta_{jk}+\partial_ju^{(m)}\cdot\partial_ku^{(m)})$ is the induced
metric for the model, because $Du^{(m)}_i\to Du^{(m)}$ and the transverse metric of the
chart structure tends to $g_0$ by \eqref{eq:delta}; the inverse charts
$(\Phi^{(m)}_i)^{-1}\to(\Phi^{(m)})^{-1}$ locally uniformly by Step~2; and $\chi_m\gamma$ is
continuous with compact support), so each term converges to
\[
r_m^n\int_{\{|x|<1/2\}}(\chi_m\gamma)\big(q^{(m)}(x),\theta^{(m)}(x,s')\big)\sqrt{\det G^{(m)}(x,s')}\,dx ,
\]
where $q^{(m)}(x)$ is the spatial component of $(\Phi^{(m)})^{-1}(\hat\Gamma^{(m)}(x,s'))$. This limit is
the corresponding term for $\cM$ at the time $s$: $\mu_s$ is the model area measure of the
embedded slice $\cM(s)$, and in the model the chart $\Phi^{(m)}$ is $D_{1/r_m}$ composed
with a left translation and a rotation, isometries of the model, so again model areas
are $r_m^n$ times areas in the chart. Summing over $m$ gives
$\int\gamma\,d\fm^i_{s_i}\to\int\gamma\,d\fm_s$. Testing with $\gamma(q,\vartheta)=\phi(q)\eta(\vartheta)$,
where $\phi\in C_c(\cU(s))$ and $\eta\in C_c(\R)$ equals one on
$[-C_\theta,C_\theta]$, gives the convergence of the spatial area measures.

The charts at the centres $X^{(m)}$, and the convergences in
them, were obtained along further subsequences depending on $s$, $(s_i)$ and $\gamma$. Every
subsequence of the subsequence of Step~1 has such a further subsequence, along which
$\int\gamma\,d\fm^i_{s_i}\to\int\gamma\,d\fm_s$, and this limit is determined by $\cM$ and $\theta$
(which is determined by the trace, by \eqref{eq:theta-from-graph}). Hence the convergence holds
along the subsequence of Step~1 itself, for every $s<T$, $s_i\to s$ and $\gamma$.
\end{proof}

\begin{lemma}[one sheet near a plane]\label{lem:onesheet}
In the situation of Proposition~\ref{prop:compact}, let the ambient coordinates be
those of the rescaled chart of a centre $p_i\in M_i$ at the scale $r_i/\lambda_i$,
$r_i\le r_{M_i}$ (so that Corollary~\ref{cor:rM} applies to the transitions to the charts
of nearby centres), possibly composed with a fixed unitary rotation $\tilde U_0$,
$\tilde U_0(w,z)=(U_0w,z)$, under which everything below is invariant. Suppose that $0\in\cM$, $T\ge0$, $\cC_4(0)\subset\cU$, and that
\[
\cM\cap\cC_4(0)=\Pi\cap\cC_4(0),\qquad \Pi:=\{(w,-4\bar\theta s,s):\ w\in P,\ s\le T\},
\]
the lift through the origin of the Lagrangian subspace $P=\{y=0\}$, with
$|\bar\theta|\le C_\theta$. Suppose also that $\cC_4(0)\subset\cU_i$ and
$\Kreg(\cM_i,\cdot)\le C_*$ on $\cM_i\cap\cC_4(0)$ for $i$ large. Let
$\bB^{3/2}:=\{|x|<\frac32,|y|<3,|z|<\frac94\bb^2\}\times(-\frac94,\frac94)$ and
$\Omega^{3/2}_{T_i}:=\{|x|<\frac32,|s|<\frac94,s\le T_i\}$. Then for $i$ large
$\cM_i\cap\bB^{3/2}$ is the Legendrian graph $\hat\Gamma_{\tilde u_i,\tilde\zeta_i}(\Omega^{3/2}_{T_i})$
of a single function $\tilde u_i\colon\Omega^{3/2}_{T_i}\to\R^n$ (in the ambient
coordinates, over the base $P$), with $\|\tilde u_i\|_{2,\alpha;\Omega^{3/2}_{T_i}}\le C$,
$\tilde u_i\to0$ in $C^{2,\alpha'}$ and $\tilde\zeta_i(x,s)\to-4\bar\theta s$ uniformly on
compact subsets (up to the upper end), and $\hat F_i$ is injective on the preimage of
$\bB^{3/2}$.
\end{lemma}
\begin{proof}
\emph{Scales and the model pieces.} Put $\varrho:=\min\{1/(2C_*),\frac18\}$ and
$\Pi':=\Pi\cap\{|x|\le2,\ |s|\le3\}$, a compact subset of $\cC_4(0)$ (its points have
$|w|\le2$ and $|z|\le12C_\theta<16\bb^2$). For $Y\in\Pi'$ and $Y_i\in\cM_i$ with
$Y_i\to Y$, Step 2 of the proof of Proposition~\ref{prop:compact} at the scale $\varrho$
(admissible for $i$ large, with $\cC_{4\varrho}(Y)\subset\cC_4(0)\subset\cU$) gives
witnesses $(U_i,u_i,\zeta_i)$ in the charts $\Phi_i:=\tilde U_i\Lambda_{Y_i,\varrho}$,
converging to a witness $(W,u,\zeta)$ of $(\cM,Y)$ in $\Phi:=\tilde W\Lambda_{Y,\varrho}$.
Here $W=\lim U_iV^{\rm c}_i$ along the subsequence of that step, where $V^{\rm c}_i\in U(n)$
are the rotations of Corollary~\ref{cor:rM}(c) used there (denoted by $V_i$ in that step).

\emph{The full transition.} Write $\Psi_i:=\Phi_i^{-1}$ for the map from the
$Y_i$-chart coordinates to the ambient coordinates. The chart of $Y_i$ is related to the
ambient chart by Corollary~\ref{cor:rM}(b), applied to the pair $(p_i,q(Y_i))$ with
$Y_i\in\cC_4(0)$ for $i$ large (as $Y_i\to Y\in\Pi'\subset\cC_4(0)$ and $\cC_4(0)$ is
open). The spatial part $\{|w|<4,\ |z|<16\bb^2\}$ of this cylinder lies in $\cB_8$, as
Corollary~\ref{cor:rM}(b) requires. Hence there are $V_i\in U(n)$ with
$\Lambda_{Y_i,\varrho}=D_{1/\varrho}\circ(\tilde V_i\ell_{Y_i}^{-1}+E'_i)$,
$\|E'_i\|_{C^{k_0}(\cB_8)}\le\delta_*\lambda_i^{-2}$, where
$D_{\varrho}(w,z,s)=(\varrho w,\varrho^2z,\varrho^2s)$. Composing with the witness
rotation $U_i$ and inverting, with $W_i:=U_iV_i\in U(n)$,
\begin{equation}\label{eq:Psi}
\Psi_i=\ell_{Y_i}\circ\tilde W_i^{-1}\circ D_\varrho+E_i,\qquad
\|E_i\|_{C^{k_0}(\bB)}\le C\lambda_i^{-2}\to0 ,
\end{equation}
and the time component of $\Psi_i$ is exactly $s\mapsto s_{Y_i}+\varrho^2s$.

The rotations $V_i$ need not be the rotations $V^{\rm c}_i$ of Step~2: there
Corollary~\ref{cor:rM}(c), with the radius $R\ge8$ fixed in that step, writes the
transition of the same pair as $\tilde V^{\rm c}_i\ell_{Y_i}^{-1}+E^{\rm c}_i$ on $\cB_R$,
with $\|E^{\rm c}_i\|_{C^{k_0}(\cB_R)}\le\delta_*(R/8)^4\lambda_i^{-2}$. For $V\in U(n)$ the
$w$-component of $\tilde V\ell_{Y_i}^{-1}$ is $w\mapsto V(w-w_{Y_i})$, with $w$-derivative $V$. Hence, on
$\cB_8\subset\cB_R$, where $\tilde V_i\ell_{Y_i}^{-1}+E'_i=\tilde V^{\rm c}_i\ell_{Y_i}^{-1}+E^{\rm c}_i$,
the matrix $V_i-V^{\rm c}_i$ is the $w$-derivative of the $w$-component of
$E^{\rm c}_i-E'_i$, so that
\[
|V_i-V^{\rm c}_i|\le\|E'_i\|_{C^1(\cB_8)}+\|E^{\rm c}_i\|_{C^1(\cB_8)}
\le\delta_*\big(1+(R/8)^4\big)\lambda_i^{-2}\longrightarrow0 .
\]
As $U_i$ is unitary, $|W_i-U_iV^{\rm c}_i|=|V_i-V^{\rm c}_i|$. Hence $W_i\to W$ along the
subsequence of that step, and the limit map is $\Psi:=\ell_Y\circ\tilde W^{-1}\circ D_\varrho$.
Since $\cM\cap\bB_{Y,\varrho}=\Pi\cap\bB_{Y,\varrho}$ is the lift of $P$ and the limit
witness graph is $\Psi(\hat\Gamma_{u,\zeta})$, the linear map $W^{-1}$ maps $P$ onto $P$
and $u\equiv0$. (The rotation $V_i$ of the chart transition and the rotation $U_i$ of
the witness are not separately controlled; only their product is constrained in the limit.) The
vertical component of $\ell_{Y}$ is
\begin{equation}\label{eq:vertical}
z\big(\ell_{Y}(w',z')\big)=z'+z_Y+\tfrac12(y_Y\cdot x'-x_Y\cdot y'),
\end{equation}
where the last term is the cross term of the group law. This term is accounted for below by working with $\Psi_i$
rather than with differences of the ambient height.

\emph{Transport to the fixed base.} Write $\tau$ for the local time of the $Y_i$-chart
and $t=s_{Y_i}+\varrho^2\tau$ for the ambient time. Here $\pi(w,z)=w$ is the horizontal
projection (on space-time points the time is ignored). The horizontal component
$\pi\circ\Psi_i$ does not depend on $z$: by Fact~\ref{fact:charts}(i) every chart
straightens the Reeb field to $2\partial_z$, so every transition $\varphi$ satisfies
$\varphi_*\partial_z=\partial_z$ and has the form $(w,z)\mapsto(\tilde w(w),z+c(w))$; the
horizontal components of $\ell_{Y_i}$, $\tilde W_i^{-1}$, $D_\varrho$ and $E_i$ are
functions of $w$ alone. Consider the map $\Xi_i\colon\Omega_{T'_i}\to P\times\R$,
$\Xi_i(x,\tau):=(f_i(x,\tau),\ s_{Y_i}+\varrho^2\tau)$, with
\begin{equation}\label{eq:Xi}
f_i(x,\tau):=\pi_P\circ\pi\circ\Psi_i(x,u_i(x,\tau),\cdot)
=x_{Y_i}+\varrho\,\pi_PW_i^{-1}(x,u_i(x,\tau))+\pi_P\pi E_i(x,u_i(x,\tau)),
\end{equation}
where $\pi_P$ is the orthogonal projection of $\C^n$ onto $P$ along $J_0P$. Thus
$f_i(x,\tau)$ is the spatial component of the ambient base projection
$(w,z,t)\mapsto(\pi_P w,t)$ of the graph point over $(x,\tau)$, and $f_i$ depends
on $(x,\tau)$ only through $(x,u_i(x,\tau))$. Hence $f_i$ is $C^{2,\alpha}$ in the parabolic
sense with norms bounded in terms of $\|u_i\|_{2,\alpha}\le1$ and of the $C^{k_0}$ maps
involved, uniformly in $i$ (same exponent $\alpha$, since a $C^{k_0}$ map composed with a
$C^{2,\alpha}$ map is $C^{2,\alpha}$). Its spatial Jacobian is
$D_xf_i=\varrho\,\pi_PW_i^{-1}(I,Du_i)+O(\lambda_i^{-2})$, which converges in
$C^{1,\alpha'}$, on compact subsets of $\Omega_{T'}$ up to the upper end, to the constant
map $\varrho\,\pi_PW^{-1}|_P$, which is $\varrho$ times an isometry of $P$ ($W$ is unitary
and $W^{-1}P=P$).

\emph{Preliminary bounds on the time window $|\tau|\le\frac{49}{64}$.} The statements of
this paragraph refer to $\{|x|<1,\ |\tau|\le\frac{49}{64},\ \tau\le T'_i\}$, which contains
all local times used below. Since $u_i\to0$ in $C^1$ on compact subsets of $\Omega_{T'}$
up to the upper end, and since on the strip $T'<\tau\le T'_i$ (when $T'_i>T'$) we have
$|u_i(x,\tau)-u_i(x,T')|\le T'_i-T'$ and $|Du_i(x,\tau)-Du_i(x,T')|\le(T'_i-T')^{\alpha/2}$ by
$\|u_i\|_{2,\alpha}\le1$, the convergence $u_i\to0$ in $C^1$ is uniform on
$\{|x|\le\frac{15}{16},\ |\tau|\le\frac{49}{64},\ \tau\le T'_i\}$. Hence $|u_i|=o(1)$ for
$|x|\le\frac{15}{16}$. Along the radial segment at the same time, using $\|Du_i\|_0\le1$, we get
$|u_i|\le\frac1{16}+o(1)$ for $\frac{15}{16}\le|x|<1$. With $|\pi_PW_i^{-1}x|\ge(1-o(1))|x|$
for $x\in P$ this gives $|f_i(x,\tau)-x_{Y_i}|\ge\varrho\big((1-o(1))|x|-|u_i(x,\tau)|\big)-o(1)$,
which exceeds $\frac34\varrho$ for $i$ large whenever $|x|\ge\frac78$ (for
$\frac78\le|x|\le\frac{15}{16}$ because $|u_i|=o(1)$ there, for $|x|\ge\frac{15}{16}$ because
$\frac{15}{16}-\frac1{16}=\frac78$). Moreover $f_i(\cdot,\tau)$ is injective on the convex
set $\{|x|<\frac{15}{16}\}$ for $i$ large, by the mean value theorem and the uniform
convergence of $D_xf_i$ to $\varrho$ times an isometry.

\emph{The fixed ambient base.} Let
\[
B^{(Y_i)}:=\big\{\xi\in P:\ |\xi-x_{Y_i}|<\tfrac34\varrho\big\}\times\big(s_{Y_i}-\tfrac9{16}\varrho^2,\ s_{Y_i}+\tfrac9{16}\varrho^2\big)\cap\{t\le T_i\}.
\]
Let $\mathrm{pr}_{P,t}\colon\Hn\times\R\to P\times\R$, $\mathrm{pr}_{P,t}(w,z,t):=(\pi_Pw,t)$. The
\emph{ambient base point} of a point $Z$ of the ambient coordinates is $\mathrm{pr}_{P,t}(Z)$.
This is not $\pi\times\mathrm{id}$: the value $\pi(Z)=w$ also records the $J_0P$-component
of $w$. Write $G_i(x,\tau):=\Psi_i(x,u_i(x,\tau),\zeta_i(x,\tau),\tau)$ for
$(x,\tau)\in\Omega_{T'_i}$, so that $\Psi_i(\hat\Gamma_{u_i,\zeta_i})=G_i(\Omega_{T'_i})$. Since
$\pi\circ\Psi_i$ does not depend on $z$ and the time component of $\Psi_i$ is exact,
\eqref{eq:Xi} gives
\[
\mathrm{pr}_{P,t}\circ G_i=\Xi_i\qquad\text{on }\Omega_{T'_i}.
\]
Let $G_i(x,\tau)$ be a point of the graph whose ambient base point $(\xi,t)=\Xi_i(x,\tau)$
lies in $B^{(Y_i)}$. Since the time component of $\Xi_i$ is exact, its local time satisfies
$|\tau|<\frac9{16}$, and then, by the preliminary bounds, $|x|<\frac78$. Over each
$(\xi,t)\in B^{(Y_i)}$ there is at most one such point, by the injectivity of
$f_i(\cdot,\tau)$ on $\{|x|<\frac{15}{16}\}$. Here $T'_i=(T_i-s_{Y_i})/\varrho^2$ (Step 2 of the
proof of Proposition~\ref{prop:compact} at the scale $\varrho$), so that $\tau\le T'_i$
exactly when $t\le T_i$.

By the inverse function theorem, for $i$
large $\Xi_i$ is a diffeomorphism from
\[
D_i:=\{|x|<\tfrac78\}\times(-\tfrac{49}{64},\tfrac{49}{64})\cap\{\tau\le T'_i\}\subset\Omega_{T'_i}
\]
onto a subset of $P\times\R$ that is relatively open in $\{t\le T_i\}$ and contains
$B^{(Y_i)}$ (the image of the ball of radius $\frac78$ under
the limit map is the ball of radius $\frac78\varrho$ about $x_{Y_i}$, and the perturbation
is $o(1)$ in $C^1$). The inverse of this diffeomorphism is $(\xi,t)\mapsto(\phi_i(\xi,t),\tau)$, $\tau=(t-s_{Y_i})/\varrho^2$,
where $f_i(\phi_i(\xi,t),\tau)=\xi$. The map $\phi_i$ is $C^{2,\alpha}$ in the parabolic
sense with bounds depending on $\varrho$ and on the bounds above: for fixed $t$,
$\phi_i(\cdot,t)$ is the inverse of the spatial $C^{2,\alpha}$ diffeomorphism $f_i(\cdot,\tau)$,
whose Jacobian is uniformly invertible, so $D_\xi\phi_i=(D_xf_i)^{-1}\circ(\phi_i,\tau)$ and
$D^2_\xi\phi_i$ are bounded in $C^{0,\alpha}$. Moreover, differentiating $f_i(\phi_i(\xi,t),\tau)=\xi$
in the ambient time,
\begin{equation}\label{eq:phi-t}
\partial_t\phi_i(\xi,t)=-\varrho^{-2}\big[(D_xf_i)^{-1}\partial_\tau f_i\big]\big(\phi_i(\xi,t),\tau\big),
\end{equation}
which is bounded in $C^{0,\alpha}$ because $\partial_\tau f_i$ is (it involves
$\partial_\tau u_i$) and $\phi_i$ is Lipschitz in $(\xi,t)$.

Next,
\[
\Psi_i(\hat\Gamma_{u_i,\zeta_i})\cap\mathrm{pr}_{P,t}^{-1}(B^{(Y_i)})
=G_i\big((\Xi_i|_{D_i})^{-1}(B^{(Y_i)})\big).
\]
Indeed, by the identity $\mathrm{pr}_{P,t}\circ G_i=\Xi_i$, a point of the left-hand side is
$G_i(x,\tau)$ with $(x,\tau)\in\Omega_{T'_i}$ and $\Xi_i(x,\tau)\in B^{(Y_i)}$. By the bounds
$|\tau|<\frac9{16}$ and $|x|<\frac78$ obtained above, $(x,\tau)\in D_i$, so
$(x,\tau)=(\Xi_i|_{D_i})^{-1}(\Xi_i(x,\tau))$. Conversely, $D_i\subset\Omega_{T'_i}$, and the
same identity shows that $G_i$ maps $(\Xi_i|_{D_i})^{-1}(B^{(Y_i)})$ into
$\mathrm{pr}_{P,t}^{-1}(B^{(Y_i)})$. A point $G_i(p)$ with $p\in(\Xi_i|_{D_i})^{-1}(B^{(Y_i)})$ has
ambient base point $\Xi_i(p)$, so this set has exactly one point over each
$(\xi,t)\in B^{(Y_i)}$, namely, in the coordinates $(x,y,z,t)$,
\[
G_i\circ(\Xi_i|_{D_i})^{-1}(\xi,t)=\big(\xi,\ v^y_i(\xi,t),\ v^z_i(\xi,t),\ t\big).
\]
(This uniqueness concerns the single piece $\Psi_i(\hat\Gamma_{u_i,\zeta_i})$; the
single-valuedness of the whole trace is proved below.) Hence this set is the
Legendrian graph, over the ambient base $B^{(Y_i)}$, of the pair $v_i=(v_i^y,v_i^z)$ (the
$J_0P$-component and the height), with
\[
\|v^y_i\|_{2,\alpha;B^{(Y_i)}}\le C,\qquad v^y_i\to0\ \text{in }C^{2,\alpha'},\qquad
v^z_i(\xi,t)\to-4\bar\theta\,t\ \text{uniformly},
\]
where $C$ depends on $\varrho$, $n$, $\bb$ and $C_\theta$ only. The bound holds because
$v^y_i=\pi_{J_0P}\circ\pi\circ\Psi_i(\cdot,u_i,\cdot)\circ(\phi_i,\tau)$ is a composition of
$C^{k_0}$ maps with $u_i$ and $\phi_i$, which are $C^{2,\alpha}$ with uniform bounds (same
exponent). The convergence $v^y_i\to0$ in $C^{2,\alpha'}(B^{(Y_i)})$, for every
$\alpha'<\alpha$, requires care at the moving upper end, since the convergence in Step~2
of the proof of Proposition~\ref{prop:compact} is on compact subsets of $\Omega_{T'}$, which do not
contain the strip $T'<\tau\le T'_i$ when $T'_i>T'$. We prove it in four steps, along the
subsequence of that step.

\emph{(a) The window.} Put $\mathsf W_i:=\{|x|\le\frac{15}{16},\ |\tau|\le\frac{49}{64},\ \tau\le T'_i\}$.
Every graph point over $B^{(Y_i)}$ has chart base point in $D_i\subset \mathsf W_i$, so the
behaviour of $u_i$ on $\mathsf W_i$ controls $v^y_i$ on $B^{(Y_i)}$. If $T'_i\le T'$, then
$\mathsf W_i=\mathsf W\cap\Omega_{T'_i}$ with $\mathsf W:=\{|x|\le\frac{15}{16},\ |\tau|\le\frac{49}{64},\ \tau\le T'\}$
a compact subset of $\Omega_{T'}$, and the convention of Step~2 gives $u_i\to0$ in
$C^{2,\alpha'}(\mathsf W_i)$ directly. If $T'>\frac{49}{64}$, then $T'_i>\frac{49}{64}$ for $i$
large (as $T'_i\to T'$) and $\mathsf W_i=\mathsf W$ is a compact subset of $\Omega_{T'}$, with the same
conclusion. The remaining case is $T'\le\frac{49}{64}$ with $T'_i>T'$, where $\mathsf W_i$ contains
the strip $\{|x|\le\frac{15}{16},\ T'<\tau\le\min\{T'_i,\frac{49}{64}\}\}$.

\emph{(b) The strip.} Let $T'\le\frac{49}{64}$ and $T'_i>T'$ (recall $T'\ge0$, as
$s_Y\le T$ for $Y\in\Pi'$). The norm $\|u_i\|_{2,\alpha;\Omega_{T'_i}}\le1$ of
Section~\ref{subsec:scale} contains the parabolic H\"older seminorms $[Du_i]_\alpha$,
$[D^2u_i]_\alpha$ and $[\partial_\tau u_i]_\alpha$, taken with respect to the parabolic
distance $\|X-Y\|=\max\{|x-y|,|s-t|^{1/2}\}$ of Appendix~\ref{app:schauder}. Hence, for
$|x|\le\frac{15}{16}$ and $T'<\tau\le T'_i$,
\begin{gather*}
|g(x,\tau)-g(x,T')|\le(\tau-T')^{\alpha/2}\le(T'_i-T')^{\alpha/2}\quad\text{for }g\in\{Du_i,D^2u_i,\partial_\tau u_i\},\\
|u_i(x,\tau)-u_i(x,T')|\le T'_i-T' .
\end{gather*}
The set $\mathsf W^0:=\{|x|\le\frac{15}{16},\ |\tau|\le\frac{49}{64},\ \tau\le T'\}$ is a compact
subset of $\Omega_{T'}$ containing the slice $\tau=T'$, so by Step~2 (convergence up to
the upper end) $u_i$, $Du_i$, $D^2u_i$ and $\partial_\tau u_i$ tend to $0$ uniformly on
$\mathsf W^0$, in particular at the time $T'$. With $T'_i-T'\to0$ this gives
\[
\sigma_i:=\sup_{\mathsf W_i}\big(|u_i|+|Du_i|+|D^2u_i|+|\partial_\tau u_i|\big)\longrightarrow0 ,
\]
the supremum being over all of $\mathsf W_i$, strip included. In the two cases of (a) the same
holds directly. No point with $\tau>T'_i$ is used.

\emph{(c) Interpolation.} Let $g$ be a function on $\mathsf W_i$ with $[g]_{\alpha;\mathsf W_i}\le1$ and
$\|g\|_{0;\mathsf W_i}\le\epsilon$, and let $\alpha'<\alpha$. For $X\ne Y$ in $\mathsf W_i$: if
$\|X-Y\|^\alpha\le2\epsilon$, then
$|g(X)-g(Y)|\le\|X-Y\|^\alpha=\|X-Y\|^{\alpha'}\|X-Y\|^{\alpha-\alpha'}\le\|X-Y\|^{\alpha'}(2\epsilon)^{1-\alpha'/\alpha}$;
otherwise
$|g(X)-g(Y)|\le2\epsilon=(2\epsilon)^{\alpha'/\alpha}(2\epsilon)^{1-\alpha'/\alpha}<\|X-Y\|^{\alpha'}(2\epsilon)^{1-\alpha'/\alpha}$.
Hence $[g]_{\alpha';\mathsf W_i}\le(2\epsilon)^{1-\alpha'/\alpha}$. Applied to
$g=Du_i,D^2u_i,\partial_\tau u_i$ (each with $[g]_\alpha\le1$ by $\|u_i\|_{2,\alpha}\le1$)
and $\epsilon=\sigma_i$, this gives
\[
\|u_i\|_{2,\alpha';\mathsf W_i}\le3\sigma_i+3(2\sigma_i)^{1-\alpha'/\alpha}\longrightarrow0 ,
\]
uniformly on all of $\mathsf W_i$, including the strip, with constants independent of $i$.

\emph{(d) Composition.} By \eqref{eq:Xi}, $f_i$ depends on $(x,\tau)$ only through
$(x,u_i(x,\tau))$, with the linear map $\varrho\,\pi_PW_i^{-1}$ and the map
$\pi_P\pi E_i$, $\|E_i\|_{C^{k_0}(\bB)}\to0$. By (c) and $W_i\to W$, on $\mathsf W_i$ we have
$D_xf_i=\varrho\,\pi_PW_i^{-1}(I,Du_i)+\pi_P\pi\,DE_i(x,u_i)(I,Du_i)\to L_0:=\varrho\,\pi_PW^{-1}|_P$
in $C^{1,\alpha'}$, $D^2_xf_i\to0$ in $C^{0,\alpha'}$ and
$\partial_\tau f_i=\varrho\,\pi_PW_i^{-1}(0,\partial_\tau u_i)+\pi_P\pi\,DE_i(x,u_i)(0,\partial_\tau u_i)\to0$
in $C^{0,\alpha'}$, where $L_0$ is $\varrho$ times an isometry of $P$. Now consider
$\phi_i$ on $B^{(Y_i)}$, where $(\phi_i(\xi,t),\tau)\in D_i\subset \mathsf W_i$. The map
$(\xi,t)\mapsto(\phi_i(\xi,t),\tau)$, $\tau=(t-s_{Y_i})/\varrho^2$, is Lipschitz for the
parabolic distances with a constant independent of $i$:
$|\tau-\tau'|^{1/2}=\varrho^{-1}|t-t'|^{1/2}$, and
$|\phi_i(\xi,t)-\phi_i(\xi',t')|\le\|D_\xi\phi_i\|_0|\xi-\xi'|+\|\partial_t\phi_i\|_0|t-t'|\le C\|(\xi,t)-(\xi',t')\|$
by the bounds on $D_\xi\phi_i$ and $\partial_t\phi_i$ above and $|t-t'|\le\frac98\varrho^2$
on $B^{(Y_i)}$. Hence composition with $(\phi_i,\tau)$ maps functions converging in
$C^{0,\alpha'}(\mathsf W_i)$ to functions converging in $C^{0,\alpha'}(B^{(Y_i)})$, and
$C^{0,\alpha'}$ is an algebra under products. Differentiating
$f_i(\phi_i(\xi,t),\tau)=\xi$ twice in $\xi$ and once in $t$,
\[
D_\xi\phi_i=(D_xf_i)^{-1}\circ(\phi_i,\tau),\qquad
D^2_\xi\phi_i=-\Big[(D_xf_i)^{-1}D^2_xf_i\big((D_xf_i)^{-1}\cdot,(D_xf_i)^{-1}\cdot\big)\Big]\circ(\phi_i,\tau),
\]
and $\partial_t\phi_i$ is \eqref{eq:phi-t}, with the factor $\varrho^{-2}=d\tau/dt$. Since
$D_xf_i\to L_0$ in $C^{0,\alpha'}(\mathsf W_i)$ with uniformly invertible values,
$(D_xf_i)^{-1}\to L_0^{-1}$ in $C^{0,\alpha'}(\mathsf W_i)$. Hence $D_\xi\phi_i\to L_0^{-1}$,
$D^2_\xi\phi_i\to0$ and $\partial_t\phi_i\to0$ in $C^{0,\alpha'}(B^{(Y_i)})$. Moreover,
$\phi_i-\ell^{\rm aff}_i\to0$ uniformly for the affine maps $\ell^{\rm aff}_i(\xi):=L_0^{-1}(\xi-x_{Y_i})$,
because $\xi=f_i(\phi_i,\tau)=x_{Y_i}+L_0\phi_i+o(1)$ uniformly. Thus $\phi_i-\ell^{\rm aff}_i\to0$
in $C^{2,\alpha'}(B^{(Y_i)})$. Next, $u_i\circ(\phi_i,\tau)\to0$ in
$C^{2,\alpha'}(B^{(Y_i)})$: its derivatives are
$D_\xi[u_i\circ(\phi_i,\tau)]=(Du_i\circ(\phi_i,\tau))D_\xi\phi_i$,
$D^2_\xi[u_i\circ(\phi_i,\tau)]=(D^2u_i\circ(\phi_i,\tau))[D_\xi\phi_i,D_\xi\phi_i]+(Du_i\circ(\phi_i,\tau))D^2_\xi\phi_i$
and
$\partial_t[u_i\circ(\phi_i,\tau)]=(Du_i\circ(\phi_i,\tau))\partial_t\phi_i+\varrho^{-2}(\partial_\tau u_i)\circ(\phi_i,\tau)$,
products of factors bounded in $C^{0,\alpha'}$ in which at least one factor tends to $0$
in $C^{0,\alpha'}$ by (c). Finally, as $\pi\circ\Psi_i$ does not depend on $z$,
\[
v^y_i=y_{Y_i}+\varrho\,\pi_{J_0P}W_i^{-1}\big(\phi_i,\ u_i\circ(\phi_i,\tau)\big)+\pi_{J_0P}\pi E_i\big(\phi_i,\ u_i\circ(\phi_i,\tau)\big),
\]
where $y_{Y_i}\to y_Y=0$ ($Y\in\Pi$); $\pi_{J_0P}W_i^{-1}(\phi_i,0)\to0$ in $C^{2,\alpha'}$,
because $\pi_{J_0P}W_i^{-1}|_P\to\pi_{J_0P}W^{-1}|_P=0$ ($W^{-1}P=P$) and $\phi_i$ is
bounded in $C^{2,\alpha'}$; $\pi_{J_0P}W_i^{-1}(0,u_i\circ(\phi_i,\tau))\to0$ in
$C^{2,\alpha'}$ by the composition estimate just proved; and the $E_i$ term tends to $0$ in $C^{2,\alpha'}$ by
the chain rule, since $E_i\to0$ in $C^{k_0}$ and the inner maps are bounded in
$C^{2,\alpha'}$. This proves $v^y_i\to0$ in $C^{2,\alpha'}(B^{(Y_i)})$, on the whole base,
including the points with $T'<\tau\le T'_i$; the flow is never extended beyond $T'_i$.

The height $v^z_i$ involves $\zeta_i$, of which only uniform convergence is used. On
the strip $T'<\tau\le T'_i$ the uniform convergence of $u_i$ and of
$\zeta_i+4\bar\theta\tau$ follows from $|\partial_\tau u_i|\le1$,
$|\partial_\tau\zeta_i|\le1+4C_\theta$ and $T'_i\to T'$: the values at the time $T'$
converge, and the increments are at most $(2+8C_\theta)(T'_i-T')\to0$.

\emph{The comparison plane and the inner box.} Let
$\Pi_{\rm ext}:=\{(w,-4\bar\theta t,t):\ w\in P,\ t\in\R\}$, the untruncated lift of $P$.
It is used only as a set with which $\cM_i$ is compared, and nothing is claimed about
$\cM$ after the time $T$. By the preliminary bounds, every graph point over $B^{(Y_i)}$
comes from a chart base point with $|x|<\frac78$, $|\tau|<\frac9{16}$, so it lies in the
inner box
\[
\bB^{\rm in}:=\{|x|<\tfrac{15}{16},\ |y|<2,\ |z|<\tfrac{\bb^2}2\}\times(-\tfrac{49}{64},\tfrac{49}{64})\subset\bB
\]
of the $Y_i$-chart. We locate the plane piece $\Pi_{\rm ext}\cap\mathrm{pr}_{P,t}^{-1}(B^{(Y_i)})$
in the $Y_i$-chart. (On $\Pi_{\rm ext}$ we have $w\in P$, so $\mathrm{pr}_{P,t}$ agrees there with
$(w,z,t)\mapsto(w,t)$, and the plane piece is the set of points of $\Pi_{\rm ext}$ whose
pair $(w,t)$ lies in $B^{(Y_i)}$.) In the limit chart $\Phi=\tilde W\Lambda_{Y,\varrho}$ the set
$\Pi_{\rm ext}$ is $\{(x,0,-4\bar\theta\tau,\tau)\}$: $Y\in\Pi$ gives $\ell_Y^{-1}(\Pi_{\rm ext})=\Pi_{\rm ext}$
(the cross term of \eqref{eq:vertical} vanishes on $P$ as $y_Y=0$), $\tilde W$ preserves
$P$, and $D_{1/\varrho}$ maps $(w,-4\bar\theta t,t)$ to $(w/\varrho,-4\bar\theta t/\varrho^2,t/\varrho^2)$.
Now $\Phi_i\to\Phi$ in $C^{k_0}$ on compact sets, and the point of $\Pi_{\rm ext}$
above $(\xi,t)\in B^{(Y_i)}$ is within $o(1)$ (ambient metric $\mathsf m$) of the graph
point above $(\xi,t)$, by $v^y_i\to0$ and $v^z_i+4\bar\theta t\to0$. Hence the $\Phi_i$-coordinates
of the plane piece over $B^{(Y_i)}$ lie within $o(1)$ of the compact set
$\{(x,0,-4\bar\theta\tau,\tau):\ |x|\le\frac78,\ |\tau|\le\frac9{16}\}$. This set has Euclidean
distance at least $\min\{\frac1{16},\ 2,\ \frac{\bb^2}2-\frac94C_\theta,\ \frac{49}{64}-\frac9{16}\}=\frac1{16}$
from the complement of $\bB^{\rm in}$ (the vertical margin is
$\frac{\bb^2}2-\frac94C_\theta=\frac74C_\theta+2\ge2$ by $\bb^2=8C_\theta+4$). Hence, for $i$
large, every point of the $Y_i$-chart within distance $\frac1{32}$ of the chart image of
the plane piece lies in $\bB^{\rm in}$.

Let $G:=\{|w|<5,\ |z|<17\bb^2\}\times(-17,17)$, a fixed region of the ambient
coordinates with $\overline G\subset\cB_8\times\R$. It contains the $1$-neighbourhood
(metric $\mathsf m$) of $\overline{\cC_4(0)}$, hence of every plane piece above (whose
points have $|\xi|<3$, $|t|<4$, $|z|\le16C_\theta$). On $G$ the chart
$\Phi_i=D_{1/\varrho}\circ(\tilde W_i\ell_{Y_i}^{-1}+\tilde U_iE'_i)$ (composed with the time
shift) is defined, is injective for $i$ large (a $C^1$-small perturbation of an affine
bijection on a convex set), and is Lipschitz for $\mathsf m$ with constant $5/\varrho^2$ for $i$ large. Indeed, $D_{1/\varrho}$ has Lipschitz constant $\varrho^{-2}$, the affine map
$\ell_{Y_i}^{-1}$ has Lipschitz constant at most $1+|w_{Y_i}|\le3+o(1)$ by
\eqref{eq:vertical}, $\tilde W_i$ is an isometry, and $\|E'_i\|_{C^1(\cB_8)}\to0$ by
Corollary~\ref{cor:rM}(b). Put
\[
\epsilon_0:=\frac{\varrho^2}{160} .
\]
If $Z\in\Hn\times\R$ is within $\mathsf m$-distance $\epsilon_0$ of the plane piece over
$B^{(Y_i)}$, then $Z\in G$, so $\Phi_i(Z)$ is within $5\epsilon_0/\varrho^2=\frac1{32}$ of the
chart image of the plane piece. Hence $\Phi_i(Z)\in\bB^{\rm in}$ and
$Z=\Psi_i(\Phi_i(Z))\in\Psi_i(\bB^{\rm in})$. Thus, for $i$ large, $\Psi_i(\bB^{\rm in})$
contains the $\epsilon_0$-neighbourhood of $\Pi_{\rm ext}\cap\mathrm{pr}_{P,t}^{-1}(B^{(Y_i)})$.

\emph{The cover.} For $Y\in\Pi'$ with ambient base point $(x_Y,s_Y)$ let
$B^{(Y)}_{1/4}:=\{|\xi-x_Y|<\frac3{16}\varrho\}\times(s_Y-\frac9{64}\varrho^2,s_Y+\frac9{64}\varrho^2)$
be the quartered base. The compact region
\[
K:=\{|\xi|\le\tfrac32+\tfrac{\varrho}{8}\}\times\big[-\tfrac94-\tfrac{\varrho^2}8,\ \min\{\tfrac94,T\}+\tfrac{\varrho^2}8\big]
\]
is covered by the open sets $B^{(Y)}_{1/4}$, $Y\in\Pi'$. Indeed, for $(\xi,t)\in K$ put
$t_0:=\min\{t,T\}$ and $Y:=(\xi,0,-4\bar\theta t_0,t_0)$. Then $Y\in\Pi'$ ($|\xi|<2$,
$|t_0|<3$) and $|t-t_0|\le\frac{\varrho^2}8<\frac9{64}\varrho^2$, so $(\xi,t)\in B^{(Y)}_{1/4}$. When $T=\infty$, $t_0=t$. The points of $K$ with $t>T$ are not base points of $\Pi'$, which is truncated at $T$, but they are covered by the quartered bases of points of $\Pi'$ at the time $T$. Choose a finite subcover by $B^{(Y^{(m)})}_{1/4}$, $m=1,\dots,N$. This subcover depends
only on $\Pi'$, $\varrho$ and $T$. Let $Y^{(m)}_i\in\cM_i$ with $Y^{(m)}_i\to Y^{(m)}$. If the
centres move by less than $\frac{3}{16}\varrho$ in $\xi$ and $\frac{9}{64}\varrho^2$ in $t$,
then the halved bases
$B^{(Y^{(m)}_i)}_{1/2}:=\{|\xi-x_{Y^{(m)}_i}|<\frac38\varrho\}\times(s_{Y^{(m)}_i}-\frac9{32}\varrho^2,s_{Y^{(m)}_i}+\frac9{32}\varrho^2)$
contain the quartered bases $B^{(Y^{(m)})}_{1/4}$, hence cover $K$. Moreover
$\Omega^{3/2}_{T_i}\subset K$ for $i$ large: its points have $|\xi|<\frac32$, $|t|<\frac94$ and
$t\le T_i$. If $T\ge\frac94$ (including $T=\infty$) this is immediate, and if $T<\frac94$ it
holds once $T_i<T+\frac{\varrho^2}{16}$, which includes the strip $T<t\le T_i$ when
$T_i>T$. Let $i$ be so large that all the above holds for every $m$ and that every point
of $\cM_i\cap\overline{\bB^{3/2}}$ is within $\mathsf m$-distance $\epsilon$ of $\Pi'$, with $2(1+4C_\theta)\epsilon<\epsilon_0$. The last condition holds for $i$ large by local set
convergence in $\cU$: otherwise a subsequence of points of $\cM_i\cap\overline{\bB^{3/2}}$
at distance $\ge\epsilon$ from $\Pi$ would accumulate at a point of the compact set
$\overline{\bB^{3/2}}\subset\cC_4(0)\subset\cU$ not in $\cM$, contradicting condition (b)
of local convergence (and points of $\Pi$ within
$\epsilon$ of $\overline{\bB^{3/2}}$ lie in $\Pi'$).

\emph{Single-valuedness.} Let $Z,Z'\in\cM_i\cap\overline{\bB^{3/2}}$ have the same
ambient base point $(\xi,t)$, $|\xi|<\frac32$. Each is within $\epsilon$ of some point
$(\xi_1,0,-4\bar\theta t_1,t_1)\in\Pi'$, so $|\xi-\xi_1|,|t-t_1|\le\epsilon$ and the height of
$Z$ differs from $-4\bar\theta t$ by at most $\epsilon+4C_\theta\epsilon$. Thus $Z$ and $Z'$
are within $2(1+4C_\theta)\epsilon<\epsilon_0$ of the point of $\Pi_{\rm ext}$ above
$(\xi,t)$ (this point lies in $\Pi$ unless $T<t\le T_i$; the comparison uses only the set
$\Pi_{\rm ext}$). Let $m$ be such that $(\xi,t)$ lies in the halved base
$B^{(Y^{(m)}_i)}_{1/2}$. Here $t\le T_i$, because $(\xi,t)$ is the ambient base point of
$Z\in\cM_i$ and $\cN_i\subset\Sigma_i\times(-\infty,T_i]$. Hence
$(\xi,t)\in B^{(Y^{(m)}_i)}_{1/2}\cap\{t\le T_i\}\subset B^{(Y^{(m)}_i)}$, by
$\frac38\varrho<\frac34\varrho$ and $\frac9{32}\varrho^2<\frac9{16}\varrho^2$. (The halved
bases are not truncated at $T_i$, since they must cover the part of $K$ above $T_i$; the
full bases are, and the inclusion holds after intersecting with $\{t\le T_i\}$.)
Then $Z,Z'\in\Psi^{(m)}_i(\bB^{\rm in})\subset\Psi^{(m)}_i(\bB)$.
By Definition~\ref{def:K}(i), within this box the trace is contained in
$\Psi^{(m)}_i(\hat\Gamma_{u^{(m)}_i,\zeta^{(m)}_i})$, whose part in
$\mathrm{pr}_{P,t}^{-1}(B^{(Y^{(m)}_i)})$ is the graph of $v^{(m)}_i$: one point above each
ambient base point. Hence $Z=Z'$.

\emph{Full base, bounds and injectivity.} Every $(\xi,t)\in\Omega^{3/2}_{T_i}$ lies in
some $B^{(Y^{(m)}_i)}_{1/2}$ (as $\Omega^{3/2}_{T_i}\subset K$), and $t\le T_i$ there, so
$(\xi,t)\in B^{(Y^{(m)}_i)}_{1/2}\cap\{t\le T_i\}\subset B^{(Y^{(m)}_i)}$. Over this base the
graph of $v^{(m)}_i$ is realised by the trace (Lemma~\ref{lem:filling} in the
$Y^{(m)}_i$-chart, transported by $\Psi^{(m)}_i$). Hence $\cM_i\cap\bB^{3/2}$ is the
Legendrian graph of the single-valued pair
$(\tilde u_i,\tilde\zeta_i):=(v^{y,(m)}_i,v^{z,(m)}_i)$ on
$B^{(Y^{(m)}_i)}_{1/2}\cap\Omega^{3/2}_{T_i}\subset B^{(Y^{(m)}_i)}$ (the pieces agree on
overlaps, being the same set, and these sets cover $\Omega^{3/2}_{T_i}$). This graph lies in
$\bB^{3/2}$ because
$|v^y_i|\to0$ and $|v^z_i+4\bar\theta t|\to0$ with $4C_\theta\cdot\frac94<\frac94\bb^2$.

The uniform $C^{2,\alpha}$ bounds and the convergence hold piecewise, hence on $\Omega^{3/2}_{T_i}$. Indeed, two points of $\Omega^{3/2}_{T_i}$ at parabolic distance less than
$\frac38\varrho$ lie in a common full base $B^{(Y^{(m)}_i)}$ (take $m$ with the first point
in the halved base; $\frac38\varrho+\frac38\varrho\le\frac34\varrho$ and
$\frac9{32}\varrho^2+\frac9{64}\varrho^2\le\frac9{16}\varrho^2$), so their H\"older quotients are
bounded by $C$. For two points at parabolic distance $\ge\frac38\varrho$, the H\"older
quotient of each $g\in\{D\tilde u_i,D^2\tilde u_i,\partial_t\tilde u_i\}$ is at most
$2\|g\|_{0;\Omega^{3/2}_{T_i}}(\frac38\varrho)^{-\alpha}$, the sup norm being the largest of
those on the pieces. For the convergence the same two estimates are used with $\alpha'$ in
place of $\alpha$, since $v^{y,(m)}_i\to0$ in $C^{2,\alpha'}(B^{(Y^{(m)}_i)})$ for each of the
finitely many $m$.

Injectivity: two parameter points with the same image
$Z\in\bB^{3/2}$ both lie in the preimage of a box $\Psi^{(m)}_i(\bB)$ containing $Z$, on
which $\hat F_i$ is injective by Definition~\ref{def:K}(ii).

\emph{Quantifiers.} The construction extracts subsequences (in Step 2 of
Proposition~\ref{prop:compact} at each of the finitely many centres $Y^{(m)}$, and for the
rotations $W^{(m)}_i$), while the lemma asserts its conclusion for all large $i$ of the
original sequence. The cover $\{B^{(Y^{(m)})}_{1/4}\}$ and the constants $\varrho$,
$\epsilon_0$, $\epsilon$ and $C$ depend only on $\Pi'$, $T$, $C_*$, $n$, $\bb$, $C_\theta$ and the
structural bounds, not on the subsequence. If the one-sheet statement, the bound
$\|\tilde u_i\|_{2,\alpha}\le C$ or the injectivity failed along a subsequence, the same
construction applied to that subsequence (which satisfies the same hypotheses and, since local
convergence passes to subsequences, has
the same limit $\cM$) would give, along a further subsequence, the conclusion with the
same $C$, a contradiction. The convergence statements follow in the same way: every
subsequence has a further subsequence along which $\tilde u_i\to0$ in $C^{2,\alpha'}$ and
$\tilde\zeta_i\to-4\bar\theta t$ uniformly on compact subsets (up to the upper end), and the
limits are determined by $\cM\cap\cC_4(0)=\Pi\cap\cC_4(0)$. Hence the whole sequence converges.
\end{proof}

\Needspace{6\baselineskip}
\subsection{The Schauder improvement}\label{subsec:schauder}

\begin{lemma}[Schauder improvement]\label{lem:schauder}
Let $u_i$ be solutions of \eqref{eq:grapheq} for structures $\varsigma_i$ with
$\delta_i\to0$ in \eqref{eq:delta}, on $\Omega^{3/2}_{T'_i}=\{|x|<\frac32,\ |s|<\frac94,\ s\le T'_i\}$,
$T'_i\to T'\in[0,\infty]$, with $\|u_i\|_{2,\alpha;\Omega^{3/2}_{T'_i}}\le C_0$ and
$u_i\to0$ in $C^2$ on compact subsets of $\Omega^{3/2}_{T'}$ (one-sided at $s=T'$; when
$T'_i<T'$, convergence on a compact subset $\cK'$ of $\Omega^{3/2}_{T'}$ means convergence
on $\cK'\cap\Omega^{3/2}_{T'_i}$, the convention of Step~2 of the proof of
Proposition~\ref{prop:compact}). Then,
for every $\mu\in[1,\frac54)$, $\|u_i\|_{2,\alpha;\Omega^{\mu}_{T'_i}}\to0$, where
$\Omega^\mu_{T'_i}=\{|x|<\mu,\ |s|<\mu^2,\ s\le T'_i\}$; in particular on the unit base
domain.
\end{lemma}
\begin{proof}
Write the system as $(\partial_s-4\Delta_E)u_i^\gamma=f_i^\gamma$ with
\[
f_i^\gamma:=\big(a_i^{kl}(x,u_i,Du_i)-4\delta^{kl}\big)D_{kl}u_i^\gamma+b_i^\gamma(x,u_i,Du_i)
=\big(4((I+Du_i^TDu_i)^{-1})_{kl}-4\delta^{kl}+e_i^{kl}\big)D_{kl}u_i^\gamma+b^\gamma_i .
\]
Let $\Omega^{5/4}:=\{|x|<\frac54,|s|<\frac{25}{16},s\le T'_i\}$. First, $u_i$ and $Du_i$
tend to $0$ uniformly on $\Omega^{5/4}$. For $0<\epsilon<\frac18$ let
$K_\epsilon:=\{|x|\le\frac54-\epsilon,\ -\frac{25}{16}+\epsilon^2\le s\le\min\{T',\frac{25}{16}\}-\epsilon^2\}$,
a compact subset of $\Omega^{3/2}_{T'}$. Every point of $\Omega^{5/4}$ is joined to a point
of $K_\epsilon$ by a path in $\Omega^{3/2}_{T'_i}$ (radial in $x$, then in $s$) of spatial
length at most $\epsilon$ and time length at most $\epsilon^2+(T'_i-T')_+$. (The convention in
the statement gives the meaning of $\sup_{K_\epsilon}$ when $T'_i<T'$. When $T'=\infty$,
$(T'_i-T')_+=0$ and $\Omega^{5/4}$ is not truncated for $i$ large, as $T'_i>\frac{25}{16}$.
When $T'_i<0$, the paths stay in $\{s\le T'_i\}$: every point of $\Omega^{5/4}$ then
has $s\le T'_i$, and $K_\epsilon\cap\{s\le T'_i\}$ is reached by a time move of length
less than $\epsilon^2$ for $i$ large.) Hence the bounds on
$\|Du_i\|_0$, $\|\partial_su_i\|_0$ and $[Du_i]_\alpha$ give
$\|u_i\|_{0;\Omega^{5/4}}\le\sup_{K_\epsilon}|u_i|+C_0(\epsilon+\epsilon^2+(T'_i-T')_+)$ and
$\|Du_i\|_{0;\Omega^{5/4}}\le\sup_{K_\epsilon}|Du_i|+C_0(\epsilon+(T'_i-T')_+^{1/2})^\alpha$. Let first $i\to\infty$, then $\epsilon\to0$.

Next, $[Du_i]_{\alpha;\Omega^{5/4}}\to0$. Let
$X,Y\in\Omega^{5/4}$ with $\|X-Y\|<\frac18$. The later point $Z$ has
$P_{1/4}(Z)\subset\Omega^{3/2}_{T'_i}$, and the argument of Lemma~\ref{lem:interp} (the
spatial mean value theorem and the time-difference argument proving \eqref{eq:interp2})
gives $|Du_i(X)-Du_i(Y)|\le C(n)(\|D^2u_i\|_0+\|\partial_su_i\|_0)\|X-Y\|\le C(n)C_0\|X-Y\|$.
Combined with $|Du_i(X)-Du_i(Y)|\le2\|Du_i\|_0$, this yields
$|Du_i(X)-Du_i(Y)|\le(C(n)C_0)^\alpha(2\|Du_i\|_{0;\Omega^{5/4}})^{1-\alpha}\|X-Y\|^\alpha$.
Pairs with $\|X-Y\|\ge\frac18$ satisfy $|Du_i(X)-Du_i(Y)|\le2\cdot8^\alpha\|Du_i\|_0\|X-Y\|^\alpha$.
Hence $[Du_i]_{\alpha;\Omega^{5/4}}\le C(n,\alpha)C_0^\alpha\|Du_i\|_{0;\Omega^{5/4}}^{1-\alpha}\to0$,
so $Du_i\to0$ in $C^{0,\alpha}(\Omega^{5/4})$. The function
$p\mapsto4((I+p^Tp)^{-1})-4I$ is smooth and vanishes at $p=0$, so its composition with
$Du_i$ tends to $0$ in $C^{0,\alpha}$. Moreover, $\|e_i\|_{0,\alpha}+\|b_i\|_{0,\alpha}\le C\delta_i\to0$
by \eqref{eq:coeff}. As $\|D^2u_i\|_{0,\alpha}\le C_0$, we get $\|f_i\|_{0,\alpha;\Omega^{5/4}}\to0$.

\emph{Zero-order finiteness.} The lemma does not assume $u_i(0,0)=0$, and $(0,0)$ does
not lie in $\Omega^{5/4}$ when $T'_i<0$. Take the base point $Z_i:=(0,\min\{0,T'_i\})$.
It lies in $\Omega^{5/4}$ for $i$ large, because $\min\{0,T'_i\}>-\frac{25}{16}$ for $i$
large. Indeed, for the indices with $T'_i\ge0$ the base point is $(0,0)$, whatever the size of
$T'_i$ (also for $T'_i>\frac{25}{16}$ or $T'_i=\infty$), and if infinitely many $T'_i$ are
negative, then $T'=0$, as $T'_i\to T'\ge0$, so these $T'_i$ tend to $0$. Every
$(x,s)\in\Omega^{5/4}$ has $|x|<\frac54$ and $|s-\min\{0,T'_i\}|<\frac{25}{16}$, and the
path from $Z_i$ to $(0,s)$ to $(x,s)$ stays in $\Omega^{5/4}$ (its time interval
$(-\frac{25}{16},\frac{25}{16})\cap(-\infty,T'_i]$ is convex), so
$\|u_i\|_{0;\Omega^{5/4}}\le|u_i(Z_i)|+\frac54\|Du_i\|_0+\frac{25}{16}\|\partial_su_i\|_0<\infty$.
Moreover, $\|u_i\|_{0;\Omega^{5/4}}\to0$ by the first paragraph. Thus $u_i$ has finite norm
$|u_i|_{2,\alpha;\Omega^{5/4}}$ in the sense of Appendix~\ref{app:schauder}.

\emph{The coefficients.} We apply Theorem~\ref{thm:schauder-var} componentwise ($N=n$)
to the operator $\partial_s-4\Delta_E$, with the constant coefficients $a=4I$, $c=d=0$,
$\lambda_0=1$ and $\Lambda=4n$. Indeed, the norm $\|a\|_{0,\alpha;Q_1}$ of the appendix is the sum
of the norms of the components of $a$, which is $4n$ for $a=4I$ (the H\"older seminorms
vanishing), and $a^{ij}\xi_i\xi_j=4|\xi|^2\ge|\xi|^2$. The constant of the theorem then
depends only on $n$, $\alpha$ and the radii.

We split the sequence into the indices with $T'_i\ge0$ and those with $T'_i<0$ (if one
of the two sets is finite there is nothing to prove for it), and show
$\|u_i\|_{2,\alpha;\Omega^\mu_{T'_i}}\to0$ along each. The whole sequence then
converges to $0$.

\emph{Indices with $T'_i\ge0$.} The domains $\Omega^{\mu}_{T'_i}\subset\Omega^{5/4}_{T'_i}$
are the cylinders $Q_\mu\subset Q_{5/4}$ of Appendix~\ref{app:schauder} with
$T=T'_i\in[0,\infty]$, the range in which the appendix is stated, and $u_i$ has finite
norm $|u_i|_{2,\alpha;Q_{5/4}}$ there. Theorem~\ref{thm:schauder-var} with the outer
domain $Q_{5/4}$ and the inner domain $Q_\mu$ (through its dilation clause) gives
\[
\|u_i\|_{2,\alpha;\Omega^{\mu}_{T'_i}}\le c\big(\|u_i\|_{0;\Omega^{5/4}_{T'_i}}+\|f_i\|_{0,\alpha;\Omega^{5/4}_{T'_i}}\big)\longrightarrow0
\]
with $c$ depending only on $n$, $\alpha$ and $\mu$, and not on the upper end $T'_i$, which
is what makes the one-sided case work.

\emph{Indices with $T'_i<0$.} If there are infinitely many of them, then $T'=0$, as
$T'_i\to T'\ge0$; so assume $T'=0$. Fix $\nu\in(\mu,\frac54)$, say
$\nu:=\frac12(\mu+\frac54)$, a function of $\mu$ (the letter $\nu$ is local to this
proof), and consider the indices with $T'_i<0$ and $|T'_i|<\frac{25}{16}-\nu^2$ (all but
finitely many of them, as $T'_i\to0$). Put
\[
v_i(x,\tau):=u_i(x,\tau+T'_i),\qquad g_i(x,\tau):=f_i(x,\tau+T'_i),
\]
wherever $(x,\tau+T'_i)\in\Omega^{5/4}_{T'_i}$. Two inclusions:
\begin{enumerate}[(1)]
\item The cylinder $Q_\nu=\{|x|<\nu,\ -\nu^2<\tau\le0\}$ of the appendix with $T=0$ is
mapped by $(x,\tau)\mapsto(x,\tau+T'_i)$ into $\Omega^{5/4}_{T'_i}$: for $-\nu^2<\tau\le0$
the time $s=\tau+T'_i$ satisfies $s\le T'_i$ and
$s>-\nu^2+T'_i>-\nu^2-(\frac{25}{16}-\nu^2)=-\frac{25}{16}$, while $|x|<\nu<\frac54$. Hence
$v_i\in C^{2,\alpha}(Q_\nu)$ with $|v_i|_{2,\alpha;Q_\nu}\le|u_i|_{2,\alpha;\Omega^{5/4}}<\infty$,
$\|v_i\|_{0;Q_\nu}\le\|u_i\|_{0;\Omega^{5/4}}\to0$ and
$\|g_i\|_{0,\alpha;Q_\nu}\le\|f_i\|_{0,\alpha;\Omega^{5/4}}\to0$, the time shift
preserving sup norms and parabolic H\"older seminorms.
\item $\Omega^\mu_{T'_i}=\{|x|<\mu,\ -\mu^2<s\le T'_i\}$ (as $T'_i<0<\mu^2$) is mapped by
$(x,s)\mapsto(x,s-T'_i)$ into the cylinder $Q_\mu=\{|x|<\mu,\ -\mu^2<\tau\le0\}$ of the
appendix with $T=0$: $\tau=s-T'_i\in(-\mu^2-T'_i,0]\subset(-\mu^2,0]$, since $-T'_i>0$.
Hence $\|u_i\|_{2,\alpha;\Omega^\mu_{T'_i}}\le|v_i|_{2,\alpha;Q_\mu}$.
\end{enumerate}
The time shift commutes with the constant-coefficient operator, so
$(\partial_\tau-4\Delta_E)v_i=g_i$ on $Q_\nu$. Theorem~\ref{thm:schauder-var} in its form
with $T=0$, applied through its dilation clause with the outer domain $Q_\nu$ and the
inner domain $Q_\mu$ (inner radius $\mu$, outer radius $\nu$; the same coefficients
$a=4I$, $\lambda_0=1$, $\Lambda=4n$, $c=d=0$ as above), gives
\[
|v_i|_{2,\alpha;Q_\mu}\le c'\big(\|v_i\|_{0;Q_\nu}+\|g_i\|_{0,\alpha;Q_\nu}\big)\longrightarrow0
\]
with $c'=c'(n,\alpha,\mu,\nu)=c'(n,\alpha,\mu)$, independent of $i$. By (2),
$\|u_i\|_{2,\alpha;\Omega^\mu_{T'_i}}\to0$ along these indices as well. No point with
$s>T'_i$, and no value of $u_i$ outside its domain, is used. The lemma is proved as
stated, without the additional hypothesis $T'_i\ge0$.
\end{proof}

\section{Model flows: growth, monotonicity and rigidity}\label{sec:model}

In this section the ambient space is the Heisenberg group with its model structure, and
the flows satisfy \eqref{eq:F} with $K=-2$. The results are the model inputs of the
blow-up argument. We prove a growth bound for flows with bounded regularity scale
(Lemma~\ref{lem:growth}), and the monotonicity identity of \cite[Proposition~3.2]{PaperI}
under a polynomial growth hypothesis in place of the area ratio bound
(Proposition~\ref{prop:mono-model}). We then prove the lower bound of the density at
points of the flow (Lemma~\ref{lem:density-point}), its lower semicontinuity under
convergence with moving centres (Lemma~\ref{lem:lsc}), and the rigidity of ancient flows
whose regularity scale is bounded on the whole flow and whose density is at most one at
every point and every scale (Proposition~\ref{prop:rigidity}). The rigidity rests on the equality case of the monotonicity identity, as outlined after
Theorem~\ref{thm:A}.

\begin{definition}\label{def:modelflow}
A \emph{model flow on $(t_0,T_*]$}, $-\infty\le t_0<T_*\le\infty$, is a smooth local
Legendrian flow (Definition~\ref{def:localflow}) for the model structure in
$\cU=\Hn\times(t_0,\infty)$ with upper end $T_*$ and $K=-2$; it is \emph{ancient} if
$t_0=-\infty$. It has \emph{uniform polynomial growth} if for every compact
$J\subset(t_0,T_*]$ there are $A_J$ and $m_J$ with
\[
\mu_t\big(F_t^{-1}(\cB_R)\big)\le A_J(1+R)^{m_J}\qquad\text{for all }t\in J,\ R\ge1 ,
\]
where $\cB_R=\{|w|<R,\ |z|<\bb^2R^2\}$ as in Section~\ref{subsec:rM}. A \emph{model flow
in trace form on $(t_0,T_*]$} is a smooth Legendrian flow in trace form
(Definition~\ref{def:traceflow}) in the same $\cU=\Hn\times(t_0,\infty)$, with upper end
$T_*$ and $K=-2$; for it, $\mu_t$ is the area measure of the embedded slice $\cM(t)$
(multiplicity one, Definition~\ref{def:traceflow}), and $\mu_t(F_t^{-1}(\cB_R))$ above and
below stands for $\mu_t(\cM(t)\cap\cB_R)$.
\end{definition}

For a model flow the slices $F_t\colon\cN_t\to\Hn$ are properly immersed, since
$\hat F$ is proper. For a model flow in trace form the slices $\cM(t)$, $t_0<t\le T_*$,
are embedded Legendrian submanifolds which are closed subsets of $\Hn$, since $\cM$ is
relatively closed in $\Hn\times(t_0,T_*]$.

The kernel of a centre $X=(q,t_X)$ is written in the coordinates $\ell_q^{-1}$, in which
$f_q=|w|^2/16$ and $z_q=z$. All identities of the model hold in these coordinates, because
left translations preserve the structure.

\subsection{Growth and parametrisation}

\begin{lemma}[growth from a bound on the regularity scale]\label{lem:growth}
Let $F$ be a model flow on $(t_0,T_*]$ (in the sense of Definition~\ref{def:localflow}
or of Definition~\ref{def:traceflow}), $L\ge1$, and let $J\subset(t_0+3L^{-2},T_*]$ be
compact. If $\Kreg(\cM,X)\le L$ for all $X\in\cM$ with $t(X)\in J$, then
\[
\mu_t\big(F_t^{-1}(\cB_R)\big)\le C(n,\bb)\,L^{n+2}(1+R)^{2n+2}\qquad(t\in J,\ R\ge1).
\]
In particular $F$ has uniform polynomial growth if $\Kreg$ is bounded on
$\cM\cap(\Hn\times J)$ for every compact $J\subset(t_0,T_*]$.
\end{lemma}
\begin{proof}
Put $r=1/(2L)$, an admissible radius for every $X\in\cM$ with $t(X)\in J$
(Lemma~\ref{lem:admissible}(a)). For such $X=(q,t)$ the cylinder $\cC_{3r}(X)$ lies in
$\cU$, since $t-9r^2>t_0$, so that Lemma~\ref{lem:filling} applies at $X$. Only the
containment of Definition~\ref{def:K}(i) is used here. By it, the slice
$\cM(t)\cap\ell_q(\cB_r)$, which lies in the box of scale $r$ at $X$, is contained in the
Legendrian graph of a function $u$ with $|Du|\le1$ over $\{|x|<r\}$ (in the original
scale). The area of this graph is at most $\int_{|x|<r}\sqrt{\det g}\,dx\le C(n)r^n$, where
$g_{ij}=\frac14(\delta_{ij}+\partial_iu\cdot\partial_ju)$. Hence
$\mu_t(F_t^{-1}(\ell_q\cB_r))\le C(n)r^n$ for every point $q\in\cM(t)$. For a flow in
the sense of Definition~\ref{def:localflow}, this uses the injectivity of $\hat F$ on the
preimage of the box (Definition~\ref{def:K}(ii)), so that the parameter measure of the
piece equals the area of the graph. For a flow in trace form, $\mu_t$ is the area measure
of the embedded slice $\cM(t)$ (multiplicity one), and $\mu_t(F_t^{-1}(\ell_q\cB_r))$
stands for $\mu_t(\cM(t)\cap\ell_q\cB_r)$ (Definition~\ref{def:modelflow}). This quantity
is at most the area of the graph piece containing that part of the slice, and
Definition~\ref{def:K}(ii) is not needed.

Now let $\{q_j\}_{j=1}^N$ be a maximal $\frac r8$-separated subset of $\cB_R$ for the gauge
$\|\cdot\|$ (which is symmetric, the gauge being invariant under the inversion
$(w,z)\mapsto(-w,-z)$). The sets $\ell_{q_j}\cB_{r/32}$ are pairwise disjoint: a common
point $Z$ would give $\|q_k-q_j\|\le2(\|Z-q_j\|+\|Z-q_k\|)<\frac r8$ by the quasi-triangle
inequality with constant $2$. By the Heisenberg group law they lie in $\cB_{2R}$, and they have Lebesgue measure
$c(n,\bb)(r/32)^{2n+2}$ (the Lebesgue measure is left invariant and $\delta_\lambda$ scales
it by $\lambda^{2n+2}$), so $N\le C(n,\bb)(R/r)^{2n+2}$. By maximality the sets
$\ell_{q_j}\cB_{r/8}$ cover $\cB_R$. If $\ell_{q_j}\cB_{r/8}$ meets $\cM(t)$ at a point
$q'$, then $\ell_{q_j}\cB_{r/8}\subset\ell_{q'}\cB_{r/2}\subset\ell_{q'}\cB_r$, again by the
quasi-triangle inequality, so the contribution of this set to $\mu_t(F_t^{-1}(\cB_R))$ is
at most $C(n)r^n$. The sets not meeting $\cM(t)$ contribute nothing. Summing, we obtain
$\mu_t(F_t^{-1}(\cB_R))\le C(n,\bb)(R/r)^{2n+2}r^n$.
\end{proof}

\begin{lemma}[global normal parametrisation]\label{lem:param}
Let $\alpha\in(0,1)$, and let $\Kreg$ be the regularity scale for this exponent. Let $F$ be
a model flow on $(t_0,T_*]$ with $\Kreg$ bounded on $\cM\cap(\Hn\times J)$
for every compact $J\subset(t_0,T_*]$, and with $\cM$ non-empty; or let $\cM$ be a model
flow in trace form (Definition~\ref{def:traceflow}) with the same bounds, also non-empty. Then $\hat F$
is an embedding (in trace form, $\cM$ is embedded by definition), every slice $\cM(t)$,
$t_0<t\le T_*$, is non-empty, and for every finite $t_1\in(t_0,T_*]$ there is a smooth
family of diffeomorphisms $\Phi_{t_1\to t}\colon\cM(t_1)\to\cM(t)$ such that $(y,t)\mapsto
\Phi_{t_1\to t}(y)$ is a normal parametrisation of $\cM$ on $\cM(t_1)\times(t_0,T_*]$. In
particular, for every such $t_1$ and every $T'$ with $t_1\le T'<T_*$, the restriction of
this parametrisation to $\cM(t_1)\times[t_1,T']$ satisfies the properness hypothesis of
Fact~\ref{fact:kernel}(i). Moreover the map $(y,t)\mapsto(\Phi_{t_1\to t}(y),t)$ from
$\cM(t_1)\times(t_0,T_*]$ to $\Hn\times(t_0,\infty)$ is proper, so that a model flow in trace
form as above is a model flow on $(t_0,T_*]$ in the sense of Definition~\ref{def:localflow},
with $\Sigma=\cM(t_1)$.
\end{lemma}
\begin{proof}
Injectivity of $\hat F$: two parameter points with the same image $X$ both lie in the
preimage of the box of an admissible radius at $X$, contradicting
Definition~\ref{def:K}(ii). So $\hat F$ is a proper injective immersion, hence an
embedding, and $\cM$ is an embedded $(n+1)$-submanifold of $\Hn\times(t_0,\infty)$,
relatively closed there and transverse to the slices. For a flow in trace form these three
properties of $\cM$ are part of Definition~\ref{def:traceflow}: $\cM$ is relatively closed
in $\Hn\times(t_0,T_*]$, hence in $\Hn\times(t_0,\infty)$, it is an embedded
$(n+1)$-submanifold with boundary $\cM\cap\{t=T_*\}$, and it is transverse to the slices
because the time function has no critical points on $\cM$. No global map $\hat F$ is used
in that case.

The vector field $W=(H-2\theta\overline\xi,1)$ on $\cM$ is smooth. Let
$J\subset(t_0,T_*]$ be compact and put
\[
L:=\max\Big\{1,\ \sup_{\cM\cap(\Hn\times J)}\Kreg,\ 2(\min J-t_0)^{-1/2}\Big\},
\]
the last term being omitted when $t_0=-\infty$. By hypothesis, $L<\infty$. For every
$X\in\cM$ with $t(X)\in J$, the radius $r:=1/(2L)$ is admissible for $(\cM,X)$, since
$r<1/\Kreg(\cM,X)$ (Lemma~\ref{lem:admissible}(a); $r_1=\infty$ in the model). Moreover,
$\cC_{3r}(X)\subset\Hn\times(t_0,\infty)$, since
$9r^2=9/(4L^2)\le\frac9{16}(\min J-t_0)<t(X)-t_0$ when $t_0>-\infty$, while there is
nothing to check when $t_0=-\infty$. Hence Lemma~\ref{lem:subord}(iii), in the form
$|B|(X)\le C_1/r$, gives $|B|\le2C_1L$ on $\cM\cap(\Hn\times J)$. With
$|H|\le\sqrt n|B|$ and $|\theta|\le C_\theta$, the spatial part of $W$ has
$\overline g$-length at most $2\sqrt n\,C_1L+2C_\theta$ there.

The metric space $(\Hn,\overline g)$ is complete, so an integral curve of $W$ through a
point of $\cM$, having finite $\overline g$-length over the compact time interval $J$,
stays in a compact $\overline g$-ball during $J$. Since $\cM\cap(\overline B\times J)$ is
compact for such a ball $B$, and $J$ is an arbitrary compact subinterval of
$(t_0,T_*]$, the curve cannot leave every compact subset of $\cM$ in finite time. Hence
the curve exists for all $t\in(t_0,T_*]$ (at $t=T_*$ by the one-sided smoothness). As
$\cM$ is non-empty, every slice $\cM(t)$ is therefore non-empty.

Now fix a finite $t_1\in(t_0,T_*]$. The maps $\Phi_{t_1\to t}$ are
the flow maps of $W$ between slices; they are diffeomorphisms with inverses
$\Phi_{t\to t_1}$. Since $W$ has unit time component and spatial part normal to the
slices, $(y,t)\mapsto\Phi_{t_1\to t}(y)$ is a normal parametrisation. For every $T'$ with
$t_1\le T'<T_*$, this parametrisation is proper on $\cM(t_1)\times[t_1,T']$, because $\cM$
is relatively closed and the speed is bounded on the compact interval $[t_1,T']$. More
generally, the map $\hat\Phi(y,t):=(\Phi_{t_1\to t}(y),t)$ is a bijection of
$\cM(t_1)\times(t_0,T_*]$ onto $\cM$, with the continuous inverse
$(q,t)\mapsto(\Phi_{t\to t_1}(q),t)$, and $\cM$ is closed in $\Hn\times(t_0,\infty)$.
Hence the preimage under $\hat\Phi$ of a compact $\cK\subset\Hn\times(t_0,\infty)$ is the
image of the compact set $\cK\cap\cM$ under this inverse, and $\hat\Phi$ is proper. The
proof uses the exponent only through the constant
$C_1$ of Lemma~\ref{lem:subord}(iii).
\end{proof}

\subsection{Monotonicity under polynomial growth}

\begin{proposition}[monotonicity identity in the model]\label{prop:mono-model}
Let $F$ be a model flow on $(t_0,T_*]$ with uniform polynomial growth and with a global
normal parametrisation as in Lemma~\ref{lem:param} (or, more generally, let $F$ be as in
Fact~\ref{fact:kernel}(i) with the area bound there replaced by uniform polynomial
growth). Let $\hat y\in\Hn$ and $t_0<T\le T_*$. Then $\int\rho^\chi_{\hat y,T}d\mu_t<\infty$
for $t_0<t<T$, and for $t_0<t_1<t_2<T$
\begin{equation}\label{eq:mono-model}
\int\rho^\chi_{\hat y,T}\,d\mu_{t_2}-\int\rho^\chi_{\hat y,T}\,d\mu_{t_1}=-\int_{t_1}^{t_2}\!\!\int D^\chi\rho_{\hat y,T}\,d\mu_t\,dt .
\end{equation}
Thus $t\mapsto\int\rho^\chi_{\hat y,T}d\mu_t$ is non-increasing on $(t_0,T)$, and it is
constant on $[t_1,t_2]$ if and only if
\[
z_{\hat y}=4\tau\theta\qquad\text{and}\qquad H=-\frac{(\overline\nabla f_{\hat y})^\perp}{\tau}
\qquad\text{on }\cM(t)\ \text{for all }t\in[t_1,t_2].
\]
\end{proposition}
\begin{proof}
By left invariance we may take $\hat y=0$. We follow the proof of
\cite[Proposition~3.2]{PaperI}, which uses the area bound in exactly two places: to
show that $\int\rho^\chi d\mu_t$ is finite, and to show that the cut-off remainders vanish
in the limit. We replace both estimates and keep the rest of that proof, in particular
its two integrations of the cut-off identity.

\emph{Finiteness.} Fix $t\in(t_0,T)$ and put $\tau=T-t$. On the shell
$\cB_{R+1}\setminus\cB_R$ either $|w|\ge R$, so that $\rho\le(4\pi\tau)^{-n/2}e^{-R^2/(16\tau)}$,
or $|z|\ge\bb^2R^2$, so that $|v|=|z-4\tau\theta|/\tau\ge\bb^2R^2/\tau-4C_\theta$ and
$\chi(v)\le\exp(-a(\bb^2R^2/\tau-4C_\theta)_+^2)$. Hence, with $\chi\le1$ and $\rho\le(4\pi\tau)^{-n/2}$,
\begin{equation}\label{eq:shell}
\rho^\chi\le(4\pi\tau)^{-n/2}\,\epsilon_R(\tau),\qquad
\epsilon_R(\tau):=\max\Big\{e^{-R^2/(16\tau)},\ \exp\big(-a(\bb^2R^2/\tau-4C_\theta)_+^2\big)\Big\}
\end{equation}
on the shell, and $\int\rho^\chi d\mu_t\le(4\pi\tau)^{-n/2}\sum_{R\ge0}\epsilon_R(\tau)A_J(2+R)^{m_J}<\infty$
for $t$ in a compact $J$, uniformly in $t\in J$.

\emph{Cut-off remainders.} Let $\mathfrak c\in C^\infty(\R)$ be non-increasing with
$\mathfrak c=1$ on $(-\infty,\frac12]$ and $\mathfrak c=0$ on $[1,\infty)$, and put
$\psi_R(w,z):=\mathfrak c(|w|/R)^2\,\mathfrak c(|z|/(\bb^2R^2))^2$, so that $\psi_R\uparrow1$
as $R\to\infty$. This is the cut-off $\psi_{R,Z}$ of the proof of
\cite[Proposition~3.2]{PaperI} with $Z=\bb^2R^2$, formed with the functions
$\mathfrak c(|\cdot|)$ on $\R^{2n}$ and on $\R$, which are smooth, radially non-increasing,
equal to $1$ on the ball of radius $\frac12$ and equal to $0$ outside the unit ball. The
cut-off identity of \cite[Proposition~3.2]{PaperI} only requires that $\psi_R\circ F_t$ be
supported in a fixed compact subset of the parameter domain for $t\in[t_1,t_2]$
(properness). With the parametrisation of Lemma~\ref{lem:param}, this follows from its compact-time properness assertion
with base time $t_1$ and $T'=t_2$, and the integrals, taken
over the embedded slices $\cM(t)$, do not depend on the base time. The identity reads
\begin{gather*}
\frac{d}{dt}\int\psi_R\rho^\chi d\mu_t=-\int\psi_RD^\chi\rho\,d\mu_t+R_4+R_5,\\
R_4=\int\chi(\partial_t-\Delta)\psi_R\,\rho,\qquad R_5=-2\int\chi'\langle\nabla\psi_R,\nabla v\rangle\rho .
\end{gather*}
The support of $\overline\nabla\psi_R$ lies in $\cB_{R}\setminus\cB_{R/2}$, where
$|\nabla\psi_R|^2/\psi_R\le C(R^{-2}+(1+R^2)/(\bb^4R^4))=:C\kappa_R$ and
$|(\partial_t-\Delta)\psi_R|\le C(\kappa_R+\bb^{-2}R^{-2})$ (the computation of
\cite[Proposition~3.2]{PaperI} with $Z=\bb^2R^2$), and $\kappa_R\to0$. On this set either
$|w|\ge R/2$ or $|z|\ge\bb^2R^2/4$. Since $\chi'^2/\chi=4a^2v^2e^{-av^2}\le C_ae^{-av^2/2}$,
with the weaker exponent $a/2$, both weighted kernels satisfy
\begin{gather*}
\rho\chi\ \text{ and }\ \rho\,\frac{\chi'^2}{\chi}\ \le\ C_a(4\pi\tau)^{-n/2}\,\tilde\epsilon_R(\tau),\\
\tilde\epsilon_R(\tau):=\max\Big\{e^{-R^2/(64\tau)},\ \exp\big(-\tfrac a2(\bb^2R^2/(4\tau)-4C_\theta)_+^2\big)\Big\}
\end{gather*}
there. Therefore, for $t\in[t_1,t_2]$ and $\tau\in[\tau_2,\tau_1]$, where $\tau_j=T-t_j$,
\begin{equation}\label{eq:remainder}
|R_4|+\int_{\cB_R\setminus\cB_{R/2}}\frac{\chi'^2}{\chi}\frac{|\nabla\psi_R|^2}{\psi_R}\rho\,d\mu_t
\le\varepsilon_R:=C_a(\kappa_R+R^{-2})(4\pi\tau_2)^{-n/2}\sup_{\tau\in[\tau_2,\tau_1]}\tilde\epsilon_R(\tau)\,A_J(2+R)^{m_J},
\end{equation}
and $\varepsilon_R\to0$ as $R\to\infty$, because $\tilde\epsilon_R(\tau)$ decays faster
than any power of $R$ uniformly for $\tau$ in the compact interval, while the other factors
grow polynomially. With \eqref{eq:remainder} in place of the bound $\varepsilon_{R,Z}$ of
\cite[Proposition~3.2]{PaperI}, and with the one-parameter family $\psi_R$ in place of
$\psi_{R,Z}$, the two integrations of that proof apply as follows, the constants referring
to the fixed interval $[t_1,t_2]$.
\begin{enumerate}[(1)]
\item \emph{Uniform mass bound.} $0\le\int\psi_R\rho^\chi d\mu_t\le\int\rho^\chi d\mu_t\le C_J$ for
$t\in[t_1,t_2]$, by the finiteness estimate.
\item \emph{Integrable dissipation.} Since $|\nabla v|=4|\cS|$ and
\[
|R_5|\le8\Big(\int\frac{\chi'^2}{\chi}\frac{|\nabla\psi_R|^2}{\psi_R}\rho\Big)^{1/2}\Big(\int\psi_R\chi|\cS|^2\rho\Big)^{1/2}\le\frac{c_\chi}8\int\psi_R\chi|\cS|^2\rho+C\varepsilon_R
\]
by Young's inequality, integrating the cut-off identity over $[t_1,t_2]$ and using
$c_\chi\chi|\cS|^2\le D^\chi$ gives $\frac78\int_{t_1}^{t_2}\!\int\psi_RD^\chi\rho\le2C_J+C(t_2-t_1)\varepsilon_R$.
Letting $R\to\infty$, we obtain $\int_{t_1}^{t_2}\!\int D^\chi\rho\,d\mu_tdt\le\frac{16}7C_J<\infty$
by monotone convergence.
\item \emph{Vanishing remainders.} With this finiteness and the Cauchy--Schwarz form of
the bound for $R_5$, $\int_{t_1}^{t_2}(|R_4|+|R_5|)\,dt\le(t_2-t_1)\varepsilon_R+8\big((t_2-t_1)\varepsilon_R\big)^{1/2}\big(\int_{t_1}^{t_2}\!\int\chi|\cS|^2\rho\big)^{1/2}\to0$
as $R\to\infty$. Integrating the cut-off identity once more then gives
\eqref{eq:mono-model} by monotone convergence.
\end{enumerate}

The equality case follows from $D^\chi\ge\frac{|v\chi'(v)|}\tau+c_\chi\chi|\cS|^2$,
in which $v\chi'(v)=0$ only at $v=0$.
\end{proof}

\subsection{Density at points of the flow}

\begin{lemma}[density at a point]\label{lem:density-point}
Let $F$ be a model flow on $(t_0,T_*]$ and $X=(q,t)\in\cM$ with
$\Kreg(\cM,X)<\infty$. Let $\psi\ge0$ be any bounded function with $\psi=1$ on
$\ell_q\cB_{r'}$ for some $r'>0$. Then
\[
\liminf_{\sigma\to0}\int\psi\,\rho^\chi_{q,t}\,d\mu_{t-\sigma}\ \ge\ 1 .
\]
\end{lemma}
\begin{proof}
Choose an admissible radius $r$ with $r<r'$ and $r<1/\Kreg(\cM,X)$, so small that
$\cC_{3r}(X)\subset\cU$ (Lemma~\ref{lem:admissible}(a)). We work in the rescaled chart
$\tilde U\Lambda_{X,r}$. In this chart the trace in the unit box is the Legendrian graph
of $u$ over $\Omega_{T'}$ with $\|u\|_{2,\alpha}\le1$, $u(0,0)=0$, $Du(0,0)=0$
(Lemma~\ref{lem:filling}), and $\theta$ has the bounds of Lemma~\ref{lem:subord}(iv) on
$\Omega'$. Put $\sigma'=\sigma/r^2$. By Fact~\ref{fact:kernel}(ii) the integral is
invariant under the rescaling, so we compute it in the rescaled chart at the time
$-\sigma'$, over the graph piece $\{|x|<R\sqrt{\sigma'}\}$, for a fixed $R\ge1$ and
$\sigma'\to0$. This piece lies in $\Omega'$ for $\sigma'$ small, and the integrand is
non-negative elsewhere.

\emph{The height.} Along the graph, by \eqref{eq:height} and $\overline\beta=-\frac14(y\,dx-x\,dy)$,
\[
D_x\zeta(x,s)=-2\overline\beta\big((e,Du\,e)\big)_{e}=\tfrac12\big(u(x,s)-x\cdot Du(x,s)\big)\quad\text{(as a covector in }e),
\]
so that $|D_x\zeta(x,s)|\le\frac12(|u|+|x||Du|)\le C(|x|^2+|s|)$. Here we used
$|u(x,s)|\le C(|x|^2+|s|)$ and $|Du(x,s)|\le C(|x|+|s|^{1/2})$, which follow from
$u(0,0)=0$, $Du(0,0)=0$ and the $C^{2,\alpha}$ bound. Hence
$|\zeta(x,-\sigma')-\zeta(0,-\sigma')|\le C|x|(|x|^2+\sigma')$. At $x=0$,
Lemma~\ref{lem:heightangle}(b) gives $\partial_s\zeta(0,s)=-2\overline\beta((0,\partial_su(0,s)))-4\theta(0,s)=-4\theta(0,s)$,
because $\overline\beta$ at a point with $x=0$ vanishes on vectors with vanishing
$x$-component. With $\zeta(0,0)=0$ and $|\partial_s\theta|\le C$, we get
$\zeta(0,-\sigma')=4\sigma'\theta(0,-\sigma')+O(\sigma'^2)$. Therefore, on
$\{|x|\le R\sqrt{\sigma'}\}$ and at the time $-\sigma'$, using
$|\theta(x,s)-\theta(0,s)|\le C|x|$ (Lemma~\ref{lem:heightangle}(c)),
\[
|v|=\frac{|\zeta-4\sigma'\theta|}{\sigma'}\le C\big(R^3\sqrt{\sigma'}+R\sqrt{\sigma'}+\sigma'\big)\longrightarrow0,
\]
so $\chi(v)\to1$ uniformly there.

\emph{The Gaussian integral.} The area element of the graph is $\sqrt{\det
g}\,dx=2^{-n}(1+O(|x|+\sqrt{\sigma'}))dx$, and $f=|w|^2/16=(|x|^2+|u|^2)/16$ with
$|u|\le C\sigma'$ on the piece. Substituting $x=\sqrt{\sigma'}y$,
\[
\int_{\{|x|<R\sqrt{\sigma'}\}}(4\pi\sigma')^{-n/2}e^{-f/\sigma'}\sqrt{\det g}\,dx
\longrightarrow\int_{\{|y|<R\}}(4\pi)^{-n/2}2^{-n}e^{-|y|^2/16}\,dy\ \ge\ 1-C(n)e^{-R^2/32},
\]
the full integral over $\R^n$ being $1$ (\cite[Example~3.3]{PaperI}). Combining the two
steps, we find that the lower limit of $\int\psi\rho^\chi d\mu_{t-\sigma}$ as $\sigma\to0$
is at least $1-C(n)e^{-R^2/32}$ for every $R$, which proves the lemma.
\end{proof}

\begin{lemma}[lower semicontinuity of the density]\label{lem:lsc}
In the situation of Proposition~\ref{prop:compact}, suppose in addition that the limit
region contains a full time slab around the evaluation time: $\Hn\times(t-\sigma-\eta_0,t-\sigma+\eta_0)\subset\cU$
for some $\eta_0>0$, where $X=(q,t)\in\cM$, $\sigma>0$. Let $X_i=(q_i,t_i)\in\cM_i$
converge to $X$, let $\sigma_i\to\sigma$, and let $\psi_i\ge0$ be bounded Borel functions
converging to $1$ locally uniformly on $\Hn$. Let $\rho^{\chi,i}_{X_i}$ be the kernel of
$X_i$ for the structure $\varsigma_i$ (with $\hat a_i\to4$ in place of $4$,
Fact~\ref{fact:kernel}(ii)); it is built from the chart functions $f_{q_i}$, $z_{q_i}$ of
$q_i$ (Section~\ref{subsec:inputs}) and is therefore defined only on the domain of the
chart of $q_i$, the set $\delta_{\lambda_i/r_i}\Upsilon_{p_i}(U_{p_i}\cap U_{q_i})$ in the
ambient coordinates of Proposition~\ref{prop:compact} (notation as there). As for the
cut-offs $\psi_q$ of Section~\ref{subsec:inputs}, the product $\psi_i\rho^{\chi,i}_{X_i}$
is extended by zero outside that domain; with this convention the right-hand side below
is defined for every such $\psi_i$ (in the model, where the charts are global, the
convention is void). Then
\[
\int\rho^\chi_{X}\,d\mu_{t-\sigma}\ \le\ \liminf_{i\to\infty}\int\psi_i\,\rho^{\chi,i}_{X_i}\,d\mu^i_{t_i-\sigma_i},
\]
where the left-hand side is the model kernel without cut-off, integrated over the slice
$\cM(t-\sigma)$ of the limit flow.
\end{lemma}
\begin{proof}
We regard both sides as integrals against the joint measures $\fm_{t-\sigma}$ and
$\fm^i_{t_i-\sigma_i}$ of Proposition~\ref{prop:compact}(iii), which converge weakly on
$\cU(t-\sigma)\times\R=\Hn\times\R$ by the slab hypothesis (the evaluation times
$t_i-\sigma_i$ lie in the slab for $i$ large). The integrands are
$k(x,\vartheta)=(4\pi\sigma)^{-n/2}e^{-f_q(x)/\sigma}\chi(z_q(x)/\sigma-4\vartheta)$ and
$k_i(x,\vartheta)=\psi_i(x)(4\pi\sigma_i)^{-n/2}e^{-\tilde f_i(x)/\sigma_i}\chi(\tilde z_i(x)/\sigma_i-\hat a_i\vartheta)$
on the domain of the chart of $q_i$ and $k_i:=0$ outside it, by the convention of the
statement. Thus $k_i$ is defined, bounded and non-negative on all of $\Hn\times\R$ (and
Borel, $\psi_i$ being Borel and the chart domain open), so that the right-hand side of the
lemma is defined. Here $\tilde f_i\to f_q$ and $\tilde z_i\to z_q$ locally uniformly.
Indeed, in the situation of Proposition~\ref{prop:compact} the manifolds $M_i$ may vary.
Corollary~\ref{cor:rM}(c) applies there to the charts of all centres in compact sets, with
constants independent of $M_i$, and shows that the transitions to the rescaled charts of
$q_i$ are $\tilde U_i\ell_{\hat y_i}^{-1}$ up to errors tending to $0$ in $C^{k_0}$ on
compact sets. As $f$ and $z$ are $U(n)$-invariant and $\hat y_i\to q$, the claim follows
(for a fixed manifold, Fact~\ref{fact:charts}(iii) gives the same). In particular, a fixed
compact $\cK\subset\Hn$ lies in a ball $\cB_R$, $R\ge8$, which also contains $\hat y_i$
for all $i$, and Corollary~\ref{cor:rM}(c) (with $\lambda=\lambda_i\ge R/8$ for $i$ large)
provides the transition to the chart of $q_i$ on all of $\cB_R$. Therefore $\cK$ lies in
the domain of the chart of $q_i$ for $i$ large, and the zero extension does not affect the
convergence on $\cK$. Hence $k_i\to k$ uniformly on compact subsets of
$\Hn\times[-C_\theta,C_\theta]$.

Let $\eta_R\in C_c(\Hn)$, $0\le\eta_R\le1$, $\eta_R\uparrow1$. Since the measures are
supported in $\Hn\times[-C_\theta,C_\theta]$,
\[
\int\eta_Rk\,d\fm_{t-\sigma}=\lim_i\int\eta_Rk_i\,d\fm^i_{t_i-\sigma_i}\le\liminf_i\int k_i\,d\fm^i_{t_i-\sigma_i},
\]
the equality by Proposition~\ref{prop:compact}(iii) (the support of $\eta_R$ being a
compact subset of the slice region) and the uniform convergence on that support, the
inequality because $k_i\ge0$. Letting $R\to\infty$, monotone convergence gives the claim.
\end{proof}

\subsection{Rigidity}

\begin{proposition}[rigidity of ancient flows with density one]\label{prop:rigidity}
Let $F$ be a non-empty ancient model flow on $(-\infty,T_*]$, $T_*\in(-\infty,\infty]$,
in the sense of Definition~\ref{def:localflow} or of Definition~\ref{def:traceflow},
with $|\theta|\le C_\theta$ and $\Kreg(\cM,X)\le L$ for all $X\in\cM$. Suppose that
\[
\int\rho^\chi_X\,d\mu_{t-\sigma}\le1\qquad\text{for all }X=(q,t)\in\cM\text{ and all }\sigma>0 .
\]
Then there are a Lagrangian subspace $P\subset\C^n$, a point $\hat y\in\Hn$ and a constant
$\bar\theta\in[-C_\theta,C_\theta]$ such that
\[
\cM(t)=\ell_{\hat y}\big(P\times\{-4\bar\theta t\}\big)\quad\text{for all }t\le T_*,\qquad
\theta\equiv\bar\theta ,
\]
each slice being an embedded Legendrian plane of multiplicity one which translates
vertically, its height $z$ changing at the constant rate $-4\bar\theta$.
\end{proposition}
\begin{proof}
\emph{Step 0: parametrisation and growth.} Replacing $L$ by $\max\{L,1\}$, which preserves
the hypothesis $\Kreg\le L$ and does not affect the conclusion, we may assume $L\ge1$.
By Lemma~\ref{lem:param} the flow has a global normal parametrisation and all slices are
non-empty, and by Lemma~\ref{lem:growth} it has uniform polynomial growth. Hence
Proposition~\ref{prop:mono-model} applies to every centre $(q,t_X)$ with
$X=(q,t_X)\in\cM$ on the interval $(-\infty,t_X)$. By Lemma~\ref{lem:subord}(iii), in the
form $|B|(X)\le C_1/r$ with $r=1/(2L)$, we have $|B|\le2C_1L$ on $\cM$. Indeed, the radius
$1/(2L)$ is admissible for $(\cM,X)$ at every $X\in\cM$ (Lemma~\ref{lem:admissible}(a),
$r_1=\infty$ in the model), the cylinders $\cC_{3/(2L)}(X)$ lie in
$\Hn\times(-\infty,\infty)=\cU$, and Lemma~\ref{lem:filling}, hence
Lemma~\ref{lem:subord}, is stated for both forms of the flow.

\emph{Step 1: the density is identically one.} Fix $X=(q,t_X)\in\cM$ and put
$\Theta_X(\sigma):=\int\rho^\chi_Xd\mu_{t_X-\sigma}$. By
Proposition~\ref{prop:mono-model}, $\Theta_X$ is non-decreasing in $\sigma\in(0,\infty)$;
by Lemma~\ref{lem:density-point} (with $\psi\equiv1$), $\liminf_{\sigma\to0}\Theta_X\ge1$;
by hypothesis $\Theta_X\le1$. Hence $\Theta_X\equiv1$ on $(0,\infty)$, and the equality
case of Proposition~\ref{prop:mono-model} gives, in the coordinates $\ell_q^{-1}$,
\begin{equation}\label{eq:equality}
z=4\tau\theta,\qquad H=-\frac{(\overline\nabla f)^\perp}{\tau}\qquad\text{on }\cM(t),\ t<t_X,\ \tau=t_X-t .
\end{equation}

\emph{Step 2: self-similarity about $X$.} Since $2\overline\nabla f=X^{\cH}$ and
$\overline\xi$ is normal to every Legendrian submanifold, \eqref{eq:equality} gives
\[
H-2\theta\overline\xi=-\frac{(X^{\cH})^\perp}{2\tau}-\frac{z}{2\tau}\overline\xi
=-\frac{(X^{\cH}+z\overline\xi)^\perp}{2\tau}=-\frac{X^\perp}{2\tau}
\qquad\text{on }\cM(t),\ t<t_X ,
\]
where $X=X^{\cH}+z\overline\xi$ is the generator of the dilations (the symbol $X$ is used
for the generator in this step only). This is the self-shrinker equation of
\cite[(1.2)]{ChangWuZhang2026ContactBlowup} centred at $(q,t_X)$.

Fix $t_1<t_X$, put
$\tau_1=t_X-t_1$, $\lambda(t)=(\tau_1/\tau)^{1/2}$ and $G(y,t):=\delta_{\lambda(t)}\ell_q^{-1}F(y,t)$
for $t<t_X$, with the normal parametrisation of Step 0. Since
$\partial_\lambda\delta_\lambda=\lambda^{-1}X\circ\delta_\lambda$,
$\lambda'/\lambda=1/(2\tau)$, and $d\delta_\lambda$ maps normal spaces of Legendrian
submanifolds to normal spaces (it preserves the contact distribution, commutes with
$\overline\Phi$ on it and maps $\overline\xi$ to $\lambda^2\overline\xi$), we get
\[
\partial_tG=d\delta_\lambda\Big(-\frac{X^\perp}{2\tau}\Big)+\frac{1}{2\tau}X\circ G
=\frac1{2\tau}\big(X-X^\perp\big)\circ G=\frac{1}{2\tau}X^\top\circ G ,
\]
a vector tangent to the slice $G_t(\Sigma)$. Consequently the space-time set
$\cG=\{(G(y,t),t):t<t_X\}$, an embedded submanifold relatively closed in
$\Hn\times(-\infty,t_X)$, contains the vector $(0,\partial_t)$ in its tangent space at
every point, so this set is invariant under time translation as long as the translated
point stays in $\Hn\times(-\infty,t_X)$. Indeed, an integral curve of a vector field
tangent to a relatively closed embedded submanifold, with bounded spatial speed, stays in
the submanifold. Hence $G_t(\Sigma)=G_{t_1}(\Sigma)$ for all $t<t_X$, that is,
\begin{equation}\label{eq:selfsim}
\ell_q^{-1}\cM(t)=\delta_{(\tau/\tau_1)^{1/2}}\big(\ell_q^{-1}\cM(t_1)\big)\qquad(t<t_X).
\end{equation}
(This is the equivalence stated after \cite[(1.2)]{ChangWuZhang2026ContactBlowup}.)

\emph{Step 3: $B\equiv0$.} The second fundamental form of a Legendrian submanifold equals
that of its Lagrangian projection \cite[Lemma~2.3]{ChangWuZhang2026ContactBlowup}, and
$\delta_\mu$ acts on the projection as the homothety $w\mapsto\mu w$ of $(\C^n,g_0)$.
Hence, by \eqref{eq:selfsim} with $\mu=(\tau/\tau_1)^{1/2}$, we have
$|B|_{\cM(t_1)}(y)=\mu\,|B|_{\cM(t)}(\delta_\mu y)\le2\mu\,C_1L$ for every $y\in\ell_q^{-1}\cM(t_1)$.
Letting $t\uparrow t_X$, so that $\mu\to0$, gives $B\equiv0$ on $\cM(t_1)$, for every
$t_1<t_X$.

\emph{Step 4: one plane.} Fix $t_1<t_X$ and let $N$ be a connected component of the
slice $\ell_q^{-1}\cM(t_1)$, an embedded submanifold of $\Hn$ which is closed as a subset
(properness of $\hat F$), hence complete in its induced metric, since the metric space
$(\Hn,\overline g)$ is complete. By \cite[Lemma~2.3]{ChangWuZhang2026ContactBlowup} the
projection $\pi|_N\colon N\to(\C^n,g_0)$ is an isometric immersion, and by Step 3 its
second fundamental form vanishes. An isometric immersion of a connected complete
Riemannian $n$-manifold into Euclidean space with vanishing second fundamental form maps
each geodesic onto a straight line, so its image lies in the affine $n$-plane $\Pi_N$
through $\pi(y_0)$ spanned by $d\pi(T_{y_0}N)$ for any $y_0\in N$. Moreover,
$\pi|_N\colon N\to\Pi_N$ is a local isometry from a complete manifold onto a connected
one, hence a covering map. Since $\Pi_N$ is simply connected, $\pi|_N$ is an isometry onto
$\Pi_N$. Thus $N$ is the Legendrian lift of the full affine Lagrangian plane $\Pi_N$
(Lagrangian because $N$ is Legendrian), and the lift is single-sheeted.

Next, by \eqref{eq:equality} and $H=0$, we have $(\overline\nabla f)^\perp=0$ on $N$.
Hence the horizontal radial field $X^{\cH}=2\overline\nabla f$, whose projection is the
position field $w\mapsto w$, is tangent to $N$, so the position vector $w$ is tangent to
$\Pi_N$ at every $w\in\Pi_N$. Writing $\Pi_N=w_0+V$ with $V$ linear, we obtain
$w_0\in V$, that is, $\Pi_N=V$ is a linear Lagrangian subspace $P_N$. (This step does not
use the self-similarity \eqref{eq:selfsim}, which relates different times.) Since $H=0$ on
$N$, we have $\nabla\theta=0$ and $\theta\equiv\bar\theta_N$ is constant on $N$. By
\eqref{eq:equality}, $N=P_N\times\{4\tau_1\bar\theta_N\}$, on which $v\equiv0$.

Each full lifted plane through the Reeb axis has Gaussian integral $\int\rho\,d\cH^n=1$
\cite[Example~3.3]{PaperI}, and the slice is embedded with the components $N$ pairwise
disjoint, so
\[
1=\Theta_X(\tau_1)=\sum_N\int_{P_N\times\{4\tau_1\bar\theta_N\}}\rho\,\chi(0)\,d\cH^n=\#\{N\} ,
\]
and $\ell_q^{-1}\cM(t_1)=P\times\{4\tau_1\bar\theta\}$ is a single plane of multiplicity
one, with $\theta\equiv\bar\theta$ on it. This holds for every $t_1<t_X$, with $P$ and
$\bar\theta$ independent of $t_1$ by \eqref{eq:selfsim}, so
$\ell_q^{-1}\cM(t)=P\times\{4\tau\bar\theta\}$ for all $t<t_X$.

\emph{Step 5: all times.} The description holds for $t<t_X$ and for every $X\in\cM$,
with $P$ and $\bar\theta$ a priori depending on $X$. For a centre $X=(q,t_X)$ put
$\hat y_X:=\ell_q(0,4t_X\bar\theta)$. By the group law $\ell_q\tau_{4t_X\bar\theta}=\ell_{\hat y_X}$,
so Step 4 reads $\cM(t)=\ell_{\hat y_X}(P\times\{-4\bar\theta t\})$ for $t<t_X$. Let
$X'$ be a second centre, with plane $P'$ and angle $\bar\theta'$, and let
$t<\min\{t_X,t_{X'}\}$. The projection of $\cM(t)$ to $\C^n$ is the affine plane
$P+(\hat y_X)_w=P'+(\hat y_{X'})_w$, so $P'=P$, and $\theta$ takes the constant value
$\bar\theta=\bar\theta'$ on $\cM(t)$. Hence $g:=\hat y_{X'}^{-1}\hat y_X$ satisfies
$\ell_g(P\times\{c\})=P\times\{c\}$ with $c=-4\bar\theta t$. At the point $(0,c)$ this gives
$\ell_g(0,c)=(g_w,c+g_z)\in P\times\{c\}$, so $g_w\in P$ and $g_z=0$. For $w\in P$, the
cross term $\frac12(b\cdot x-a\cdot y)$ of the group law, with $g_w=(a,b)$ and $w=(x,y)$,
vanishes because $P$ is Lagrangian. Thus $\ell_g$ acts on every $P\times\{c'\}$ as the
translation by $g_w\in P$, and $\ell_{\hat y_X}(P\times\{c'\})=\ell_{\hat y_{X'}}(P\times\{c'\})$
for all $c'\in\R$. Fix a centre $X_0$ and put $\hat y:=\hat y_{X_0}$. Given $t<T_*$, there
is a centre $X$ with $t_X>t$, since the slices are non-empty at all times (Step 0). Hence
$\cM(t)=\ell_{\hat y_X}(P\times\{-4\bar\theta t\})=\ell_{\hat y}(P\times\{-4\bar\theta t\})$.
This holds for all $t<T_*$, also when $T_*=\infty$, and at $t=T_*$ (if finite) by
continuity of the trace.
\end{proof}

\begin{remark}
The hypotheses of Proposition~\ref{prop:rigidity} are those produced by the point-picking
argument of Section~\ref{subsec:proofA}: a bound on the regularity scale on all of $\cM$ from the
normalisation, and the density bound from Lemma~\ref{lem:lsc} with the rescaled cut-offs
$\psi_i\to1$. No Type~I hypothesis and no classification of self-shrinkers is used: the
equality case of the monotonicity identity, the bound on $|B|$ and the density count
replace them. The flow may stop at a finite upper end $T_*$. For a centre
$X=(q,t_X)\in\cM$, where $t_X\le T_*$, the self-similarity \eqref{eq:selfsim} is used only
at times $t<t_X$, and the slice $t=T_*$ (if $T_*$ is finite) is reached in Step~5 by
continuity.
\end{remark}

\section{Proof of the local regularity theorem}\label{sec:proof}

In this section we prove Theorem~\ref{thm:A}, following the plan given after its statement.
Of the three lemmas of
Section~\ref{subsec:finite}, the exhaustion of Lemma~\ref{lem:exhaust}(b),(c) is the contact
counterpart of \cite[Proposition~2.8]{wh}, and the finiteness of Lemma~\ref{lem:finite} is part
of White's notion of a $C^{2,\alpha}$ flow \cite[\S2.5]{wh}. The three lemmas concern, in
order, the finiteness of the regularity scale on compact sets, the density at a point of a
smooth flow, and the exhaustion of $\cU$ by relatively compact open sets. For the first,
note that the charts of nearby centres differ by a transition, so that the comparison of
the regularity scale at nearby points is not immediate. No Type~I hypothesis is used,
and the density hypothesis enters only through Lemma~\ref{lem:density-M}, which excludes
double points, and through Lemma~\ref{lem:lsc}.

\subsection{Finiteness of the regularity scale, double points and exhaustion}\label{subsec:finite}

\begin{lemma}[finiteness on compact sets]\label{lem:finite}
Let $F$ be a smooth local Legendrian flow in $\cU$ (Definition~\ref{def:localflow})
whose map $\hat F$ is injective. Then $\Kreg(\cM,X)<\infty$ for every $X\in\cM$, and for
every compact $\cK\subset\cU$
\[
\sup\{\Kreg(\cM,X):\ X\in\cM\cap\cK\}<\infty .
\]
\end{lemma}
\begin{proof}
\emph{The trace.} The map $\hat F=(F,t)$ is an immersion (each $F_t$ is an immersion and
the time component is the projection), injective and proper, hence a smooth embedding of
$\cN$ onto $\cM$. Therefore $\cM$ is an embedded $(n+1)$-dimensional submanifold of
$M\times\R$ with boundary $\cM\cap\{t=T_*\}$, relatively closed in $\cU$; its slices
$\cM(t)=F_t(\cN_t)$ are embedded Legendrian submanifolds, and the time function has no
critical points on $\cM$.

\emph{Step 1: a local graph at one point.} Fix $X=(q,t)\in\cM$ and use the chart
$\Lambda_{X,1}$, i.e.\ $\Upsilon_q$ with the time $s=t'-t$. At the origin the structure of
Fact~\ref{fact:charts}(i) coincides with the model ($h_q(0)=g_0$, $\beta_q(0)=0$), so the
tangent space of the slice $\cM(t)$ at the origin is a horizontal Lagrangian subspace
$U_0P$ of $\C^n$, $U_0\in U(n)$. In the rotated chart $\tilde U_0^{-1}\Lambda_{X,1}$ the
projection to $P$ restricted to $\cM(t)$ has invertible differential at the origin, and
$\cM$ is an embedded manifold transverse to the slices. By the implicit function theorem
there are $\rho_0>0$ and smooth functions $u,\zeta$ (with one-sided derivatives at the
upper end) on the base set $Q_{\rho_0}$, where for $r>0$
\[
Q_r:=\{|x|<r,\ |s|<r^2,\ s\le T_*-t\}
\]
(we write $Q_r$ rather than $D_r$ for this set: in this proof $D_\lambda$, with a positive
number as subscript, denotes the parabolic dilation
$D_\lambda(q,t)=(\delta_\lambda q,\lambda^2t)$ of Lemma~\ref{lem:admissible}(c), whereas $D$
and $D^2$, without subscript or with the subscript $x$ or $\xi$, denote spatial derivatives),
with $u(0,0)=0$, $Du(0,0)=0$, $\zeta(0,0)=0$, such that
$\cM\cap V=\hat\Gamma_{u,\zeta}(Q_{\rho_0})\cap V$ for the neighbourhood
$V:=\{|w|<\rho_0,\ |z|<\bb^2\rho_0^2\}\times(-\rho_0^2,\rho_0^2)$ of the origin. Here
the sheet through the origin is the only part of the embedded manifold $\cM$ in a small
neighbourhood, and $\rho_0$ is chosen accordingly. The graphs $\Gamma_{u(\cdot,s)}$ are
Lagrangian for the transverse structure, $\zeta$ being the height of the Legendrian lift.

For $0<r\le\rho_0/3$ consider the chart $\tilde U_0^{-1}\Lambda_{X,r}=D_{1/r}\circ\tilde U_0^{-1}\Lambda_{X,1}$.
The box $\bB$ of this chart is the dilated box $D_r\bB\subset V$ in the previous coordinates, so
$\tilde U_0^{-1}\Lambda_{X,r}(\cM)\cap\bB$ is contained in the Legendrian graph of
\[
u_r(x,s):=r^{-1}u(rx,r^2s),\qquad \zeta_r(x,s):=r^{-2}\zeta(rx,r^2s),
\]
over $\Omega_{T'}$, $T'=(T_*-t)/r^2$, and
\begin{multline*}
\|u_r\|_{2,\alpha;\Omega_{T'}}=\|Du\|_{0;Q_r}+r\big(\|D^2u\|_{0;Q_r}+\|\partial_su\|_{0;Q_r}\big)\\
+r^\alpha[Du]_{\alpha;Q_r}+r^{1+\alpha}\big([D^2u]_{\alpha;Q_r}+[\partial_su]_{\alpha;Q_r}\big)\longrightarrow0
\quad(r\to0),
\end{multline*}
with $Q_r$ as above (the map $(x,s)\mapsto(rx,r^2s)$ takes $\Omega_{T'}$ onto $Q_r$ and
multiplies the parabolic distance by $r$), because $Du(0,0)=0$, $Du$ is continuous and $u$
is smooth on $\overline{Q_{\rho_0/3}}$. The injectivity condition of
Definition~\ref{def:K}(ii) holds since $\hat F$ is injective. Hence every sufficiently
small $r$ is admissible for $(\cM,X)$ and $\Kreg(\cM,X)<\infty$.

\emph{Step 2: uniformity on compact sets.} Suppose that $X_i\in\cM\cap\cK$ and
$\Kreg(\cM,X_i)\to\infty$. Since $\cM\cap\cK$ is compact, we may assume $X_i\to X=(q,t)\in\cM\cap\cK$.
Let $\rho_0$, $U_0$, $u$, $\zeta$, $V$ be as in Step 1 for $X$, and for $0<r\le\rho_0/9$ let
$u_r$, $\zeta_r$ be as there, now on the base
\[
\Omega^4_{T'}:=\{|x|<4,\ |s|<16,\ s\le T'\},\qquad T'=(T_*-t)/r^2 ,
\]
and put $\eta(r):=\|u_r\|_{2,\alpha;\Omega^4_{T'}}$ and $\omega(r):=\eta(r)+C^{\rm tr}r^2$. By the
computation of Step 1, $\eta(r)\to0$ as $r\to0$. Moreover $\zeta_r\to-4\theta(X)\,s$
uniformly on the truncated set $\{|x|\le4,\ |s|\le16,\ s\le T'\}$. Indeed, $\zeta$ is
smooth on $Q_{\rho_0}$ (one-sided at the upper end) with $\zeta(0,0)=0$,
$D_x\zeta(0,0)=0$ and $D^2_x\zeta(0,0)=0$. For the last two, note that by
\eqref{eq:height}, for every unit vector $e\in\R^n$,
$D_x\zeta(x,0)\,e=-2\beta_\varsigma\big((e,Du(x,0)\,e)\big)$ with $\beta_\varsigma$ evaluated at
$(x,u(x,0))$, where $\beta_\varsigma=\overline\beta+O(|w|^3)$,
$\overline\beta_{(x,u)}((e,Du\,e))=-\frac14(u\cdot e-x\cdot Du\,e)$, $u(\cdot,0)=O(|x|^2)$ and
$Du(\cdot,0)=O(|x|)$, so that $D_x\zeta(\cdot,0)=O(|x|^2)$. Finally,
$\partial_s\zeta(0,0)=-4\theta(X)$ by \eqref{eq:theta-from-graph}. Hence the one-sided
Taylor expansion gives
$\zeta(x,s)=-4\theta(X)s+O(|x|^3+|x||s|+s^2)$ on $Q_{\rho_0}$ and $\zeta_r(x,s)=-4\theta(X)s+O(r)$
on $\Omega^4_{T'}$.

\emph{The order of the choices.} Two constants $C_3=C_3(n,\alpha)$ and $C_4=C_4(n,\alpha)$
are fixed below, in this order and before $r$; neither depends on $r$ or $i$. First, under
the preliminary conditions
\begin{equation}\label{eq:omega-small}
r\le\min\{\rho_0/9,\ r_1\},\qquad \omega(r)\le\tfrac1{600},\qquad
\sup_{\Omega^4_{T'}}|\zeta_r|\le64C_\theta+1 ,
\end{equation}
which involve neither $C_3$ nor $C_4$, the estimate \eqref{eq:step2-a} holds for all large
$i$ with a constant $C_3=C_3(n,\alpha)$ produced by the computation below, and
\eqref{eq:step2-b}, \eqref{eq:step2-c} follow from it with the explicit constants written
there. Next, under \eqref{eq:omega-small} and the additional condition
$(6+2C_3)\,\omega(r)\le\frac1{100}$, which involves only the already fixed $C_3$, the
inversion of the base map and the bound $\|\tilde v_i\|_{2,\alpha}\le C_4\,\omega(r)$ on the
graph function below hold with a constant $C_4=C_4(n,\alpha)$ (depending on $C_3$). Finally,
since $\omega(r)\to0$ and $\zeta_r\to-4\theta(X)s$ uniformly with
$\sup_{|s|\le16}|4\theta(X)s|\le64C_\theta$, we may and do fix $r$ with
\begin{equation}\label{eq:r-choice}
r\le\min\{\rho_0/9,\ r_1\},\qquad (6+2C_3)\,\omega(r)\le\tfrac1{100},\qquad C_4\,\omega(r)\le\tfrac12,\qquad
\sup_{\Omega^4_{T'}}|\zeta_r|\le64C_\theta+1 ;
\end{equation}
the second condition implies $\omega(r)\le\frac1{600}$ (as $C_3\ge0$), so
\eqref{eq:omega-small} holds. Each estimate below names the condition it uses. We show that
$r/18$ is admissible for $(\cM,X_i)$ for all large $i$, which contradicts
$\Kreg(\cM,X_i)\to\infty$.

\emph{Box inclusion.} For $i$ large, $\|X_i-X\|_X<r/6$, and the point
$\hat y_i:=\delta_{1/r}\Upsilon_q(q_i)$ lies in $\cB_1$; both statements hold because
$X_i\to X$ and the chart of $q$ is a homeomorphism. By Lemma~\ref{lem:quasi}, whose
right-hand side $3(\|X_i-X\|_X+\|Z-X_i\|_{X_i})$ is $<3(r/6+r/6)=r\le r_M$ (as
$r\le r_1\le r_M$ by \eqref{eq:omega-small} and \eqref{eq:r1}), every $Z\in\cC_{r/6}(X_i)$ satisfies
$\|Z-X\|_X<r$, so that
\[
(\tilde U\Lambda_{X_i,r/18})^{-1}(\bB)\subset\cC_{r/6}(X_i)\subset\cC_r(X)\subset V
\qquad\text{for every }U\in U(n),
\]
using $(\tilde U\Lambda_{Y,r'})^{-1}(\bB)\subset\cC_{3r'}(Y)$ (Section~\ref{subsec:scale}).
Hence the trace in such a box consists of points of the sheet $S:=\hat\Gamma_{u,\zeta}(Q_{\rho_0})$
whose coordinates in $\Lambda^0:=\tilde U_0^{-1}\Lambda_{X,r}$ satisfy $|x|<1$, $|s|<1$, and
the injectivity condition of Definition~\ref{def:K}(ii) holds in the box.

\emph{The sheet lies in the controlled region.} In the chart $\Lambda_{X,r}$ the sheet point
over $(x,s)\in\Omega^4_{T'}$ has coordinates $\big(U_0(x,u_r(x,s)),\ \zeta_r(x,s)\big)$ at the
time $s$, and this point lies in $\cB_8=\{|w|<8,|z|<64\bb^2\}$. Indeed, since $u_r(0,0)=0$, the
path from $(0,0)$ to $(0,s)$ to $(x,s)$ gives
$|u_r(x,s)|\le|s|\,\|\partial_su_r\|_0+|x|\,\|Du_r\|_0\le20\,\eta(r)$ on $\Omega^4_{T'}$ (the
norm $\eta(r)$ has no zero-order term), so $|w|^2<16+400\,\eta(r)^2<17$ because
$\eta(r)\le\omega(r)\le\frac1{600}$ by \eqref{eq:omega-small}. Moreover, $|z|\le64C_\theta+1<64\bb^2$
by \eqref{eq:omega-small} and $\bb^2=8C_\theta+4$. Only \eqref{eq:omega-small} is used here.

\emph{The transition.} We have $\Lambda_{X_i,r}\circ\Lambda_{X,r}^{-1}=(\Theta_i,\ s\mapsto s-s_i)$
with $\Theta_i:=\delta_{1/r}\Upsilon_{q_i}\Upsilon_q^{-1}\delta_r$ and $s_i:=(t_i-t)/r^2$,
$|s_i|<\frac1{36}$. By Lemma~\ref{lem:rM}(b), as $r\le r_1\le r_*$ and $\hat y_i\in\cB_8$,
there are $V_i\in U(n)$ with $\Theta_i=\tilde V_i\ell_{\hat y_i}^{-1}+E_i$ on $\cB_8$ and
$\|E_i\|_{C^{k_0}(\cB_8)}\le C^{\rm tr}r^2$. By Fact~\ref{fact:charts}(i) both charts
straighten the Reeb field to $2\partial_z$, so $\Theta_i$ and $\tilde V_i\ell_{\hat y_i}^{-1}$
have the form $(w,z)\mapsto(\tilde w(w),z+c(w))$, and the horizontal component
\[
\vartheta_i(w):=\pi\circ\Theta_i(w,z)=V_i(w-\hat y_{i,w})+\pi E_i(w)
\]
is a function of $w$ alone, with $\|\pi E_i\|_{C^{k_0}(\{|w|<8\})}\le C^{\rm tr}r^2$.

\emph{The sheet in the chart of $X_i$.} The horizontal position in $\Lambda_{X_i,r}$ of the
sheet point over $(x,s)\in\Omega^4_{T'}$ is $g_i(x,s):=\vartheta_i\big(U_0(x,u_r(x,s))\big)$.
The point $X_i$ is the sheet point over some $(x_i,s_i)$ with $|x_i|<\frac16$, and
$g_i(x_i,s_i)=0$. Let $B^0_i:=V_iU_0(I,0)\colon\R^n\to\C^n$, $x\mapsto V_iU_0(x,0)$, an
isometric embedding with image $V_iU_0P$. Write $w_r(x,s):=U_0(x,u_r(x,s))$, so that
$g_i=\vartheta_i\circ w_r$, $D_xw_r=U_0(I,Du_r)$, $D^2_xw_r=U_0(0,D^2u_r)$,
$\partial_sw_r=U_0(0,\partial_su_r)$, and $D\vartheta_i=V_i+D(\pi E_i)$,
$D^2\vartheta_i=D^2(\pi E_i)$. The chain rule gives
\begin{gather*}
D_xg_i=\big(V_i+D(\pi E_i)\big)U_0(I,Du_r),\qquad
\partial_sg_i=\big(V_i+D(\pi E_i)\big)U_0(0,\partial_su_r),\\
D^2_xg_i=D^2(\pi E_i)\big[U_0(I,Du_r),U_0(I,Du_r)\big]+\big(V_i+D(\pi E_i)\big)U_0(0,D^2u_r),
\end{gather*}
the derivatives of $\pi E_i$ being evaluated at $w_r(x,s)$.

We bound the two kinds of
terms separately, using only \eqref{eq:omega-small}. \emph{Terms without $\pi E_i$.} Since
$V_iU_0$ is a linear isometry of $\C^n\cong\R^{2n}$, the $C^{0,\alpha}$ norms (supremum
plus parabolic H\"older seminorm) of $V_iU_0(0,D^2u_r)$ and $V_iU_0(0,\partial_su_r)$ equal
$\|D^2u_r\|_{0,\alpha}$ and $\|\partial_su_r\|_{0,\alpha}$, which are at most $\eta(r)$;
likewise $[V_iU_0(I,Du_r)]_\alpha=[Du_r]_\alpha\le\eta(r)$ and
$\|V_iU_0(0,Du_r)\|_0\le\eta(r)$. \emph{Terms containing $\pi E_i$.} The map $w_r$ takes
$\Omega^4_{T'}$ into the ball $\{|w|<\sqrt{17}\}\subset\{|w|<8\}$ (shown above), and it is
$2$-Lipschitz for the parabolic distance $\max\{|x-x'|,|s-s'|^{1/2}\}$ on the convex set
$\Omega^4_{T'}$: $|w_r(x,s)-w_r(x',s')|\le(1+\|Du_r\|_0)|x-x'|+\|\partial_su_r\|_0|s-s'|$,
and $|s-s'|\le\sqrt{32}\,|s-s'|^{1/2}$ there, so the Lipschitz constant is at most
$1+(1+\sqrt{32})\eta(r)<2$ by \eqref{eq:omega-small}. Hence, by
$\|\pi E_i\|_{C^{k_0}(\{|w|<8\})}\le C^{\rm tr}r^2$ with $k_0\ge3$ and the mean value
theorem in the ball, every derivative of $\pi E_i$ of order at most two, composed with $w_r$,
has supremum at most $C^{\rm tr}r^2$ and parabolic H\"older seminorm at most
$2\cdot8^{1-\alpha}C^{\rm tr}r^2$ on $\Omega^4_{T'}$ (whose parabolic diameter is at most $8$).
The factors $U_0(I,Du_r)$, $U_0(0,D^2u_r)$ and $U_0(0,\partial_su_r)$ have $C^{0,\alpha}$
norms at most $1+2\eta(r)\le2$, and the $C^{0,\alpha}$ norm of a product is at most the
product of the $C^{0,\alpha}$ norms; so each term containing $\pi E_i$ has $C^{0,\alpha}$ norm
at most $C(n,\alpha)\,C^{\rm tr}r^2$. Collecting, $\|D_xg_i-B^0_i\|_0\le\|V_iU_0(0,Du_r)\|_0
+\|D(\pi E_i)\circ w_r\|_0\,\|(I,Du_r)\|_0\le\eta(r)+\sqrt2\,C^{\rm tr}r^2\le2\omega$, and
on $\Omega^4_{T'}$
\begin{equation}\label{eq:step2-a}
\|D_xg_i-B^0_i\|_0\le2\omega,\qquad
[D_xg_i]_\alpha+\|D^2_xg_i\|_{0,\alpha}+\|\partial_sg_i\|_{0,\alpha}\le C_3\,\omega ,
\end{equation}
where $\omega=\omega(r)=\eta(r)+C^{\rm tr}r^2$ and $C_3=C_3(n,\alpha)$ collects the
constants above (the terms without $\pi E_i$ contribute at most $3\eta(r)$).

The slice at the time $0$ of $\Lambda_{X_i,r}(\cM)$ is Legendrian
for the structure $\varsigma_{1/r}$ of $q_i$, which at the origin is the model. Hence the
horizontal projection $\Pi_i:=D_xg_i(x_i,s_i)(\R^n)$ of its tangent space at the origin is
a Lagrangian subspace of $(\C^n,\omega_0)$, $\Pi_i=U_iP$ with $U_i\in U(n)$. By
\eqref{eq:step2-a}, the isometric embedding $U_i^{-1}B^0_i$ differs by at most $2\omega$
from $U_i^{-1}D_xg_i(x_i,s_i)$, whose image is $P$; hence
\begin{equation}\label{eq:step2-b}
\big\|\pi_{J_0P}\,U_i^{-1}V_iU_0|_P\big\|\le2\omega,\qquad
A_i:=\pi_P\,U_i^{-1}V_iU_0|_P\ \text{ satisfies }\ \|A_i^{\rm T}A_i-I\|\le4\omega^2 ,
\end{equation}
so that $A_i$ is invertible ($4\omega^2<1$ by \eqref{eq:omega-small}) and, by polar
decomposition, $A_i=O_i(A_i^{\rm T}A_i)^{1/2}$ with $O_i$ a
linear isometry of $P$ and
$\|A_i-O_i\|=\|(A_i^{\rm T}A_i)^{1/2}-I\|\le\|A_i^{\rm T}A_i-I\|\le4\omega^2\le4\omega$. In the chart
$\tilde U_i^{-1}\Lambda_{X_i,r}$ the sheet point over $(x,s)$ has horizontal coordinates
$U_i^{-1}g_i(x,s)$; write
\[
\beta_i:=\pi_P\circ U_i^{-1}g_i,\qquad \gamma_i:=\pi_{J_0P}\circ U_i^{-1}g_i
\qquad\text{on }\Omega^4_{T'} .
\]
By \eqref{eq:step2-a}, \eqref{eq:step2-b} and $D_xg_i(x_i,s_i)(\R^n)=U_iP$,
\begin{equation}\label{eq:step2-c}
\begin{gathered}
\|D_x\beta_i-O_i\|_0\le6\omega,\qquad \|D_x\gamma_i\|_0\le4\omega,\qquad
\beta_i(x_i,s_i)=0,\quad \gamma_i(x_i,s_i)=0,\quad D_x\gamma_i(x_i,s_i)=0,\\
[D_x\beta_i]_\alpha+\|D^2_x\beta_i\|_{0,\alpha}+\|\partial_s\beta_i\|_{0,\alpha}
+[D_x\gamma_i]_\alpha+\|D^2_x\gamma_i\|_{0,\alpha}+\|\partial_s\gamma_i\|_{0,\alpha}\le2C_3\,\omega .
\end{gathered}
\end{equation}

\emph{Inversion of the base map.} From here on we use, besides \eqref{eq:omega-small}, the
condition $(6+2C_3)\,\omega\le\frac1{100}$ of \eqref{eq:r-choice}, which involves only the
constant $C_3$ already fixed; it gives $1-6\omega\ge\frac{99}{100}$ and $2C_3\omega\le\frac1{100}$.
Fix $s$ with $|s|<16$, $s\le T'$. By
\eqref{eq:step2-c} and the mean value theorem on the convex set $\{|x|<4\}$,
$|\beta_i(x,s)-\beta_i(x',s)-O_i(x-x')|\le6\omega|x-x'|$, so $\beta_i(\cdot,s)$ is
injective on $\{|x|<4\}$ with $|\beta_i(x,s)-\beta_i(x',s)|\ge(1-6\omega)|x-x'|$. Its image contains the ball of radius $(1-6\omega)\frac{23}6$ about $\beta_i(x_i,s)$. Indeed, the ball $\{|x-x_i|<\frac{23}6\}$ lies in $\{|x|<4\}$, and for $|\xi-\beta_i(x_i,s)|<(1-6\omega)\frac{23}{6}$
the map $x\mapsto x_i+O_i^{-1}\big(\xi-\beta_i(x_i,s)-(\beta_i(x,s)-\beta_i(x_i,s)-O_i(x-x_i))\big)$
is a contraction of a closed ball of radius $<\frac{23}6$ about $x_i$, whose fixed point solves $\beta_i(x,s)=\xi$. Since $|\beta_i(x_i,s)|=|\beta_i(x_i,s)-\beta_i(x_i,s_i)|\le\|\partial_s\beta_i\|_0|s-s_i|\le2C_3\omega\cdot17\le\frac{17}{100}<\frac12$
and $(1-6\omega)\frac{23}6-\frac12\ge\frac{99}{100}\cdot\frac{23}6-\frac12>3$,
the image contains $\{|\xi|<3\}$. Let $\phi_i(\cdot,s)$ be the inverse of $\beta_i(\cdot,s)$
on $\{|\xi|<3\}$, and put $Q':=\{|\xi|<3,\ |s|<16,\ s\le T'\}$ (a base set in the
coordinates $(\xi,s)$). On $Q'$,
\begin{gather*}
D_\xi\phi_i=(D_x\beta_i)^{-1}\circ(\phi_i,s),\qquad
\partial_s\phi_i=-(D_x\beta_i)^{-1}\partial_s\beta_i\circ(\phi_i,s),\\
D^2_\xi\phi_i=-D_\xi\phi_i\,D^2_x\beta_i(\phi_i,s)[D_\xi\phi_i,D_\xi\phi_i],
\end{gather*}
so $\|D_\xi\phi_i\|_0\le(1-6\omega)^{-1}\le\frac{100}{99}<2$, the map $(\xi,s)\mapsto(\phi_i(\xi,s),s)$ is Lipschitz for the
parabolic distance on the bounded set $Q'$ with a constant depending on $n$ only, and
(compositions of the $C^{0,\alpha}$ functions in \eqref{eq:step2-c} with this Lipschitz map,
and products with the bounded $D_\xi\phi_i$)
\[
[D_\xi\phi_i]_\alpha+\|D^2_\xi\phi_i\|_{0,\alpha}+\|\partial_s\phi_i\|_{0,\alpha}\le C(n,\alpha)\,C_3\,\omega\qquad\text{on }Q' .
\]

\emph{The graph function.} Put $\tilde v_i(\xi,s):=\gamma_i(\phi_i(\xi,s),s)$ on $Q'$. The
sheet points over $(x,s)\in\Omega^4_{T'}$ with $|\beta_i(x,s)|<3$ are exactly the points
$(\xi,\tilde v_i(\xi,s),\cdot,s)$, $|\xi|<3$, of the Legendrian graph of $\tilde v_i$ in the
chart $\tilde U_i^{-1}\Lambda_{X_i,r}$ (with the time $s-s_i$). By the chain rule,
$D_\xi\tilde v_i=D_x\gamma_i\,D_\xi\phi_i$,
$D^2_\xi\tilde v_i=D^2_x\gamma_i[D_\xi\phi_i,D_\xi\phi_i]+D_x\gamma_i\,D^2_\xi\phi_i$ and
$\partial_s\tilde v_i=D_x\gamma_i\,\partial_s\phi_i+\partial_s\gamma_i$, all evaluated at
$(\phi_i,s)$, so with \eqref{eq:step2-c} and the bounds on $\phi_i$ (which used only
\eqref{eq:omega-small} and $(6+2C_3)\omega\le\frac1{100}$),
\begin{gather*}
\|\tilde v_i\|_{2,\alpha;Q'}\le C_4\,\omega\quad\text{with }C_4=C_4(n,\alpha)\text{ (through }C_3\text{)},\qquad\text{hence }\|\tilde v_i\|_{2,\alpha;Q'}\le\tfrac12\text{ by \eqref{eq:r-choice}},\\
\tilde v_i(0,s_i)=\gamma_i(x_i,s_i)=0,\qquad D_\xi\tilde v_i(0,s_i)=D_x\gamma_i(x_i,s_i)D_\xi\phi_i(0,s_i)=0 ,
\end{gather*}
using $\phi_i(0,s_i)=x_i$.

\emph{The witness.} In the chart $\tilde U_i^{-1}\Lambda_{X_i,r}$ the time is $\tau=s-s_i$; put
$v_i(\xi,\tau):=\tilde v_i(\xi,\tau+s_i)$ on $Q'_i:=\{|\xi|<3,\ |\tau|<9,\ \tau\le T'-s_i\}$,
which corresponds to a subset of $Q'$ since $|s_i|<\frac1{36}$. Then
$\|v_i\|_{2,\alpha;Q'_i}\le\|\tilde v_i\|_{2,\alpha;Q'}\le\frac12$, $v_i(0,0)=0$ and
$D_\xi v_i(0,0)=0$. In the chart $\tilde U_i^{-1}\Lambda_{X_i,r/18}=D_{18}\circ\tilde U_i^{-1}\Lambda_{X_i,r}$
the sheet over the base $\Omega_{T''}$, $T''=324(T'-s_i)=(T_*-t_i)/(r/18)^2$, is the
Legendrian graph of $v^{(i)}(\xi,\tau):=18\,v_i(\xi/18,\tau/324)$, and by the scaling property
of the norm $\|v^{(i)}\|_{2,\alpha;\Omega_{T''}}\le\|v_i\|_{2,\alpha;\frac1{18}\Omega}\le\frac12$,
with $v^{(i)}(0,0)=0$, $Dv^{(i)}(0,0)=0$; the graph is Lagrangian for the rotated structure,
being a piece of the sheet, and its height is $C^1$ and vanishes at the origin.
We now check the containment. By the box inclusion,
$\tilde U_i^{-1}\Lambda_{X_i,r/18}(\cM)\cap\bB$ consists of
sheet points over $(x,s)$ with $|x|<1$, $|s|<1$, whose base coordinates in
$\tilde U_i^{-1}\Lambda_{X_i,r}$ are $(\beta_i(x,s),s-s_i)$ with
$|\beta_i(x,s)|\le|\beta_i(x_i,s)|+(1+6\omega)|x-x_i|<\frac12+\frac{101}{100}\cdot\frac76<3$.
So these are points of the graph of $\tilde v_i$. Their base points lie in
$\Omega_{T''}$ because the points lie in $\bB$, and hence they are points of the graph of
$v^{(i)}$. The injectivity
condition holds in the box, and $r/18\le r_1$ by \eqref{eq:r-choice}. Hence $r/18$ is
admissible for $(\cM,X_i)$ for all large $i$, i.e.\ $\Kreg(\cM,X_i)\le18/r$, which is the
announced contradiction.
\end{proof}

\begin{lemma}[density at a point; no double points]\label{lem:density-M}
Let $F$ be a smooth local Legendrian flow in $\cU$ and $X\in\cM$. Then
\[
\liminf_{\sigma\to0}\Theta^\chi_X(\sigma)\ \ge\ \#\hat F^{-1}(X).
\]
In particular, if $F$ satisfies $(\mathrm H_\eps)$ with $\eps<1$, then $\hat F$ is
injective.
\end{lemma}
\begin{proof}
\emph{Reduction to one branch.} The set $\hat F^{-1}(X)$ is compact (properness) and
discrete (immersion), hence finite; let $(x_1,t),\dots,(x_m,t)$ be its points. Choose
disjoint open neighbourhoods $W_j$ of $(x_j,t)$ in $\cN$ with compact closures in $\cN$
on which $\hat F$ is injective (an immersion is locally injective). Then
$\Theta^\chi_X(\sigma)\ge\sum_j\int_{W_j\cap\cN_{t-\sigma}}\psi_q\rho^\chi_{q,t}d\mu_{t-\sigma}$,
the integrand being non-negative, and the $j$-th term is the density at $X$ of the
\emph{branch} $F^{(j)}:=F|_{W_j}$. The branch is a smooth local Legendrian flow in
$\cU_j:=(M\times\R)\setminus\hat F(\partial W_j)$, an open set containing $X$: the set
$\hat F(\partial W_j)$ is compact and does not contain $X$ (injectivity on
$\overline{W_j}$), $\hat F(W_j)\subset\cU_j$ for the same reason, and for a compact
$\cK\subset\cU_j$ the preimage $\hat F^{-1}(\cK)\cap W_j=\hat F^{-1}(\cK)\cap\overline{W_j}$
is compact. The map $\hat F^{(j)}$ of the branch is injective. So it suffices to prove
$\liminf_{\sigma\to0}\Theta^\chi_X(\sigma)\ge1$ for a flow with injective $\hat F$.

\emph{A flow with injective $\hat F$.} Let $X=(q,t)\in\cM$ and fix a radius
$\mathsf{s}\in(0,d_1(X,\cU)/2)$ (the letter $\rho$ is kept for the kernels $\rho_{q,t}$,
$\rho^\chi_{q,t}$); then $\overline{\cC_{\mathsf{s}}(X)}$ lies in $U_q\times\R$ and is a
compact subset of $\cU$. Indeed, $\mathsf{s}<r_1\le r_M\le\min\{r_0/64,(z_0/(128\bb^2))^{1/2}\}$
by \eqref{eq:r1} and Definition~\ref{def:rM}, so $\mathsf{s}<r_0$ and
$\bb^2\mathsf{s}^2<z_0/128$: the closed box $\{|w|\le\mathsf{s},\ |z|\le\bb^2\mathsf{s}^2\}$
lies in the image $\{|w|<r_0,\ |z|<z_0\}$ of $\Upsilon_q$, hence
$\{Y:\|Y-X\|_X\le\mathsf{s}\}=\Upsilon_q^{-1}\{|w|\le\mathsf{s},|z|\le\bb^2\mathsf{s}^2\}\times[t-\mathsf{s}^2,t+\mathsf{s}^2]$
is a compact subset of $U_q\times\R$ containing $\overline{\cC_{\mathsf{s}}(X)}$, and its
points $Y$ have $\|Y-X\|_X\le\mathsf{s}<d(X,\cU)$, so they lie in $\cU$.
Hence by Lemma~\ref{lem:finite}, $\Kreg(\cM,\cdot)\le L$ on
$\cM\cap\overline{\cC_{\mathsf{s}}(X)}$. Let $\sigma_i\to0$ realise the $\liminf$ and put
$\lambda_i:=\sigma_i^{-1/2}\to\infty$. Consider the restriction of $F$ to $\cC_{\mathsf{s}}(X)$
(Remark~\ref{rem:localflow}(3)) in the rescaled chart $\Lambda_{X,1/\lambda_i}$: it is a
smooth local Legendrian flow $\cM_i$ in $\cU_i:=\cC_{\lambda_i\mathsf{s}}(0)$, and
$\cU_i^c\to\emptyset$. Indeed, the radii $\lambda_i\mathsf{s}\to\infty$, so $\cU_i^c$
eventually leaves every compact set. Hence no subsequence of points of $\cU_i^c$ converges,
which is condition (b) of the convergence of sets used in Proposition~\ref{prop:compact}
with limit $\emptyset$, condition (a) being vacuous. The ambient structure of these flows is
$\varsigma_{\lambda_i}$, which by
Corollary~\ref{cor:rM}(a) (with $r=r_M$, $\lambda=\lambda_ir_M$) satisfies the closeness
hypothesis of Proposition~\ref{prop:compact} with $\lambda_ir_M\to\infty$. The constant
in the angle equation is $(K+2)\lambda_i^{-2}\to0$. The upper end is
$T_i=\lambda_i^2(T_*-t)\in[0,\infty]$, with $T_i\to T\in\{0,\infty\}$ along a subsequence. By Lemma~\ref{lem:admissible}(b),(c), $\Kreg(\cM_i,\cdot)\le L/\lambda_i$ on
$\cM_i\cap\cC_{\lambda_i\mathsf{s}}(0)$. Thus Proposition~\ref{prop:compact} applies with
$\cU=\Hn\times\R$: along a subsequence $\cM_i\to\cM_\infty$ locally, a model flow in
trace form with upper end $T$, and $0\in\cM_\infty$.

\emph{The limit is the tangent plane.} Let $U_0$, $u$, $\zeta$ be as in Step 1 of the
proof of Lemma~\ref{lem:finite} for $X$. For every $R\ge1$ and $i$ large, the radius $R$
is admissible for $(\cM_i,0)$, with the witness $(U_0^{-1},u_{R/\lambda_i},\zeta_{R/\lambda_i})$
in the notation of that step (the rescaling of the witness at the radius $R/\lambda_i$ of
$(\cM,X)$, Lemma~\ref{lem:admissible}(c)). By Proposition~\ref{prop:compact}(ii) at
$X=0$ with this $r=R$, $\cM_\infty\cap\bB_{0,R}=\hat\Gamma_{u^{(R)},\zeta^{(R)}}(\Omega_{T'})$
in the chart $\tilde U_0^{-1}\Lambda_{0,R}$, where $u^{(R)}=\lim u_{R/\lambda_i}=0$ (the
norms tend to $0$ by Step 1) and
$\zeta^{(R)}(x,s)=\lim\lambda_i^2R^{-2}\zeta(Rx/\lambda_i,R^2s/\lambda_i^2)=\partial_s\zeta(0,0)\,s=-4\theta(X)\,s$
(Taylor expansion of the smooth $\zeta$ with $\zeta(0,0)=0$, $D_x\zeta(0,0)=0$,
$D^2_x\zeta(0,0)=0$ by \eqref{eq:height} and $Du(0,0)=0$, and \eqref{eq:theta-from-graph}
at the origin). Since the boxes $\bB_{0,R}$ exhaust $\Hn\times\R$,
\[
\cM_\infty=\Pi:=\{(w,-4\bar\theta s,s):\ w\in U_0P,\ s\le T\},\qquad \bar\theta:=\theta(X),
\]
the lift of the tangent plane through the origin, translating vertically, its height $z$
changing at the constant rate $-4\bar\theta$.

\emph{The density.} Apply Lemma~\ref{lem:lsc} to the flows $\cM_i$, with both points of
that lemma at the origin and both of its scale parameters equal to $1$. The scales
$\sigma_i\to0$ fixed above are not the parameters of that lemma: in the chart
$\Lambda_{X,1/\lambda_i}$ the time is $s=\lambda_i^2(t'-t)$, so the kernel scale $1$ and the
evaluation time $s=-1$ correspond to the original scale $\lambda_i^{-2}=\sigma_i$ and the
original time $t-\sigma_i$. We take the full slab $\Hn\times\R$ and $\psi_i:=\psi_q$
written in the rescaled chart; this cut-off equals
$1$ on $\{|w|<\lambda_ir_0/2,\ |z|<\lambda_i^2z_0/2\}$ and hence tends to $1$ locally
uniformly. By Fact~\ref{fact:kernel}(ii) the right-hand side of Lemma~\ref{lem:lsc} is the
density of the restricted flow at $X$ at the original scale $\lambda_i^{-2}=\sigma_i$, which is at
most $\Theta^\chi_X(\sigma_i)$. The left-hand side is the density of $\Pi$ at the origin at
the scale $1$: on the slice $\Pi(-1)=U_0P\times\{4\bar\theta\}$ the caloric height is
$u=z-4\tau\theta=4\bar\theta-4\bar\theta=0$, so $\chi(v)=1$, and
\begin{multline*}
\int_{\Pi(-1)}\rho^\chi_{0,0}\,d\mu_{-1}=\int_{U_0P}(4\pi)^{-n/2}e^{-|w|^2/16}\,d\mu_{g_0}\\
=(4\pi)^{-n/2}\,2^{-n}\int_{\R^n}e^{-|x|^2/16}\,dx=(4\pi)^{-n/2}2^{-n}(16\pi)^{n/2}=1 ,
\end{multline*}
the area element of a Lagrangian plane for $g_0=\frac14|dw|^2$ being $2^{-n}$ times the
Lebesgue measure. Hence $1\le\liminf_i\Theta^\chi_X(\sigma_i)=\liminf_{\sigma\to0}\Theta^\chi_X(\sigma)$.

\emph{No double points.} If $(\mathrm H_\eps)$ holds with $\eps<1$ and $X\in\cM$ had
$\#\hat F^{-1}(X)\ge2$, then $\Theta^\chi_X(\sigma)\ge2-o(1)>1+\eps$ for small
$\sigma<d_1(X,\cU)^2$, contradicting $(\mathrm H_\eps)$.
\end{proof}

\begin{lemma}[exhaustion]\label{lem:exhaust}
Let $F$ be a smooth local Legendrian flow in $\cU$ satisfying $(\mathrm H_\eps)$ with
$\eps<1$, and let $\cU'\subset\cU$ be open with compact closure in $\cU$. Then:
\begin{enumerate}[(a)]
\item the restriction $F|_{\cU'}$ satisfies $(\mathrm H_\eps)$ in $\cU'$;
\item $\Kreg_{\cU'}(\cM\cap\cU')=\sup_{X\in\cM\cap\cU'}d_1(X,\cU')\Kreg(\cM\cap\cU',X)<\infty$;
\item if $C\ge3$ and $\Kreg_{\cU'}(\cM\cap\cU')\le C$ for every such $\cU'$, then
$\Kreg_\cU(\cM)\le C$.
\end{enumerate}
\end{lemma}
\begin{proof}
(a) For $X\in\cM\cap\cU'$ the density of the restriction is at most that of $F$ (the
integrand is non-negative) and $d_1(X,\cU')\le d_1(X,\cU)$, so the range of scales only
shrinks. (b) By Lemma~\ref{lem:density-M}, $\hat F$ is injective; by Lemma~\ref{lem:finite}
applied to the compact set $\overline{\cU'}\subset\cU$, $\Kreg(\cM,X)\le L$ for
$X\in\cM\cap\overline{\cU'}$; by Lemma~\ref{lem:admissible}(b), $\Kreg(\cM\cap\cU',X)\le\Kreg(\cM,X)$;
and $d_1\le r_1$.

(c) Fix $X\in\cM$ and let $\cU'_k\uparrow\cU$ be an exhaustion of
$\cU$ by open sets with compact closures in $\cU$. For $\delta<d_1(X,\cU)$ the closed set
$\{\|Y-X\|_X\le\delta\}=\Upsilon_q^{-1}\{|w|\le\delta,|z|\le\bb^2\delta^2\}\times[t-\delta^2,t+\delta^2]$
is a compact subset of $\cU$. Indeed, it lies in $U_q\times\R$, and the closed box lies in the
image $\{|w|<r_0,|z|<z_0\}$ of $\Upsilon_q$, because
$\delta<r_1\le r_M\le\min\{r_0/64,(z_0/(128\bb^2))^{1/2}\}$ by \eqref{eq:r1} and
Definition~\ref{def:rM}, so $\delta<r_0$ and $\bb^2\delta^2<z_0/128$; the time interval
is closed and bounded; and its points $Y$ satisfy $\|Y-X\|_X\le\delta<d(X,\cU)$, so they
lie in $\cU$. Hence this compact set lies in $\cU'_k$ for $k$ large and
$d(X,\cU'_k)\ge\delta$. It follows that $d_1(X,\cU'_k)\to d_1(X,\cU)$. For
$r<d(X,\cU'_k)/3$ the box $\bB_{X,r}\subset\cC_{3r}(X)$ lies in $\cU'_k$, so the trace of
$F$ and of its restriction in the box coincide, as do the preimages of the box under
$\hat F$; thus $r$ is admissible for $(\cM,X)$ if and only if it is admissible for
$(\cM\cap\cU'_k,X)$. By hypothesis $\Kreg(\cM\cap\cU'_k,X)\le C/d_1(X,\cU'_k)$, so every
$r<d_1(X,\cU'_k)/C\le d(X,\cU'_k)/3$ is admissible for $(\cM\cap\cU'_k,X)$, hence for
$(\cM,X)$, and $\Kreg(\cM,X)\le C/d_1(X,\cU'_k)$. Letting $k\to\infty$ gives
$d_1(X,\cU)\Kreg(\cM,X)\le C$.
\end{proof}

\subsection{Point picking and the normalised sequence}

\begin{lemma}[normalised sequence]\label{lem:pick}
Let $F$ be a smooth local Legendrian flow in $\cU\subset M\times\R$ with $\hat F$
injective and $s:=\Kreg_\cU(\cM)\in(0,\infty)$. Then there is $X_0=(q_0,t_0)\in\cM$ with
$d_1(X_0,\cU)\Kreg(\cM,X_0)\ge s/2$; put $\lambda:=\Kreg(\cM,X_0)$ and
$d_\star:=d_1(X_0,\cU)$. Let
$F^0$ be the restriction of $F$ to $\cU\cap(U_{q_0}\times\R)$, the domain of the chart of
$q_0$ (a smooth local Legendrian flow by Remark~\ref{rem:localflow}(3)), with trace
$\cM^0=\cM\cap(U_{q_0}\times\R)$, and let $\cM'$ be the trace of $F^0$ written in the
rescaled chart $\Lambda_{X_0,1/\lambda}$ (i.e.\ $\delta_\lambda\Upsilon_{q_0}$ with the time
$\lambda^2(t-t_0)$): a smooth local Legendrian flow in the \emph{rescaled region}
\[
\cU':=\Lambda_{X_0,1/\lambda}\big(\cU\cap(U_{q_0}\times\R)\big)\subset\Hn\times\R
\]
for the structure $\varsigma_\lambda$, with $r_1$ replaced by $\lambda r_1$
(Lemma~\ref{lem:admissible}(c)). Then $\lambda d_\star\ge s/2$, $\lambda r_M\ge s/2$,
$1/\lambda\le r_1$, and
\begin{enumerate}[(a)]
\item $\Kreg(\cM',0)=1$ and $\Kreg(\cM',Y)\le12$ for all $Y\in\cM'\cap\cC_{\lambda d_\star/6}(0)$;
\item $\cU'\supset\cC_{\lambda d_\star}(0)$, and the upper end of $\cM'$ is
$T'=\lambda^2(T_*-t_0)\ge0$;
\item for $Y=(q,t)\in\cM'\cap\cC_{\lambda d_\star/6}(0)$ and $0<\sigma<(\lambda d_\star/6)^2$, the
rescaled density $I'_Y(\sigma):=\int\psi'\rho^{\chi\prime}_{Y}d\mu'_{t-\sigma}$ of $\cM'$ at
$Y$ (with the cut-off $\psi'$ of the centre $Y$, written in the rescaled chart) satisfies
\[
I'_Y(\sigma)=\Theta^\chi_{\hat Y}\big(\sigma\lambda^{-2};\,F^0\big)
\le\Theta^\chi_{\hat Y}\big(\sigma\lambda^{-2};\,F\big),\qquad \sigma\lambda^{-2}<d_1(\hat Y,\cU)^2 ,
\]
where $\hat Y\in\cM^0$ is the point corresponding to $Y$.
\end{enumerate}
\end{lemma}
\begin{proof}
The point $X_0$ exists because $s$ is a finite supremum over the non-empty set $\cM$,
and $\lambda<\infty$ by Lemma~\ref{lem:finite}. Since $d_\star\le r_1\le r_M\le1$,
$\lambda\ge s/(2d_\star)\ge s/2$ and $\lambda r_M\ge\lambda d_\star\ge s/2$; and
$\lambda\ge1/r_1$ by Definition~\ref{def:K}. For $Y\in\cM$ with $\|Y-X_0\|_{X_0}<d_\star/6$:
if $Z\notin\cU$ had $\|Z-Y\|_Y<d_\star/6$, then Lemma~\ref{lem:quasi}, whose right-hand side
$3(\|Y-X_0\|_{X_0}+\|Z-Y\|_Y)$ is $<3(d_\star/6+d_\star/6)=d_\star\le r_M$, would give
$\|Z-X_0\|_{X_0}<d_\star\le d(X_0,\cU)$, which is impossible. Therefore
$d(Y,\cU)\ge d_\star/6$ and $d_1(Y,\cU)\ge d_\star/6$. Hence
$\Kreg(\cM,Y)\le s/d_1(Y,\cU)\le6s/d_\star\le12\lambda$, using $s\le2d_\star\lambda$.

The restriction $F^0$ of $F$ to $\cU\cap(U_{q_0}\times\R)$ defined in the statement,
with trace $\cM^0$, has the same regularity scale as $F$ at every
$Y\in\cM\cap\cC_{d_\star/6}(X_0)$. Indeed, let $Y=(q,t)$, work at the scale $r_M$, and put
$h:=\delta_{1/r_M}\Upsilon_{q_0}(q)$; then $|h_w|<\frac16$ and $|h_z|<\bb^2/36$, since
$\|Y-X_0\|_{X_0}<d_\star/6\le r_M/6$. By Lemma~\ref{lem:rM}(b) for the pair $(q_0,q)$, the
transition $S:=\delta_{1/r_M}\Upsilon_q\Upsilon_{q_0}^{-1}\delta_{r_M}$ is $\delta_*$-close in
$C^{k_0}(\cB_8)$ to $\tilde U\ell_h^{-1}$; in particular $q_0\in U_q$, and
$h':=S(0)=\delta_{1/r_M}\Upsilon_q(q_0)$ satisfies $|h'_w|\le|h_w|+\delta_*<\frac14$ and
$|h'_z|\le\bb^2/36+\delta_*$, so $h'\in\cB_1$. By Lemma~\ref{lem:rM}(b) for the pair
$(q,q_0)$, the transition in the other direction,
$S':=\delta_{1/r_M}\Upsilon_{q_0}\Upsilon_q^{-1}\delta_{r_M}$ from the rescaled chart of $q$ to
that of $q_0$, is defined on $\cB_8$ and $\delta_*$-close there to $\tilde U'\ell_{h'}^{-1}$.
For $\hat z\in\cB_3$ the group law gives $|(\ell_{h'}^{-1}\hat z)_w|\le3+\frac14$ and
$|(\ell_{h'}^{-1}\hat z)_z|\le9\bb^2+\frac{\bb^2}{36}+\delta_*+\frac12\cdot\frac14\cdot3$, so that
$|S'(\hat z)_w|<3+\frac14+\delta_*<4$ and $|S'(\hat z)_z|<9\bb^2+\frac{\bb^2}{36}+\frac38+2\delta_*<16\bb^2$
(using $\bb^2\ge4$, $\delta_*\le\frac1{16}$). Thus $S'(\cB_3)\subset\cB_4$. For $r\le r_1$ the
cylinder $\cC_{3r}(Y)$ has rescaled $q$-coordinates in $\cB_3$ ($3r\le3r_M$,
$9\bb^2r^2\le9\bb^2r_M^2$), hence lies in $\Upsilon_{q_0}^{-1}(\delta_{r_M}\cB_4)\subset U_{q_0}$
($4r_M<r_0$ and $16\bb^2r_M^2<z_0$ by Definition~\ref{def:rM}). Every box
$(\tilde U\Lambda_{Y,r})^{-1}(\bB)$, $r\le r_1$, lies in $\cC_{3r}(Y)$
(Section~\ref{subsec:scale}), so the traces of $F$ and $F^0$ in each such box coincide, as
do the preimages of the box under $\hat F$; the admissible radii of $(\cM,Y)$ and
$(\cM^0,Y)$ are therefore the same, and $\Kreg(\cM^0,Y)=\Kreg(\cM,Y)$.

The rescaled flow $\cM'$ is $F^0$ written in the chart $\Lambda_{X_0,1/\lambda}$, as
defined in the statement, so
Lemma~\ref{lem:admissible}(c) gives (a): $\Kreg(\cM',0)=\lambda^{-1}\Kreg(\cM^0,X_0)=1$ and
$\Kreg(\cM',Y')=\lambda^{-1}\Kreg(\cM^0,\hat Y')\le12$ for $Y'\in\cM'\cap\cC_{\lambda d_\star/6}(0)$.
For (b), every $Y$ with $\|Y-X_0\|_{X_0}<d_\star\le d(X_0,\cU)$ lies in $\cU$, and
$\cC_{d_\star}(X_0)\subset U_{q_0}\times\R$ because $d_\star\le r_1\le r_M$ gives
$d_\star<r_0$ and $\bb^2d_\star^2<z_0$ (Definition~\ref{def:rM}). Hence
$\cC_{d_\star}(X_0)\subset\cU\cap(U_{q_0}\times\R)$, and since
$\Lambda_{X_0,1/\lambda}(\cC_{d_\star}(X_0))=\cC_{\lambda d_\star}(0)$ by the homogeneity
of the gauge, $\cC_{\lambda d_\star}(0)\subset\cU'$. The upper end $T_*$ of $F^0$ becomes
$T'=\lambda^2(T_*-t_0)\ge0$ in the chart. For (c), Fact~\ref{fact:kernel}(ii) applied to $F^0$ gives
$I'_Y(\sigma)=\Theta^\chi_{\hat Y}(\sigma\lambda^{-2};F^0)$; the density of the restriction
$F^0$ is at most that of $F$, the integrand being non-negative and the parameter domain of
$F^0$ being contained in that of $F$; and $\sigma\lambda^{-2}<(d_\star/6)^2\le d_1(\hat Y,\cU)^2$.
\end{proof}

\subsection{Proof of Theorem~\ref{thm:A}}\label{subsec:proofA}

\begin{proof}[Proof of Theorem~\ref{thm:A}]
\emph{Step 1: the contradiction sequence.} Suppose the theorem is false for the given
$n,a,\alpha,C_\theta$. Then for every $i=1,2,\dots$ there are a compact Sasakian manifold $M_i$
(with its constants $K_i$, $r_{M_i}$, $r_1^{(i)}$), an open set $\cU_i\subset M_i\times\R$
and a smooth local Legendrian flow $F_i$ in $\cU_i$ satisfying $(\mathrm H_{\eps_i})$ with
$\eps_i:=\min\{1/i,\frac12\}$ and $\Kreg_{\cU_i}(\cM_i)>\max\{i,3\}$. By
Lemma~\ref{lem:exhaust}(c), applied with $C=\max\{i,3\}$, there is an open
$\cV_i\subset\cU_i$ with compact closure in $\cU_i$ such that the restriction of $F_i$ to
$\cV_i$, which satisfies $(\mathrm H_{\eps_i})$ by (a), has $\Kreg_{\cV_i}>\max\{i,3\}$.
By (b) this supremum is finite, and by Lemma~\ref{lem:density-M} the map $\hat F_i$ is
injective. Replacing $(\cU_i,F_i)$ by $(\cV_i,F_i|_{\cV_i})$ we may therefore assume
\[
s_i:=\Kreg_{\cU_i}(\cM_i)\in(i,\infty),\qquad \hat F_i\ \text{injective},\qquad (\mathrm H_{\eps_i}),\ \eps_i\to0 .
\]

\emph{Step 2: normalisation.} We apply Lemma~\ref{lem:pick} to each $F_i$, with
$X_i=(q_i,t_i)$, $\lambda_i=\Kreg(\cM_i,X_i)$ and $d_{\star,i}=d_1(X_i,\cU_i)$. The rescaled flows
$\cM'_i$ (the restrictions of $F_i$ to $\cU_i\cap(U_{q_i}\times\R)$ written in the charts
$\delta_{\lambda_i}\Upsilon_{q_i}$, time $\lambda_i^2(t-t_i)$) in the rescaled regions
\[
\cU'_i:=\Lambda_{X_i,1/\lambda_i}\big(\cU_i\cap(U_{q_i}\times\R)\big)\subset\Hn\times\R
\]
of Lemma~\ref{lem:pick}(b) satisfy
$\Kreg(\cM'_i,0)=1$, $\Kreg(\cM'_i,\cdot)\le12$ on $\cC_{R_i}(0)$ with
$R_i:=\lambda_id_{\star,i}/6\ge s_i/12\to\infty$, and $\cU'_i\supset\cC_{6R_i}(0)$. Hence
$\cU'^c_i\to\emptyset$: $\cU'^c_i$ eventually leaves every compact set, so the hypothesis of
Proposition~\ref{prop:compact} on the complements holds in its subsequence form.

The ambient coordinates are the rescaled charts of the centres
$q_i$ at the scales $1/\lambda_i=r_{M_i}/\Lambda_i$ with $\Lambda_i:=\lambda_ir_{M_i}\ge s_i/2\to\infty$,
so by Corollary~\ref{cor:rM} the structures $\varsigma_i$ satisfy the closeness hypothesis
of Proposition~\ref{prop:compact} with $r_i=r_{M_i}$ and $\Lambda_i$ in place of $\lambda_i$ (only the normalised bounds of
Corollary~\ref{cor:rM} enter; no lower bound on the scales $r_{M_i}$ is used). The constants in the angle equation are $\kappa_{\theta,i}=(K_i+2)\lambda_i^{-2}$ (the constant $\kappa_\theta$ of Section~\ref{sec:graphs} for these charts), with
$|\kappa_{\theta,i}|\le|K_i+2|\lambda_i^{-2}\le(\lambda_ir^{(i)}_1)^{-2}\le4s_i^{-2}\to0$
by \eqref{eq:r1}. The upper ends $T_i\ge0$ converge, along a subsequence, to
$T\in[0,\infty]$, and $|\theta_i|\le C_\theta$. Hence Proposition~\ref{prop:compact}
applies with $\cU=\Hn\times\R$ and $L=12$ on every compact set (for $i$ large). Along a subsequence, $\cM'_i\to\cM$ locally, where $\cM$ is a smooth model flow in trace form in
$\Hn\times\R$ with upper end $T$, i.e.\ an ancient model flow on $(-\infty,T]$ in the sense of Definition~\ref{def:modelflow}. Moreover, $0\in\cM$ since $0\in\cM'_i$ for all $i$, and by
Proposition~\ref{prop:compact}(ii$'$) (no margin, $\cU$ being the whole space),
$\Kreg(\cM,X)\le12$ for every $X\in\cM$.

\emph{Step 3: the limit has density at most one.} Let $X=(q,t)\in\cM$ and $\sigma>0$, and
let $Z_i=(\hat y_i,t'_i)\in\cM'_i$ with $Z_i\to X$ (the centres $X_i=(q_i,t_i)$ of Step 2
keep their meaning). Let $\hat Y_i=(p'_i,\hat t_i)\in\cM_i$ be the point of the original
flow corresponding to $Z_i$; thus $p'_i\in M_i$, $\hat y_i=\delta_{\lambda_i}\Upsilon_{q_i}(p'_i)\in\Hn$
and $t'_i=\lambda_i^2(\hat t_i-t_i)$. Let $I'_{Z_i}(\sigma)=\int\psi_i\rho^{\chi,i}_{Z_i}d\mu^i_{t'_i-\sigma}$
be the integral on the right of Lemma~\ref{lem:lsc}, with $\rho^{\chi,i}_{Z_i}$ the kernel of
the centre $Z_i$ for the structure $\varsigma_i$ and $\psi_i:=\psi_{p'_i}$ written in the
ambient chart. For $i$ large, $Z_i\in\cC_{R_i}(0)$ and $\sigma<R_i^2$, so by
Lemma~\ref{lem:pick}(c) and $(\mathrm H_{\eps_i})$ for $F_i$,
\[
I'_{Z_i}(\sigma)=\Theta^\chi_{\hat Y_i}\big(\sigma\lambda_i^{-2};F_i|_{\cU_i\cap(U_{q_i}\times\R)}\big)
\le\Theta^\chi_{\hat Y_i}\big(\sigma\lambda_i^{-2};F_i\big)\le1+\eps_i .
\]
Next we show that $\psi_i\to1$ locally uniformly on $\Hn$. Let $\cK\subset\Hn$ be compact and
$R\ge8$ with $\cK\cup\{\hat y_i\}\subset\cB_R$ for all $i$. By Corollary~\ref{cor:rM}(c)
(with $p=q_i$, $q=p'_i$, $r=r_{M_i}$ and $\lambda=\Lambda_i\ge R/8$ for $i$ large, so that
$\delta_{\lambda/r}=\delta_{\lambda_i}$ and $\hat y=\hat y_i\in\cB_R$), the transition from the
ambient chart $\delta_{\lambda_i}\Upsilon_{q_i}$ to the rescaled chart $\delta_{\lambda_i}\Upsilon_{p'_i}$ of the moving centre
is $\delta_*(R/8)^4\Lambda_i^{-2}$-close in $C^{k_0}(\cB_R)$ to $\tilde U_i\ell_{\hat y_i}^{-1}$.
By the group law, $\ell_{\hat y_i}^{-1}(\cB_R)\subset\cB_{3R}$,
so for $i$ large the coordinates in the chart of $p'_i$ of all points of $\cK$ lie in
$\cB_{4R}$, a set independent of $i$. In the chart of $p'_i$ the cut-off
$\psi_{p'_i}=(\zeta_*\circ\Upsilon_{p'_i})^2$ equals $1$ where $|w|<\lambda_ir_0^{(i)}/2$ and
$|z|<\lambda_i^2z_0^{(i)}/2$ (rescaled coordinates), and
$\lambda_ir_0^{(i)}/2\ge32\Lambda_i\to\infty$, $\lambda_i^2z_0^{(i)}/2\ge64\bb^2\Lambda_i^2\to\infty$
by Definition~\ref{def:rM} ($r_0\ge64r_M$, $z_0\ge128\bb^2r_M^2$); hence $\psi_i=1$ on
$\cK$ for $i$ large. Lemma~\ref{lem:lsc} (with the points $Z_i\to X$, with $\sigma_i=\sigma$ and the full slab) now gives
\[
\int\rho^\chi_X\,d\mu_{t-\sigma}\le\liminf_iI'_{Z_i}(\sigma)\le1\qquad\text{for all }X\in\cM,\ \sigma>0 .
\]

\emph{Step 4: the limit is a plane.} By Step 2, $\cM$ is a non-empty ancient model flow
on $(-\infty,T]$ in trace form, with $|\theta|\le C_\theta$ and $\Kreg\le12$ on $\cM$; by
Step 3 its density is at most one at every point and scale. Proposition~\ref{prop:rigidity}
gives a Lagrangian subspace $P_0=U_0P$, a point $\hat y$ and $\bar\theta\in[-C_\theta,C_\theta]$
with $\cM(t)=\ell_{\hat y}(P_0\times\{-4\bar\theta t\})$ for $t\le T$ and $\theta\equiv\bar\theta$.
Since $0\in\cM(0)=\ell_{\hat y}(P_0\times\{0\})$, we have $\hat y\in P_0\times\{0\}$, and
$P_0\times\{0\}$ is an abelian subgroup of $\Hn$ (for $w,w'\in P_0$ the cross term
$\frac12(b\cdot x'-a\cdot y')$ of the group law is $\omega_0$-type and vanishes on a
Lagrangian subspace). Hence $\cM(t)=P_0\times\{-4\bar\theta t\}$ and
\[
\cM=\Pi:=\{(w,-4\bar\theta t,t):\ w\in P_0,\ t\le T\}.
\]
Compose all the rescaled charts with the fixed unitary rotation $\tilde U_0^{-1}$: the
ambient charts $\delta_{\lambda_i}\Upsilon_{q_i}$ and likewise the charts
$\delta_{\lambda_i}\Upsilon_p$ of all other centres $p\in M_i$. This changes neither the
hypotheses nor the conclusions of Steps 2--3. Indeed, the model structure, the gauge, the
functions $f$ and $z$, and the convergence are invariant under $\tilde U_0$; so are the
closeness bounds, the norms of Section~\ref{subsec:rM} being invariant under unitary
rotations. The admissible radii do not change, the rotation of a witness in
Definition~\ref{def:K} being arbitrary. Finally, $\tilde U_0$ is an automorphism of the
Heisenberg group, because $U_0$ preserves $\omega_0$ and hence the cross term of the group
law; so the transitions of Corollary~\ref{cor:rM}, which are now conjugated by
$\tilde U_0$, keep their form, with $U_0^{-1}UU_0$ and $\tilde U_0^{-1}\hat y$ in place of
the rotation $U$ and the centre $\hat y$ of that corollary. After this rotation,
$P_0=P=\{y=0\}$.

\emph{Step 5: one sheet and the Schauder improvement.} The hypotheses of
Lemma~\ref{lem:onesheet} hold with $C_*=12$: $0\in\cM$, $T\ge0$, $\cC_4(0)\subset\cU=\Hn\times\R$,
$\cM\cap\cC_4(0)=\Pi\cap\cC_4(0)$, and for $i$ large $\cC_4(0)\subset\cC_{R_i}(0)\subset\cU'_i$
with $\Kreg(\cM'_i,\cdot)\le12$ there. Hence, for $i$ large, $\cM'_i\cap\bB^{3/2}$ is the
Legendrian graph $\hat\Gamma_{\tilde u_i,\tilde\zeta_i}(\Omega^{3/2}_{T_i})$ of a single function
with $\|\tilde u_i\|_{2,\alpha;\Omega^{3/2}_{T_i}}\le C$, $\tilde u_i\to0$ in $C^{2,\alpha'}$ on
compact subsets of $\Omega^{3/2}_T$ (up to the upper end); moreover $\hat F_i$ is injective
on the preimage of $\bB^{3/2}$.

The bound and the convergence on compact subsets imply that
$\tilde u_i$ and $D\tilde u_i$ tend to $0$ uniformly on all of $\Omega^{3/2}_{T_i}$. For $0<\epsilon<\frac18$ let $K_\epsilon:=\{|x|\le\frac32-\epsilon,\ -\frac94+\epsilon^2\le s\le\min\{T,\frac94\}-\epsilon^2\}$,
a compact subset of $\Omega^{3/2}_T$. Every point of $\Omega^{3/2}_{T_i}$ is joined to a
point of $K_\epsilon$ by a path in $\Omega^{3/2}_{T_i}$ (radial in $x$, then in $s$) of spatial
length at most $\epsilon$ and time length at most $\epsilon^2+(T_i-T)_+$, so
$\|D\tilde u_i\|_{0;\Omega^{3/2}_{T_i}}\le\sup_{K_\epsilon}|D\tilde u_i|+C\big(\epsilon+(T_i-T)_+^{1/2}\big)^\alpha$
and $\|\tilde u_i\|_{0;\Omega^{3/2}_{T_i}}\le\sup_{K_\epsilon}|\tilde u_i|+C\big(\epsilon+\epsilon^2+(T_i-T)_+\big)$,
by the bounds on $[D\tilde u_i]_\alpha$, $\|D\tilde u_i\|_0$ and $\|\partial_s\tilde u_i\|_0$. Letting first $i\to\infty$ and then $\epsilon\to0$ gives the claim. In particular
$|\tilde u_i|<2$ and $|D\tilde u_i|\le1$ on $\Omega^{3/2}_{T_i}$ for $i$ large, which is the
range on which the coefficients of Lemma~\ref{lem:grapheq} are defined and estimated.

The projection of this Legendrian graph moves by the mean curvature flow of the transverse
metric (\cite[Lemma~2.3]{ChangWuZhang2026ContactBlowup}, whose pointwise proof applies to
local projections as explained in \cite[\S4.6]{ChangWuZhang2026ContactBlowup}), so by
Lemma~\ref{lem:grapheq} (with $\delta_i=\delta_*\Lambda_i^{-2}\to0$ in \eqref{eq:delta})
$\tilde u_i$ solves \eqref{eq:grapheq} on $\Omega^{3/2}_{T_i}$. The hypotheses of
Lemma~\ref{lem:schauder} are $\delta_i\to0$, the bound $\|\tilde u_i\|_{2,\alpha;\Omega^{3/2}_{T_i}}\le C$ and
the convergence $\tilde u_i\to0$ in $C^2$ on compact subsets. Applied with $\mu=\frac65$,
the lemma gives
\[
\|\tilde u_i\|_{2,\alpha;\Omega^{6/5}_{T_i}}\longrightarrow0 .
\]

\emph{Step 6: the contradiction.} We show that the radius $\frac98$ is admissible for
$(\cM'_i,0)$ for $i$ large, so that $\Kreg(\cM'_i,0)\le\frac89<1$, contradicting
$\Kreg(\cM'_i,0)=1$. Let $A_i:=D\tilde u_i(0,0)$, a symmetric $n\times n$ matrix (the slice
$\cM'_i(0)$ is Lagrangian at the origin for $\omega_0$, the structure at the origin of
the chart being the model) with $A_i\to0$, let $S_i:=(I+A_i^2)^{-1/2}$ (symmetric,
commuting with $A_i$), and define $U_i$ on $\C^n=\R^n\times\R^n\ni(x,y)$ by the real block
formula
\begin{gather*}
U_i(x,y):=\big(S_i(x-A_iy),\ S_i(A_ix+y)\big),\\
\text{i.e.}\quad U_i=(I-J_0\mathbf A_i)S_i\ \text{ with }\ \mathbf A_i:=\mathrm{diag}(A_i,A_i),
\end{gather*}
where $J_0(x,y)=(y,-x)$ is the complex structure of Section~\ref{sec:setting} (for which
$\zeta_j=\frac12(x_j-\sqrt{-1}\,y_j)$ are holomorphic; in the coordinates $\zeta$ the map is
$\zeta\mapsto(I-\sqrt{-1}\,A_i)S_i\zeta$). Since $\mathbf A_i$ and $S_i$ are symmetric and
commute with each other and with $J_0$, and $J_0^{\rm T}=-J_0$, we have
$(I-J_0\mathbf A_i)^{\rm T}(I-J_0\mathbf A_i)=(I+J_0\mathbf A_i)(I-J_0\mathbf A_i)=I-J_0^2\mathbf A_i^2=I+\mathbf A_i^2$,
hence $U_i^{\rm T}U_i=S_i(I+\mathbf A_i^2)S_i=I$ and $U_iJ_0=J_0U_i$, i.e.\ $U_i\in U(n)$;
moreover
\begin{gather*}
\|U_i-I\|\le\|S_i-I\|+\|J_0\mathbf A_iS_i\|\le C|A_i|\to0,\\
U_i(x,0)=(S_ix,\ A_iS_ix),\qquad\text{so}\qquad U_iP=\{(x',A_ix')\}=T_0\big(\pi\cM'_i(0)\big),
\end{gather*}
where the projected tangent space is taken on the sheet in $\bB^{3/2}$.
This sheet projects at time zero to the graph of $\tilde u_i(\cdot,0)$, with
$D\tilde u_i(0,0)=A_i$, and $S_i$ is invertible.

In the chart $\tilde U_i^{-1}$ (composed with the ambient chart) the sheet
$\{(x,\tilde u_i(x,s),\tilde\zeta_i(x,s),s):(x,s)\in\Omega^{6/5}_{T_i}\}$ has the base map
$b_i(x,s):=\pi_PU_i^{-1}(x,\tilde u_i(x,s))=x+\pi_P(U_i^{-1}-I)(x,\tilde u_i(x,s))$, with
$\|b_i-\mathrm{id}\|_{2,\alpha;\Omega^{6/5}_{T_i}}\le C\|U_i-I\|(1+\|\tilde u_i\|_{2,\alpha})\to0$,
and the $J_0P$-component $c_i(x,s):=\pi_{J_0P}U_i^{-1}(x,\tilde u_i(x,s))=\tilde u_i+\pi_{J_0P}(U_i^{-1}-I)(x,\tilde u_i)$,
with $\|c_i\|_{2,\alpha;\Omega^{6/5}_{T_i}}\to0$. For $i$ large, $b_i(\cdot,s)$ is a
diffeomorphism of $\{|x|<\frac65\}$ onto a set containing $\{|\xi|<\frac{23}{20}\}$, with
inverse $\mathrm{id}+o(1)$ in $C^{2,\alpha}$ (parabolic sense, the time derivative of the
inverse being $-(D_xb_i)^{-1}\partial_sb_i$ composed with the inverse). Hence over
$\{|\xi|<\frac{23}{20}\}\times(-\frac{36}{25},\frac{36}{25})\cap\{s\le T_i\}$ the sheet is, in the
chart $\tilde U_i^{-1}$, the Legendrian graph of $u'_i:=c_i\circ(b_i^{-1},\cdot)$ with
$\|u'_i\|_{2,\alpha}\to0$, $u'_i(0,0)=0$ and $Du'_i(0,0)=0$ (the tangent plane at the
origin is $U_iP$, which $\tilde U_i^{-1}$ maps to $P$). In the chart
$\tilde U_i^{-1}D_{8/9}$, which is $\tilde U_i^{-1}\Lambda_{0,9/8}$ for the rescaled flow,
the witness function is $u^{(i)}(\xi,s):=\frac89u'_i(\frac98\xi,\frac{81}{64}s)$ on
$\Omega_{T''}$, $T''=\frac{64}{81}T_i$, with
\[
\|u^{(i)}\|_{2,\alpha;\Omega_{T''}}\le(\tfrac98)^{1+\alpha}\|u'_i\|_{2,\alpha;\Omega^{9/8}_{T_i}}\longrightarrow0
\]
(each term of the norm is multiplied by a power $\mu^j$, $j\le1+\alpha$, under
$u\mapsto\mu^{-1}u(\mu x,\mu^2s)$ with $\mu=\frac98$). The box of this chart is
$D_{9/8}\tilde U_i(\bB)\subset\bB^{3/2}$ for $i$ large ($U_i\to I$), so
$\tilde U_i^{-1}D_{8/9}(\cM'_i)\cap\bB$ consists of points of the graph of $\tilde u_i$ with
$|x|<\frac98(1+o(1))<\frac65$, which are points of the graph of $u^{(i)}$ over $\Omega_{T''}$. The injectivity condition holds by Lemma~\ref{lem:onesheet}. The graph is Lagrangian for the rotated structure, with a $C^1$ height vanishing at the origin. Finally, $\frac98\le\lambda_ir_1^{(i)}$,
the bound on admissible radii in the rescaled picture (Lemma~\ref{lem:admissible}(c)).
Thus $\frac98$ is admissible for $(\cM'_i,0)$ for $i$ large, which gives the announced contradiction.
Hence there are $\eps$ and $C$ as claimed; they depend only on $n,a,\alpha,C_\theta$, since
the sequence was arbitrary apart from these data.
\end{proof}

\begin{remark}
(1) The constants are not explicit; as in \cite{wh} the proof is by compactness. The
exhaustion of Lemma~\ref{lem:exhaust} requires $C\ge3$, which is no restriction.
(2) The density hypothesis is used in Step 1, through Lemma~\ref{lem:exhaust}(a),(b) and
Lemma~\ref{lem:density-M} (finiteness and no double points), and in Step 3, through
Lemma~\ref{lem:pick}(c) and Lemma~\ref{lem:lsc}. The normalisation of Step 2 uses
$d_1=\min\{d,r_1\}$ so that
$\lambda_ir_{M_i}\to\infty$ even when the manifolds $M_i$ vary.
(3) Applying Theorem~\ref{thm:A} to the flow truncated at the time of the
centre gives the past estimate, as shown in the proof of Corollary~\ref{cor:past}.
\end{remark}

\subsection{The model version and the gap theorem}\label{subsec:gap}

The proofs of Sections~\ref{sec:setting}--\ref{sec:proof} use the ambient manifold only
through Fact~\ref{fact:charts}, Fact~\ref{fact:kernel}, the constants
$r_0,z_0,r_*,C^{\rm tr},C_\sharp,k_0$ entering Definition~\ref{def:rM}, the cut-off
$\zeta_*$, and the smoothness and properness of the flows. The compactness of $M$ enters
only through the proofs of these facts in \cite{PaperI}. The inputs for closed flows,
recalled at the beginning of Section~\ref{sec:terminal}, are not used before that section.

For the model $(\Hn,\overline\eta,\overline\xi,\overline\Phi,\overline g)$ we use the charts
$\Upsilon_p:=\ell_p^{-1}$, which are defined on all of $\Hn$, so that $r_0=z_0=\infty$. In
these charts the identities and expansions of Fact~\ref{fact:charts}(i)--(iii) hold exactly, with
$h_p=g_0$, $\beta_p=\overline\beta$ and vanishing errors; Fact~\ref{fact:kernel}(i)
concerns the model itself, and (ii) holds with $\hat a\equiv4$. The transitions of Lemma~\ref{lem:rM}(b) are the left translations $\ell_{\hat y}^{-1}$
themselves, so $C^{\rm tr}=0$ with any $r_*$, and $C_\sharp=0$. We take $\zeta_*\equiv1$, no compact support being needed for global charts, so that $\psi_q\equiv1$. Definition~\ref{def:rM} then allows $r_M:=1$, and \eqref{eq:r1}
gives $r_1=1$ since $K=-2$.

However, not all the uniform bounds of Fact~\ref{fact:charts}(i) and of the paragraph
following it extend to the whole of these charts. Those on $\beta_p$ itself and on the
coefficients of the metric and of its inverse do not: the
coefficients of $\overline\beta$ are linear in $w$, and those of $\overline g$ and of its
inverse are unbounded (for instance
$\overline g(\partial_{x_j},\partial_{x_j})=\frac14+\frac1{16}y_j^2$ and
$\overline g^{-1}(dz,dz)=|w|^2+4$). These coefficients are the same in every chart
$\ell_p^{-1}$, and on each $\cB_R$ they are bounded, with all their derivatives, by
constants depending only on $n$ and $R$; thus the coefficient bounds hold for the model
on bounded chart regions only, while the bounds on $\overline\nabla\zeta_q$ and
$\overline\nabla^2\zeta_q$ derived from them are trivial, as $\zeta_q\equiv1$. This
suffices for the statements of Sections~\ref{sec:setting}--\ref{sec:proof} that are used
in the proof of Theorem~\ref{thm:A} and in the rest of this section: their proofs use the
charts only through the exact identities above and the group law, or on bounded chart
regions. With these choices they hold for the model, with
$d_1=\min\{d,1\}$ and the density $\int\rho^\chi_{q,t}d\mu_{t-\sigma}$ without cut-off;
in particular the proof of Theorem~\ref{thm:A} applies to local flows in open subsets of
$\Hn\times\R$ with this $r_1=1$. Since that proof is by contradiction, it gives its own constants:

\begin{claim}
Let $0<\alpha<1$. There are $\eps_{\rm m}=\eps_{\rm m}(n,a,\alpha,C_\theta)\in(0,\frac12]$ and
$C_{\rm m}=C_{\rm m}(n,a,\alpha,C_\theta)\in[3,\infty)$ (any smaller $\eps_{\rm m}$ and any larger
$C_{\rm m}$ then have the same property) such that every smooth local Legendrian flow $F$ in an
open set $\cU\subset\Hn\times\R$ for the model structure, with any upper end and
$|\theta|\le C_\theta$, satisfying $\int\rho^\chi_{q,t}\,d\mu_{t-\sigma}\le1+\eps_{\rm m}$ for all
$X=(q,t)\in\cM$ and $0<\sigma<\min\{d(X,\cU),1\}^2$, has
$\min\{d(X,\cU),1\}\,\Kreg^{(1)}(\cM,X)\le C_{\rm m}$ for all $X\in\cM$, where $\Kreg^{(1)}$ is the
regularity scale with $r_1=1$.
\end{claim}
\begin{proof}
We run Steps 1--6 of the proof of Theorem~\ref{thm:A} for a contradiction sequence of model flows
$F_i$ in $\cU_i\subset\Hn\times\R$ with the chart data above. Step 1 uses
Lemma~\ref{lem:exhaust}(a)--(c) and Lemma~\ref{lem:density-M}. Their proofs use the ambient
manifold only through the inputs listed before the Claim (Fact~\ref{fact:charts}, the
charts $\Upsilon_q$ with their domains $U_q$, Fact~\ref{fact:kernel}, the constants
$r_0,z_0,r_*,C^{\rm tr},C_\sharp,k_0$ of Definition~\ref{def:rM} and \eqref{eq:r1}, the
cut-off $\zeta_*$, and the smoothness and properness of the flows), all of which are available in the
model in the form described there, and through Lemma~\ref{lem:quasi},
Lemma~\ref{lem:rM}(b), Corollary~\ref{cor:rM}(a), Lemma~\ref{lem:admissible},
Lemma~\ref{lem:finite}, Proposition~\ref{prop:compact} and Lemma~\ref{lem:lsc}, which are
statements of Sections~\ref{sec:setting}--\ref{sec:proof} and hold for the model as noted
above. In Step 2,
Lemma~\ref{lem:pick} with $r_M=1$ gives $\lambda_i\ge s_i/2\to\infty$; the ambient charts are the
model charts $\delta_{\lambda_i}\ell_{q_i}^{-1}$, so the closeness hypothesis of
Proposition~\ref{prop:compact} holds with zero errors and $\kappa_{\theta,i}=0$. In Step 3,
Lemma~\ref{lem:pick}(c) and Lemma~\ref{lem:lsc} apply with $\psi_i\equiv1$. Steps 4 and 5 use
Proposition~\ref{prop:rigidity}, which concerns the model, and Lemmas~\ref{lem:onesheet} and~\ref{lem:schauder}, which
hold for the model as noted above. Step 6 is unchanged: its last condition $\frac98\le\lambda_ir_1$,
the bound on admissible radii in the rescaled picture (Lemma~\ref{lem:admissible}(c)), reads
$\frac98\le\lambda_i$ and holds for $i$ large, so $\Kreg(\cM'_i,0)\le\frac89<1$, contradicting
Lemma~\ref{lem:pick}(a). The argument does not use Theorem~\ref{thm:A} or
Theorems~\ref{thm:A0}--\ref{thm:alpha}.
\end{proof}

\noindent For $0<\alpha<1$ let $\eps_{\rm cpt}=\eps_{\rm cpt}(n,a,\alpha,C_\theta)\in(0,\frac12]$ and
$C_{\rm cpt}=C_{\rm cpt}(n,a,\alpha,C_\theta)\in[3,\infty)$ be constants with which
Theorem~\ref{thm:A} holds. From now on $\eps$ and $C$ denote the constants
$\min\{\eps_{\rm cpt},\eps_{\rm m}\}$ and $\max\{C_{\rm cpt},C_{\rm m}\}$, and this pair is
meant whenever the constants of Theorem~\ref{thm:A} are referred to. By the monotonicity
clauses, Theorem~\ref{thm:A} and the Claim both hold with them, and $\eps\le\frac12$,
$C\ge3$. In the model, however, we use Definition~\ref{def:K} with $r_1=\infty$, and the
natural statement is scale invariant.

\begin{theorem}[local regularity in the model]\label{thm:A0}
Let $0<\alpha<1$ and let $\eps=\eps(n,a,\alpha,C_\theta)$ and $C=C(n,a,\alpha,C_\theta)$
be the constants fixed above, $\eps=\min\{\eps_{\rm cpt},\eps_{\rm m}\}$ and
$C=\max\{C_{\rm cpt},C_{\rm m}\}$. Let $\cU\subset\Hn\times\R$ be open and let $F$
be a smooth local Legendrian flow in $\cU$ for the model structure ($K=-2$), with any
upper end and $|\theta|\le C_\theta$, such that
\[
\Theta^\chi_X(\sigma):=\int\rho^\chi_{q,t}\,d\mu_{t-\sigma}\le1+\eps\qquad\text{for all }X=(q,t)\in\cM\text{ and all }0<\sigma<d(X,\cU)^2 .
\]
Then, with $\Kreg$ the regularity scale of Definition~\ref{def:K} for the exponent
$\alpha$ and $r_1=\infty$: $d(X,\cU)\,\Kreg(\cM,X)\le C$ for every $X\in\cM$ with
$d(X,\cU)<\infty$, and every radius is admissible for $(\cM,X)$ when $d(X,\cU)=\infty$.
\end{theorem}
\begin{proof}
Let $\Kreg^{(1)}$ be the regularity scale with $r_1=1$, i.e.\ with the admissible radii
restricted to $r\le1$; by Lemma~\ref{lem:admissible}(a), $\Kreg^{(1)}=\max\{\Kreg,1\}$. By the
Claim above, which applies since $\eps\le\eps_{\rm m}$ and $\min\{d(X,\cU),1\}\le d(X,\cU)$,
and by $C_{\rm m}\le C$, we have
$\min\{d(X,\cU),1\}\Kreg^{(1)}(\cM,X)\le C$ for all $X\in\cM$; moreover $C\ge3$.

Let $X\in\cM$ with $0<d(X,\cU)<\infty$ and put $\lambda:=1/d(X,\cU)$. The dilated flow
$D_\lambda F$ in $D_\lambda\cU$, $D_\lambda(q,t)=(\delta_\lambda q,\lambda^2t)$, is a smooth
local Legendrian flow for the model structure with the same angle bound. Its densities
are those of $F$ (Fact~\ref{fact:kernel}(ii) with the model structure, where
$\hat a\equiv4$), and $d(D_\lambda X,D_\lambda\cU)=\lambda d(X,\cU)=1$ by the homogeneity of
the gauge. Hence the dilated flow satisfies the hypothesis, and
$\Kreg^{(1)}(D_\lambda\cM,D_\lambda X)\le C$.
By Lemma~\ref{lem:admissible}(c), $\Kreg(D_\lambda\cM,D_\lambda X)=\lambda^{-1}\Kreg(\cM,X)$,
hence $\max\{\lambda^{-1}\Kreg(\cM,X),1\}\le C$ and $d(X,\cU)\Kreg(\cM,X)\le C$.

If $d(X,\cU)=\infty$, let $\varrho>0$ and consider the restriction of $F$ to
$\cC_\varrho(X)$ (Remark~\ref{rem:localflow}(3)). It satisfies the hypothesis in
$\cC_\varrho(X)$, its densities being at most those of $F$, and $d(X,\cC_\varrho(X))=\varrho$. By the first part, $\Kreg(\cM\cap\cC_\varrho(X),X)\le C/\varrho$, so every $r<\varrho/C$ is
admissible for the restriction (Lemma~\ref{lem:admissible}(a)), and then also for
$(\cM,X)$, since for $r<\varrho/3$ the box $\bB_{X,r}\subset\cC_{3r}(X)$ lies in
$\cC_\varrho(X)$, where the traces coincide (as in the proof of Lemma~\ref{lem:exhaust}(c)).
Letting $\varrho\to\infty$ shows that every radius is admissible.
\end{proof}

\begin{theorem}[gap theorem for ancient model flows]\label{thm:gap}
Put $\bar\eps:=\eps(n,a,\frac12,C_\theta)$, the constant $\eps$ fixed before
Theorem~\ref{thm:A0}, for the exponent $\alpha=\frac12$. Let $F$ be a non-empty ancient model flow on $(-\infty,T_*]$,
$T_*\in(-\infty,\infty]$, in the sense of Definition~\ref{def:localflow}, or in trace form
with $\Kreg$, for some exponent in $(0,1)$, bounded on compact time intervals (then, by
Lemma~\ref{lem:param}, it is a flow in the sense of Definition~\ref{def:localflow}), with
$|\theta|\le C_\theta$ and
\[
\int\rho^\chi_X\,d\mu_{t-\sigma}\le1+\bar\eps\qquad\text{for all }X=(q,t)\in\cM\text{ and all }\sigma>0 .
\]
Then there are a Lagrangian subspace $P\subset\C^n$, a point $\hat y\in\Hn$ and
$\bar\theta\in[-C_\theta,C_\theta]$ with
\[
\cM(t)=\ell_{\hat y}\big(P\times\{-4\bar\theta t\}\big)\quad\text{for all }t\le T_*,\qquad\theta\equiv\bar\theta :
\]
the flow is a single Legendrian plane of multiplicity one, translating vertically, its
height $z$ changing at the constant rate $-4\bar\theta$.
\end{theorem}
\begin{proof}
Throughout, $\Kreg$ denotes the regularity scale for the exponent $\frac12$, and $C$ is the
constant fixed before Theorem~\ref{thm:A0}, taken for the same exponent; thus
Theorem~\ref{thm:A0} holds for the exponent $\frac12$ with the constants $\bar\eps$ and $C$.

\emph{Step 1: every radius is admissible.} Fix a finite $t_*\le T_*$ and let $F^{\le t_*}$ be the
truncation at $t_*$ (Remark~\ref{rem:localflow}(2); if $T_*$ is finite, $t_*=T_*$ is
allowed, and then $F^{\le t_*}=F$), a smooth local Legendrian flow in
$\Hn\times\R$ with upper end $t_*$ whose densities at its points are those of $F$. For
$X\in\cM^{\le t_*}$ and $\varrho>0$, Theorem~\ref{thm:A0} applied to the restriction of
$F^{\le t_*}$ to $\cC_\varrho(X)$ (densities at most those of $F$, $d(X,\cC_\varrho(X))=\varrho$)
gives $\Kreg(\cM^{\le t_*}\cap\cC_\varrho(X),X)\le C/\varrho$, and, as in the proof of
Theorem~\ref{thm:A0}, every $r<\varrho/\max\{C,3\}$ is admissible for $(\cM^{\le t_*},X)$.
As $\varrho$ is arbitrary, every radius is admissible for $(\cM^{\le t_*},X)$, and by
Lemma~\ref{lem:density-M} the map $\hat F^{\le t_*}$ is injective.

\emph{Step 2: the slices are unions of planes.} Fix $X=(q,t)\in\cM$ with $t\le t_*$. For
$R\ge1$ let $(U_R,u^{(R)},\zeta^{(R)})$ be a witness of $(\cM^{\le t_*},X)$ at the radius
$R$, and let $u^R(x,s):=R\,u^{(R)}(x/R,s/R^2)$, $\zeta^R(x,s):=R^2\zeta^{(R)}(x/R,s/R^2)$ be
the graph functions in the chart $\tilde U_R\Lambda_{X,1}$, over the base
$\{|x|<R,\ |s|<R^2,\ s\le t_*-t\}$. From $\|u^{(R)}\|_{2,1/2}\le1$, $u^{(R)}(0,0)=0$ and
$Du^{(R)}(0,0)=0$ we get on the unit base $\Omega_{t_*-t}$
\begin{gather*}
|Du^R(x,s)|=\big|Du^{(R)}(x/R,s/R^2)-Du^{(R)}(0,0)\big|\le\|(x/R,s/R^2)\|^{1/2}\le R^{-1/2},\\
|u^R(x,s)|\le|s|\,\|\partial_su^R\|_0+|x|\,\|Du^R\|_0\le R^{-1}+R^{-1/2},
\end{gather*}
since $\partial_su^R=R^{-1}(\partial_su^{(R)})(x/R,s/R^2)$; moreover
$|D_x\zeta^R|=\sup_{|e|=1}|2\overline\beta((e,Du^Re))|\le\frac12(|u^R|+|x||Du^R|)\to0$ by
\eqref{eq:height} with the model form $\overline\beta$, and $|\partial_s\zeta^R|\le1+4C_\theta$
by Lemma~\ref{lem:heightangle}. The cylinder $\cC_1(X)$ lies in the box of every chart
$\tilde U_R\Lambda_{X,1}$, so by Definition~\ref{def:K}(i) and Lemma~\ref{lem:filling}
(applicable since $\cC_{3R}(X)\subset\Hn\times\R$),
\[
\cM^{\le t_*}\cap\cC_1(X)=\hat\Gamma_{u^R,\zeta^R}\cap\cC_1(X)\qquad\text{for every }R\ge1 .
\]

By Arzel\`a--Ascoli there is a sequence $R\to\infty$ with $U_R\to U_\infty$ and
$\zeta^R\to\zeta_\infty$ uniformly on the unit base. The bound $|\zeta^R|\le4C_\theta+o(1)$
on the unit base is obtained as follows. By \eqref{eq:theta-from-graph} in the model, where
$\overline\beta_{(x,u)}((0,\partial_su))=\frac14x\cdot\partial_su$ vanishes at $x=0$, we have
$\partial_s\zeta^R(0,s)=-4\theta(0,s)$ with $\theta$ the Legendrian angle at the sheet point
over $(0,s)$. Hence $|\zeta^R(0,s)|\le4C_\theta|s|<4C_\theta$ by $\zeta^R(0,0)=0$ and
$|\theta|\le C_\theta$, and $|\zeta^R(x,s)-\zeta^R(0,s)|\le|x|\sup|D_x\zeta^R|=o(1)$ by the
bound on $D_x\zeta^R$ above. Along this sequence the graphs converge uniformly on the unit base to
$\hat\Gamma_{0,\zeta_\infty}$ with $\zeta_\infty=\zeta_\infty(s)$ and $|\zeta_\infty|\le4C_\theta<\bb^2$.
Here the graphs $\hat\Gamma_{u^R,\zeta^R}$ are written in the charts $\tilde U_R\Lambda_{X,1}$
and are compared in the fixed coordinates $\Lambda_{X,1}$, i.e.\ pulled back by
$\tilde U_R^{-1}$: their points $(U_R^{-1}(x,u^R(x,s)),\zeta^R(x,s),s)$ converge uniformly
to $(U_\infty^{-1}(x,0),\zeta_\infty(s),s)$, and $\hat\Gamma_{0,\zeta_\infty}$ is
accordingly read in the chart $\tilde U_\infty\Lambda_{X,1}$ (the set
$\Lambda_{X,1}(\cC_r(X))=\{|w|<r,\ |z|<\bb^2r^2\}\times(-r^2,r^2)$ is invariant under every
$\tilde U$, so the cylinders below have the same description in all these charts). This
limit is the lift of the plane $U_\infty^{-1}P$ through the origin in the coordinates
$\Lambda_{X,1}$ (equivalently, of $P$ in the chart $\tilde U_\infty\Lambda_{X,1}$),
translating vertically.

We now restrict to the smaller cylinder $\cC_{1/2}(X)$, whose points have base coordinates
in a compact subset of the unit base. After this restriction, the left-hand side of the
displayed identity above, which does not depend on $R$, equals
$\hat\Gamma_{0,\zeta_\infty}\cap\cC_{1/2}(X)$. Indeed, a point of $\cM^{\le t_*}\cap\cC_{1/2}(X)$ lies
on every graph and hence on the uniform limit. Conversely, a point of
$\hat\Gamma_{0,\zeta_\infty}\cap\cC_{1/2}(X)$ is the limit of the points of the graphs with
the same base coordinates, which belong to the fixed trace $\cM^{\le t_*}$ and lie in
$\cC_1(X)$ for large $R$, and $\cM^{\le t_*}$ is relatively closed in the open cylinder
$\cC_1(X)$ (the trace being closed in $\Hn\times\R$). Thus every slice
$\cM(t)$, $t\le t_*$, is a closed embedded Legendrian submanifold of $\Hn$ (properness and
Step 1) which coincides, near each of its points, with an affine horizontal Legendrian
plane; each of its connected components is therefore open and closed in such a plane,
i.e.\ a whole affine Legendrian plane $\ell_{\hat y_j}(P_j\times\{z_j\})$.

\emph{Step 3: constant angle and vertical translation.} On an affine Legendrian plane we
have $H=0$, so $\nabla\theta=-\overline\Phi H=0$ and $\theta$ is constant on each component of
each slice. The angle equation $(\partial_t-\Delta)\theta=0$ then gives $\partial_t\theta=0$
along the normal parametrisation, and $\partial_tF=H-2\theta\overline\xi=-4\theta\partial_z$
shows that each plane keeps its horizontal position and translates vertically, its height
$z$ changing at the constant rate $-4\bar\theta_j$ determined by its angle. The components correspond across times
through the global normal parametrisation of Lemma~\ref{lem:param}, available since Step 1
bounds the regularity scale on $\cM^{\le t_*}$. Hence, for $t\le t_*$,
$\cM(t)=\bigcup_j\ell_{\hat y_j}(P_j\times\{z_j-4\bar\theta_jt\})$.

\emph{Step 4: a single plane.} Fix $X=(q,t)\in\cM$ on the plane $j=1$ and $\sigma>0$. On
the plane $k$ at the time $t-\sigma$, the function $z_X=z\circ\ell_X^{-1}$ equals its
value at the time $t$ plus $4\bar\theta_k\sigma$ (vertical translation; $z\circ\ell_X^{-1}$
differs from $z$ by a function of $w$), so the caloric height is
$u=z_X-4\sigma\theta=z_X(\cdot,t)=c_k+\ell_k(w)$, a constant plus a linear function of
$w$ (the cross term of the group law), and $v=u/\sigma$. As $\sigma\to\infty$,
$\chi(v)\to1$ uniformly on $\{|w|\le A\sqrt\sigma\}$ for every $A$, while
$\int_{\text{plane }k}\rho_X\,d\mu=e^{-\mathrm{dist}^2/(16\sigma)}\to1$, where $\mathrm{dist}$
is the Euclidean distance in the $w$-coordinates from the centre to the affine Lagrangian
plane (the Gaussian integral over such a plane is explicit, the area element being
$2^{-n}$ times the Lebesgue measure). Hence $\int_{\text{plane }k(t-\sigma)}\rho^\chi_Xd\mu\to1$ for every $k$.
By Fatou's lemma for sums, the number of planes is at most
$\liminf_{\sigma\to\infty}\int\rho^\chi_Xd\mu_{t-\sigma}\le1+\bar\eps<2$, by the
normalisation $\bar\eps\le\frac12$ of Theorem~\ref{thm:A}. So there is one
plane, of multiplicity one by Step 1: $\cM(t)=\ell_{\hat y}(P\times\{z_0-4\bar\theta t\})$
for $t\le t_*$, and $z_0$ is absorbed into $\hat y$ by a vertical translation. Since the
finite $t_*\le T_*$ was arbitrary and the slices at a common time agree, this holds for every
finite $t\le T_*$. If $T_*$ is finite, the slice $t=T_*$ is covered by the choice $t_*=T_*$,
for which $F^{\le t_*}=F$ and Step~1 gives the injectivity of $\hat F$ on the whole
parameter domain, the slice $t=T_*$ included; so the multiplicity is one there as well.
\end{proof}

\begin{theorem}[a threshold independent of $\alpha$]\label{thm:alpha}
Let $\bar\eps$ be as in Theorem~\ref{thm:gap}. For every $\alpha\in(0,1)$ there is
$C_\alpha=C(n,a,\alpha,C_\theta)$ such that Theorems~\ref{thm:A} and~\ref{thm:A0} hold
with $\eps:=\bar\eps/2$ and $C:=C_\alpha$, the regularity scale being the one of
Definition~\ref{def:K} for the exponent $\alpha$.
\end{theorem}
\begin{proof}
Suppose the statement for Theorem~\ref{thm:A} fails for some $\alpha$. Then there are
flows $F_i$ in $\cU_i\subset M_i\times\R$ satisfying $(\mathrm H_{\bar\eps/2})$ with
$\Kreg_{\cU_i}(\cM_i)\to\infty$ (for the exponent $\alpha$). We run the proof of
Theorem~\ref{thm:A} for these flows with the fixed threshold $\eps_i:=\bar\eps/2$ in place
of $\eps_i\to0$. Steps 1 and 2 are unchanged and give an ancient model flow $\cM$ in
trace form with $0\in\cM$ and $\Kreg(\cM,\cdot)\le12$ for the exponent $\alpha$, which by
Lemma~\ref{lem:param}, applied with this exponent, is a model flow in the sense of
Definition~\ref{def:localflow}. Step 3 gives
$\int\rho^\chi_Xd\mu_{t-\sigma}\le1+\bar\eps/2\le1+\bar\eps$ for all $X\in\cM$ and
$\sigma>0$. Theorem~\ref{thm:gap} shows that $\cM$ is a single translating Legendrian plane
of multiplicity one, which passes through the origin at the time $0$; exactly as in Step 4
this gives $\cM=\Pi$ after a fixed rotation. Steps 5 and 6 (Lemma~\ref{lem:onesheet} with
$C_*=12$ and Lemma~\ref{lem:schauder} for the exponent $\alpha$) then yield the
contradiction $\Kreg(\cM'_i,0)\le\frac89$. The same argument, applied to sequences of
model flows with $\Kreg^{(1)}$ and $d_1=\min\{d,1\}$ as in the proof of the Claim of
Section~\ref{subsec:gap} (whose Step~6 uses $\lambda_i\to\infty$), gives a constant for the model
with $r_1=1$ and the threshold $\bar\eps/2$. With $C_\alpha$ the larger of the two constants,
both statements hold, and the dilation argument of the proof of Theorem~\ref{thm:A0} gives
Theorem~\ref{thm:A0} with $(\bar\eps/2,C_\alpha)$.
\end{proof}

\begin{remark}
The threshold $\bar\eps$ is not claimed to be optimal, and $C_\alpha$ may depend on
$\alpha$. The gap theorem uses the truncations at finite $t_*\le T_*$ and the exhaustion by the
cylinders $\cC_\varrho(X)$ in place of a global statement for flows in $\Hn\times\R$. For flows
in the sense of Definition~\ref{def:localflow}, no growth or regularity hypothesis is assumed
beyond the density bound, since Theorem~\ref{thm:A0} provides the regularity scale. For flows in
trace form, the theorem assumes in addition that $\Kreg$, for some exponent, is bounded on
compact time intervals, which by Lemma~\ref{lem:param} makes them flows in the sense of
Definition~\ref{def:localflow}.
\end{remark}

\section{Terminal regularity for closed flows}\label{sec:terminal}

At a point $(p,T)$ of the first singular time, Theorem~\ref{thm:B} starts with a
density bound at that one point. To apply Theorem~\ref{thm:A} we need density bounds
at all nearby flow points. Following \cite[Remark~3.3, Theorem~3.5]{wh}, we apply the
local estimate to flows truncated at times $T_i\uparrow T$, in a cylinder about
$(p,T)$ which crosses time $T$. In the Euclidean setting, continuity gives the bound
at nearby centres at one fixed scale. Monotonicity then gives it at smaller scales. Here the density at a centre is defined through
the contact normal chart of that centre, and we use no continuity of these charts in the
centre. Lemma~\ref{lem:compare-centres} replaces this continuity by a comparison, at a
fixed small scale, with an auxiliary density built from the chart of a fixed point. Under the Type~I
and ambient hypotheses of Corollary~\ref{cor:Bprime}, the terminal density at
a curvature singularity is the total multiplicity of the lifted tangent planes.

Throughout the section, $F\colon\Sigma\times[0,T)\to M$ is a closed flow satisfying (S) and (F) of
\cite[\S2.2]{PaperI}: $M$ is a compact Sasakian manifold, $\Sigma$ a closed $n$-manifold,
$F$ a smooth normally parametrised family of Legendrian immersions satisfying
\eqref{eq:F} with $T<\infty$, and $C_\theta=\sup|\theta|<\infty$. For such flows we use the
following inputs from \cite{PaperI}. The first is the projected area ratio bound of
\cite[(D) and Lemma~2.4]{PaperI}: there is $D_0$, depending on $(M,\eta,g)$, $T$,
$C_\theta$ and $F_0$, such that
\begin{equation}\label{eq:D}
\mu_t\big(F_t^{-1}(\Upsilon_q^{-1}\{|w|<R\})\big)\le D_0R^n\qquad\text{for all }q\in M,\ 0<R\le r_0,\ t\in[0,T).
\end{equation}

\emph{Almost monotonicity} \cite[Proposition~3.5, Lemma~3.6]{PaperI}. Put
\begin{equation}\label{eq:tau0}
\tau_0:=\min\big\{1,\ |K+2|^{-1},\ z_0/(32C_\theta+1)\big\}\in(0,1],
\end{equation}
where the term $|K+2|^{-1}$ is omitted if $K=-2$. The corresponding constant in
\cite[\S2.6]{PaperI} is $\min\{\tau_0,T\}$. For a centre $(q,t_*)$ with $0<t_*\le T$,
the window $\min\{\tau_0,t_*\}$ in \eqref{eq:almost-mono} and \eqref{eq:C5} therefore
agrees with that of \cite{PaperI}. The tail estimate \eqref{eq:tail} below also allows
centre times beyond $T$. There is
$C_2$, depending only on $n,a,K,C_\theta$, the constant $D_0$ of \eqref{eq:D}, $\mu_0(\Sigma)$
and $(M,\eta,g)$, such that for every centre $(q,t_*)$, $0<t_*\le T$, and
$t_*-\min\{\tau_0,t_*\}\le t<t_*$,
\begin{equation}\label{eq:almost-mono}
\frac{d}{dt}\int\psi_q\rho^\chi_{q,t_*}\,d\mu_t\le-\frac34\int\psi_qD^\chi\rho_{q,t_*}\,d\mu_t+C_2\Big(1+\tau^2\int\psi_q|H|^2\rho_{q,t_*}\,d\mu_t\Big),
\end{equation}
the last term being integrable in time; hence
$\Theta^\chi(q,t_*):=\lim_{t\to t_*}\int\psi_q\rho^\chi_{q,t_*}d\mu_t$ exists. Moreover,
with $C_5$ independent of the centre,
\begin{equation}\label{eq:C5}
\Theta^\chi_{q,t_*}(\sigma')\le\Theta^\chi_{q,t_*}(\sigma)+C_5\sigma\qquad\text{for }0<\sigma'\le\sigma\le\min\{\tau_0,t_*\},
\end{equation}
where $\Theta^\chi_{q,t_*}(\sigma)=\int\psi_q\rho^\chi_{q,t_*}d\mu_{t_*-\sigma}$.

\emph{Tails} \cite[Lemma~3.7]{PaperI}. For $q\in M$, $t>0$ and $0<\sigma\le\min\{\tau_0,t\}$
with $t-\sigma<T$ (the centre time $t$ may be $\ge T$; the slice $t-\sigma$ lies in $[0,T)$,
and $\Theta^\chi_{q,t}(\sigma)=\int\psi_q\rho^\chi_{q,t}\,d\mu_{t-\sigma}$ as in \eqref{eq:C5}),
and for $R,Z>0$, the part of $\Theta^\chi_{q,t}(\sigma)$ coming from
$\{f_q\ge R^2\}\cup\{|z_q|\ge Z\}$ is at most
\begin{multline}\label{eq:tail}
\tail(R,Z;\sigma):=C(n)D_0e^{-R^2/(2\sigma)}\\
+(4\pi\sigma)^{-n/2}D_0(4R)^n\exp\big(-a(Z/\sigma-8C_\theta)_+^2\big).
\end{multline}
For $\sigma\le T$ this is \cite[Lemma~3.7]{PaperI}, whose window is
$0<\sigma\le\min\{\tau_0,T,t\}$ with $t-\sigma<T$. For $T<\sigma\le\min\{\tau_0,t\}$,
which can occur only when $t>T$ and $T<\tau_0$, the bound follows from the proof of
\cite[Lemma~3.7]{PaperI}: that proof uses only the area ratio bound (D) of
\cite[Lemma~2.4]{PaperI} on the slice $t-\sigma\in[0,T)$, the bound $|\theta|\le C_\theta$
and $a_K(\sigma)\le8\sigma$, the last because $|K+2|\sigma\le1$ for $\sigma\le\tau_0$; it
uses neither $t<T$ nor $\sigma\le T$.

The trace of the flow is $\cM:=\{(F_t(x),t):x\in\Sigma,\ 0<t<T\}\subset M\times(0,T)$ and
$\overline\cM$ is its closure in $M\times\R$. The \emph{terminal density} at $(p,T)$ is
$\Theta^\chi(p,T)=\lim_{\sigma\to0}\Theta^\chi_{p,T}(\sigma)$, with
$\Theta^\chi_{p,T}(\sigma)=\int\psi_p\rho^\chi_{p,T}\,d\mu_{T-\sigma}$
(see \eqref{eq:almost-mono}). A point $p\in M$ is a \emph{curvature singularity} of $F$ at
$T$ if for every neighbourhood $V$ of $p$ in $M$ the quantity
$\sup\{|B|(x,t):\ F_t(x)\in V\}$ is unbounded as $t\to T$.

\subsection{Densities at nearby centres at a fixed scale}\label{subsec:compare-centres}

The density $\Theta^\chi_{q,t}(\sigma)$ is defined through the chart of $q$, and we do not
use any continuity of the charts in their centre. Instead we compare it, at scales which
are small compared with $r_M$ but fixed, with an auxiliary density built from the chart of
a fixed point $p$. Fix $p\in M$ and $0<r\le r_M$. For $q\in\Upsilon_p^{-1}(\delta_r\cB_1)$,
i.e.\ $|w_p(q)|<r$ and $|z_p(q)|<\bb^2r^2$, put $\hat y:=\Upsilon_p(q)$ and
\begin{gather*}
\tilde f_q:=f_{\hat y}\circ\Upsilon_p,\qquad \tilde z_q:=z_{\hat y}\circ\Upsilon_p,\qquad
\tilde\rho_q:=(4\pi\sigma)^{-n/2}e^{-\tilde f_q/\sigma},\qquad
\tilde\rho^\chi_{q,t}:=\tilde\rho_q\,\chi\Big(\frac{\tilde z_q-a_K(\sigma)\theta}{\sigma}\Big),\\
\tilde\Theta_{q,t}(\sigma):=\int\psi_p\,\tilde\rho^\chi_{q,t}\,d\mu_{t-\sigma}\qquad(\sigma\le t<T+\sigma),
\end{gather*}
where the fixed chart pulls back the model potentials
$f_{\hat y}=f\circ\ell_{\hat y}^{-1}$ and $z_{\hat y}=z\circ\ell_{\hat y}^{-1}$ of
Section~\ref{subsec:inputs}. The weight uses the ambient coefficient $a_K(\sigma)$ and the
cut-off $\psi_p$. For $t=\sigma$ the integral is taken over the initial slice, against $\mu_0$.
For $q=p$ the auxiliary density is the
density itself, $\tilde\Theta_{p,t}(\sigma)=\Theta^\chi_{p,t}(\sigma)$. For fixed $\sigma$ the map
$(q,t)\mapsto\tilde\Theta_{q,t}(\sigma)$ is continuous on
$\{(q,t):\ q\in\Upsilon_p^{-1}(\delta_r\cB_1),\ \sigma<t<T+\sigma\}$: the integrand is a
bounded continuous function of $(\hat y,t,x)\in\delta_r\cB_1\times(\sigma,T+\sigma)\times\Sigma$,
the flow being smooth at the times $t-\sigma<T$, so the integral over the
compact $\Sigma$ is continuous.

\begin{lemma}[comparison at a fixed scale]\label{lem:compare-centres}
There is $C_6=C_6(n)$ with the following property. Let $p\in M$, $\Lambda\ge1$ and
$0<r\le r_M$ with $2C^{\rm tr}r^2\Lambda^2\le1$, put $\sigma:=r^2/\Lambda^2$ and assume
$\sigma\le\tau_0$. Then for every $q\in\Upsilon_p^{-1}(\delta_r\cB_1)$ and every $t$ with
$\sigma\le t\le T$,
\[
\Theta^\chi_{q,t}(\sigma)\le e^{\epsilon_1}\,\tilde\Theta_{q,t}(\sigma)+C_6D_0\,\epsilon_1
+C(n)D_0\Big(e^{-\Lambda^2/8}+\Lambda^ne^{-a(\bb^2\Lambda^2-8C_\theta)_+^2}\Big),\qquad
\epsilon_1:=2C^{\rm tr}r^2\Lambda^2 .
\]
\end{lemma}
\begin{proof}
Put $G:=\Upsilon_p^{-1}(\delta_r\cB_4)$, $R:=r/2$ and $Z:=\bb^2r^2$, and write
$\hat y_r:=\delta_{1/r}\hat y\in\cB_1$.

\emph{The two kernels on $G$.} By Lemma~\ref{lem:rM}(b) for the pair $(p,q)$ at the scale
$r\le r_M\le r_*$, the transition $S:=\delta_{1/r}\Upsilon_q\Upsilon_p^{-1}\delta_r$ satisfies
$S=\tilde U\ell_{\hat y_r}^{-1}+E$ on $\cB_8$ with $\|E\|_{C^{k_0}(\cB_8)}\le C^{\rm tr}r^2$.
In the rescaled coordinates $\hat w=\delta_{1/r}\Upsilon_p(\cdot)$ on $\cB_8$ we have
$r^{-2}f_q=f\circ S$ and $r^{-2}\tilde f_q=f_{\hat y_r}=f\circ\tilde U\ell_{\hat y_r}^{-1}$
(homogeneity of $f$ and $U(n)$-invariance). Since $f=|w|^2/16$ is Lipschitz with
constant $\frac{10}8$ on $\{|w|\le10\}$, it follows that
\[
|f_q-\tilde f_q|\le2C^{\rm tr}r^4,\qquad |z_q-\tilde z_q|\le C^{\rm tr}r^4\qquad\text{on }\Upsilon_p^{-1}(\delta_r\cB_8)\supset G .
\]
Both cut-offs equal $1$ on $G$: $\psi_p=1$ on $\Upsilon_p^{-1}\{|w|<r_0/2,|z|<z_0/2\}\supset G$
by Definition~\ref{def:rM}. Moreover, $S(\cB_4)\subset\cB_6$ (horizontal part at most
$4+1+\delta_*$, vertical part at most $16\bb^2+\bb^2+\frac12\cdot4+\delta_*$), so $G$ lies
in $\Upsilon_q^{-1}\{|w|<6r,|z|<36\bb^2r^2\}\subset\Upsilon_q^{-1}\{|w|<r_0/2,|z|<z_0/2\}$,
where $\psi_q=1$ ($6r_M\le r_0/8$, $36\bb^2r_M^2\le36z_0/128$). Hence on $G$, with
$v=(z_q-a_K(\sigma)\theta)/\sigma$ and $\tilde v=(\tilde z_q-a_K(\sigma)\theta)/\sigma$,
\[
|v-\tilde v|\le C^{\rm tr}r^4/\sigma=\tfrac12\epsilon_1,\qquad
e^{-f_q/\sigma}\le e^{\epsilon_1}e^{-\tilde f_q/\sigma},\qquad
\psi_q\rho^\chi_{q,t}\le e^{\epsilon_1}\psi_p\,\tilde\rho_q\big(\chi(\tilde v)+\tfrac12\epsilon_1\big),
\]
using $\|\chi'\|_\infty=\sqrt{2a/e}\le1$.

\emph{Tails of the true density.} By the tail estimate \eqref{eq:tail} (centre $(q,t)$,
$\sigma\le\min\{\tau_0,t\}$, $t-\sigma<T$), the part of $\Theta^\chi_{q,t}(\sigma)$ coming
from the complement of $A:=\{f_q<R^2\}\cap\{|z_q|<Z\}$ is at most
$\tail(R,Z;\sigma)=C(n)D_0e^{-\Lambda^2/8}+(4\pi)^{-n/2}D_0(2\Lambda)^ne^{-a(\bb^2\Lambda^2-8C_\theta)_+^2}$.
Moreover, $A\subset G$. Indeed, in the rescaled chart of $q$, $A$ is $\{|w|<2,|z|<\bb^2\}$. By
Lemma~\ref{lem:rM}(b) for the pair $(q,p)$ (the point $p$ lies in $U_q$ with rescaled
coordinates $h':=S(0)$, $|h'_w|\le1+\delta_*$, $|h'_z|\le\bb^2+\delta_*$, so $h'\in\cB_8$),
the transition from the rescaled chart of $q$ to that of $p$ is $\delta_*$-close on $\cB_8$ to $\tilde U'\ell_{h'}^{-1}$. For $|w|<2$, $|z|<\bb^2$ the group law then gives a
horizontal part of modulus $<2+1+2\delta_*<4$ and a vertical part of modulus
$<\bb^2+\bb^2+\frac12(1+\delta_*)\cdot2+2\delta_*<16\bb^2$. Hence $A\subset\Upsilon_p^{-1}(\delta_r\cB_4)=G$.

\emph{Gaussian mass of the auxiliary kernel on $G$.} We claim
\[
M_G:=\int_{F_{t-\sigma}^{-1}(G)}\tilde\rho_q\,d\mu_{t-\sigma}\le C_{\rm mass}(n)\,D_0 ,
\]
and we put $C_6:=2C_{\rm mass}$. Every point of $G$ has rescaled $p$-coordinates in
$\cB_4\subset\cB_8$, where $S$ is controlled, and $G\subset U_q$ by the inclusion
$S(\cB_4)\subset\cB_6$ above. For $R'\ge r$ the set $G\cap\{|w_p-\hat y_w|<R'\}$ lies in
$\{|w_p|<2R'\}$, whose preimage under $F_{t-\sigma}$ has area at most $D_0(2R')^n$ by
\eqref{eq:D} if $2R'\le r_0$ and at most $D_0r_0^n\le D_0(2R')^n$ otherwise (all of
$F^{-1}_{t-\sigma}(U_p)$ having area at most $D_0r_0^n$). For $\sqrt\sigma\le R'<r$ the set
$G\cap\{|w_p-\hat y_w|<R'\}$ lies in $\{|w_q|<2R'\}$. Indeed, for its points, in the rescaled
coordinates and with $\hat w_w$ the horizontal component of the point $\hat w$, we have
$|S(\hat w)_w-U(\hat w_w-\hat y_{r,w})|\le\|E\|_{C^0(\cB_8)}\le C^{\rm tr}r^2\le\frac1{2\Lambda^2}\le\frac1\Lambda\le\frac{R'}r$,
so $|\hat w_w-\hat y_{r,w}|<R'/r$ gives $|S(\hat w)_w|<2R'/r$. Hence the preimage of this set
has area at most $D_0(2R')^n$ by \eqref{eq:D} in the chart of $q$.

Decompose $G$ into
$G\cap\{|w_p-\hat y_w|<\sqrt\sigma\}$ and the shells
$G\cap\{2^k\sqrt\sigma\le|w_p-\hat y_w|<2^{k+1}\sqrt\sigma\}$, $k\ge0$, on which
$\tilde\rho_q\le(4\pi\sigma)^{-n/2}e^{-4^k/16}$. By the two cases (with $R'=2^{k+1}\sqrt\sigma$),
their preimages have area at most $D_0(2^{k+2}\sqrt\sigma)^n$, and summing gives
$M_G\le D_0(4\pi)^{-n/2}\big(2^n+\sum_{k\ge0}2^{(k+2)n}e^{-4^k/16}\big)=:C_{\rm mass}(n)D_0$.

\emph{Conclusion.} Since the integrand is non-negative and $A\subset G$,
\begin{align*}
\Theta^\chi_{q,t}(\sigma)&\le\int_{F^{-1}(A)}\psi_q\rho^\chi_{q,t}\,d\mu_{t-\sigma}+\tail(R,Z;\sigma)\\
&\le e^{\epsilon_1}\int_{F^{-1}(G)}\psi_p\tilde\rho_q\chi(\tilde v)\,d\mu_{t-\sigma}+e^{\epsilon_1}\tfrac{\epsilon_1}2M_G+\tail(R,Z;\sigma),
\end{align*}
where the first integral is at most $\tilde\Theta_{q,t}(\sigma)$ and, as $\epsilon_1\le1$,
$e^{\epsilon_1}\frac{\epsilon_1}2M_G\le2\epsilon_1C_{\rm mass}D_0=C_6D_0\epsilon_1$. Every use
of the transition $S$ above concerns points of $G$, whose rescaled coordinates lie in
$\cB_4\subset\cB_8$.
\end{proof}

\subsection{Terminal regularity}

\begin{theorem}[terminal regularity: empty or one sheet]\label{thm:B}
Assume (S) and (F), and let $\eps=\eps(n,a,\frac12,C_\theta)\in(0,\frac12]$ be the constant
of Theorem~\ref{thm:A} for $\alpha=\frac12$ (by Theorem~\ref{thm:alpha}, $\bar\eps/2$ with
any exponent $\alpha$ may be used instead). Let $p\in M$ satisfy $\Theta^\chi(p,T)<1+\eps$.
Then there is $\varrho>0$ such that, with $V:=\Upsilon_p^{-1}\{|w|<\varrho,\ |z|<\bb^2\varrho^2\}$
and $\delta:=\varrho^2$, one of the following holds.
\begin{enumerate}[(a)]
\item $(p,T)\notin\overline\cM$, $F_t(\Sigma)\cap V=\emptyset$ for all $t\in[T-\delta,T)$, and
$\Theta^\chi(p,T)=0$.
\item $(p,T)\in\overline\cM$, and $\cM\cap(V\times[T-\delta,T))$ is a single embedded
Legendrian sheet, smooth up to the time $T$. Precisely, there are $Y_0\in\cM$, a radius
$r_{\mathrm{adm}}>0$ (the \emph{admissible radius} of the terminal sheet, chosen in Step~3 of
the proof; it is unrelated to the transition scale $r_*$ of Lemma~\ref{lem:rM}(b)),
a chart $\tilde U\Lambda_{Y_0,r_{\mathrm{adm}}}$ in which $\overline{\cC_\varrho(p,T)}$ lies in
$\{|x|<\frac14,|y|<2,|z|<\bb^2\}\times(-\frac1{16},\frac1{16})$, and functions $u,\theta,\zeta$
on $Q:=\{|x|<\frac12,\ |s|<\frac14,\ s\le T'\}$, $T'=(T-t(Y_0))/r_{\mathrm{adm}}^2$, smooth up to the
upper end $s=T'$ (one-sided), such that $\cM\cap\overline{\cC_\varrho(p,T)}$ is the part
lying in $\overline{\cC_\varrho(p,T)}$ of the Legendrian graph of $u$ over $Q\cap\{s<T'\}$,
with height $\zeta$ and Legendrian angle $\theta$; the slices $\cM(t)\cap V$ converge
smoothly on compact subsets of $V$, as $t\to T$, to the embedded Legendrian sheet given by
$(u,\zeta)(\cdot,T')$, which passes through $p$; $|B|$ is bounded on
$\{(x,t):F_t(x)\in V,\ T-\delta\le t<T\}$; $\Theta^\chi(p,T)\ge1$; and the map
$\hat F=(F,t)$ is injective on $\hat F^{-1}(V\times[T-\delta,T))$, that is, $F_t$ is
injective on $F_t^{-1}(V)$ for every $t\in[T-\delta,T)$, so that the sheet is covered
exactly once by the flow.
\end{enumerate}
In particular $p$ is not a curvature singularity of $F$ at $T$; equivalently,
$\Theta^\chi(p,T)\ge1+\eps$ at every curvature singularity.
\end{theorem}
\begin{proof}
Let $C$ be the constant of Theorem~\ref{thm:A} for $\alpha=\frac12$; we may assume
$C\ge3$. Put $\eta:=\frac12(1+\eps-\Theta^\chi(p,T))>0$.

\emph{Step 1: the scale.} Choose $\Lambda\ge1$ so large that
$C(n)D_0(e^{-\Lambda^2/8}+\Lambda^ne^{-a(\bb^2\Lambda^2-8C_\theta)_+^2})\le\eta/8$, with the
constant of Lemma~\ref{lem:compare-centres}. Then choose $0<r\le r_M$ so small that
$\epsilon_1=2C^{\rm tr}r^2\Lambda^2$ satisfies $\epsilon_1\le1$, $e^{\epsilon_1}\le1+\eta/8$
and $C_6D_0\epsilon_1\le\eta/8$, and that $\sigma_*:=r^2/\Lambda^2$ satisfies
$\sigma_*\le\min\{\tau_0,T/2\}$, $C_5\sigma_*\le\eta/8$ and
$\Theta^\chi_{p,T}(\sigma_*)<1+\eps-\tfrac74\eta$. The last condition can be met because
$\Theta^\chi_{p,T}(\sigma)\to\Theta^\chi(p,T)=1+\eps-2\eta$ as $\sigma\to0$
(see \eqref{eq:almost-mono}).

\emph{Step 2: the density hypothesis near $(p,T)$.} Since
$\tilde\Theta_{p,T}(\sigma_*)=\Theta^\chi_{p,T}(\sigma_*)<1+\eps-\frac74\eta$ and
$(q,t)\mapsto\tilde\Theta_{q,t}(\sigma_*)$ is continuous at $(p,T)$, there is
$\rho_1\in(0,\min\{r,\sqrt{\sigma_*}\}]$ such that $\tilde\Theta_{q,t}(\sigma_*)<1+\eps-\frac32\eta$ for all
$(q,t)\in N:=\Upsilon_p^{-1}(\delta_{\rho_1}\cB_1)\times(T-\rho_1^2,T]$. For such $(q,t)$ with
$t\le T$, Lemma~\ref{lem:compare-centres} gives
\[
\Theta^\chi_{q,t}(\sigma_*)\le(1+\tfrac\eta8)(1+\eps-\tfrac32\eta)+\tfrac\eta8+\tfrac\eta8\le1+\eps-\eta ,
\]
using $1+\eps-\frac32\eta\le2$, which holds as $\eps\le\frac12$. The bound \eqref{eq:C5}
is valid for the centres $(q,t)$, $0<t\le T$, provided $\sigma_*\le\min\{\tau_0,t\}$, which
holds as $t>T-\rho_1^2\ge T-\sigma_*\ge\sigma_*$ by the choice of $\rho_1$ and
$\sigma_*\le T/2$. It gives, for $0<\sigma'\le\sigma_*$,
\begin{equation}\label{eq:Heps-near}
\Theta^\chi_{q,t}(\sigma')\le\Theta^\chi_{q,t}(\sigma_*)+C_5\sigma_*\le1+\eps-\tfrac78\eta<1+\eps\qquad\text{for all }(q,t)\in N,\ t\le T .
\end{equation}

\emph{Step 3: Theorem~\ref{thm:A} on the truncations.} Put
$\rho_2:=\min\{\rho_1,(\sigma_*/2)^{1/2}\}$ and $\cU:=\cC_{\rho_2}(p,T)=\Upsilon_p^{-1}\{|w|<\rho_2,|z|<\bb^2\rho_2^2\}\times(T-\rho_2^2,T+\rho_2^2)$,
an open set of $M\times\R$ which crosses the time $T$. For $T_i\uparrow T$ the truncations
$F^{\le T_i}$, restricted to $\cU$, are smooth local Legendrian flows in $\cU$ with upper
end $T_i$ (Definition~\ref{def:localflow}): the map $\hat F$ is proper on
$\Sigma\times[0,T_i]$ because $\Sigma$ is closed. Every point $Y=(q,t)$ of their traces
lies in $N$ with $t\le T$, and $d_1(Y,\cU)^2\le d(Y,\cU)^2<2\rho_2^2\le\sigma_*$ (the point
$(q,T-\rho_2^2)\notin\cU$ has quasi-distance $<\sqrt2\rho_2$ from $Y$). The density of the restriction to $\cU$ is at most the density of the closed flow (the
integrand is non-negative and the parameter domain is smaller), so \eqref{eq:Heps-near}
gives the hypothesis $(\mathrm H_\eps)$ for each $F^{\le T_i}$ in $\cU$. Theorem~\ref{thm:A} then gives
\[
d_1(Y,\cU)\,\Kreg(\cM^{\le T_i},Y)\le C\qquad\text{for all }Y\in\cM^{\le T_i}\cap\cU,\ i\in\N .
\]
For $Y\in\cM\cap\cC_{\rho_2/6}(p,T)$ we have $d(Y,\cU)\ge\rho_2/6$ by Lemma~\ref{lem:quasi},
exactly as in the proof of Lemma~\ref{lem:pick} ($\cU$ is the cylinder of radius $\rho_2$
about $(p,T)$ and $\rho_2\le r_M$). Moreover $\rho_2\le r_1$: if $K\ne-2$, then
$\rho_2\le(\sigma_*/2)^{1/2}\le(\tau_0/2)^{1/2}<|K+2|^{-1/2}$, since $\tau_0\le|K+2|^{-1}$
(see \eqref{eq:tau0}). Hence $d_1(Y,\cU)\ge\rho_2/6$ and, with
\[
r_{\mathrm{adm}}:=\rho_2/(12C)<d_1(Y,\cU)/C\qquad(\text{so also }r_{\mathrm{adm}}\le\rho_2/36\le r_1/36<r_1/2),
\]
the radius $r_{\mathrm{adm}}$ (the \emph{admissible radius} of the terminal sheet; it is not
the transition scale $r_*$ of Lemma~\ref{lem:rM}(b), which enters this section only
through $r_M\le r_*$) is admissible for $(\cM^{\le T_i},Y)$ for every $Y\in\cM\cap\cC_{\rho_2/6}(p,T)$
with $t(Y)\le T_i$ (Lemma~\ref{lem:admissible}(a)), uniformly in $i$.

\emph{Step 4: the dichotomy and the uniform graph.} Put $\varrho:=r_{\mathrm{adm}}/100$, and let $V$
and $\delta=\varrho^2$ be as in the statement, so that $V\times(T-\delta,T+\delta)=\cC_\varrho(p,T)$.
If $(p,T)\notin\overline\cM$, then $\cM\cap\cC_{\varrho'}(p,T)=\emptyset$ for some
$\varrho'\in(0,\varrho]$, and (a) holds with $\varrho'/2$ in place of $\varrho$: the set
$V\times[T-\delta,T)$ formed with $\varrho'/2$, including its closed lower end, lies in the
open cylinder $\cC_{\varrho'}(p,T)$. Then also $\Theta^\chi(p,T)=0$: for
$0<\sigma<\min\{\delta,\tau_0,T\}$ the slice $F_{T-\sigma}(\Sigma)$ misses
$V=\Upsilon_p^{-1}\{|w|<\varrho'/2,\ |z|<\bb^2\varrho'^2/4\}$, and $\psi_p$ vanishes outside $U_p$, so
$\Theta^\chi_{p,T}(\sigma)$ comes only from $\{f_p\ge(\varrho'/8)^2\}\cup\{|z_p|\ge\bb^2\varrho'^2/4\}$
(recall $f_p=|w|^2/16$). By the tail estimate \eqref{eq:tail} with $R=\varrho'/8$ and $Z=\bb^2\varrho'^2/4$,
\[
\Theta^\chi_{p,T}(\sigma)\le C(n)D_0e^{-\varrho'^2/(128\sigma)}
+(4\pi\sigma)^{-n/2}D_0(\varrho'/2)^n\exp\big(-a(\bb^2\varrho'^2/(4\sigma)-8C_\theta)_+^2\big),
\]
which tends to $0$ as $\sigma\to0$ (for $\sigma\le\bb^2\varrho'^2/(4(16C_\theta+1))$ the bracket is at
least $\bb^2\varrho'^2/(8\sigma)$).

Now assume
$(p,T)\in\overline\cM$ and fix $Y_0\in\cM\cap\cC_\varrho(p,T)$. The quasi-symmetry
$\|(p,T)-Y_0\|_{Y_0}\le2\|Y_0-(p,T)\|_{(p,T)}$ holds by Lemma~\ref{lem:rM}(b) at the scale
$\|Y_0-(p,T)\|_{(p,T)}\le r_M$: the coordinates of $(p,T)$ in the chart of $Y_0$ are
$\tilde U(\hat y^{-1})+E(0)$ with $\|E\|_{C^0}\le\delta_*$ in rescaled units, which changes
the gauge by at most $2\delta_*^{1/2}\le\frac12$. Together with Lemma~\ref{lem:quasi}, the quasi-symmetry gives, for
$Z\in\overline{\cC_\varrho(p,T)}$,
\[
\|Z-Y_0\|_{Y_0}\le3(2\varrho+\varrho)=9\varrho<\tfrac14r_{\mathrm{adm}} .
\]
Hence in every chart $\tilde U\Lambda_{Y_0,r_{\mathrm{adm}}}$, $U\in U(n)$, the closed cylinder
$\overline{\cC_\varrho(p,T)}$ lies in $\{|x|<\frac14,|y|<2,|z|<\bb^2\}\times(-\frac1{16},\frac1{16})$,
with base coordinates in $\{|x|<\frac14,\ |s|<\frac1{16}\}$. The closure
$\{|x|\le\frac14,\ |s|\le\frac1{16}\}$ of this base set is a compact subset of the inner base
$\{|x|<\frac12,|s|<\frac14\}$, with margins $\frac14$ in $x$ and $\frac3{16}$ in $s$. Moreover,
$T'=(T-t(Y_0))/r_{\mathrm{adm}}^2<\varrho^2/r_{\mathrm{adm}}^2=10^{-4}$.

By Step 3, for every $i$ with $T_i\ge t(Y_0)$ there is a witness $(U_i,u_i,\zeta_i)$ of
$(\cM^{\le T_i},Y_0)$ at the radius $r_{\mathrm{adm}}$: $\|u_i\|_{2,1/2;\Omega_{T'_i}}\le1$,
$T'_i=(T_i-t(Y_0))/r_{\mathrm{adm}}^2\uparrow T'$, $u_i(0,0)=0$, $Du_i(0,0)=0$, and
$\cM^{\le T_i}\cap\bB_{Y_0,r_{\mathrm{adm}}}\subset\hat\Gamma_{u_i,\zeta_i}(\Omega_{T'_i})$ in the chart
$\tilde U_i\Lambda_{Y_0,r_{\mathrm{adm}}}$, $\bB_{Y_0,r_{\mathrm{adm}}}:=(\tilde U_i\Lambda_{Y_0,r_{\mathrm{adm}}})^{-1}(\bB)\supset\cC_{r_{\mathrm{adm}}}(Y_0)\supset\overline{\cC_\varrho(p,T)}$.
Since $\cC_{3r_{\mathrm{adm}}}(Y_0)\subset\cU$ (as $C\ge3$, $3r_{\mathrm{adm}}\le\rho_2/(4C)\le\rho_2/12<\rho_2/6\le d(Y_0,\cU)$),
Lemma~\ref{lem:filling} shows that the graph over $\Omega_{T'_i}$ is realised by the trace,
$u_i$ solves \eqref{eq:grapheq}, and the Legendrian angle $\theta_i$ of the graph solves
\eqref{eq:angleeq}.

\emph{Estimates of all orders, uniform in $i$.} The structure in the chart
$\tilde U_i\Lambda_{Y_0,r_{\mathrm{adm}}}$ is the fixed smooth structure $\varsigma_{1/r_{\mathrm{adm}}}$ of the chart
of $Y_0$ rotated by $U_i\in U(n)$. Since $U(n)$ is compact, the coefficients of
\eqref{eq:grapheq}, as functions of $(x,u,Du)$ on $\{|x|<2,\ |u|<2,\ |Du|\le1\}$, and the
coefficients $G^{ij}$, $G^{ij}\,{}^G\Gamma^k_{ij}$ and $H^k$ of \eqref{eq:angleeq}, as
functions of $(x,u,Du,D^2u)$ on $\{|x|<2,\ |u|<2,\ |Du|\le1,\ |D^2u|\le1\}$ (the last two
involve $D^2u$; see the proof of Lemma~\ref{lem:angleeq}), have derivatives of
every order bounded independently of $i$, by constants depending on the flow, on the chart
of $Y_0$ and on $r_{\mathrm{adm}}$. For every $i$ and every $(x,s)\in\Omega_{T'_i}$ the
arguments $(x,u_i,Du_i)$ and $(x,u_i,Du_i,D^2u_i)$ lie in these sets, since $|x|<1$,
$|u_i|<\frac32$ (Lemma~\ref{lem:heightangle}) and
$|Du_i|,|D^2u_i|\le\|u_i\|_{2,1/2;\Omega_{T'_i}}\le1$.

Put $Q^{(i)}_\rho:=\{|x|<\rho,\ |s|<\rho^2,\ s\le T'_i\}$, and let
\[
r_m:=\tfrac{17}{32}+2^{-m-5},\qquad \hat r_m:=\tfrac12+2^{-m-5}\qquad(m=0,1,2,\dots),
\]
so that $r_0=\frac9{16}$, $r_m\downarrow\frac{17}{32}$, $\hat r_0=\frac{17}{32}$ and
$\hat r_m\downarrow\frac12$. By the last assertion of Lemma~\ref{lem:filling}, applied as
above with $\cC_{3r_{\mathrm{adm}}}(Y_0)\subset\cU$, the functions $u_i$, $\zeta_i$ and $\theta_i$ are smooth
on $\Omega_{T'_i}$ up to the upper end, so the norms below are finite on every
$Q^{(i)}_\rho$, $\rho<1$; only their bounds matter.

(1) \emph{The graph function.} We claim by induction on $m\ge0$ that
$\|u_i\|_{C^{m+4,1/2}(Q^{(i)}_{r_m})}\le C_m$ with $C_m$ independent of $i$, where the
H\"older norms are parabolic (a time derivative counts twice). For $m=0$ this is
Lemma~\ref{lem:subord}(i) on $\Omega^{(9/16)}_{T'_i}=Q^{(i)}_{r_0}$. Assume it for $m$. For
every multi-index $\beta$ of length $m+3$, $w:=D^\beta u_i$ solves on $Q^{(i)}_{r_m}$ a
system $\partial_sw^\gamma=\bar a^{ij}D_{ij}w^\gamma+g^\gamma$, obtained by differentiating
\eqref{eq:grapheq} $m+3$ times, where $g$ is a polynomial in the derivatives of the
coefficients $a,b$ of order at most $m+3$ (evaluated at $(x,u_i,Du_i)$) and in the
derivatives of $u_i$ of order at most $m+4$. By the induction hypothesis and the
coefficient bounds, $\|g\|_{0,1/2;Q^{(i)}_{r_m}}\le C$. Theorem~\ref{thm:schauder-var} with
the outer domain $Q^{(i)}_{r_m}$ and the inner domain $Q^{(i)}_{r_{m+1}}$ bounds
$|w|_{2,1/2;Q^{(i)}_{r_{m+1}}}$, i.e.\ the spatial derivatives of $u_i$ of order $m+5$ and
$\partial_sD^\beta u_i$, in $C^{0,1/2}$. The remaining derivatives (time derivatives of
higher order and mixed derivatives) are recovered from the equation and its
derivatives. This gives the claim for $m+1$. Consequently every $C^k$ norm of $u_i$ on
$Q^{(i)}_{17/32}$ is bounded independently of $i$.

(2) \emph{The angle.} The coefficients of the linear equation \eqref{eq:angleeq} are
smooth functions of $(x,u_i,Du_i,D^2u_i)$, hence bounded in every $C^k(Q^{(i)}_{17/32})$ by (1).
We claim $\|\theta_i\|_{C^{m+2,1/2}(Q^{(i)}_{\hat r_m})}\le C'_m$. For $m=0$ this is
Lemma~\ref{lem:subord}(iv), whose proof applies Theorem~\ref{thm:schauder-var} with the
inner domain $\Omega^{(9/16)}_{T'_i}\supset Q^{(i)}_{\hat r_0}$, using only the bound
$|\theta_i|\le C_\theta$.
Assume the claim for $m$. Differentiating \eqref{eq:angleeq} $m+1$ times gives for each
$D^\beta\theta_i$, $|\beta|=m+1$, a linear equation of the same type on $Q^{(i)}_{\hat r_m}$
with a right-hand side bounded in $C^{0,1/2}$ by the claim for $m$ and the coefficient
bounds, and Theorem~\ref{thm:schauder-var} with the outer domain $Q^{(i)}_{\hat r_m}$ and the
inner domain $Q^{(i)}_{\hat r_{m+1}}$ gives the claim for $m+1$.

(3) \emph{The height.} By \eqref{eq:height} and \eqref{eq:theta-from-graph},
$D_x\zeta_i\cdot e=-2\beta_\varsigma((e,Du_i\,e))$ for unit vectors $e$, and $\partial_s\zeta_i=-2\beta_\varsigma((0,\partial_su_i))-4\theta_i$,
with $\zeta_i(0,0)=0$; hence every $C^k$ norm of $\zeta_i$ on $Q^{(i)}_{1/2}$ is bounded
independently of $i$, by (1) and (2).

In each step the outer domain carries the bounds of the required order and the inner
domain has a positive margin; all constants are independent of $i$ and of $T'_i$, and
every $Q^{(i)}_{\hat r_m}$ contains $Q^{(i)}_{1/2}$.

\emph{The limit and its extension to $s=T'$.} The functions $u_i,\theta_i,\zeta_i$ are
defined for $s\le T'_i$ only, and $T'_i\uparrow T'$. For a margin $0<\nu<\frac1{16}$ let
$K_\nu:=\{|x|\le\frac12-\nu,\ -\frac14+\nu\le s\le T'-\nu\}$, a compact subset of
$Q^-:=\{|x|<\frac12,\ -\frac14<s<T'\}$ on which all $u_i,\theta_i,\zeta_i$ with
$T'_i\ge T'-\nu$ are defined. By the uniform bounds and Arzel\`a--Ascoli, a diagonal
subsequence (along which also $U_i\to U$) converges in $C^\infty$ on compact subsets of
$Q^-$ to functions $u,\theta,\zeta$. Every derivative of $u,\theta,\zeta$ is bounded on $Q^-$
by the corresponding constant, so each derivative is uniformly continuous on
$\{|x|\le\frac12-\nu,\ -\frac14+\nu\le s<T'\}$ and extends continuously to $s=T'$, with
compatible limits. Thus $u,\theta,\zeta$ are smooth up to the upper end on
$Q:=\{|x|<\frac12,\ |s|<\frac14,\ s\le T'\}$ (one-sided at $s=T'$). We pass to the limit in
\eqref{eq:grapheq}, \eqref{eq:height}, \eqref{eq:theta-from-graph} and \eqref{eq:angleeq}
on compact subsets of $Q^-$. (In \eqref{eq:angleeq} the zero-order constant is
$\kappa_\theta=(K+2)r_{\mathrm{adm}}^2$, the structure of the chart being
$\varsigma_{1/r_{\mathrm{adm}}}$, i.e.\ $r=r_{\mathrm{adm}}$ and $\lambda=1$ in the notation of
Section~\ref{sec:graphs}; it is the same fixed constant for all $i$ and does not involve
the margin $\nu$.) The limiting equations show that the limit graph is Lagrangian for the rotated structure, with
height $\zeta$ and Legendrian angle $\theta$, and it moves by the projected flow up to
$s=T'$ by continuity.

\emph{The trace in the target cylinder.} Let $Z\in\cM\cap\overline{\cC_\varrho(p,T)}$. Then
$t(Z)<T$, so $t(Z)\le T_i$ for $i$ large. By the inclusion above and
Definition~\ref{def:K}(i), $Z$ lies on the graph of $u_i$ for all such $i$, with base
coordinates $(x_i,s_i)$ in $\{|x|<\frac14,|s|<\frac1{16}\}$ and $s_i\le T'_i$. The charts
$\tilde U_i\Lambda_{Y_0,r_{\mathrm{adm}}}$ converge to $\tilde U\Lambda_{Y_0,r_{\mathrm{adm}}}$, so $(x_i,s_i)\to(x,s)$,
the base coordinates of $Z$ in the limit chart, with $s<T'$ (as $t(Z)<T$). Hence
$(x_i,s_i)\in K_\nu$ for some fixed $\nu>0$ and all large $i$, and the uniform
convergence $u_i\to u$ on $K_\nu$ shows that $Z$ lies on the graph of $u$. Conversely,
let $(x,s)$ with $|x|<\frac14$, $|s|<\frac1{16}$, $s<T'$, be such that the point
$Z:=(x,u(x,s),\zeta(x,s))$ of the limit graph at the time $s$ lies in $\overline{\cC_\varrho(p,T)}$,
and put $t:=t(Y_0)+r_{\mathrm{adm}}^2s<T$. For $i$ large, $(x,s)\in K_\nu$ with $s\le T'_i$, and the
point $Z_i:=(x,u_i(x,s),\zeta_i(x,s))$ of the graph of $u_i$ at the time $s$ lies on the
trace of $F$ at the time $t$ (Lemma~\ref{lem:filling}): $Z_i\in F_t(\Sigma)$. By the
convergence of the charts and of $u_i,\zeta_i$ on $K_\nu$, $Z_i\to Z$. Since the slice
$F_t(\Sigma)$ is compact, $Z\in F_t(\Sigma)$, i.e.\ $Z\in\cM$. Thus
$\cM\cap\overline{\cC_\varrho(p,T)}$ is exactly the part of the graph of $u$ over
$Q\cap\{s<T'\}$ lying in $\overline{\cC_\varrho(p,T)}$.

By the smoothness of $u,\zeta$ up to
$s=T'$, the slices $\cM(t)\cap V$ converge smoothly on compact subsets of $V$, as $t\to T$,
to the embedded Legendrian sheet $\{(x,u(x,T'),\zeta(x,T'))\}$ (in the chart). Since
$(p,T)\in\overline\cM$, the point $(p,T)$, whose chart coordinates have $|x|<\frac14$ and
time $T'$, is a limit of points of the graphs and lies on this sheet. The curvature bound
on $V\times[T-\delta,T)\subset\overline{\cC_\varrho(p,T)}$ follows from Lemma~\ref{lem:subord}(iii)
applied to each $F^{\le T_i}$ at the admissible radius $r_{\mathrm{adm}}$, uniformly in $i$. This is (b)
except for the last two claims, the density bound and the injectivity, which are proved in
Step~5.

\emph{Step 5: the local flow at $(p,T)$, $\Theta^\chi(p,T)\ge1$, and the injectivity in (b).} In the chart
$\tilde U\Lambda_{Y_0,r_{\mathrm{adm}}}$ let $\Gamma:=\hat\Gamma_{u,\zeta}(Q_\Gamma)$ with the base domain
$Q_\Gamma:=\{|x|<\frac38,\ -\frac3{16}<s\le T'\}$ (recall $T'<10^{-4}$). It is an embedded
$(n+1)$-dimensional submanifold with boundary $\Gamma\cap\{s=T'\}$, smooth up to the
boundary, Legendrian in each slice, on which the time function has no critical points.
It is relatively closed in the open set
$\cU'':=(\tilde U\Lambda_{Y_0,r_{\mathrm{adm}}})^{-1}\big(\{|x|<\tfrac38,|y|<2,|z|<\bb^2\}\times(-\tfrac3{16},\tfrac3{16})\big)$;
the graph stays inside this set ($|u|<\frac32$, $|\zeta|<\bb^2/2$ by Lemma~\ref{lem:heightangle},
passed to the limit), and $(p,T)\in\cU''$. The normal velocity $H-2\theta\xi$ of the graph is
smooth up to $s=T'$. Integrating it as in Step 3 of the proof of
Proposition~\ref{prop:compact}, we find that $(p,T)$ has an open neighbourhood $\cU'''\subset\cU''$ in which $\Gamma\cap\cU'''$ is the
trace of a smooth local Legendrian flow $G$ (Definition~\ref{def:localflow}) with upper end $T$ and injective $\hat G$. The embedding $\hat G$ is proper as a map into $\cU'''$, since its image
is relatively closed there. Replacing $\cU'''$ by
$\cU'''\cap\cC_\varrho(p,T)\subset\cU''\cap\cC_\varrho(p,T)$, still an open neighbourhood of
$(p,T)$, keeps $G$ a smooth
local Legendrian flow with the same properties (Remark~\ref{rem:localflow}(3)). For $t<T$, its slices now lie in the portion of the graph identified with the
trace of $F$ in Step~4. Lemma~\ref{lem:density-M} applied
to $G$ gives $\liminf_{\sigma\to0}\Theta^\chi_{(p,T)}(\sigma;G)\ge1$. For $t<T$ the slices of
$G$ are pieces of the slices of $F$ in $\overline{\cC_\varrho(p,T)}$, each covered exactly
once by $F_t$: the maps $\hat F^{\le T_i}$ are injective on $(\hat F^{\le T_i})^{-1}(\cU)$ by
Lemma~\ref{lem:density-M} applied to $F^{\le T_i}$ in $\cU$, which satisfies $(\mathrm H_\eps)$
with $\eps\le\frac12<1$. Hence $\Theta^\chi_{p,T}(\sigma;G)\le\Theta^\chi_{p,T}(\sigma)$ for
small $\sigma$, and $\Theta^\chi(p,T)=\lim_{\sigma\to0}\Theta^\chi_{p,T}(\sigma)\ge1$.

\emph{The injectivity in (b).} First, $V\times[T-\delta,T)\subset\cU=\cC_{\rho_2}(p,T)$: indeed
$\varrho=r_{\mathrm{adm}}/100=\rho_2/(1200C)<\rho_2$, so that
$\{|w|<\varrho,\ |z|<\bb^2\varrho^2\}\subset\{|w|<\rho_2,\ |z|<\bb^2\rho_2^2\}$, and
$\delta=\varrho^2<\rho_2^2$ gives $[T-\delta,T)\subset(T-\rho_2^2,T+\rho_2^2)$. Now let
$(x,t),(x',t')\in\Sigma\times[0,T)$ with $F_t(x)\in V$, $F_{t'}(x')\in V$,
$t,t'\in[T-\delta,T)$ and $\hat F(x,t)=\hat F(x',t')$. Then $t=t'<T$, and we choose $i$
with $T_i\ge t$. Both $(x,t)$ and $(x',t)$ lie in $\Sigma\times[0,T_i]$ and in
$\hat F^{-1}(\cU)$, so the injectivity of $\hat F^{\le T_i}$ on $(\hat F^{\le T_i})^{-1}(\cU)$ established
above gives $x=x'$. Hence $\hat F$ is injective on $\hat F^{-1}(V\times[T-\delta,T))$, which
is the last claim of (b). Equivalently, each $F_t$, $T-\delta\le t<T$, is injective on
$F_t^{-1}(V)$.

\emph{Step 6: curvature singularities.} In case (a) no point of the flow lies in $V$ for
$t\in[T-\delta,T)$, and in case (b) $|B|$ is bounded there; so $p$ is not a curvature
singularity. Contrapositively, $\Theta^\chi(p,T)\ge1+\eps$ at every curvature singularity.
\end{proof}

\begin{remark}
The proof uses the density hypothesis only in the form \eqref{eq:Heps-near} on the
truncated flows in the open set $\cU$, which crosses the terminal time. The argument is
the truncation device of \cite[Remark~3.3, Theorem~3.5]{wh}. Here
Lemma~\ref{lem:compare-centres} gives the comparison at a fixed scale,
and \eqref{eq:C5} propagates the bound to smaller scales. No information on the flow after the time $T$ is used or obtained.
\end{remark}

\subsection{Type~I singularities}

\begin{corollary}[total multiplicity at a Type~I curvature singularity]\label{cor:Bprime}
Assume (S), (F), the Type~I hypothesis (I) of \cite[\S2.2]{PaperI} ($|B|^2\le C_I/(T-t)$)
and that $M$ is $\eta$-Einstein with the constant $K$ of (F). Let $p$ be a curvature
singularity of $F$ at $T$. Then $\Theta^\chi(p,T)$ is an integer at least $2$, and for
every sequence $\lambda_i\to\infty$ a subsequence of the contact parabolic rescalings at
$(p,T)$ converges, as in \cite[Theorem~5.5]{PaperI}, to a flow whose time-$s$ slice is a
finite union of horizontal Legendrian planes $P_k\times\{4|s|\bar\theta'_k\}$, $k=1,\dots,N'$,
with multiplicities $m'_k\in\N$ (the index $k$ counts the planes; distinct planes may have
the same Legendrian angle), of total multiplicity $\sum_{k=1}^{N'}m'_k=\Theta^\chi(p,T)\ge2$.
\end{corollary}
\begin{proof}
By Theorem~\ref{thm:B}, $\Theta^\chi(p,T)\ge1+\eps$. Under (I) the flow extends
continuously to the time $T$ and $p\in F_T(\Sigma)$ (\cite[Lemma~5.1]{PaperI}, the point
$p$ being a curvature singularity). Hence, by \cite[Theorem~5.5]{PaperI}, a subsequence of the
rescalings at $(p,T)$ converges smoothly on compact subsets of $\Hn\times(-\infty,0)$ to
a non-empty flow $F^\infty$ whose slices are finite unions of horizontal Legendrian planes
$P_k\times\{4|s|\bar\theta'_k\}$, $k=1,\dots,N'$, counted with multiplicities $m'_k\in\N$.
Here the index $k$ runs over the planes (the components of $F^\infty_s$, each $P_k$ a Lagrangian
subspace of $\C^n$), and the angles $\bar\theta'_k$ need not be distinct. Under the same
reparametrisations the Legendrian angles converge locally uniformly, the angle on the
plane $P_k\times\{4|s|\bar\theta'_k\}$ tending to $\bar\theta'_k$ (\cite[Theorem~5.5(a)]{PaperI},
whose index of the components is our $k$, and the hypothesis (C) of \cite[\S5]{PaperI},
which holds here by \cite[Proposition~5.4]{PaperI}). We pass to a further subsequence along which the conclusions of
\cite[Proposition~4.2(c)]{PaperI} hold with $S_0=1,2,3,\dots$ (successive further
subsequences and a diagonal sequence). The smooth limit and the convergence of the angles
are unchanged. Along this subsequence, for every $s<0$, both the area measures
$\mu^i_s=(F^i_s)_\#\mu$ and the joint measures $\fm^i_s=(F^i_s,\theta^i_s)_\#\mu$ converge
weakly, the latter to $(F^\infty_s,\theta^\infty_s)_\#\mu^\infty$, which carries each plane
$P_k\times\{4|s|\bar\theta'_k\}$, with its multiplicity $m'_k$, at the angle $\bar\theta'_k$. This is
the argument of the proof of \cite[Corollary~5.3]{PaperI} from (C), on every compact
interval of negative times as in the proof of \cite[Theorem~5.5]{PaperI}. Hence the limit
is a weak contact tangent flow for all negative times at $(p,T)$ in the sense of
\cite[Definition~4.1]{PaperI}, along the subsequence of \cite[Proposition~4.2(c)]{PaperI}.
Its measure at the time $s$, the weak limit of the $\mu^i_s$, is the area measure of
$F^\infty_s$ counted with multiplicities,
$\mu_s=\sum_{k=1}^{N'}m'_k\,\cH^n\llcorner\big(P_k\times\{4|s|\bar\theta'_k\}\big)$.

Let $\bar\theta_1<\dots<\bar\theta_N$ be the distinct values among
$\bar\theta'_1,\dots,\bar\theta'_{N'}$, so that $N\le N'$. Comparing with the decomposition
$\mu_s=\sum_j\|\cC_j\|\times\delta_{4|s|\bar\theta_j}$ of \cite[Theorem~4.5(v)]{PaperI}
according to the distinct heights $4|s|\bar\theta_j$ shows that the angles of that theorem
are these $\bar\theta_j$, and that the cone $\cC_j$ carrying the angle $\bar\theta_j$ is the
union of the planes $P_k$ with $\bar\theta'_k=\bar\theta_j$, counted with their
multiplicities: $\|\cC_j\|=\sum_{k:\,\bar\theta'_k=\bar\theta_j}m'_k\,\cH^n\llcorner P_k$.
Each $P_k$ is a linear subspace, of density $1$ at the origin, so the density $m_j$ of
$\cC_j$ at its vertex is $m_j=\sum_{k:\,\bar\theta'_k=\bar\theta_j}m'_k$. By
\cite[Corollary~4.7]{PaperI} with $\varpi\equiv1$,
\[
\Theta^\chi(p,T)=\lim_{t\to T}\int\psi_p\rho^\chi_{p,T}\,d\mu_t=\sum_{j=1}^Nm_j=\sum_{k=1}^{N'}m'_k ,
\]
the total multiplicity, an integer; being at least $1+\eps>1$, it is at least $2$.
\end{proof}

\begin{remark}\label{rem:Bprime}
(1) The conclusion allows a single plane of multiplicity at least two as well as several
planes, with equal or different Legendrian angles; it does not assert that two smooth
sheets or a double point have formed. Total multiplicity at least two is a necessary
condition for a curvature singularity, not a sufficient one.
(2) In the regular case the density of \cite[\S4.2]{ChangWuZhang2026ContactBlowup} is the
transverse one; Corollary~\ref{cor:Bprime} is stated with the height-localised density,
whose terminal value counts the lifted planes at their own heights (compare
\cite[Corollary~1.3, Remark~4.13]{ChangWuZhang2026ContactBlowup}). We do not claim more
than the integrality and the bound.
\end{remark}

\appendix
\section{Interior Schauder estimates on one-sided cylinders}\label{app:schauder}

This appendix provides the interior estimate used in Lemma~\ref{lem:subord}. The estimate
concerns linear parabolic systems with a common scalar principal part, H\"older
coefficients and lower-order terms, on cylinders truncated at an upper end, and its
constants are independent of the truncation time. For the heat operator
$\partial_s-\Delta$ the estimate is \cite[Theorem~8.2]{wh}. Lemma~\ref{lem:const} deduces
from it, by parabolic scaling and a linear change of the space variables, the estimate
\eqref{eq:const} on backward cylinders for the operators $\partial_s-A^{ij}D_{ij}$ with a
constant symmetric positive definite matrix $A$. The variable-coefficient case is obtained
from that lemma by freezing the coefficients on small cylinders. We follow the elliptic
argument of \cite[Theorem~6.2]{gt}, with the interior weighted seminorms of
\cite[\S6.1, Lemma~6.32]{gt} adapted to the parabolic setting.

Every local cylinder used below is \emph{backward} in time, so that it never crosses the
upper end. This is why the constants are independent of the truncation.

\emph{Notation.} Points are $X=(x,s)$, $Y=(y,t)$ in $\R^n\times\R$, the parabolic
distance is $\|X-Y\|:=\max\{|x-y|,|s-t|^{1/2}\}$, and for $T\in[0,\infty]$ and $r>0$
\[
Q_r:=\{(x,s):\ |x|<r,\ |s|<r^2,\ s\le T\},\qquad
P_r(X_0):=\{(x,s):\ |x-x_0|<r,\ s_0-r^2<s\le s_0\}
\]
are the truncated cylinder centred at the origin and the \emph{backward cylinder} with
top centre $X_0=(x_0,s_0)$. Note that $X_0\in P_r(X_0)$ and that $P_r(X_0)$ is convex.
The H\"older seminorm and norm are
$[f]_{\alpha;Q}:=\sup\{|f(X)-f(Y)|/\|X-Y\|^\alpha:\ X\ne Y\in Q\}$ and
$\|f\|_{0,\alpha;Q}:=\|f\|_{0;Q}+[f]_{\alpha;Q}$. For $u\in C^{2,\alpha}(Q)$ (parabolic
sense: $u,Du,D^2u,\partial_su$ continuous with finite seminorms) we write
\[
|u|_{2,\alpha;Q}:=\|u\|_{0;Q}+\|Du\|_{0;Q}+\|D^2u\|_{0;Q}+\|\partial_su\|_{0;Q}+[Du]_{\alpha;Q}+[D^2u]_{\alpha;Q}+[\partial_su]_{\alpha;Q};
\]
for scalar $u$ this is bounded by the parabolic $C^{2,\alpha}$ norm of \cite[\S8]{wh}, and it
contains the norm of Section~\ref{subsec:scale}. For a vector-valued $u$ all norms are taken
componentwise and summed, so that $|u|_{2,\alpha;Q}$ is bounded by $C(n,N)$ times the
parabolic $C^{2,\alpha}$ norm of \cite[\S8]{wh}, where $N$ is the number of components of $u$
($C(n,N)=n^2N$ suffices). Likewise, for the matrix coefficient $a=(a^{ij})$ and the
vector coefficient $c=(c^k)$ of Theorem~\ref{thm:schauder-var}, the norms are taken
componentwise and summed: $\|a\|_{0,\alpha;Q}:=\sum_{i,j}\|a^{ij}\|_{0,\alpha;Q}$ and
$\|c\|_{0,\alpha;Q}:=\sum_k\|c^k\|_{0,\alpha;Q}$. Correspondingly, we put
$\|a\|_{0;Q}:=\sum_{i,j}\|a^{ij}\|_{0;Q}$ and $|a(X)-a(Z)|:=\sum_{i,j}|a^{ij}(X)-a^{ij}(Z)|$
for the pointwise quantities in the proof of that theorem.

\subsection{Interpolation on backward cylinders}

\begin{lemma}[interpolation]\label{lem:interp}
There is $C_1=C_1(n,\alpha)$ such that for every backward cylinder $P_r(X_0)$, every
$u\in C^{2,\alpha}(P_r(X_0))$ and every $0<\epsilon\le1$, with the norms on the left
taken over $P_{r/2}(X_0)$ and those on the right over $P_r(X_0)$,
\begin{align}
r\|Du\|_0+r^2\|D^2u\|_0+r^2\|\partial_su\|_0
&\le\epsilon\,r^{2+\alpha}\big([D^2u]_\alpha+[\partial_su]_\alpha\big)+C_1\epsilon^{-2/\alpha}\|u\|_0,\label{eq:interp1}\\
r^{1+\alpha}[Du]_\alpha&\le C_1\big(r^2\|D^2u\|_0+r^2\|\partial_su\|_0\big).\label{eq:interp2}
\end{align}
\end{lemma}
\begin{proof}
The map $(x,s)\mapsto(x_0+rx,s_0+r^2s)$ sends $P_1(0)$ onto $P_r(X_0)$ and $P_{1/2}(0)$
onto $P_{r/2}(X_0)$, and multiplies $\|D^ju\|_0$ by $r^j$, $\|\partial_su\|_0$ by $r^2$,
and the seminorms of $D^ju$, $\partial_su$ by $r^{j+\alpha}$, $r^{2+\alpha}$. Hence it
suffices to take $X_0=0$ and $r=1$; all norms below are over $P_1(0)$ unless indicated
otherwise. Fix
$X=(x,s)\in P_{1/2}(0)$, a unit vector $e$ and $0<\mu\le\frac12$. Since $|x|<\frac12$ and
$s\in(-\frac14,0]$, the segment $\{(x+\tau e,s):|\tau|\le\mu\}$ and the time segment
$\{(x,s-\tau):0\le\tau\le\mu^2\}$ lie in $P_1(0)$. We write $u_{ee}:=D^2_{ee}u$ and
suppress the time $s$ in the first two lines. Taylor's formula with integral remainder
gives
\begin{gather*}
\begin{aligned}
u(x+\mu e)-2u(x)+u(x-\mu e)={}&\mu^2u_{ee}(x)\\
&+\int_0^\mu(\mu-\tau)\big[u_{ee}(x+\tau e)+u_{ee}(x-\tau e)-2u_{ee}(x)\big]d\tau,
\end{aligned}\\
u(x+\mu e)-u(x-\mu e)=2\mu D_eu(x)+\int_0^\mu(\mu-\tau)\big[u_{ee}(x+\tau e)-u_{ee}(x-\tau e)\big]d\tau,\\
u(x,s)-u(x,s-\mu^2)=\mu^2\partial_su(X)+\int_0^{\mu^2}\big[\partial_su(x,s-\tau)-\partial_su(X)\big]d\tau ,
\end{gather*}
and the integrals are bounded by $[D^2u]_\alpha\mu^{2+\alpha}$ (using
$\int_0^\mu(\mu-\tau)\tau^\alpha d\tau\le\mu^{2+\alpha}/2$), by $2^{\alpha-1}[D^2u]_\alpha\mu^{2+\alpha}$
and by $[\partial_su]_\alpha\mu^{2+\alpha}$. Hence
\begin{gather*}
|D^2_{ee}u(X)|\le4\mu^{-2}\|u\|_0+\mu^\alpha[D^2u]_\alpha,\qquad
|D_eu(X)|\le\mu^{-1}\|u\|_0+\mu^{1+\alpha}[D^2u]_\alpha,\\
|\partial_su(X)|\le2\mu^{-2}\|u\|_0+\mu^\alpha[\partial_su]_\alpha .
\end{gather*}
Polarisation ($D_{ij}u=D^2_{ee}u-\frac12D^2_{e_ie_i}u-\frac12D^2_{e_je_j}u$ for
$e=(e_i+e_j)/\sqrt2$) gives $\|D^2u\|_{0;P_{1/2}(0)}\le C(n)(\mu^{-2}\|u\|_0+\mu^\alpha[D^2u]_\alpha)$.
Choosing $\mu:=\min\{\frac12,(\epsilon/C(n))^{1/\alpha}\}$ gives \eqref{eq:interp1}
with $C_1\epsilon^{-2/\alpha}\ge C(n)\mu^{-2}$.

For \eqref{eq:interp2} let $X=(x,s)$, $Y=(y,t)\in P_{1/2}(0)$, say $t\le s$, and put
$h:=(s-t)^{1/2}<\frac12$. Then $|Du(x,s)-Du(y,s)|\le\|D^2u\|_0|x-y|$, since the segment
lies in $B_{1/2}\times\{s\}\subset P_1(0)$. For the time increment, we use the second
Taylor formula above at $(y,s)$ and at $(y,t)$ with the step $h$; the points
$(y\pm he,s)$, $(y\pm he,t)$ and the time segments between them lie in $P_1(0)$, since
$|y|+h<1$ and $t>-\frac14$. With the remainder now bounded by $h^2\|D^2u\|_0$, this gives
\[
2h\,|D_eu(y,s)-D_eu(y,t)|\le\sum_\pm|u(y\pm he,s)-u(y\pm he,t)|+2h^2\|D^2u\|_0
\le2h^2\|\partial_su\|_0+2h^2\|D^2u\|_0 ,
\]
so $|D_eu(y,s)-D_eu(y,t)|\le h(\|\partial_su\|_0+\|D^2u\|_0)$ when $h>0$. Since
$\|X-Y\|<1$ on $P_{1/2}(0)$, we get $|Du(X)-Du(Y)|\le C(n)(\|D^2u\|_0+\|\partial_su\|_0)\|X-Y\|^\alpha$,
which is \eqref{eq:interp2}.
\end{proof}

\subsection{Constant coefficients on backward cylinders}

\begin{lemma}[constant coefficients]\label{lem:const}
Let $A$ be a constant symmetric matrix with $\lambda_0I\le A\le\Lambda I$,
$0<\lambda_0\le1\le\Lambda$, let $L_A:=\partial_s-A^{ij}D_{ij}$, and put
$\gamma:=\frac12(\lambda_0/\Lambda)^{1/2}\le\frac12$. There is $C_0=C_0(n,\alpha,\lambda_0,\Lambda)$
such that for every backward cylinder $P_\sigma(X_0)$ and every $v\in C^{2,\alpha}(P_\sigma(X_0))$,
\begin{equation}\label{eq:const}
\sigma^{2+\alpha}\big([D^2v]_{\alpha;P_{\gamma\sigma}(X_0)}+[\partial_sv]_{\alpha;P_{\gamma\sigma}(X_0)}\big)
\le C_0\big(\|v\|_{0;P_\sigma(X_0)}+\sigma^2\|L_Av\|_{0;P_\sigma(X_0)}+\sigma^{2+\alpha}[L_Av]_{\alpha;P_\sigma(X_0)}\big).
\end{equation}
\end{lemma}
\begin{proof}
The map $(x,s)\mapsto(x_0+\sigma x,s_0+\sigma^2s)$ sends $P_1(0)$ onto $P_\sigma(X_0)$,
$P_\gamma(0)$ onto $P_{\gamma\sigma}(X_0)$, and transforms $L_A$ into $\sigma^2L_A$. This
produces the powers of $\sigma$ in \eqref{eq:const}, so we may take $X_0=0$ and
$\sigma=1$. The linear change of variables $x=A^{1/2}\tilde x$ transforms $L_A$ into the
heat operator $\partial_s-\Delta_{\tilde x}$: for $\tilde v(\tilde x,s):=v(A^{1/2}\tilde x,s)$
we have $D^2_{\tilde x}\tilde v=A^{1/2}(D^2v)A^{1/2}$ and
$\Delta_{\tilde x}\tilde v=A^{ij}D_{ij}v$.
Since $|A^{1/2}\tilde x|\le\Lambda^{1/2}|\tilde x|$, the set $\{|\tilde x|<\Lambda^{-1/2}\}$
is mapped into $B_1$, so $\tilde v$ is defined on the backward cylinder
$\tilde P:=\{|\tilde x|<\Lambda^{-1/2}\}\times(-\Lambda^{-1},0]$, whose image in the $x$
variables lies in $P_1(0)$. We apply Theorem~8.2 of \cite{wh} with $T=0$; its domains
$\{\|X\|<r,\ \tau(X)\le0\}$ are exactly the backward cylinders $P_r(0)$. With the radii
$r'=\Lambda^{-1/2}$ and $r=\frac12\Lambda^{-1/2}$, it gives
\[
\|\tilde v\|_{2,\alpha;P_r(0)}\le c\big(\|\tilde v\|_{0;\tilde P}+\|(\partial_s-\Delta)\tilde v\|_{0,\alpha;\tilde P}\big),\qquad c=c(n,\alpha,\Lambda),
\]
the norm on the left being the parabolic $C^{2,\alpha}$ norm of \cite[\S8]{wh}, which
contains $[D^2_{\tilde x}\tilde v]_\alpha+[\partial_s\tilde v]_\alpha$. We now return to
the $x$ variables. Since $|A^{-1/2}x|\le\lambda_0^{-1/2}|x|$, the cylinder $P_r(0)$ in the $\tilde x$
variables contains, in the $x$ variables, the set $\{|x|<\lambda_0^{1/2}r\}\times(-r^2,0]\supset P_\gamma(0)$,
because $\gamma=\lambda_0^{1/2}r$ and $\gamma^2\le r^2$. The change of variables is
linear with $\|A^{1/2}\|\le\Lambda^{1/2}$ and $\|A^{-1/2}\|\le\lambda_0^{-1/2}$, so
parabolic distances and the H\"older seminorms of $D^2v$, $\partial_sv$ and $L_Av$ in the
two sets of variables are comparable with constants depending on $\lambda_0,\Lambda,\alpha$
only. This gives \eqref{eq:const}.
\end{proof}

\subsection{Weighted seminorms and the variable-coefficient estimate}

Let $0<\sigma_0\le1$ be a constant, fixed in the proof of Theorem~\ref{thm:schauder-var}.
For $X\in Q_1$ put
\[
\mathsf d(X):=\min\big\{\sigma_0,\ 1-|x|,\ (1+s)^{1/2}\big\},
\]
the parabolic distance from $X$ to the lateral and lower boundary of $Q_1$, capped by
$\sigma_0$. The upper end $s=T$ does not enter. Two properties are used throughout:
\begin{enumerate}[(a)]
\item $P_{\mathsf d(X)}(X)\subset Q_1$ for every $X\in Q_1$: on this cylinder
$|y|<|x|+\mathsf d(X)\le1$, $t>s-(1+s)=-1$, and $t\le s\le T$, $t\le s<1$;
\item $\mathsf d$ is $1$-Lipschitz for the parabolic distance, since $1-|x|$ is
$1$-Lipschitz in $x$ and $|(1+s)^{1/2}-(1+t)^{1/2}|\le|s-t|^{1/2}$. In particular
$\mathsf d\ge\mathsf d(X)-\sigma$ on $P_\sigma(X)$.
\end{enumerate}
For $u\in C^{2,\alpha}(Q_1)$ and $X,Y\in Q_1$ let $\mathsf d_{X,Y}:=\min\{\mathsf d(X),\mathsf d(Y)\}$.
We define the interior weighted quantities as in \cite[(4.17), (6.10)]{gt}, adapted to the
parabolic setting: the distance to the boundary in \cite{gt} is replaced by the capped
parabolic distance $\mathsf d$, and the Euclidean distance by the parabolic distance. Only
the weighted structure is borrowed, and all estimates for these quantities are proved
below. The quantities are
\begin{gather*}
[u]^*_1:=\sup_{X\in Q_1}\mathsf d(X)|Du(X)|,\qquad
[u]^*_2:=\sup_{X\in Q_1}\mathsf d(X)^2|D^2u(X)|,\qquad
[u]^*_\partial:=\sup_{X\in Q_1}\mathsf d(X)^2|\partial_su(X)|,\\
[Du]^*_{1,\alpha}:=\sup_{X\ne Y}\mathsf d_{X,Y}^{1+\alpha}\frac{|Du(X)-Du(Y)|}{\|X-Y\|^\alpha},\qquad
[D^2u]^*_{2,\alpha}:=\sup_{X\ne Y}\mathsf d_{X,Y}^{2+\alpha}\frac{|D^2u(X)-D^2u(Y)|}{\|X-Y\|^\alpha},\\
[\partial_su]^*_{2,\alpha}:=\sup_{X\ne Y}\mathsf d_{X,Y}^{2+\alpha}\frac{|\partial_su(X)-\partial_su(Y)|}{\|X-Y\|^\alpha},
\end{gather*}
and the two sums are
\[
S(u):=[D^2u]^*_{2,\alpha}+[\partial_su]^*_{2,\alpha},\qquad W(u):=[u]^*_1+[u]^*_2+[u]^*_\partial+[Du]^*_{1,\alpha}.
\]
If $|u|_{2,\alpha;Q_1}<\infty$, then $S(u)+W(u)\le4|u|_{2,\alpha;Q_1}$, since $\mathsf d\le\sigma_0\le1$.

\begin{lemma}[weighted interpolation]\label{lem:winterp}
For $u\in C^{2,\alpha}(Q_1)$ with $|u|_{2,\alpha;Q_1}<\infty$ and $0<\epsilon\le1$,
\[
[u]^*_1+[u]^*_2+[u]^*_\partial\le\epsilon\,S(u)+C_1'\epsilon^{-2/\alpha}\|u\|_{0;Q_1},\qquad
[Du]^*_{1,\alpha}\le C_1'\big([u]^*_1+[u]^*_2+[u]^*_\partial\big),
\]
with $C_1'=C_1'(n,\alpha)$.
\end{lemma}
\begin{proof}
Fix $X\in Q_1$ and write $\delta_X:=\mathsf d(X)$ for the value of the weight at $X$. This
is a local abbreviation, used in this proof and in the proof of Lemma~\ref{lem:compare};
the letter $d$ is reserved for the zero-order coefficient of
Theorem~\ref{thm:schauder-var}. We apply \eqref{eq:interp1} on
$P_{\delta_X/2}(X)\subset Q_1$ (radius $r=\delta_X/2$), whose inner cylinder
$P_{\delta_X/4}(X)$ contains $X$. On $P_{\delta_X/2}(X)$ we
have $\mathsf d\ge\delta_X/2$ by (b), so $\mathsf d_{Y,Z}\ge\delta_X/2$ for $Y,Z\in P_{\delta_X/2}(X)$ and
$(\delta_X/2)^{2+\alpha}[D^2u]_{\alpha;P_{\delta_X/2}(X)}\le[D^2u]^*_{2,\alpha}$, and likewise for
$\partial_su$. Hence
\[
\tfrac{\delta_X}2|Du(X)|+\tfrac{\delta_X^2}4|D^2u(X)|+\tfrac{\delta_X^2}4|\partial_su(X)|\le\epsilon'\,S(u)+C_1\epsilon'^{-2/\alpha}\|u\|_{0;Q_1}
\]
for every $0<\epsilon'\le1$. Taking the supremum over $X$ with $\epsilon'=\epsilon/12$
gives the first inequality.

For the second inequality, let $X\ne Y\in Q_1$, let $Z$ be the later of the two points (the
one with the larger time; either one if the times agree), and put
$\delta_Z:=\mathsf d(Z)\ge\mathsf d_{X,Y}$.
If $\|X-Y\|<\delta_Z/4$, both points lie in $P_{\delta_Z/4}(Z)$ (the earlier point has time in
$(s_Z-\delta_Z^2/16,s_Z]$ and spatial distance $<\delta_Z/4$ from the spatial component of $Z$), so by \eqref{eq:interp2} on
$P_{\delta_Z/2}(Z)$, where $\mathsf d\ge\delta_Z/2$,
\begin{multline*}
\mathsf d_{X,Y}^{1+\alpha}\frac{|Du(X)-Du(Y)|}{\|X-Y\|^\alpha}\le\delta_Z^{1+\alpha}[Du]_{\alpha;P_{\delta_Z/4}(Z)}\\
\le2^{1+\alpha}C_1\,\tfrac{\delta_Z^2}4\big(\|D^2u\|_{0;P_{\delta_Z/2}(Z)}+\|\partial_su\|_{0;P_{\delta_Z/2}(Z)}\big)
\le2^{1+\alpha}C_1\big([u]^*_2+[u]^*_\partial\big).
\end{multline*}
If $\|X-Y\|\ge\delta_Z/4\ge\mathsf d_{X,Y}/4$, then
$\mathsf d_{X,Y}^{1+\alpha}|Du(X)-Du(Y)|/\|X-Y\|^\alpha\le4^\alpha\mathsf d_{X,Y}(|Du(X)|+|Du(Y)|)\le2\cdot4^\alpha[u]^*_1$.
\end{proof}

\begin{lemma}[local and weighted seminorms]\label{lem:compare}
Let $0<\gamma\le\frac12$ and $u\in C^{2,\alpha}(Q_1)$ with $|u|_{2,\alpha;Q_1}<\infty$, and put
\[
S_\gamma(u):=\sup_{Z\in Q_1}\mathsf d(Z)^{2+\alpha}\Big([D^2u]_{\alpha;P_{\gamma\mathsf d(Z)/2}(Z)}+[\partial_su]_{\alpha;P_{\gamma\mathsf d(Z)/2}(Z)}\Big).
\]
Then $S(u)\le2S_\gamma(u)+4(2/\gamma)^\alpha\big([u]^*_2+[u]^*_\partial\big)$, and for every
$Z\in Q_1$ and $0<\sigma\le\mathsf d(Z)/2$,
\begin{gather*}
\sigma^{2+\alpha}\big([D^2u]_{\alpha;P_\sigma(Z)}+[\partial_su]_{\alpha;P_\sigma(Z)}\big)\le S(u),\qquad
\sigma^{1+\alpha}[Du]_{\alpha;P_\sigma(Z)}\le[Du]^*_{1,\alpha},\\
\sigma^j\|D^ju\|_{0;P_\sigma(Z)}\le[u]^*_j\ (j=1,2),\qquad \sigma^2\|\partial_su\|_{0;P_\sigma(Z)}\le[u]^*_\partial .
\end{gather*}
\end{lemma}
\begin{proof}
For the second group of inequalities, note that $\mathsf d\ge\mathsf d(Z)-\sigma\ge\sigma$
on $P_\sigma(Z)$ by (b), so every weight $\mathsf d(X)$ or $\mathsf d_{X,Y}$ there is at
least $\sigma$. For the first inequality, let $X\ne Y\in Q_1$, $f\in\{D^2u,\partial_su\}$, $Z$ the
later point and $\delta_Z:=\mathsf d(Z)\ge\mathsf d_{X,Y}$, as in the proof of Lemma~\ref{lem:winterp}.
If $\|X-Y\|<\gamma\delta_Z/2$, both points lie
in $P_{\gamma\delta_Z/2}(Z)$ and $\mathsf d_{X,Y}^{2+\alpha}|f(X)-f(Y)|/\|X-Y\|^\alpha\le\delta_Z^{2+\alpha}[f]_{\alpha;P_{\gamma\delta_Z/2}(Z)}\le S_\gamma(u)$.
If $\|X-Y\|\ge\gamma\delta_Z/2\ge\gamma\mathsf d_{X,Y}/2$, then
$\mathsf d_{X,Y}^{2+\alpha}|f(X)-f(Y)|/\|X-Y\|^\alpha\le(2/\gamma)^\alpha\mathsf d_{X,Y}^2(|f(X)|+|f(Y)|)\le2(2/\gamma)^\alpha([u]^*_2+[u]^*_\partial)$.
Summing over the two choices of $f$ gives the bound stated.
\end{proof}

\begin{theorem}[interior Schauder estimate, variable coefficients, one-sided]\label{thm:schauder-var}
Let $0<\lambda_0\le1\le\Lambda$, $0<\alpha<1$, $N\ge1$, $0<\rho<1$ and $T\in[0,\infty]$.
Let $u\in C^{2,\alpha}(Q_1;\R^N)$ with $|u|_{2,\alpha;Q_1}<\infty$ satisfy
\[
\partial_su^\nu=a^{ij}(x,s)D_{ij}u^\nu+c^k(x,s)D_ku^\nu+d(x,s)u^\nu+f^\nu(x,s),\qquad\nu=1,\dots,N,
\]
on $Q_1$, with $a^{ij}=a^{ji}$, $a^{ij}\xi_i\xi_j\ge\lambda_0|\xi|^2$ and
$\|a\|_{0,\alpha;Q_1}+\|c\|_{0,\alpha;Q_1}+\|d\|_{0,\alpha;Q_1}\le\Lambda$. Then
\[
|u|_{2,\alpha;Q_\rho}\le C\big(\|u\|_{0;Q_1}+\|f\|_{0,\alpha;Q_1}\big),
\]
with $C=C(n,N,\alpha,\lambda_0,\Lambda,\rho)$ independent of $T$. By parabolic dilation
the same holds for $Q_{\rho}\subset Q_{\rho'}$, $0<\rho<\rho'$, in place of
$Q_\rho\subset Q_1$, with $C$ depending also on $\rho'$; and if $u$ is only known to be
$C^{2,\alpha}$ with finite norm on the sets $Q_{\rho'}$, $\rho'<1$, the estimate on
$Q_\rho$ follows by applying it on $Q_{\rho'}$ for one $\rho'\in(\rho,1)$.
\end{theorem}
\begin{proof}
The system is diagonal in the component index $\nu$ with the same scalar operator, so it
suffices to treat $N=1$. In what follows, $C$ denotes constants depending on
$n,\alpha,\lambda_0,\Lambda$ only, until $\rho$ enters in the last step. Let
$\gamma=\frac12(\lambda_0/\Lambda)^{1/2}$ and $C_0$ be as in Lemma~\ref{lem:const}, and let
$S_\gamma(u)$, $S(u)$, $W(u)$ be the quantities defined above for a $\sigma_0\in(0,1]$ to
be chosen. All of them are finite.

\emph{Step 1: freezing at one point.} Fix $Z\in Q_1$ and put $\sigma:=\mathsf d(Z)/2$,
so that $P_\sigma(Z)\subset Q_1$ and $\mathsf d\ge\sigma$ on $P_\sigma(Z)$. On $P_\sigma(Z)$
we write the equation as $L_{a(Z)}u=g$ with
\[
g:=f+\big(a^{ij}-a^{ij}(Z)\big)D_{ij}u+c^kD_ku+du .
\]
Since $|a(X)-a(Z)|\le\Lambda\|X-Z\|^\alpha\le\Lambda\sigma^\alpha$ on $P_\sigma(Z)$, and
$[u]_{\alpha;P_\sigma(Z)}\le(2\sigma)^{1-\alpha}(\|Du\|_{0;P_\sigma(Z)}+2\sigma\|\partial_su\|_{0;P_\sigma(Z)})$
by the mean value theorem on the convex set $P_\sigma(Z)$ (where $\|X-Y\|<2\sigma$), we
have, with all norms over $P_\sigma(Z)$,
\begin{gather*}
\sigma^2\|g\|_0\le\sigma^2\|f\|_0+\Lambda\sigma^{2+\alpha}\|D^2u\|_0+\Lambda\sigma^2\|Du\|_0+\Lambda\sigma^2\|u\|_0,\\
\sigma^{2+\alpha}[g]_\alpha\le\sigma^{2+\alpha}[f]_\alpha+\Lambda\sigma^{2+2\alpha}[D^2u]_\alpha
+\Lambda\sigma^{2+\alpha}\big(\|D^2u\|_0+[Du]_\alpha+\|Du\|_0\big)\\
\qquad\qquad+\Lambda\sigma^{2+\alpha}\big([u]_\alpha+\|u\|_0\big).
\end{gather*}
By Lemma~\ref{lem:compare} (second group, with this $\sigma=\mathsf d(Z)/2$) and
$\sigma\le\frac12$,
\[
\sigma^2\|g\|_0+\sigma^{2+\alpha}[g]_\alpha\le\|f\|_{0,\alpha;Q_1}+\Lambda\sigma_0^\alpha S(u)+C\Lambda\big(W(u)+\|u\|_{0;Q_1}\big),
\]
where the term $\Lambda\sigma^{2+2\alpha}[D^2u]_{\alpha;P_\sigma(Z)}\le\Lambda\sigma^\alpha S(u)\le\Lambda\sigma_0^\alpha S(u)$
is the only one carrying the top seminorm.

The frozen matrix $a(Z)$ is symmetric and satisfies $a(Z)\ge\lambda_0I$ by hypothesis. It
also satisfies $a(Z)\le\Lambda I$, because
the largest eigenvalue of a symmetric matrix is at most the sum of the absolute values of
its entries ($a^{ij}(Z)\xi_i\xi_j\le\sum_{i,j}|a^{ij}(Z)|$ for $|\xi|=1$, as $|\xi_i\xi_j|\le1$)
and $\sum_{i,j}|a^{ij}(Z)|\le\|a\|_{0;Q_1}\le\Lambda$ by the componentwise convention.
Hence Lemma~\ref{lem:const} applies on $P_\sigma(Z)$ with $A=a(Z)$ and the stated
$\lambda_0$, $\Lambda$, $\gamma$, and gives
\begin{multline*}
\mathsf d(Z)^{2+\alpha}\big([D^2u]_{\alpha;P_{\gamma\mathsf d(Z)/2}(Z)}+[\partial_su]_{\alpha;P_{\gamma\mathsf d(Z)/2}(Z)}\big)
=2^{2+\alpha}\sigma^{2+\alpha}\big([D^2u]_{\alpha;P_{\gamma\sigma}(Z)}+[\partial_su]_{\alpha;P_{\gamma\sigma}(Z)}\big)\\
\le C\big(\Lambda\sigma_0^\alpha S(u)+W(u)+\|u\|_{0;Q_1}+\|f\|_{0,\alpha;Q_1}\big).
\end{multline*}

\emph{Step 2: supremum and absorption.} Taking the supremum over $Z\in Q_1$ gives
$S_\gamma(u)\le C(\Lambda\sigma_0^\alpha S(u)+W(u)+\|u\|_{0;Q_1}+\|f\|_{0,\alpha;Q_1})$, and the
first inequality of Lemma~\ref{lem:compare} turns this into
\[
S(u)\le C_2\Lambda\sigma_0^\alpha\,S(u)+C_2\big(W(u)+\|u\|_{0;Q_1}+\|f\|_{0,\alpha;Q_1}\big),
\qquad C_2=C_2(n,\alpha,\lambda_0,\Lambda).
\]
By Lemma~\ref{lem:winterp}, $W(u)\le(1+C_1')([u]^*_1+[u]^*_2+[u]^*_\partial)\le(1+C_1')(\epsilon S(u)+C_1'\epsilon^{-2/\alpha}\|u\|_{0;Q_1})$
for every $0<\epsilon\le1$. Now fix $\sigma_0\in(0,1]$ with $C_2\Lambda\sigma_0^\alpha\le\frac14$
and $\epsilon\in(0,1]$ with $C_2(1+C_1')\epsilon\le\frac14$. Both choices depend only on
$n,\alpha,\lambda_0,\Lambda$. Since $S(u)<\infty$, the term $\frac12S(u)$ can be
absorbed, and we obtain
\[
S(u)\le C\big(\|u\|_{0;Q_1}+\|f\|_{0,\alpha;Q_1}\big),\qquad\text{hence also}\qquad W(u)\le C\big(\|u\|_{0;Q_1}+\|f\|_{0,\alpha;Q_1}\big).
\]

\emph{Step 3: the unweighted norm on $Q_\rho$.} On $Q_\rho$ we have
$\mathsf d\ge\mathsf d_\rho:=\min\{\sigma_0,1-\rho,(1-\rho^2)^{1/2}\}>0$, so every weight in
$S(u)$ and $W(u)$ restricted to points of $Q_\rho$ is at least $\mathsf d_\rho$. Thus
$|u|_{2,\alpha;Q_\rho}\le\|u\|_{0;Q_1}+\mathsf d_\rho^{-2-\alpha}(S(u)+W(u))$, which gives the
theorem with $C=C(n,\alpha,\lambda_0,\Lambda,\rho)$. For $N>1$ the estimate is applied to
each component.

\emph{Independence of $T$.} The cylinders $P_\sigma(Z)$, $Z\in Q_1$, are backward and lie
in $Q_1$, so Theorem~8.2 of \cite{wh} is used only in its $T=0$ form (Lemma~\ref{lem:const}).
The weight $\mathsf d$ does not involve the upper end, and the interpolation inequalities
are proved on backward cylinders. No constant depends on $T$, and no point with time
above $T$ is ever used.
\end{proof}

\begin{remark}
The estimate is applied in Lemma~\ref{lem:subord} with $a=\bar a(x,s)$, $c=0$, $d=0$
(for the spatial derivatives $D_ku$ and $D_{kl}u$, with the terms produced by differentiation of the coefficients
on the right-hand sides) and with $a=G^{ij}$, $c^k=-G^{ij}\,{}^G\Gamma^k_{ij}-H^k$, $d=\kappa_\theta$,
$f=0$ (for $\theta$). The coefficients are bounded in $C^{0,\alpha}$ by constants
depending on $n$, $\delta_*$ and $\|u\|_{2,\alpha}\le1$ only. The solutions to which the
theorem is applied ($D_ku$, $D_{kl}u$, $\theta$) have finite parabolic $C^{2,\alpha}$ norms on
the outer domains used there ($\Omega^{(7/8)}_{T'}$, $\Omega^{(5/8)}_{T'}$, $\Omega^{(3/4)}_{T'}$),
whose closures are compact subsets of the base domain up to its closed upper end. These
norms are finite because the witness $u$ and the angle $\theta(x,s)$ are smooth up to the
upper end by the last assertion of Lemma~\ref{lem:filling}. Only this finiteness is used,
never the size of the norms. The domains
$\Omega^{(\rho)}_{T'}$ of Lemma~\ref{lem:subord} are the cylinders $Q_\rho$ with $T=T'$, and the inner
domain $\Omega^{(9/16)}_{T'}$ contains $\Omega'_{T'}$.
\end{remark}

\section*{AI Declaration}
ChatGPT (OpenAI) and Claude (Anthropic) assisted with language editing and
exploratory calculations. The authors independently checked all mathematical
arguments and computations.

\bibliographystyle{amsalpha}
\bibliography{references}
\end{document}